\documentclass[review]{elsarticle} 

\usepackage{hyperref}

\journal{Applied Mathematical Modelling}

\usepackage{amssymb, amsmath, amsthm}
\usepackage{amsfonts} 
\usepackage{color}
\usepackage{mathrsfs}
\usepackage{comment}
\usepackage{dcolumn}
\usepackage{bm} 
\usepackage{listings}

\usepackage[hang,small,bf]{caption}
\usepackage[subrefformat=parens]{subcaption}
\numberwithin{equation}{section}
\newcommand{\e}{\varepsilon}
\renewcommand{\d}{\mathrm{d}}
\newcommand{\R}{\mathbb R}
\newcommand{\N}{\mathbb N}

\usepackage{verbatim}
\usepackage{subcaption}
\usepackage{multicol}
\newtheorem{thm}{Theorem}[section]

\newtheorem{remark}[thm]{Remark}

\newtheorem{definition}[thm]{Definition}

\newtheorem{example}[thm]{Example}

\begin{document}

\begin{frontmatter}

\title{Nesterov's acceleration 
for level set-based topology optimization 
using reaction-diffusion equations\tnoteref{mytitlenote}}
\tnotetext[mytitlenote]{This study is partially supported 
by a project JPNP20004 subsidized by the New Energy and Industrial Technology Development Organization (NEDO) and JSPS KAKENHI Grant Number JP22K20331.}

\author[mymainaddress]{Tomoyuki Oka}
\ead{tomoyuki-oka@g.ecc.u-tokyo.ac.jp}

\author[mymainaddress2]{Ryota Misawa}

\author[mymainaddress]{Takayuki Yamada\corref{mycorrespondingauthor}}
\cortext[mycorrespondingauthor]{Corresponding author}
\ead{t.yamada@mech.t.u-tokyo.ac.jp}

\address[mymainaddress]{Graduate School of Engineering, The University of Tokyo, \\Yayoi 2-11-16, Bunkyo-ku, Tokyo 113-8656, Japan}

\address[mymainaddress2]{Graduate School of Science and Engineering, Saitama University, \\
255 Shimo-Okubo, Sakura-ku, Saitama City, Saitama 338-8570, Japan
}

\begin{abstract}
This paper discusses level set-based structural optimization. 
Level set-based structural optimization is a method used to determine an optimal configuration for minimizing objective functionals by updating level set functions characterized as solutions to partial differential equations (PDEs) (e.g.,~Hamilton-Jacobi and reaction-diffusion equations). 
In this study, based on Nesterov's accelerated gradient method, a nonlinear (damped) wave equation will be derived as a PDE satisfied by level set functions and applied to minimum mean compliance problems. 
Numerically, the method developed in this study will yield convergence to an optimal configuration faster than methods using only a reaction-diffusion equation, and moreover, its FreeFEM++ code will also be described.
\end{abstract}

\begin{keyword} 
topology optimization, level set method, reaction-diffusion equation, nonlinear (damped) wave equation, Nesterov's accelerated gradient method
\MSC[2010] {\emph{Primary}: 80M50; \emph{Secondary}: 35Q93, 47J35}
\end{keyword}
\end{frontmatter}


\section{Introduction}\label{sec1}

\emph{Topology optimization} is a structural optimization with the highest degree of design freedom.
 It is a method for determining an optimal material configuration (denoted by a set $\Omega_{\rm opt}\subset \R^{d}$ ($d=2$ or $3$) below) to minimize objective functionals and is currently being developed in various fields. 
Generally, the following objective functional $F_\Omega:H^1(\Omega)^d\to \R$ is treated in usual variational methods\/{\rm :} 
\begin{align*}
F_\Omega(\boldsymbol{v})=\int_{\Omega} f_{\rm d}(x,\boldsymbol{v},\boldsymbol{\nabla} \boldsymbol{v})\, \d x + \int_{\partial \Omega} f_{\rm b}(x,\boldsymbol{v},\boldsymbol{\nabla} \boldsymbol{v})\, \d\sigma,
\end{align*}
where $f_{\rm d}$ and $f_{\rm b}$ are functions defined on $\Omega\times \R^d\times \R^{d\times d}$ and $\partial\Omega\times \R^d\times \R^{d\times d}$, respectively.
On the other hand, the following set function 
$F:\mathcal{U}_{\rm ad}\to\R$ is treated in topology optimization\/{\rm:}
\begin{align*}
F(\Omega)=\int_{\Omega} f_{\rm d}(x,\boldsymbol{u}_\Omega,\boldsymbol{\nabla} \boldsymbol{u}_\Omega)\, \d x + \int_{\partial \Omega} f_{\rm b}(x,\boldsymbol{u}_\Omega,\boldsymbol{\nabla} \boldsymbol{u}_\Omega)\, \d\sigma.
\end{align*}
Here, $\mathcal{U}_{\rm ad}:=$\{$\Omega \subset D\colon$constraint conditions\} and $D\subset \R^{d}$ denote 
a family of sets for admissible domains and a fixed design domain such that $\partial\Omega\cap \partial D\neq \emptyset$, respectively, and the state variable $\boldsymbol{u}_\Omega\in H^1(\Omega)^d$ is a (vector-valued) function determined by fixing $\Omega\subset D$; for instance, 
a solution to the Euler-Lagrange equation (i.e.,~its Fr\'echet derivative $F_{\Omega}'$ is zero). 
Then consider the following minimization problem\/{\rm :}
\begin{align}\label{min-prob1}
\inf_{\Omega\in \mathcal{U}_{\rm ad}} F(\Omega).
\end{align}
Therefore, topology optimization can be regarded as a minimization problem with a set $\Omega\in \mathcal{U}_{\rm ad}$ as a variable, and \eqref{min-prob1} can be replaced by the following distribution problem of materials\/{\rm :} 
\begin{align}\label{min-prob2}
\inf_{\chi_\Omega\in L^{\infty}(D;\{0,1\})} F(\chi_\Omega)
\quad \text{subject to constraint conditions},
\end{align}
where $F:L^{\infty}(D;\{0,1\})\to \R$ is a functional given by
\begin{align*}
F(\chi_\Omega)=\int_{D} f_{\rm d}(x,\boldsymbol{u}_\Omega,\boldsymbol{\nabla} \boldsymbol{u}_\Omega)\chi_\Omega(x)\, \d x 
+ \int_{\partial D} f_{\rm b}(x,\boldsymbol{u}_\Omega,\boldsymbol{\nabla} \boldsymbol{u}_\Omega)\, \d\sigma
\end{align*}
and $\chi_\Omega\in L^{\infty}(D)$ is a characteristic function defined by 
\begin{align*}
\chi_\Omega(x)=
\begin{cases}
1,\quad x\in \overline{\Omega}:=\Omega\cup\partial\Omega,\\
0,\quad x\in D\setminus \overline{\Omega}.
\end{cases}
\end{align*}
Hence, topology optimization is the minimization problem for $\chi_\Omega\in L^{\infty}(D)$. Moreover, $\Omega_{\rm opt}\in \mathcal{U}_{\rm ad}$ may be characterized as the domain of the minimizer $\chi_\Omega\in L^{\infty}(D)$ 
in \eqref{min-prob2} such that $\chi_\Omega=1$.
Thus, topology optimization implies that changes in the shape of $\partial\Omega$ and the topology of $\Omega\in \mathcal{U}_{\rm ad}$, such as an increase or decrease in the number of holes, can be allowed in the optimization procedure.
However, various issues remain on how to determine $\Omega_{\rm opt} \in \mathcal{U}_{\rm ad}$ even if such configurations exist.  

\subsection{Homogenization-based topology optimization.}

The existence of $\Omega_{\rm opt}\in \mathcal{U}_{\rm ad}$
for generalized problems is obtained by the \emph{homogenization theory} based on \emph{$H$ {\rm(}or $G${\rm)}-convergence}  (see, e.g.,~\cite[Theorem 3.2.1]{A02} and \cite{MT97}), and therefore, it forms the basis for numerical analysis. 
Homogenization is a method for replacing heterogeneous materials, which possess many microstructures, with an equivalent homogeneous material. 
As for the numerical analysis, the so-called \emph{homogenization design method} was first developed in \cite{BK88}.
Moreover, its simplified version, the so-called \emph{
Solid Isotropic Material with Penalization {\rm (}SIMP{\rm )} method} \cite{B88}, is frequently used by replacing the characteristic functions with density functions, but there are certain issues\/{\rm :}
~(i)\,An optimized configuration $\Omega_{\rm opt} \in \mathcal{U}_{\rm ad}$ typically includes grayscale domains since the density function takes a value in $[0,1]$. 
In other words, $\partial\Omega_{\rm opt}$ is not clearly expressed unless such domains are removed.
(ii)\,True material properties of composite materials with microstructures are not devised due to oversimplification; indeed,
it holds only material density (see, e.g.,~\cite{ACMOY19} for the resurrection of the homogenization method \cite{S07,W11} for filtering  and \cite{ad1,ad6} for meshless).

\subsection{Level set-based structural optimization} 
To overcome the aforementioned issues, a level set method, which was first proposed in a previous study \cite{OS88} to implicitly represent the evolution of interfaces, was introduced, and the following level set function $\phi\in H^1(D)$ was employed\/{\rm:}
\begin{align}\label{levelset}
\phi(x)
\begin{cases}
>0,
\quad &x\in\Omega,\\
=0,\quad &x\in\partial\Omega,\\
<0,\quad &x\in D\setminus \overline{\Omega}.\\
\end{cases}
\end{align}
Thus, $[\phi>0]:=\{x\in D\colon \phi(x)>0\}$, $[\phi<0]$ and $[\phi=0]$ represent \emph{material domains}, \emph{void domains} and \emph{structural boundaries}, respectively, and $\chi_{\Omega}\in L^{\infty}(D;\{0,1\})$ in \eqref{min-prob2} can be replaced by  
\begin{align*}
\chi_\phi(x)=
\begin{cases}
1\quad &\text{ if } \phi(x)\ge 0,\\  
0\quad &\text{ if } \phi(x)< 0.
\end{cases}
\end{align*}
Moreover, $\phi\in H^1(D)$ is determined by combining  Lagrange's method of undetermined multipliers with the Karush-Kuhn-Tucker (KKT) conditions since $F(\chi_\Omega)$ can also be replaced with $F(\phi)$.
However, the direct derivation of $\phi\in H^1(D)$ is nearly impossible in general. 

Alternatively, some methods that can be used to update 
\eqref{levelset} by introducing a fictitious time variable $t\in (0,+\infty)$ have been devised. 
As a typical example for \eqref{levelset}, the following signed distance function is known\/{\rm :}
\begin{align}\label{distfunc}
\phi(x)=
\begin{cases}
d(x),\quad &x\in\Omega,\\
0,\quad \quad &x\in\partial\Omega,\\  
-d(x),\quad &x\in D\setminus \overline\Omega,
\end{cases}
\end{align}
where $d(x)=\inf_{y\in \partial\Omega}|x-y|$. 
By noting that \eqref{distfunc} solves an eikonal equation, the following Hamilton-Jacobi equation is derived by differentiating for \eqref{distfunc} regarding $t\in (0,\infty)$\/{\rm :} 
\begin{align}\label{HJ}
\partial_t \phi(x,t)+v(x,t) |\nabla \phi (x,t)| =0,
\quad (x,t)\in D\times (0,+\infty),
\end{align}
where $\partial_t=\partial/\partial t$ and $v$ is a given function. 
In previous studies \cite{AJT02, AJT04, WWG03}, the equivalence of transporting the solution to \eqref{HJ} and moving $\partial\Omega$ along the descent gradient direction of functionals was used, and the shape derivative (based on the Fr\'echet derivative) was applied; therefore, optimized shape $\partial\Omega_{\rm opt}$ can be expressed by  employing the solution to \eqref{HJ} and the topology can also be changed by reducing the number of holes, but changing the topologies by generating new holes is not; in other words, $\Omega_{\rm opt}\in \mathcal{U}_{\rm ad}$ deeply depends on initial configurations
 (see also for other level set methods \cite{ad2,ad3,ad4,ad5,ad7}). 
With the aid of the bubble method \cite{EKS94}, a concept of \emph{topological derivative} was introduced in \cite{AGJ05}, and the issue of the dependence on initial configurations might have been resolved. The topological derivative \cite{Td1,Td2} represents an influence when a sufficiently small ball $B_\e(x)\subset \Omega$ is created and is defined as follows\/{\rm :} 
 \begin{definition}[Topological derivative]
 Let $A=\{\Omega\subset D\colon \Omega \text{ is open in } D\}$. 
 A function $J$ defined in $A$ is said to be topologically differentiable at $\Omega_0$ and at point $x\in \Omega_0$ if the following limit exists\/{\rm :}
\begin{align*}
\d_{\rm T}J(\Omega_0,x):=\lim_{\e\to 0_+}\frac{J(\Omega_0\setminus \overline{B_\e(x)})-J(\Omega_0)}{|B_\e(x)|}.
\end{align*}
Here, $B_\e(x):=\{ y\in \Omega_0\colon |x-y|<\e \}$ and $|B_\e(x)|$ denotes the Lebesgue measure of $B_\e(x)$. 
 \end{definition}
 \begin{example}\label{EX}
 \rm
Let $J(\Omega)=|\Omega |$. Then $\d_{\rm T}J(\Omega)=-1$ and 
 $J(\Omega)$ can be regarded as the functional for $\chi_\phi\in L^{\infty}(D)$ by noting that
 $$
J(\Omega)=\int_D \chi_\phi(x)\, \d x=:J(\chi_\phi),
 $$
which implies that $\d_{\rm T} J(\Omega)$ can be identified with $-J'(\chi_\phi)$.
\end{example}
Generally, the solution to \eqref{HJ} does not coincide with \eqref{distfunc}, and note that the bubble method is heuristic. Thus the issue of the dependence on initial configurations still remains (e.g.,~\cite{S99} for the reinitialization of level set functions). 

\subsection{Level set-based topology optimization using a reaction-diffusion equation}
To avoid the dependence on initial configurations, 
in a previous study \cite{YINT10}, 
\eqref{levelset} was modified as follows\/{\rm :}
\begin{align}\label{levelset2}
\phi(x)
\begin{cases}
\in (0,1],
\quad &x\in\Omega,\\
=0,\quad &x\in\partial\Omega,\\
\in [-1,0),
\quad &x\in D\setminus \overline{\Omega}.\\
\end{cases}
\end{align}
Notably, \eqref{distfunc} can no longer be taken as \eqref{levelset2}.
Furthermore, based on the concept of the gradient descent method, the following time evolution equation was employed\/{\rm :} 
\begin{align}\label{ee1}
\partial_t \phi =\rho\d_{\rm T} F,
\end{align}
where $\rho>0$. 
In terms of the regularity, the following (reaction) diffusion equation was applied\/{\rm :}
\begin{align}\label{ee2}
\partial_t\phi=\rho\d_{\rm T}F+\tau\Delta\phi
\quad \text{ in } D\times (0,+\infty),
\end{align}
where $\tau>0$.
Thus $\tau\Delta\phi$ plays a role in the regularization term since \eqref{ee2} is expected to have the smoothing effect, and the solution to \eqref{ee1} may be approximated by the solution to \eqref{ee2} for $\tau>0$ small enough. 
Actually, $\Omega_{\rm opt}\in\mathcal{U}_{\rm ad}$ exhibits complex configurations if $\tau>0$ is small enough and vice versa. Hence, an appropriate value $\tau>0$ prevents the generation of excessively geometrically complex configurations, and $\Omega_{\rm opt}\in\mathcal{U}_{\rm ad}$ with the geometric simplicity is obtained. On the other hand, by employing $\d_{\rm T}F$, the dependence on initial configurations is entirely removed since generating holes in $D$ can be allowed. Therefore, in terms of practicality, level set-based topology optimization with \eqref{ee2} may be valid for various problems since the issues (grayscale domains, reinitialization for level set functions and dependence on initial configurations) are all removed. 
However, the mathematical justification of replacing the topological derivative for the Fr\'echet derivative and optimality remain future issues.
 
\subsection{Setting of the optimization problem}
This paper concerns an optimal design that minimizes the following volume-constrained mean compliance\/{\rm :}
\begin{align}\label{op1}
\inf_{\Omega\in \mathcal{U}_{\rm ad}}
\left\{F(\Omega):=\langle \boldsymbol{t}, \boldsymbol{u}_\Omega\rangle_{H^{1/2}(\Gamma_t)^d}\right\},
\end{align}
where 
$$
\mathcal{U}_{\rm ad}=\{\Omega \subset D \colon |\Omega|\le G_{\text{max}}|D|   \},
$$
$G_{\text{max}}>0$, 
the vector-valued function ${\boldsymbol u}_\Omega\in C^2(\Omega)^d\cap C(\overline{\Omega})^d$ is a unique (classical) solution to the following linearized elasticity system\/{\rm:}
\begin{align*}
\begin{cases}
-\mathrm{\bf div}[\mathbb{D}{\boldsymbol \varepsilon}({\boldsymbol u}_\Omega)]=0&\text{ in } \Omega,\\ 
{\boldsymbol u}_\Omega =0 &\text{ on } \Gamma_D,\\  
-\mathbb{D}{\boldsymbol \varepsilon}({\boldsymbol u}_\Omega)\cdot \boldsymbol{ n}=\boldsymbol{ t}  &\text{ on } \Gamma_t,\\ 
-\mathbb{D}{\boldsymbol \varepsilon}({\boldsymbol u}_\Omega)\cdot \boldsymbol{ n}=0  &\text{ on } \Gamma_N:=\partial\Omega\setminus (\Gamma_D\cup \Gamma_t). 
\end{cases}
\end{align*}
Here and henceforth, $\Gamma_D\cap\Gamma_t=\emptyset$, $|\Omega|$ represents the Lebesgue measure of $\Omega\subset \R^d$. 
The forth-order elastic tensor $\mathbb{D}=\mathbb{D}_{ijk\ell}e_i\otimes e_j\otimes e_k\otimes e_\ell$
and the strain tensor ${\boldsymbol \varepsilon}({\boldsymbol u}_\Omega)$ are given by
\begin{align*}
\mathbb{D}_{ijk\ell}&=\frac{E \nu}{(1+\nu)(1-2\nu)}\delta_{ij}\delta_{k\ell}
+\frac{E}{2(1+\nu)}(\delta_{ik}\delta_{j\ell}+\delta_{i\ell}\delta_{jk})
\end{align*}
for some $E, \nu>0$ and 
$$
{\boldsymbol \varepsilon}({\boldsymbol u}_\Omega)=\frac{1}{2}\left({\boldsymbol \nabla} {\boldsymbol u}_\Omega+({\boldsymbol \nabla} {\boldsymbol u}_\Omega)^{\bf T}\right), \quad
\boldsymbol \nabla {\boldsymbol u}_\Omega=\partial_{x_i}({\boldsymbol{u}_\Omega})_{j}e_i\otimes e_j,  
$$
respectively. The traction $\boldsymbol{ t}\in \R^d$ is a constant vector.  
In particular, $\boldsymbol{n}$, $\delta_{ij}$ and $e_k$ stand for the outer unit normal vector, Kronecker delta and 
$k$-th vector of the canonical basis of $\R^d$, respectively.

The optimization problem \eqref{op1} can be replaced by the following minimization problem with the level set function $\phi\in H^1(D)$ given by \eqref{levelset2}\/{\rm :}  
\begin{align}\label{op2}
\inf_{\phi\in H^1(D;[-1,1])} \left\{F(\phi):=\int_{\Gamma_t} \boldsymbol{ t}\cdot \boldsymbol{ u}_{\phi}(x)\, \d\sigma \right\}
\end{align}
subject to
\begin{align}
G(\phi):=\int_D (\chi_\phi(x)- G_{\text{max}})\, \d x\le 0,
\label{constraint}
\end{align}
where $D$ is a bounded domain in $\R^d$ such that $\Omega\subset D$ and  $\partial D$ is the boundary of $D$ with  $\Gamma_D$ and $\Gamma_t$ such that 
$\Gamma_D\cap \Gamma_t=\emptyset$, 
the state variable ${\boldsymbol u}_\phi\in V$ satisfies
\begin{align}
\int_D \boldsymbol{A}_\phi(x){\boldsymbol \varepsilon}(\boldsymbol{ u}_\phi)(x)\colon {\boldsymbol \varepsilon}({\boldsymbol{ v}})(x)\, \d x
=
\int_{\Gamma_t} {\bf t}\cdot {\boldsymbol{ v}}(x) \, \d\sigma
\quad \text{ for all }\ {\boldsymbol{ v}}\in V, 
\label{gov1}
\end{align}
$\boldsymbol{A}_\phi=\mathbb{D}\chi_\phi+\e {\mathbb{I}}(1-\chi_\phi)$
and 
$
V=\{ \boldsymbol{ v}\in H^1(D)^d\colon \boldsymbol{ v}=0\ \text{ on } \Gamma_D\}
$.
Furthermore, by Lagrange's method of undetermined multipliers, the objective functional $F(\phi)$ in \eqref{op2} is replaced with
\begin{align}
\overline{F}(\phi,\boldsymbol{u}_\phi,\boldsymbol{\tilde{u}},\lambda)
&=\int_{\Gamma_t}\boldsymbol{ t}\cdot (\boldsymbol{ u}_\phi(x)+\tilde{\boldsymbol{ u}}(x))\, \d\sigma \nonumber\\
&\qquad -\int_D \boldsymbol{A}_\phi(x) {\boldsymbol \varepsilon}(\boldsymbol{ u}_\phi)(x)\colon {\boldsymbol \varepsilon}(\tilde{\boldsymbol{ u}})(x)\, \d x
+\lambda G(\phi), \label{min-prob}
\end{align}
which implies that the optimization problem \eqref{op2}--{\eqref{constraint}} may also be replaced with the unconstrained minimization problem for \eqref{min-prob}.  
Here and henceforth, $\overline{F}$ denotes the Lagrangian of $F$ and $\lambda\ge 0$ and $\tilde{\boldsymbol{ u}}\in V$ stand for the Lagrange multiplier. 

\subsection{Aims and plan of this paper}
In this paper, we shall provide a method for level set-based topology optimization that converges to optimized configurations faster than the method based on the (reaction) diffusion equation \eqref{ee2}. 
To this end, instead of the usual gradient descent method, 
 Nesterov's accelerated gradient method \cite{N83} will be introduced, and a nonlinear (damped) wave equation will be applied as a partial differential equation (PDE) to update the level set function (see the next section below for a derivation).

This paper is organized as follows. In the next section, we shall set up time evolution equations to update the level set functions.
In particular, we shall show that the level set function satisfies \eqref{ee2} and a nonlinear damped wave equation, according to the usual gradient descent method and Nesterov's accelerated gradient method \cite{N83}, respectively.
Thus, Section \ref{S:form} offers a new idea and is the most contributing.
Section \ref{S:algo} describes the numerical algorithm for the minimization problem of \eqref{min-prob} and Section \ref{S:reslut} deals with the main results of this paper. 
Furthermore, we shall emphasize that convergence to optimal configurations for the minimization problem of \eqref{min-prob} is improved through typical numerical examples.
The FreeFEM++~\cite{H12} code will be described in 
\ref{Ap} (see also \cite{AP06}). The final section will conclude this paper.

\section{Formulation of nonlinear hyperbolic-parabolic equations}\label{S:form}

To find $\Omega_{\rm opt}\in \mathcal{U}_{\rm ad}$ for the minimization problem of \eqref{min-prob}, we shall formulate the equations satisfied by the level set functions.
To this end, we shall first derive a reaction-diffusion equation. 
Noting that, for any functional $\mathscr{F}:H^1(D)\to\R$, it holds that
\begin{align}\label{perturb}
\mathscr{F}(\phi)\le \mathscr{F}(\phi)+\frac{\tau}{2}\int_D|\nabla \phi(x)|^2\, \d x 
\quad \text{for all $\tau>0$.}
\end{align}
Let $\tilde{\mathscr{F}}:H^1(D)\to\R$ be the right-hand side in \eqref{perturb}. Then, by replacing $\mathscr{F}$ with $\tilde{\mathscr{F}}$,  the gradient descent method such as \eqref{ee1} yields  
\begin{align}\label{rde}
\partial_t\phi=-\rho\mathscr{F}'(\phi)+\tau \Delta \phi,
\end{align}
which implies that, for $\tau>0$ small enough, one can choose \eqref{rde} as the equation which the level set function in \eqref{min-prob} satisfies by noting that
$$
\tilde{\mathscr{F}}(\phi)-\mathscr{F}(\phi)\to 0
\quad \text{ as }\ \tau\to 0_+.
$$
\begin{remark}\label{R:derivation}
\rm

We note that \eqref{rde} does not coincide with \eqref{ee2}. 
However, it is well-known that 
the replacement of $\mathscr{F}'$ with $-\rm{d}_{\rm T}\mathscr{F}$ is adequate in various optimization problems; indeed, if $\varphi\in H^1(D)$ is a 
critical point for $\mathscr{F}$, then $\partial_{\chi_\phi}\mathscr{F}{|_{\phi=\varphi}}= 0$ follows since $\mathscr{F}'= \partial_{\chi_\phi}\mathscr{F}\partial_{\phi}\chi_\phi$ formally. As in Example \ref{EX}, 
identifying $\partial_{\chi_{\phi}}\mathscr{F}$ with $-\rm{d}_{\rm T}\mathscr{F}$, we observe that,
for any $\psi\in H^1(D)$,
$\mathscr{F}(\varphi+k\psi)
=\mathscr{F}(\varphi)+k\langle \mathscr{F}'(\varphi),\psi \rangle_{H^1(D)}+o(k)
= \mathscr{F}(\varphi)-k\langle {\rm d}_{\rm T}\mathscr{F},\psi \rangle_{H^1(D)}+o(k)$
as $k\to 0_+$. Hence, setting $\psi=\rm{d}_{\rm T}\mathscr{F}$, we obtain \eqref{ee1}
for a suitable initial level set function $\phi_0\in L^{\infty}(D)$, which implies that \eqref{rde} can be replaced with \eqref{ee2} under this setting, and therefore, the replacement of $\mathscr{F}'$ with $-\rm{d}_{\rm T}\mathscr{F}$ and a perturbation of the Dirichlet energy yield \eqref{ee2}. 

\end{remark}

\subsection{Nonlinear damped wave equation}
In this subsection, we shall establish another equation satisfied by \eqref{levelset2} to converge faster to an optimized configuration. 
Similar to the gradient descent method, we recall that the following improved gradient descent method, so-called \emph{Nesterov's accelerated gradient method}, was developed in \cite{N83}\/{\rm :}
\begin{align}
&\psi_{n+1}(x)=\phi_{n}(x)-k\mathscr{F}'(\phi_n),\label{AGM}\\
&\phi_n(x)=\psi_{n}(x)+\frac{n-1}{n+2}(\psi_n-\psi_{n-1})\label{inertia}
\end{align}
for $n\in \mathbb{N}$.
Here, $k>0$, $\phi_0$ stands for the initial level set function and $\psi_0=\phi_0$. Therefore, the gradient descent method and \eqref{AGM}-\eqref{inertia} are equivalent until $n=1$ (i.e., $\phi_{m}=\psi_{m}$ for $m=0,1$). 
Conversely, as for $n\ge 2$, the second term of the right-hand side in \eqref{inertia} plays a role in the inertia term.

Now, we consider another equation satisfied by \eqref{levelset2}.
Let $n\in \N$ be large enough and identified with $n+2$ (i.e.,~$n\simeq n+2$). Then, setting $k=(\varDelta t)^2\rho$ and noting that
$$
\psi_{n+1}=\psi_n+(1-3/n)(\psi_n-\psi_{n-1})-k\mathscr{F}'(\phi_n),
$$ 
one can derive that
$$
\frac{\psi_{n+1}-2\psi_n+\psi_{n-1}}{(\varDelta t)^2}+\frac{3}{n\varDelta t}\frac{\psi_n-\psi_{n-1}}{\varDelta t}=-\rho \mathscr{F}'(\phi_n).
$$
Hence, by setting $\psi_{n\pm i}=\psi(x,(n\pm i)\varDelta t)$ and $i=0,1$, we have
\begin{align}\label{nconti}
\partial_{t}^2\psi+(3/t) \partial_t \psi =-\rho\mathscr{F}'(\psi) 
\end{align}
formally (see \cite{SBC14} for justification).
Thus, combining \eqref{perturb} with \eqref{nconti}, we obtain \begin{align}\label{nldw}
\partial_{t}^2\psi+(3/t) \partial_t \psi
&=
-\rho\mathscr{F}'(\psi)+\tau \Delta \psi.
\end{align}

\begin{remark}
\rm
Notably, two initial conditions are required to solve uniquely \eqref{nldw}, which implies that 
it is necessary to have an initial data and another data updated by employing it. In this study, \eqref{rde} will be applied the first few times to construct the initial data and get the same regularity in \eqref{rde} (see also Remarks \ref{R:reg2} and \ref{R:reg} below).
Actually, hyperbolic equations do not have smoothing effects in general. Hence, 
we note that the scheme \eqref{AGM}-\eqref{inertia} does not mean a regularization scheme.
\end{remark}

\subsection{Update of level set functions}
Based on \eqref{rde} and \eqref{nldw}, we shall set up equations which \eqref{levelset2} for the minimization problem of \eqref{min-prob} satisfies. 
By \cite[Appendix B]{OYIN15}, it holds that 
\begin{align*}
\d_{\rm T}\overline{F}(\phi)
=
\mathbb{A}\chi_{\phi}{\boldsymbol \varepsilon}({\boldsymbol u})\colon {\boldsymbol \varepsilon}({\boldsymbol u})-\lambda,
\end{align*} 
where $\mathbb{A}=\mathbb{A}_{ijk\ell}e_i\otimes e_j\otimes e_k\otimes e_\ell$ is given by 
\begin{align*}
\mathbb{A}_{ijk\ell}&=\frac{-3(1-\nu)}{2(1+\nu)(7-5\nu)}\Bigl[\frac{(1-14\nu+15\nu^2)E}{(1-2\nu)^2}\delta_{ij}\delta_{k\ell}
+5E(\delta_{ik}\delta_{j\ell}+\delta_{i\ell}\delta_{jk})\Bigl].
\end{align*}
As in Remark \ref{R:derivation}, 
we replace $\overline{F}'$ with $
- \d_{\rm T}\overline{F},
$
which implies that \eqref{rde} coincides with \eqref{ee2}. 
Here we put $\mathscr{F}=\overline{F}$.
Then we set $\phi\in L^{\infty}(0,+\infty;H^1_0(D))$ as a unique weak solution to 
\begin{align}\label{NLDW}
\begin{cases}
\partial_{t}^2\phi+(3/t)\partial_t \phi-\tau \Delta \phi=\rho \d_{\rm T}\overline{F}(\phi)  \text{ in } D\times (0,+\infty),\\
\phi\lvert_{\partial D}=0,\ \phi\lvert_{t=0}=\phi_{0},\ \partial_t\phi\lvert_{t=0}=\phi_{1},
\end{cases}
\end{align}
where $\phi_0\in H^2(D)\cap H^1_0(D)$ and $\phi_1\in H^1_0(D)$.

\begin{remark}[Well-posedness and boundary conditions] \label{R:reg2}
\rm 
The characteristic function $\chi_\phi\in L^{\infty}(D;\{0,1\})$ is replaced by an approximated Lipchitz continuous function in terms of numerical analysis (see \ref{Ap} below). Thus standard general theories 
ensure well-posedness for \eqref{NLDW} (see, e.g.,~\cite{CH} for details). Here, the homogeneous Dirichlet boundary condition for \eqref{NLDW} is imposed, but only for the uniqueness of solutions. Thus, other boundary conditions can also be allowed.
\end{remark}

\begin{remark} \label{R:reg}
\rm
Since $(3/t)\to 0$ as $t\to +\infty$, the damping term $(3/t)\partial_t \phi$ 
may be ignored in \eqref{NLDW} for simplicity.
Indeed, in order to construct $\phi_0\in H^2(D)\cap H^1_0(D)$ and $\phi_1\in H^1_0(D)$, let $\psi\in L^{\infty}(s,+\infty;H^1_0(D))$ be a unique weak solution to
\begin{align*}
\begin{cases}
\partial_{t}^2\psi+(3/t)\partial_t\psi-\tau \Delta \psi=\rho \d_{\rm T}\overline{F}(\psi) \ \text{ in } D\times (s,+\infty),\\
\psi\lvert_{\partial D}=0,\quad \psi\lvert_{t=s}=\varphi\lvert_{t=s},\quad \partial_t\psi\lvert_{t=s}= 
\partial_t\varphi\lvert_{t=s}
\end{cases}
\end{align*}
for some $s\in (0,+\infty)$ and $\varDelta t>0$. Here, $\varphi\in L^2((0,s];H^1_0(D))$ is a unique weak solution to 
\begin{align*}
\begin{cases}
\partial_{t}\varphi-\tau \Delta \varphi= \rho \d_{\rm T}\overline{F}(\varphi) \ \text{ in } D\times (0,s],\\
\varphi\lvert_{\partial D}=0,\quad \varphi\lvert_{t=0}=\varphi_0\in L^{\infty}(D). 
\end{cases}
\end{align*}
Since the reaction term is numerically treated as a Lipchitz continuous function, there exists $T_{\rm max}>0$ such that $\varphi\in C(0,T_{\rm max};H^2(D)\cap H^1_0(D))$. 
In particular, by \eqref{levelset2}, one can choose $T_{\rm max} \in (0,+\infty)$ as a large number, and so is $s\in (0,+\infty)$; in other words, $(3/t)\partial_t \psi$ is small enough to be negligible. 
 
Furthermore, in terms of the regularity of solutions,
the control of geometric complexity will be expected, and hence, the optimized configuration $\Omega_{\rm opt}\in\mathscr{U}_{\rm ad}$ will also be as smooth as that for the reported method \cite{YINT10}.
\end{remark}

\section{Numerical algorithm}\label{S:algo} 
In this section, we shall describe a numerical algorithm to solve the optimization problem for \eqref{min-prob} by updating the level set function (see 
\ref{Ap} below for technical details).

\vspace{2mm}
\noindent
{\bf Step\,1.}
Set the fixed design domain $D\subset \R^d$, boundary conditions for \eqref{gov1} and 
the initial level set function ($\emph{ophi}:=\phi_0$ in the code). 

\vspace{2mm}
\noindent
{\bf Step\,2.}
Determine the state value $\boldsymbol{u}_\phi\in V$. To this end, discretizing $D\subset \R^d$ with finite elements, we solve \eqref{gov1} using the finite element method.

\vspace{2mm}
\noindent
{\bf Step\,3.}
Compute the functionals $F(\phi)$ and $G(\phi)$ (\emph{obj} and \emph{Gv} in the cade, respectively). 
Here, we note that $G(\phi)$ is normalized in the code.

\vspace{2mm}
\noindent
{\bf Step\,4.}
Check for convergence. In the code,  based on the gradient descent method,
 we define convergence conditions as follows\/{\rm :}
\begin{align}\label{condi-conv}
\| \emph{LsfDiff}\,\|_{L^{\infty}(D)}
<\emph{eps0pt},\quad \emph{Gv}\le 0.
\end{align}
Here, $\emph{eps0pt}\in\R$ is the criterion for convergence (see Step\,7 below for \emph{LsfDiff}). 
If the conditions in \eqref{condi-conv} are all satisfied, then we terminate the optimization. 
Otherwise, we proceed to the next step. 

\vspace{2mm}
\noindent
{\bf Step\,5.}
Compute the topological derivative $\d_{\rm T} F$ and the Lagrange multiplier $\lambda\ge 0$ ($\emph{Td1}(\cdot,\cdot)$ and \emph{LagGV} in the code, respectively). In particular, we set the following normalizer for dimensionless in the code\/{\rm:}
$$
\emph{AbsTd1}:=\frac{\int_{D} |\d_{\rm T} F|\, \d x}{|D|}.
$$
On the other hand, as for the Lagrange multiplier $\emph{LagGV}$, we employ \emph{augmented Lagrangian's method} as follows\/{\rm :}
$$
\emph{LagGV}=\emph{LagGVp}+
\emph{LagGVD},
$$
where $\emph{LagGVp}$ is the previous version of $\emph{LagGV}$ and $\emph{LagGVD}$ is some normalized volume functional (see 
\ref{Ap} below for details). 

\vspace{2mm}
\noindent
{\bf Step\,6.}
Solve PDEs. In the code, 
let $\emph{Iter}\in\mathbb{Z}$ be an iteration number.
Choosing $\emph{StatIt}\in\mathbb{Z}$ to which \eqref{NLDW} applies, we solve the following either (i) or (ii) using the finite difference method discretized in the time direction\/{\rm :} 
\begin{itemize}
\item[(i)]
In case  $\emph{Iter}\le \emph{StatIt}$, for all $\tilde{\Phi}\in \tilde{V}:=\{ \tilde{\Phi}\in H^1(D)\colon \tilde{\Phi}=0 \text{ on } \Gamma_D\}$, 
\begin{align*}
&0=\int_D \left(\frac{\emph{phi}}{\emph{dt}}\right)\tilde{\Phi}(x)\, \d x
 +\int_D \emph{tau} \nabla (\emph{phi})\cdot \nabla\tilde{\Phi}(x)\, \d x\\
&\quad-
\int_D 
\emph{CdF}\left(\frac{\emph{Td1}(u,\emph{ophi})}{\emph{AbsTd1}}-{\emph{LagGv}}\right) \tilde{\Phi}(x)\, \d x
 -
\int_D \left(\frac{\emph{ophi}}{{\emph{dt}}}\right) \tilde{\Phi}(x)\, \d x \nonumber
\end{align*}
Here $\emph{dt}:=\varDelta t>0$, $\emph{tau}:=\tau>0$ and $\emph{CdF}:=\rho>0$ are given parameters.

\item[(ii)]
In case $\emph{Iter} >  \emph{StatIt}$, for all $\tilde{\Phi}\in \tilde{V}$,  
\begin{align*}
&0=\int_D \left(\frac{\emph{phi}}{\emph{dt}}\right)\tilde{\Phi}(x)\, \d x
+\int_D \emph{tau} \nabla (\emph{phi})\cdot \nabla\tilde{\Phi}(x)\, \d x\\
&-
\int_D 
\emph{CdF}\left(\frac{\emph{Td1}(u,\emph{ophi})}{\emph{AbsTd1}}-\emph{LagGv}\right) \tilde{\Phi}(x)\, \d x
 -
\int_D\left(\frac{2\emph{ophi}-\emph{oophi}}{\emph{dt}}\right)\tilde{\Phi}(x)\, \d x.
\nonumber
\end{align*}
Here we used the fact in Remark \ref{R:reg} for simplicity since the damping term eventually becomes negligibly small. 
\end{itemize}

\vspace{2mm}
\noindent
{\bf Step\,7.}
Normalize \emph{phi} as follows\/{\rm :}
$$
\emph{phi}=
\begin{cases}
\text{sgn}(\emph{phi})&\text{ if } |\emph{phi}|>1,\\
\emph{phi}&\text{ otherwise}
\end{cases}
$$
and set \emph{LsfDiff}$:=$\emph{phi}-\emph{ophi}. Return to {Step\,2} after setting the next initial level set functions as \emph{ophi}=\emph{phi} and \emph{oophi}=\emph{ophi}.

\section{Main results}\label{S:reslut}

In this section, we shall describe numerical examples for the two-dimensional case mainly and numerically show that the method based on \eqref{NLDW} and \eqref{ee2} converges to an optimized configuration faster than the method based on only \eqref{ee2}. 

Let $D\subset \R^2 $ be a rectangle and set
Young's modulus $\bm{E}>0$, Poisson's ratio ${\nu}>0$, 
and the Lame coefficients ${\tilde\lambda}>0$ and $\mu>0$ as follows{\rm :}
$$
E=2.1\times 10^{11}, \quad \nu=0.3, \quad 
\tilde{\lambda}=\frac{E\nu}{(1+\nu)(1-2\nu)},\quad 
\mu=\frac{E}{2(1+\nu)}. 
$$
In order to solve \eqref{gov1} in Step\,2,
the elasticity tensor and traction vector replaced by $D$ and $g$ according to the code, respectively, are set as follows\/{\rm :}
$$
D=
\begin{pmatrix}
\tilde\lambda+2\mu & \tilde\lambda& 0 \\
\tilde\lambda & \tilde\lambda+2\mu& 0 \\
0&0&\mu
\end{pmatrix},\quad
g=(0,-1.0\times 10^3).
$$
Here, we note that the elasticity tensor is rewritten as the matrix in terms of the finite element method.
Moreover, 
we set the topological derivative $\d_{\rm T}\overline{F}$ as follows\/{\rm :}  
\begin{align*}
\d_{\rm T}\overline{F}
&=
A\chi _{\phi}\boldsymbol{\epsilon}(\boldsymbol{u}_\phi)\cdot \boldsymbol{\epsilon}(\boldsymbol{u}_\phi)-\lambda,
\end{align*} 
where, $\boldsymbol{\epsilon}(\boldsymbol{u}_\phi)=(\partial_{x_1}u_1, \partial_{x_2}u_2, \partial_{x_2}u_1+ \partial_{x_1}u_2)$,
$\boldsymbol{u}_\phi=(u_1,u_2)\in V$ and 
$$
A=
\begin{pmatrix}
A_1+2A_2 &  A_1&0\\
A_1 & A_1+2A_2&0\\
0&0&A_2
\end{pmatrix}.
$$
Here, $A_1$ and $A_2$ are given by
$$
A_1=-\frac{3(1-\nu)(1-14\nu+15\nu^2)}{2(1+\nu)(7-5\nu)(1-2\nu)^2}E
\quad \text{ and }\quad
A_2=\frac{15(1-\nu)}{2(1+\nu)(7-5\nu)}E,
$$
respectively.
We consider the three  models 
 (see Figure~\ref{fig:initial} below).  
For simplicity, the method with  \eqref{ee2} and the method with \eqref{NLDW} (i.e.,~nonlinear hyperbolic-parabolic equations) are described as (RD) and (NLHP), respectively.

In this paper, we define (NLHP) converging faster than (RD), if (NLHP) has fewer iteration numbers that satisfy all convergence conditions than (RD) (see Step\,4 in \S \ref{S:algo}); indeed, the convergence condition for the level set functions mentioned in the previous section is standard in the gradient descent method, and moreover, if 
$\|\phi_{n+1}-\phi_{n}\|_{L^{\infty}(D)}< \e $ for $\e>0$ small enough, then
the configurations $\Omega_{\phi_{n}}=\{x\in D\colon \chi_{\phi_{n}}(x)=1\}$ and $\Omega_{\phi_{n+1}}$ can be (almost) identified.

\begin{figure*}[htbp]
    \begin{tabular}{ccc}
      \hspace*{-5mm} 
      \begin{minipage}[t]{0.32\hsize}
        \centering
        \includegraphics[keepaspectratio, scale=0.07]{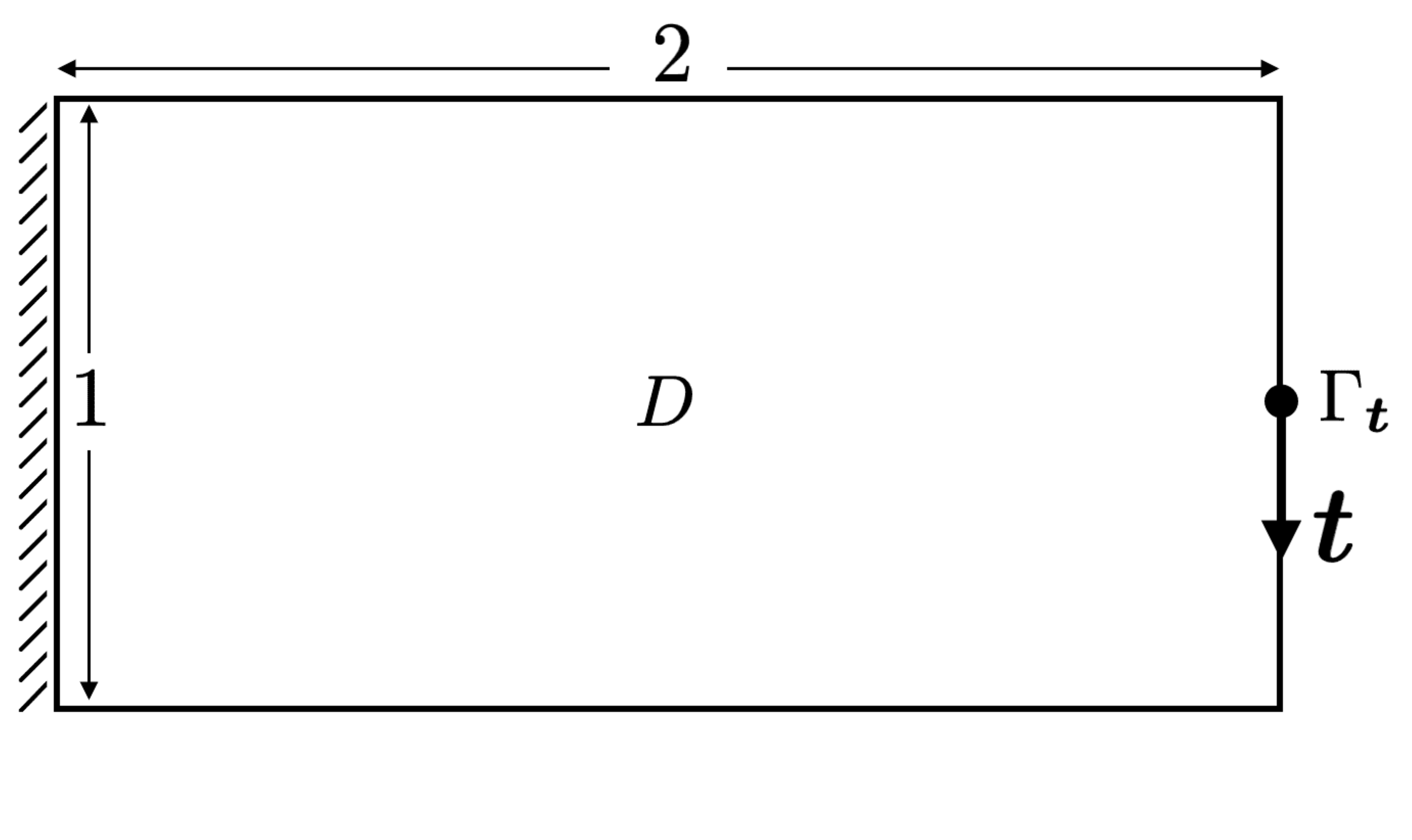}
        \subcaption{Cantilever}
        \label{ica}
      \end{minipage} 
      \quad
      \begin{minipage}[t]{0.32\hsize}
        \centering
        \includegraphics[keepaspectratio, scale=0.07]{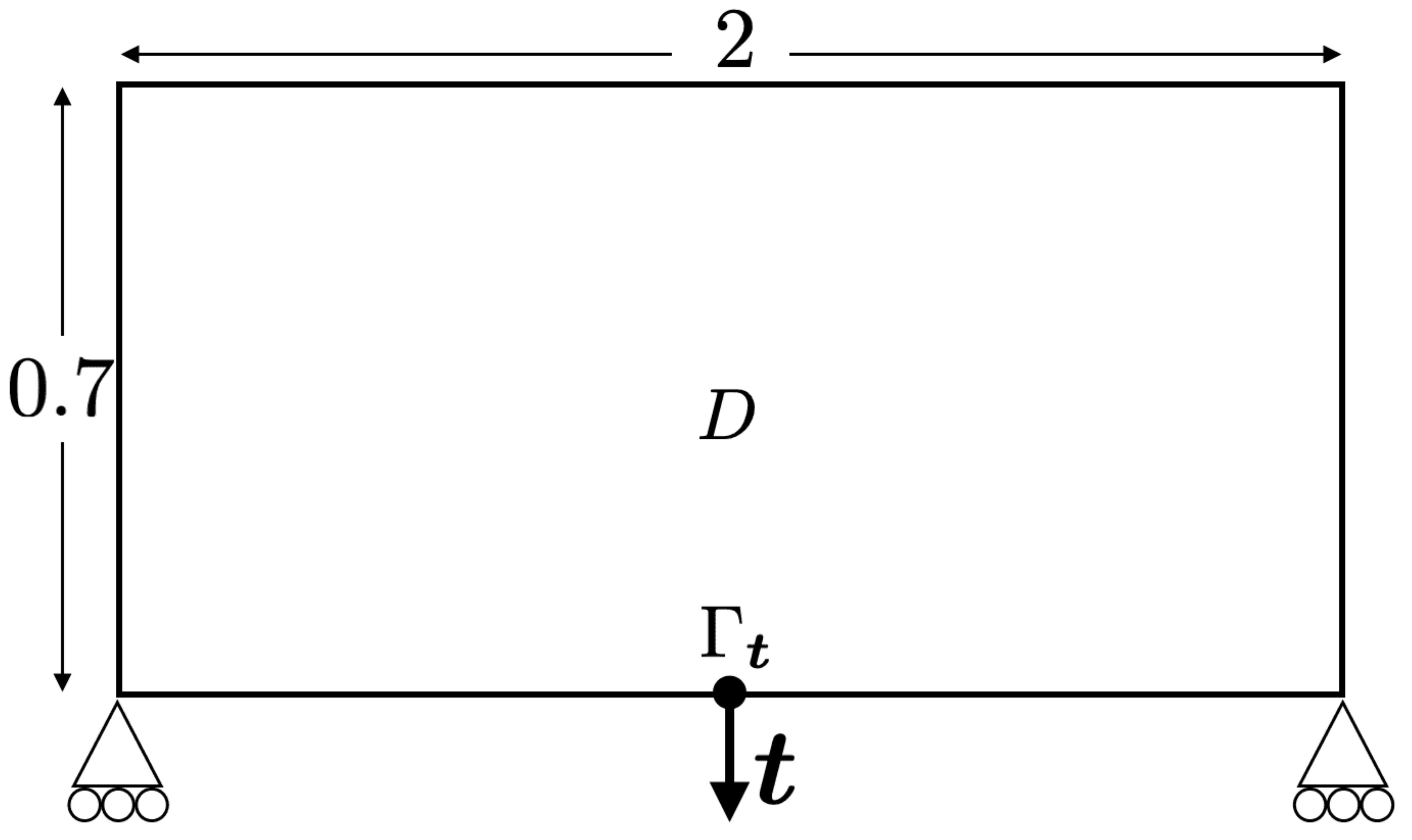}
        \subcaption{Bridge}
        \label{ibr}
      \end{minipage} 
      \quad 
       \begin{minipage}[t]{0.32\hsize}
        \centering
        \includegraphics[keepaspectratio, scale=0.07]{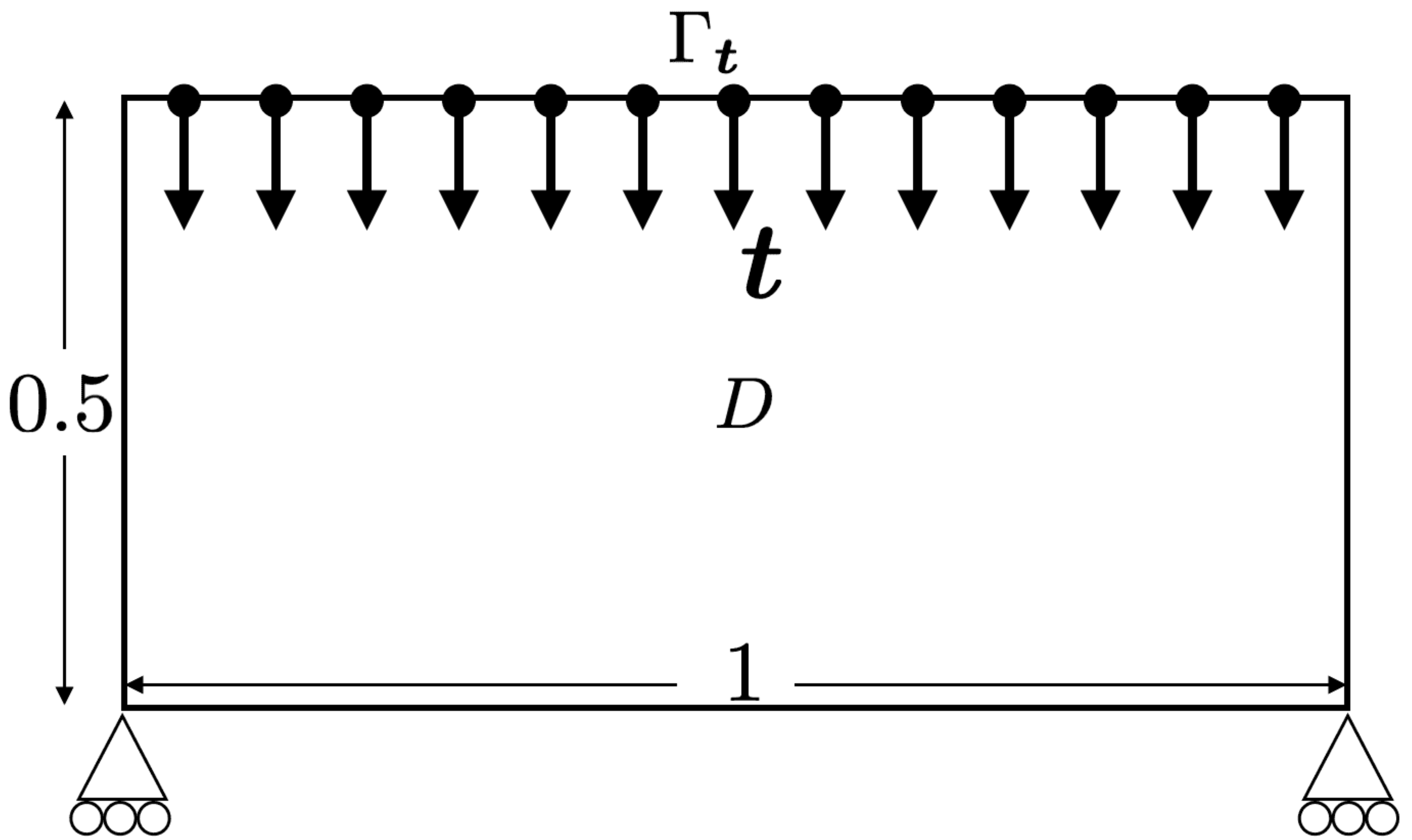}
        \subcaption{Radiator}
        \label{ira}
      \end{minipage} 
      \end{tabular}
       \caption{ Fixed design domain $D$ and boundary conditions.}
    \label{fig:initial}
  \end{figure*}

\subsection{Cantilever model}\label{S:ca}

Based on Figure \ref{ica}, we consider the so-called \emph{cantilever model}.
The given parameters are the same as in the code (see \ref{Ap}).
In particular, we set $(\tau, G_{\rm max})=(5.0\times 10^{-4}, 0.45)$, and 
the number of triangles $n_t$ and  maximum edge size $h_{\rm max}$ are set to 
$(n_t,h_{\rm max})=(38400, 0.0144
)$.
As for the convergence criterion, we choose $\emph{eps0pt}=1.0\times 10^{-3}$ in the code.

\begin{figure*}[htbp]
    \hspace*{-5mm} 
    \begin{tabular}{ccccc}
          \begin{minipage}[t]{0.2\hsize}
        \centering
        \includegraphics[keepaspectratio, scale=0.09]{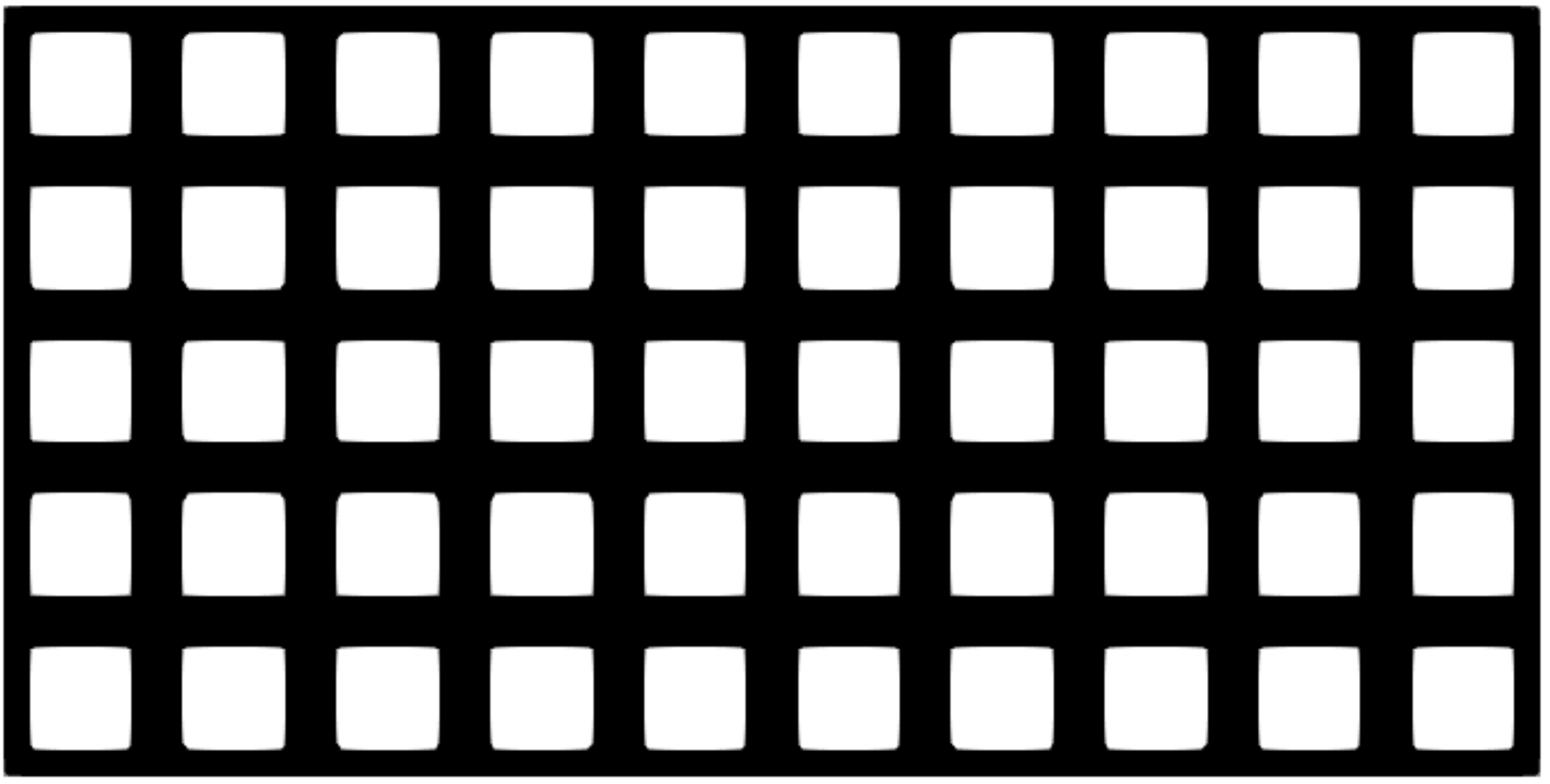}
        \subcaption{Step\,0}
        \label{2-a}
      \end{minipage} 
      \begin{minipage}[t]{0.2\hsize}
        \centering
        \includegraphics[keepaspectratio, scale=0.09]{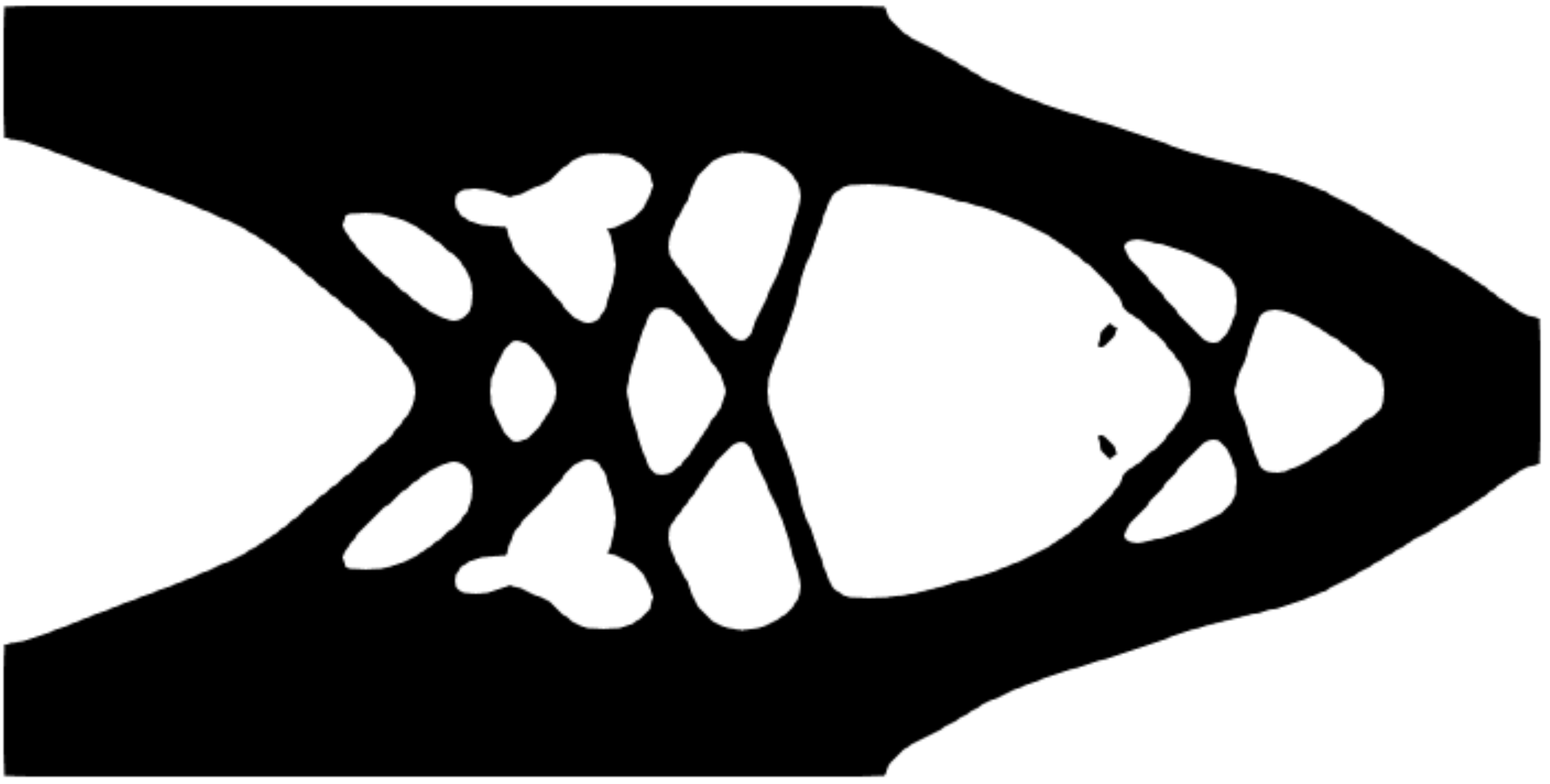}
        \subcaption{Step\,10}
        \label{2-b}
      \end{minipage} 
      \begin{minipage}[t]{0.2\hsize}
        \centering
        \includegraphics[keepaspectratio, scale=0.09]{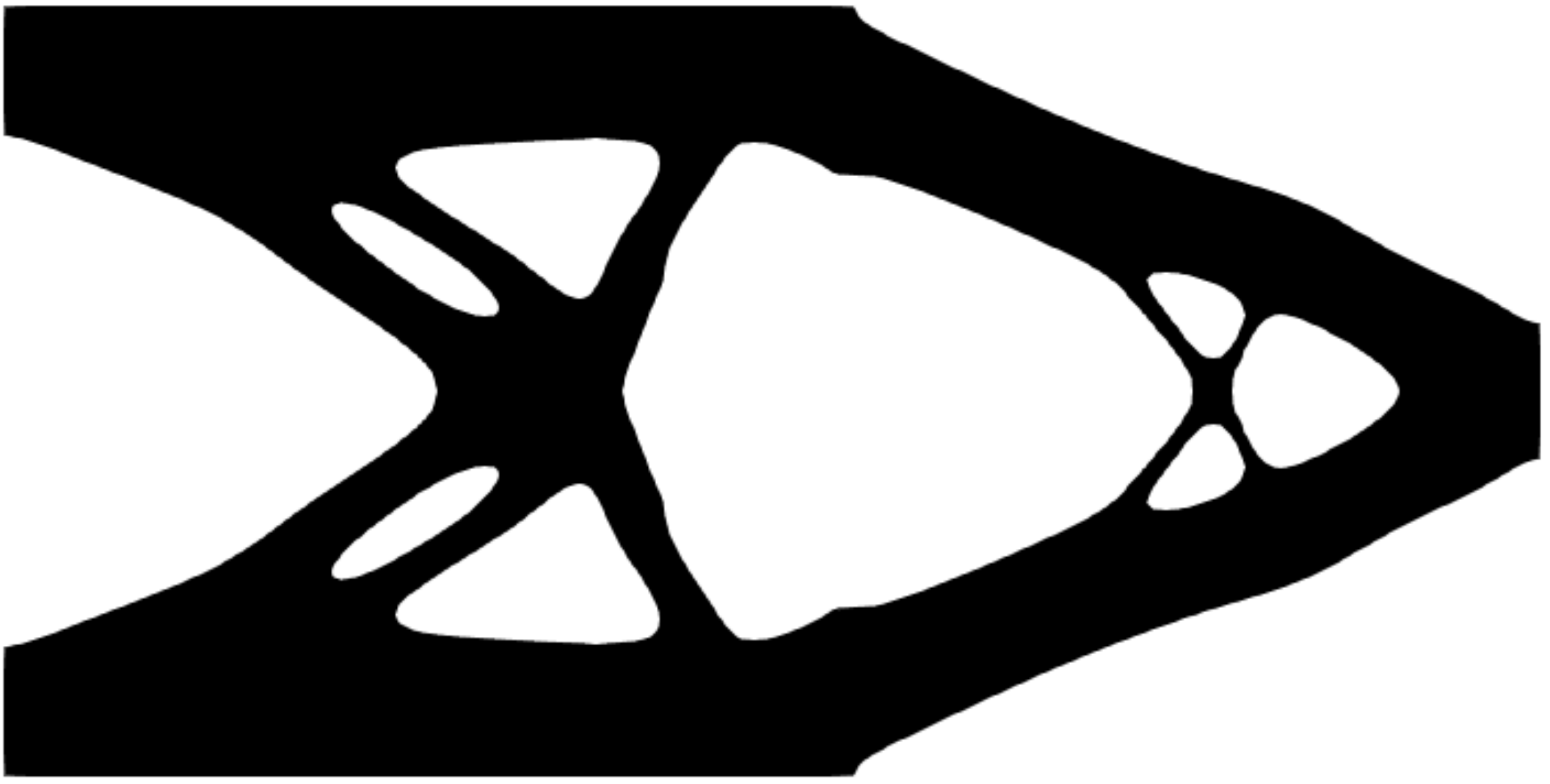}
        \subcaption{Step\,15}
        \label{2-c}
      \end{minipage} 
         \begin{minipage}[t]{0.2\hsize}
        \centering
        \includegraphics[keepaspectratio, scale=0.09]{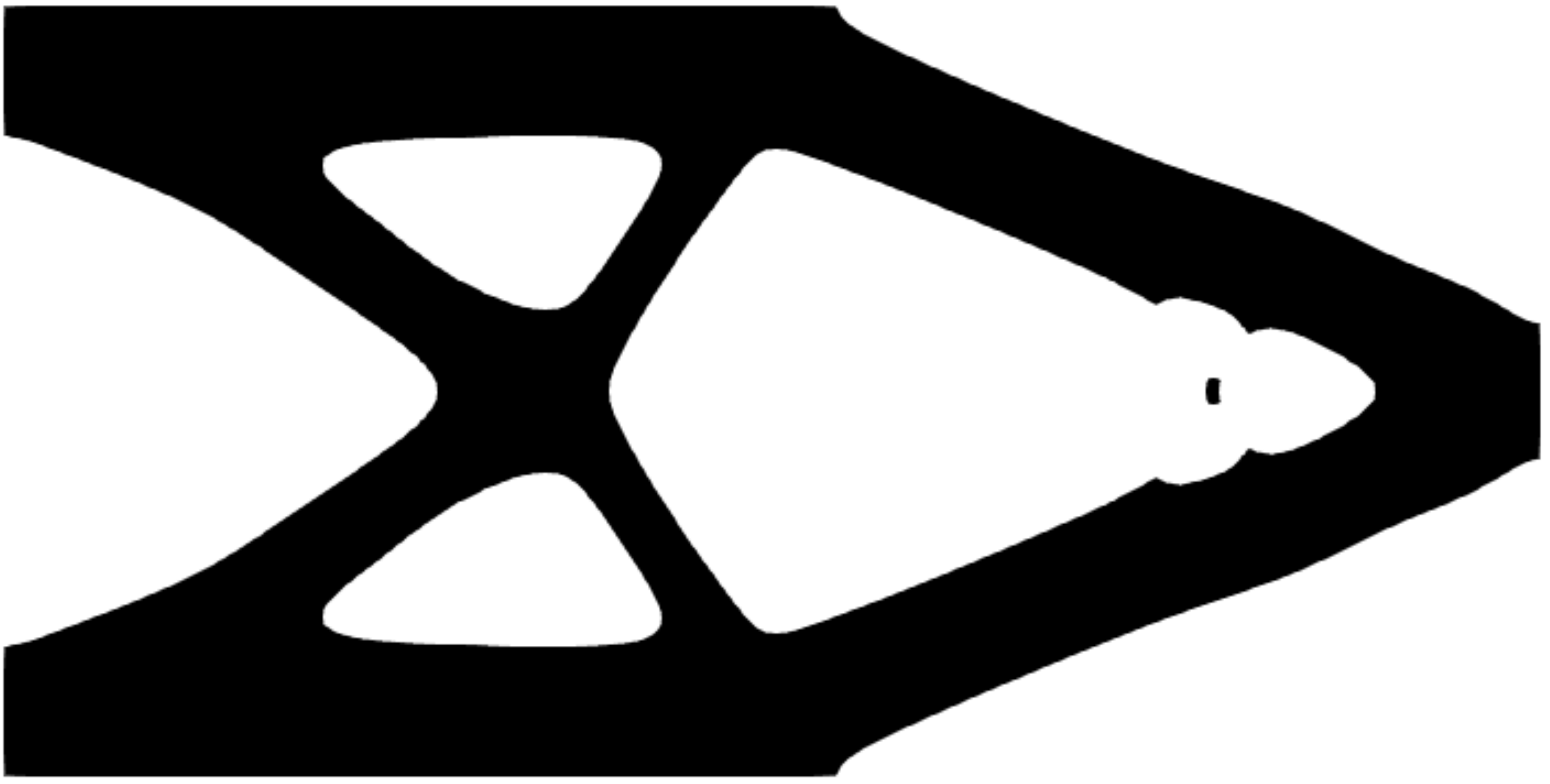}
        \subcaption{Step\,20}
      \end{minipage}
       \begin{minipage}[t]{0.2\hsize}
        \centering
        \includegraphics[keepaspectratio, scale=0.09]{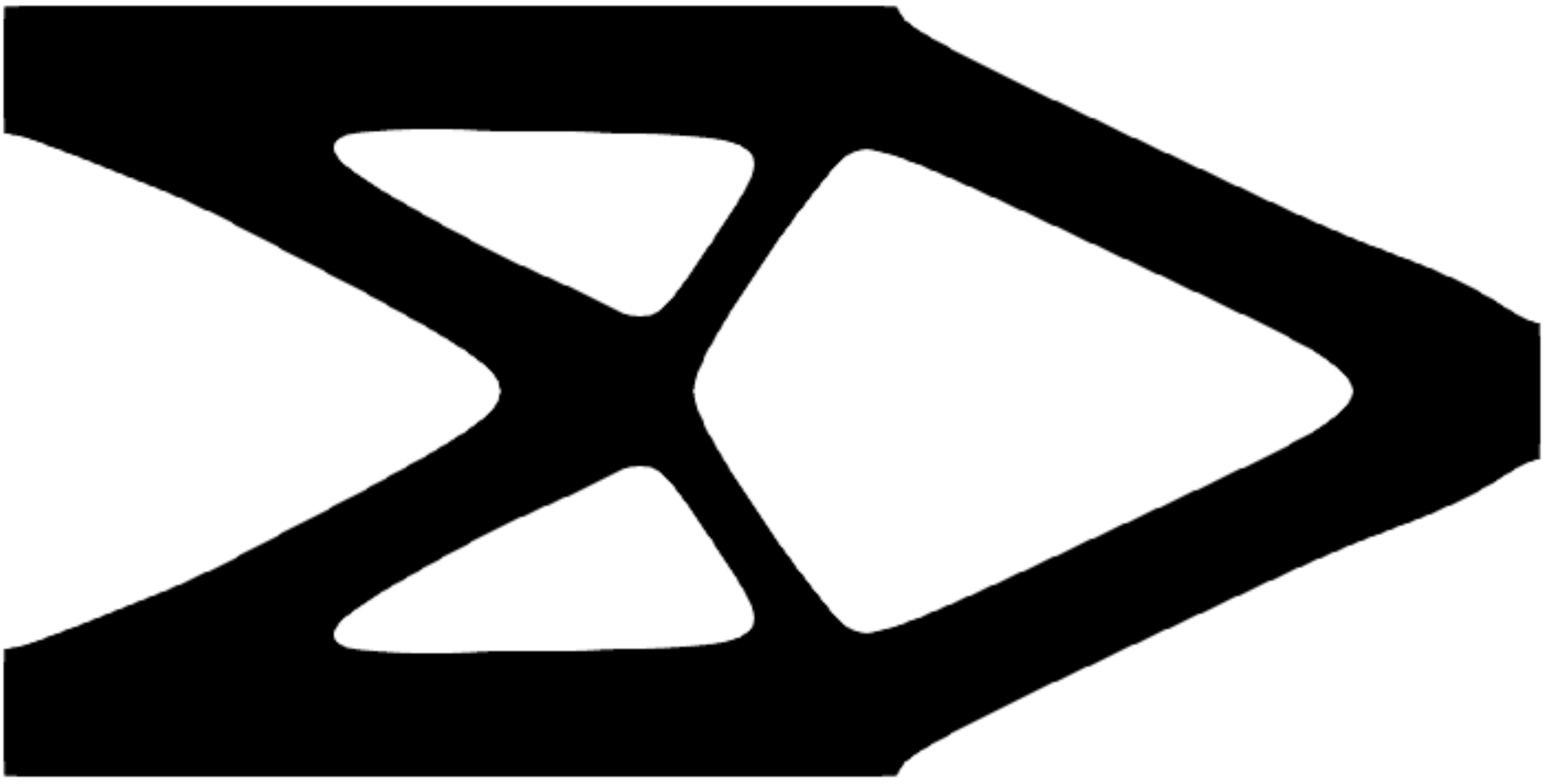}
        \subcaption{Step\,704$^{\#}$}
      \end{minipage}  
      \\ 
    \begin{minipage}[t]{0.2\hsize}
        \centering
        \includegraphics[keepaspectratio, scale=0.09]{cap0.pdf}
        \subcaption{Step\,0}
        \label{2-e}
      \end{minipage} 
      \begin{minipage}[t]{0.2\hsize}
        \centering
        \includegraphics[keepaspectratio, scale=0.09]{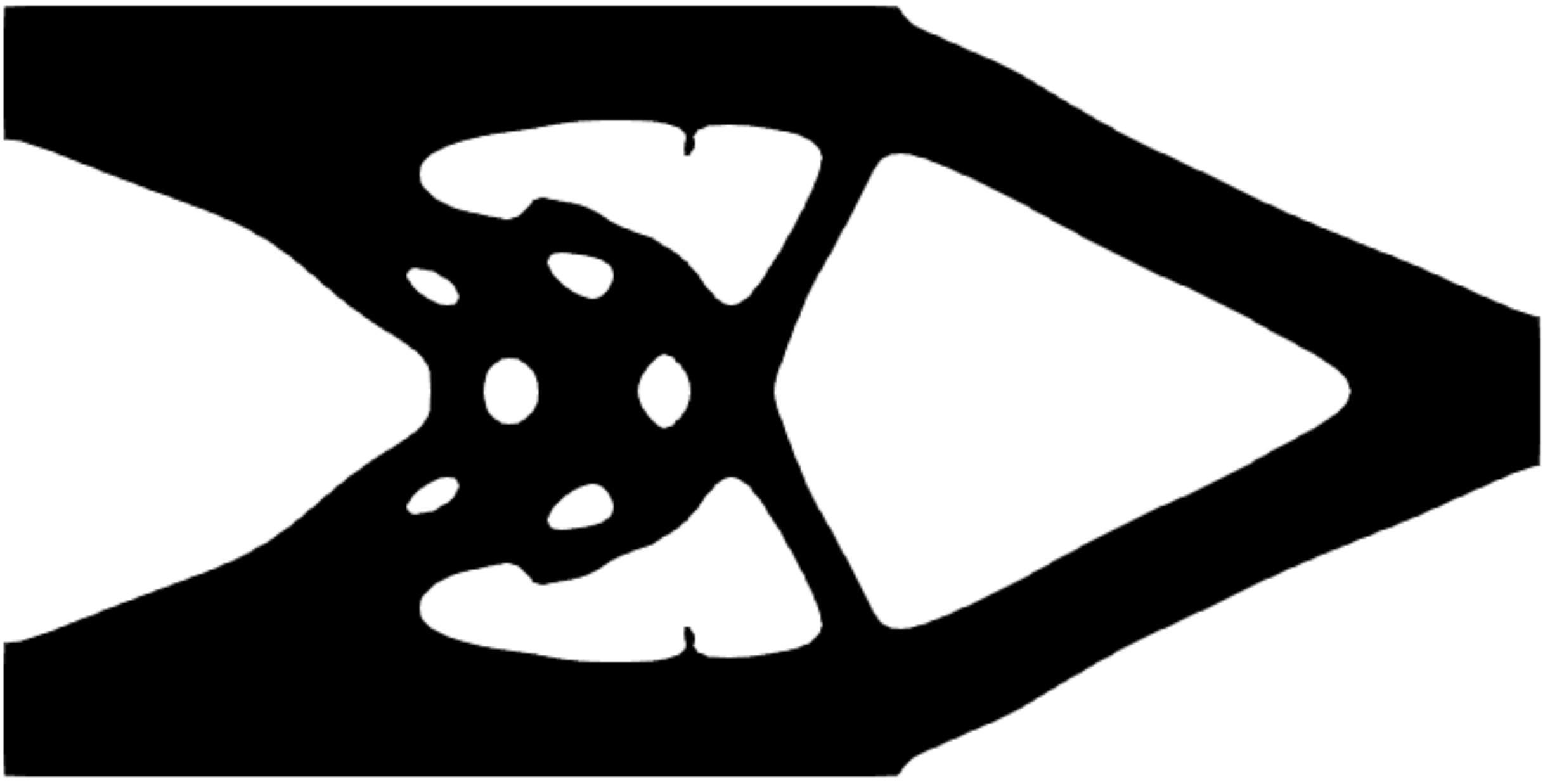}
        \subcaption{Step\,10}
        \label{2-f}
      \end{minipage} 
      \begin{minipage}[t]{0.2\hsize}
        \centering
        \includegraphics[keepaspectratio, scale=0.09]{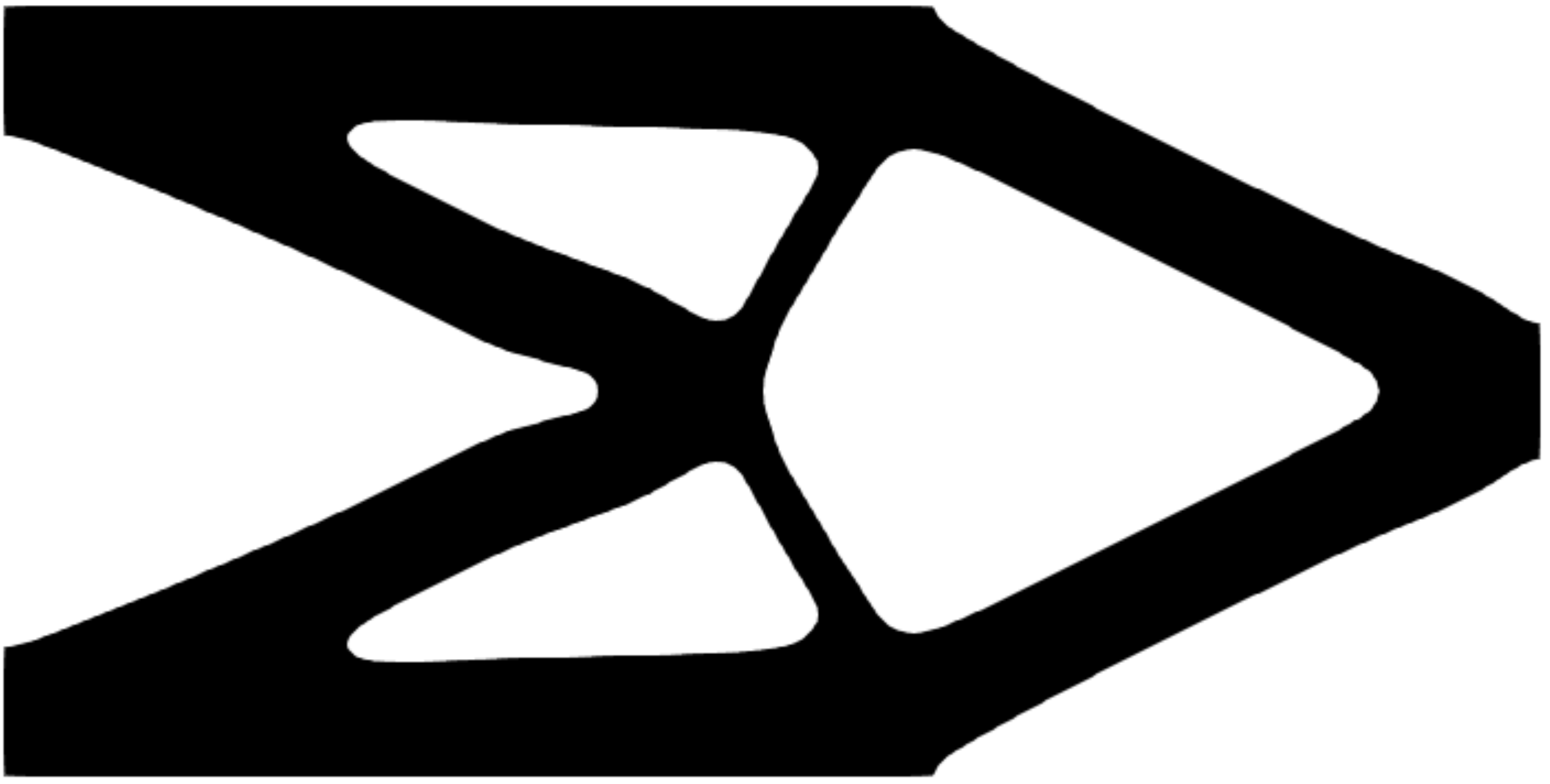}
        \subcaption{Step\,15}
        \label{2-g}
      \end{minipage} 
       \begin{minipage}[t]{0.2\hsize}
        \centering
        \includegraphics[keepaspectratio, scale=0.09]{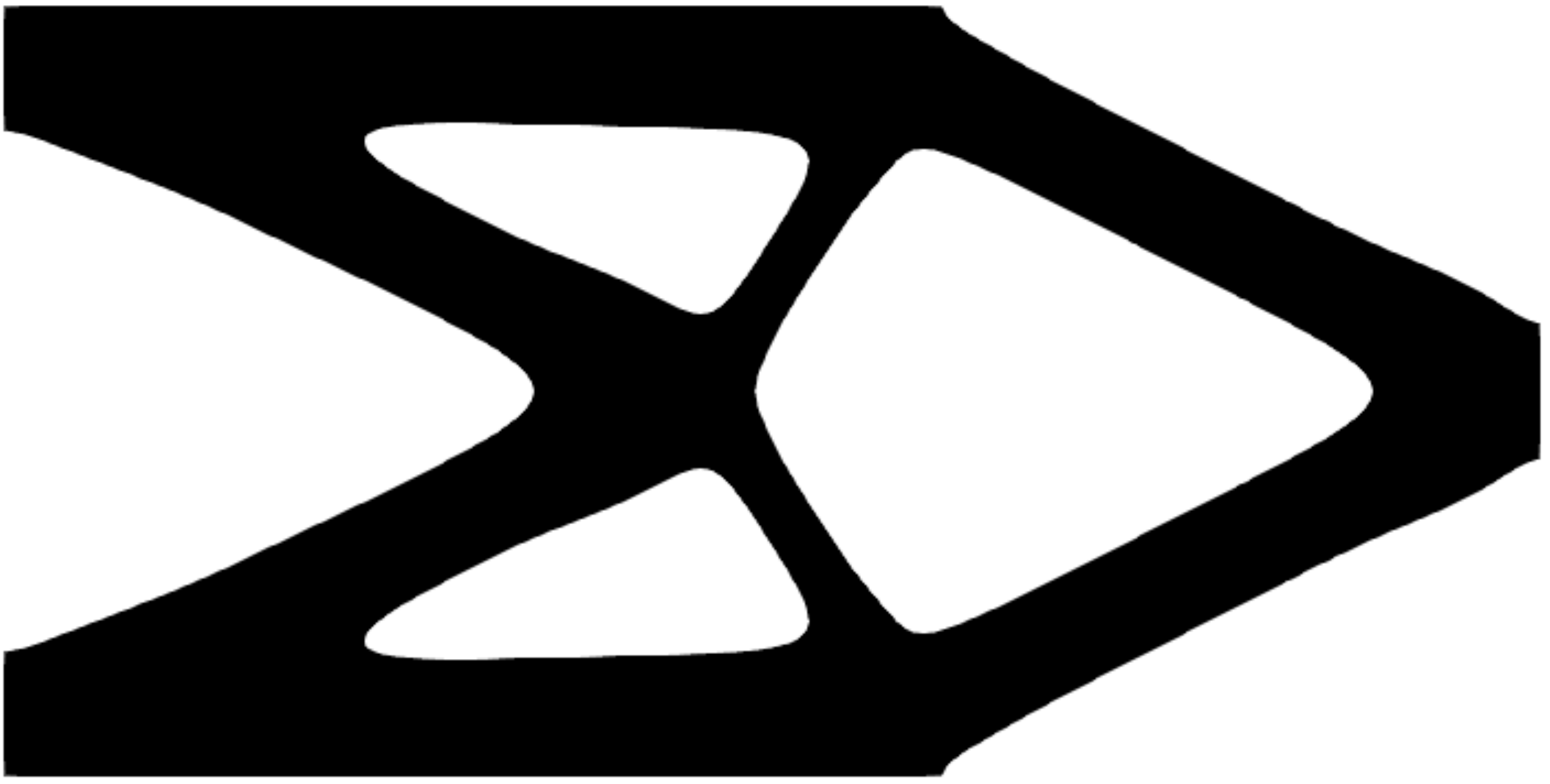}
        \subcaption{Step\,20}
        \label{2-h}
      \end{minipage} 
     \begin{minipage}[t]{0.2\hsize}
        \centering
        \includegraphics[keepaspectratio, scale=0.09]{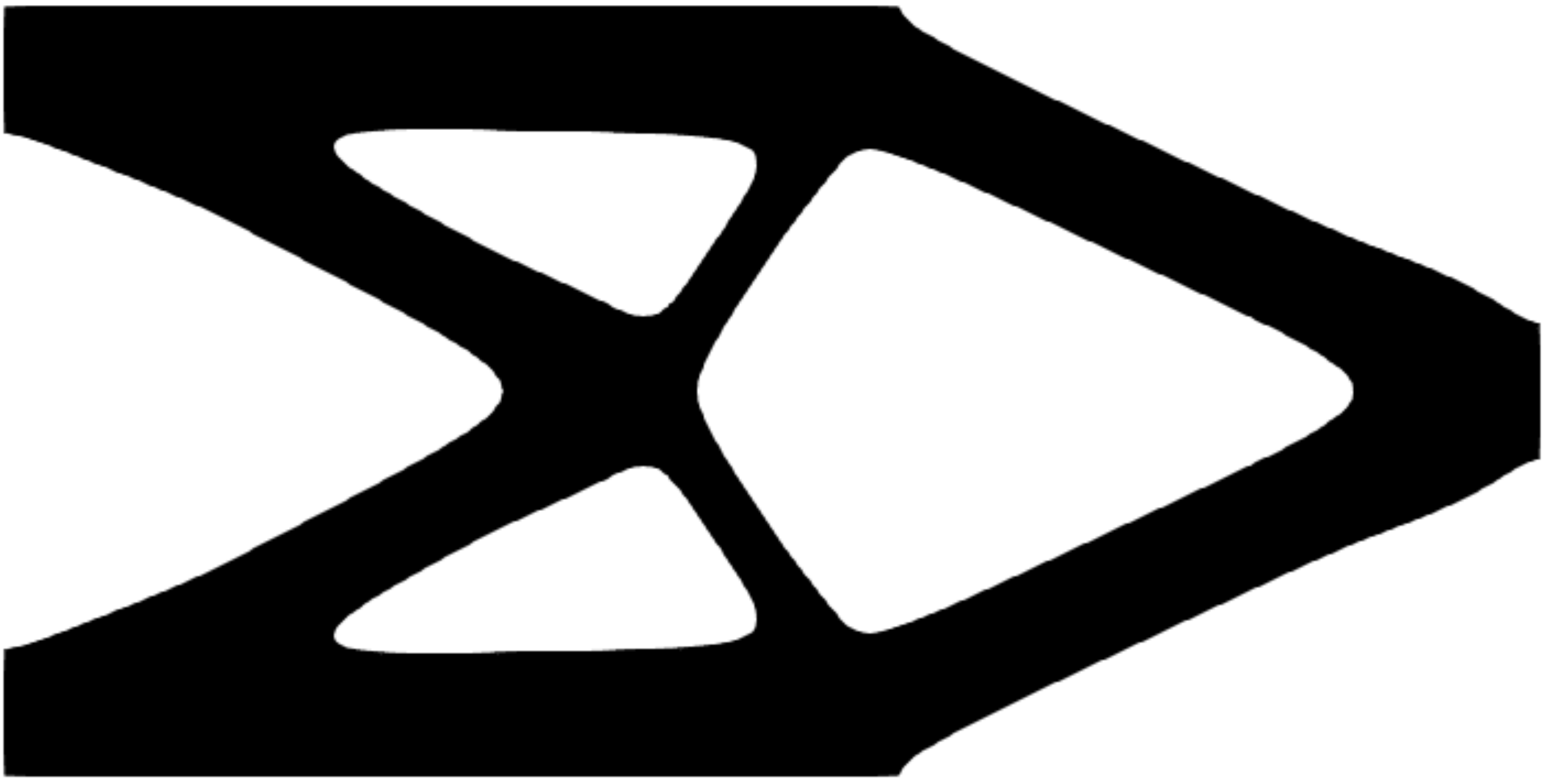}
        \subcaption{Step\,389$^\#$}
        \label{2-h}
      \end{minipage} 
    \end{tabular}
     \caption{Configuration $\Omega_{\phi_n}\subset D$ for the case where the initial configuration is the periodically perforated domain. Figures (a)--(e) and (f)--(j) represent the results of (RD) and (NLHP), respectively.     
The symbol $^{\#}$ implies the final step.}
     \label{fig:ppd}
  \end{figure*}

\begin{figure*}[htbp]
    \begin{tabular}{ccc}
      \hspace*{-5mm} 
      \begin{minipage}[t]{0.49\hsize}
        \centering
        \includegraphics[keepaspectratio, scale=0.33]{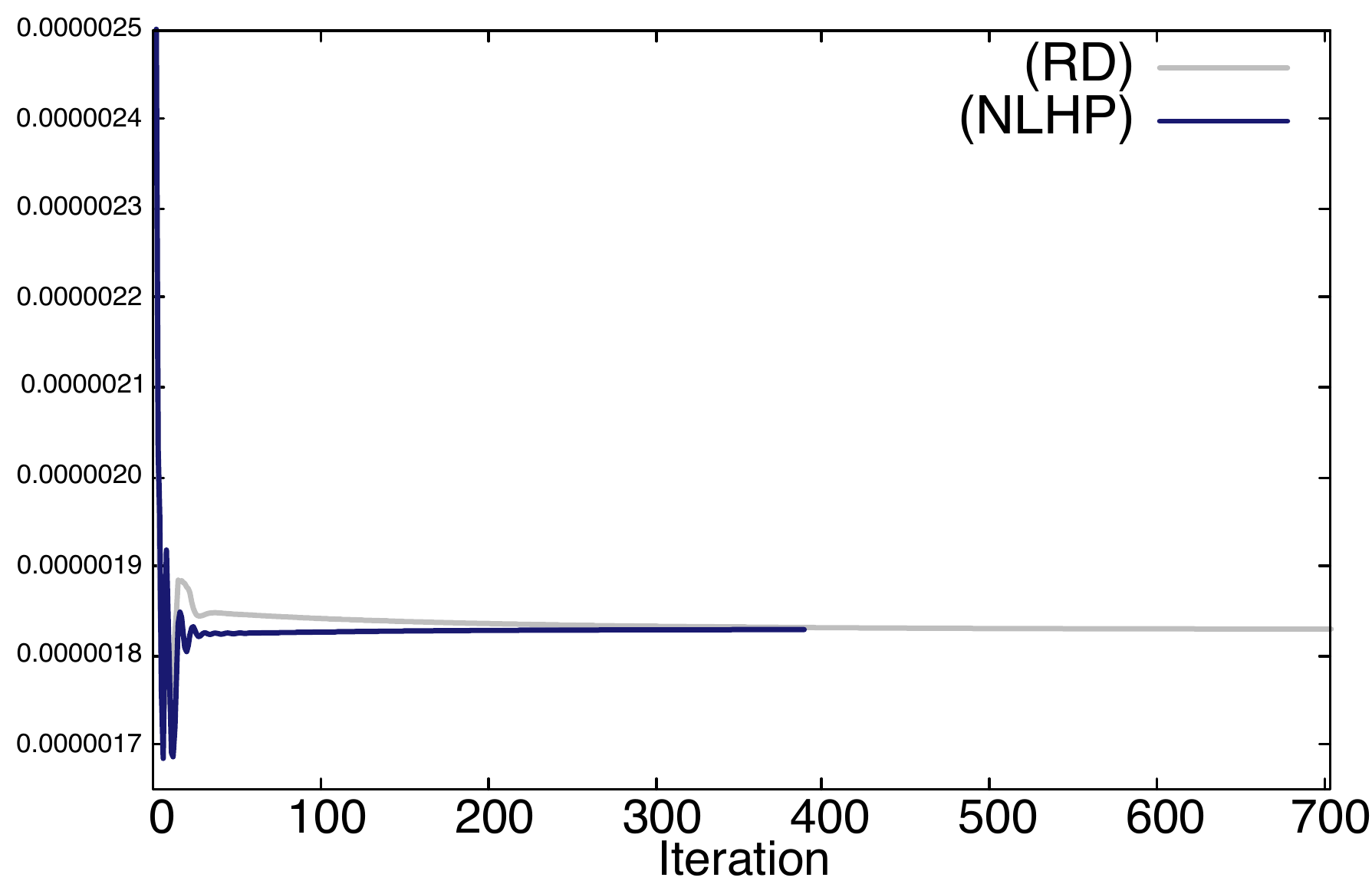}
        \subcaption{$F(\phi_n)$}
        \label{ca1-a}
      \end{minipage} 
      \begin{minipage}[t]{0.49\hsize}
        \centering
        \includegraphics[keepaspectratio, scale=0.33]{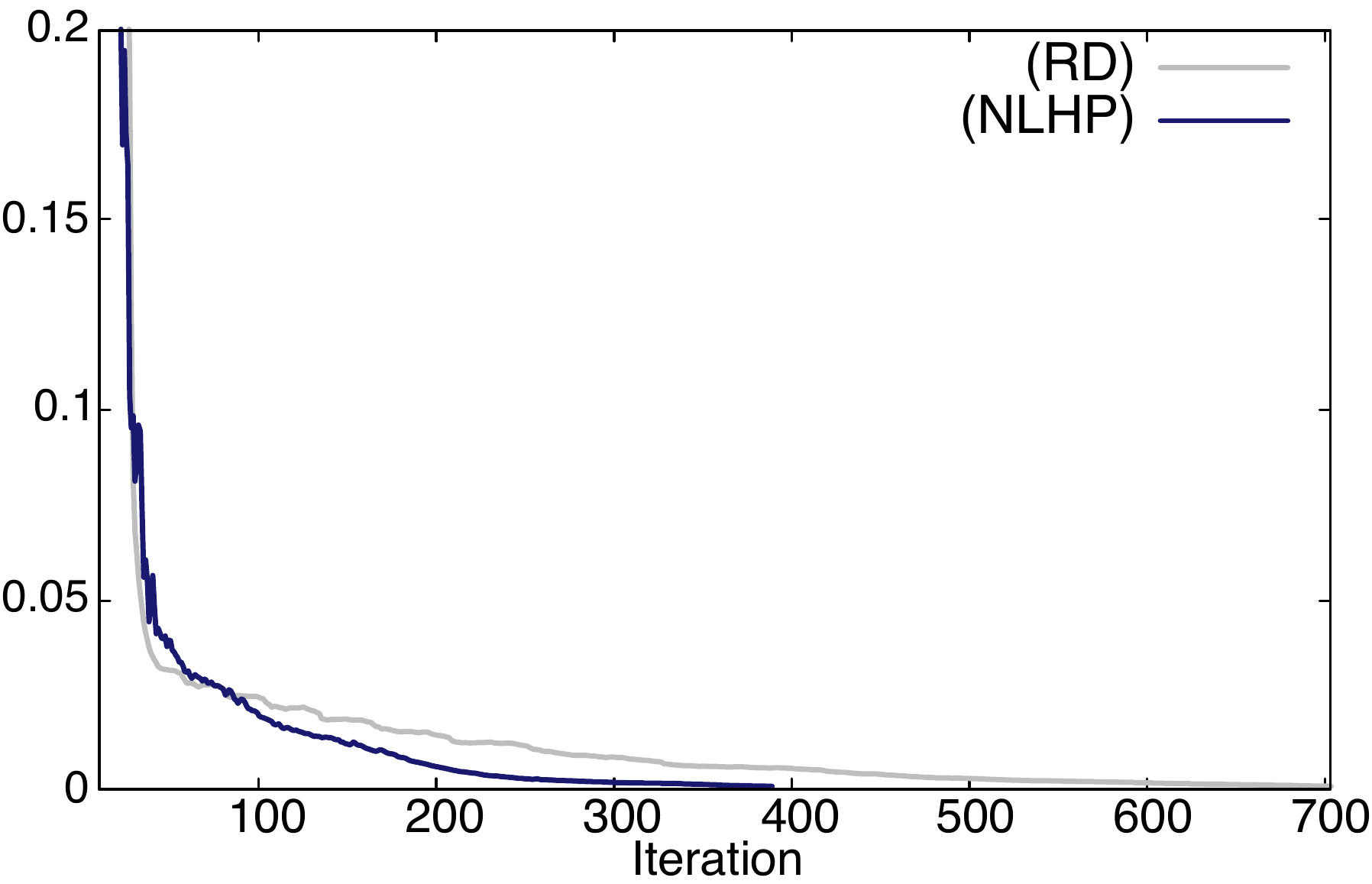}
        \subcaption{$\|\phi_{n+1}-\phi_n\|_{L^{\infty}(D)}$}
        \label{ca1-b}
      \end{minipage} &
      \end{tabular}
       \caption{ Objective functional and convergence condition for \S \ref{S:ca}-(i).}
    \label{fig:cag1}
  \end{figure*}

\vspace{2mm}
\noindent
{\bf Case (i) (Periodically perforated domain).\,} We first consider the case where 
the initial configuration is a periodically perforated domain.
Then Figures \ref{fig:ppd} and \ref{fig:cag1} are obtained, and
one can confirm that (NLHP) satisfies the convergence condition in Figure \ref{ca1-b} faster than (RD).
In particular, it is noteworthy that (NLHP) optimizes the topology in only 15 steps (see Figure \ref{2-g}).

\begin{figure*}[htbp]
   \hspace*{-5mm} 
    \begin{tabular}{ccccc}
      \begin{minipage}[t]{0.2\hsize}
        \centering
        \includegraphics[keepaspectratio, scale=0.09]{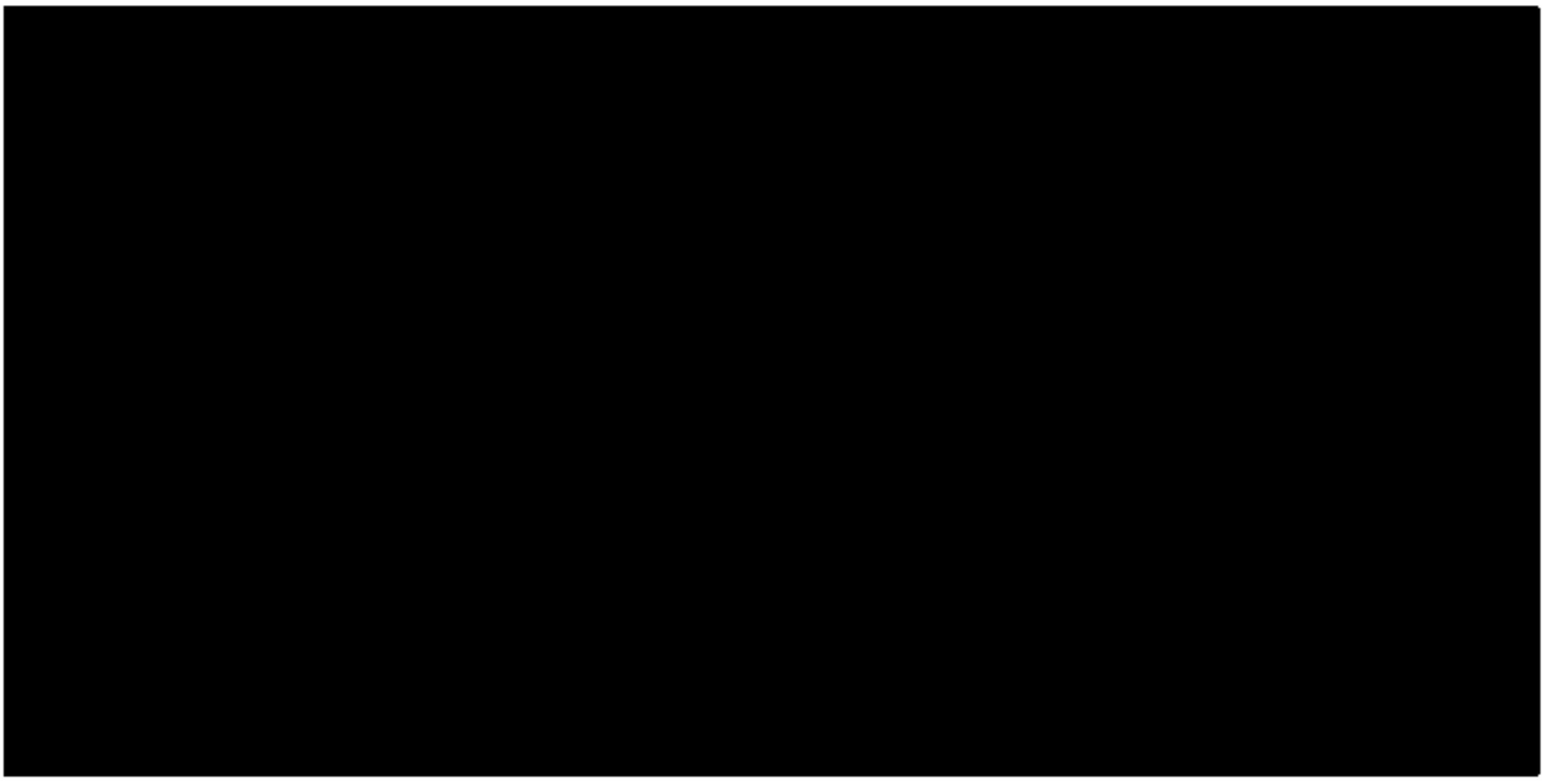}
        \subcaption{Step\,0}
        \label{1-a}
      \end{minipage} 
      \begin{minipage}[t]{0.2\hsize}
        \centering
        \includegraphics[keepaspectratio, scale=0.09]{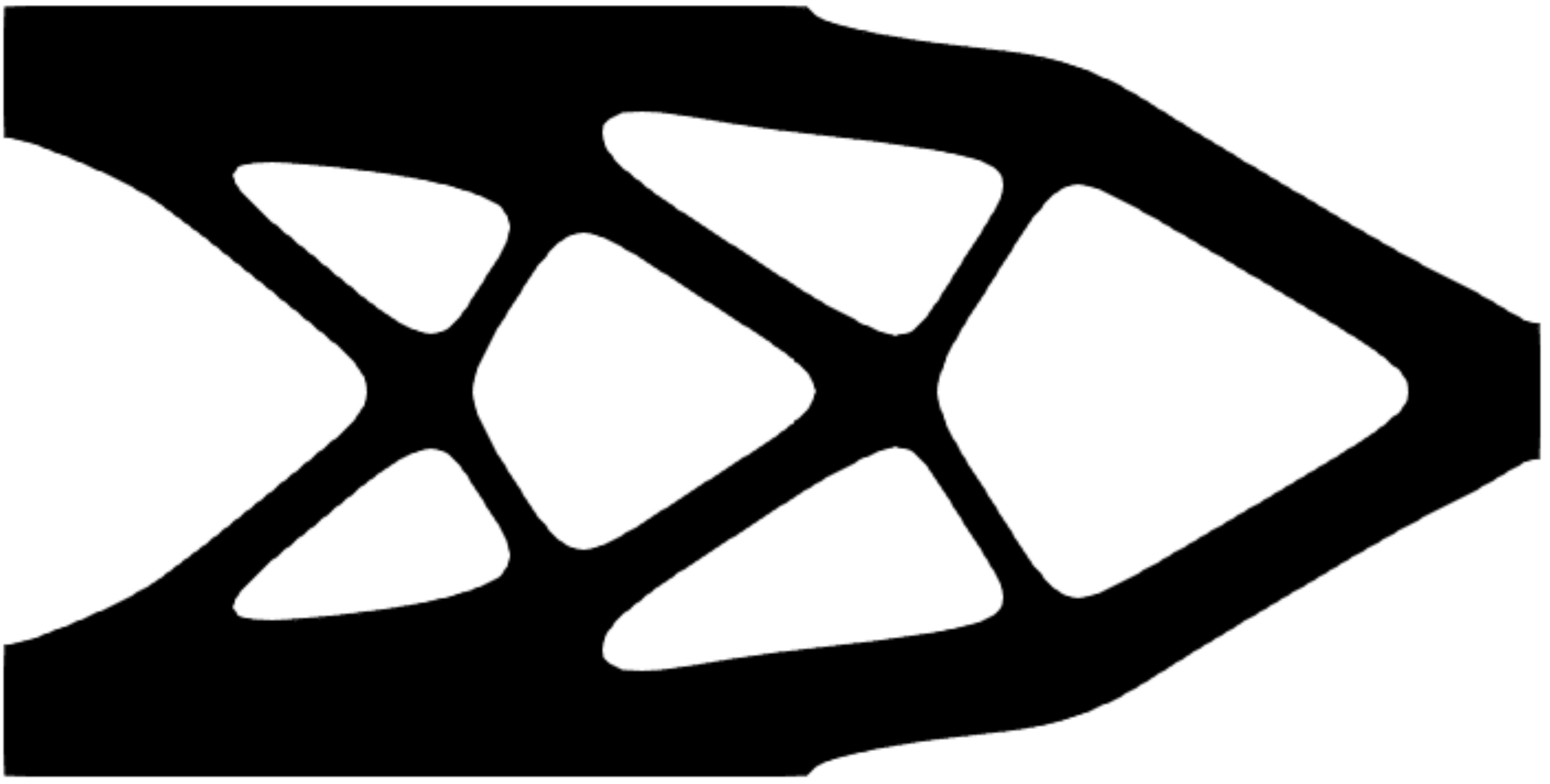}
        \subcaption{Step\,50}
        \label{1-b}
      \end{minipage} 
      \begin{minipage}[t]{0.2\hsize}
        \centering
        \includegraphics[keepaspectratio, scale=0.09]{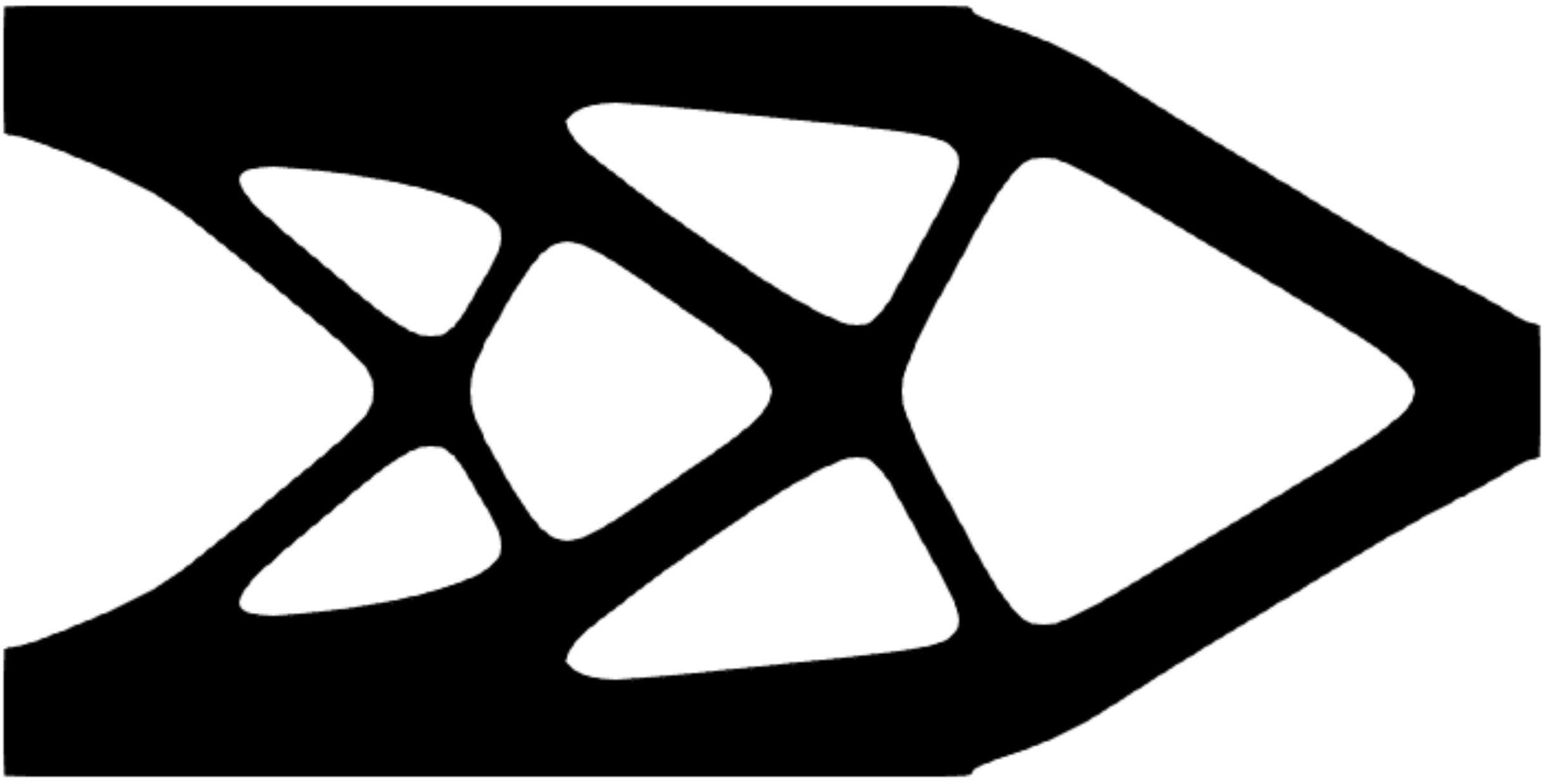}
        \subcaption{Step\,150}
        \label{1-c}
      \end{minipage} 
         \begin{minipage}[t]{0.2\hsize}
        \centering
        \includegraphics[keepaspectratio, scale=0.09]{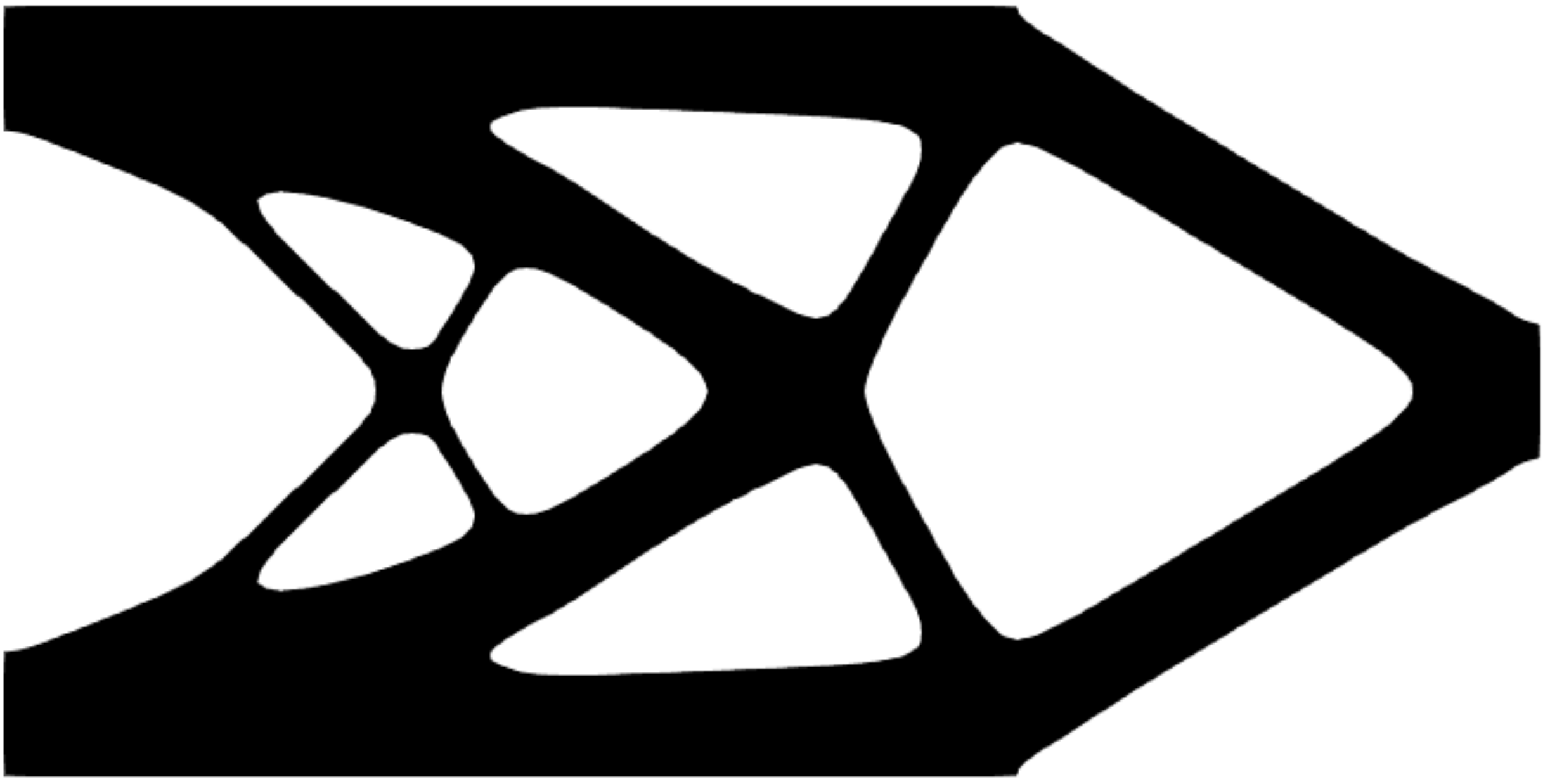}
        \subcaption{Step\,250}
        \label{1-d}
      \end{minipage}
           \begin{minipage}[t]{0.2\hsize}
        \centering
        \includegraphics[keepaspectratio, scale=0.09]{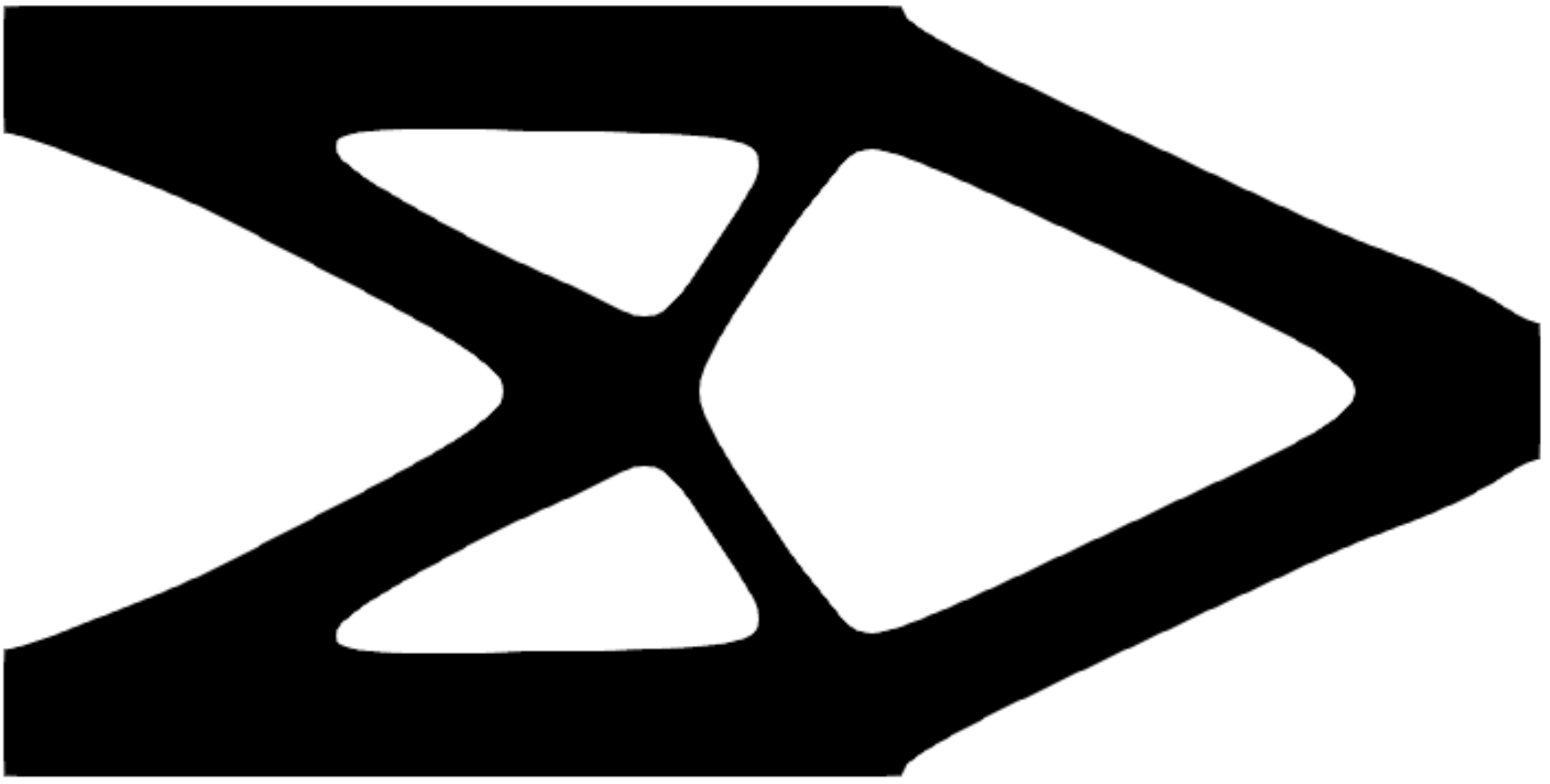}
        \subcaption{Step\,1076$^{\#}$}
        \label{2-h}
      \end{minipage}  
      \\ 
    \begin{minipage}[t]{0.2\hsize}
        \centering
        \includegraphics[keepaspectratio, scale=0.09]{caw0.pdf}
        \subcaption{Step\,0}
        \label{1-e}
      \end{minipage} 
      \begin{minipage}[t]{0.2\hsize}
        \centering
        \includegraphics[keepaspectratio, scale=0.09]{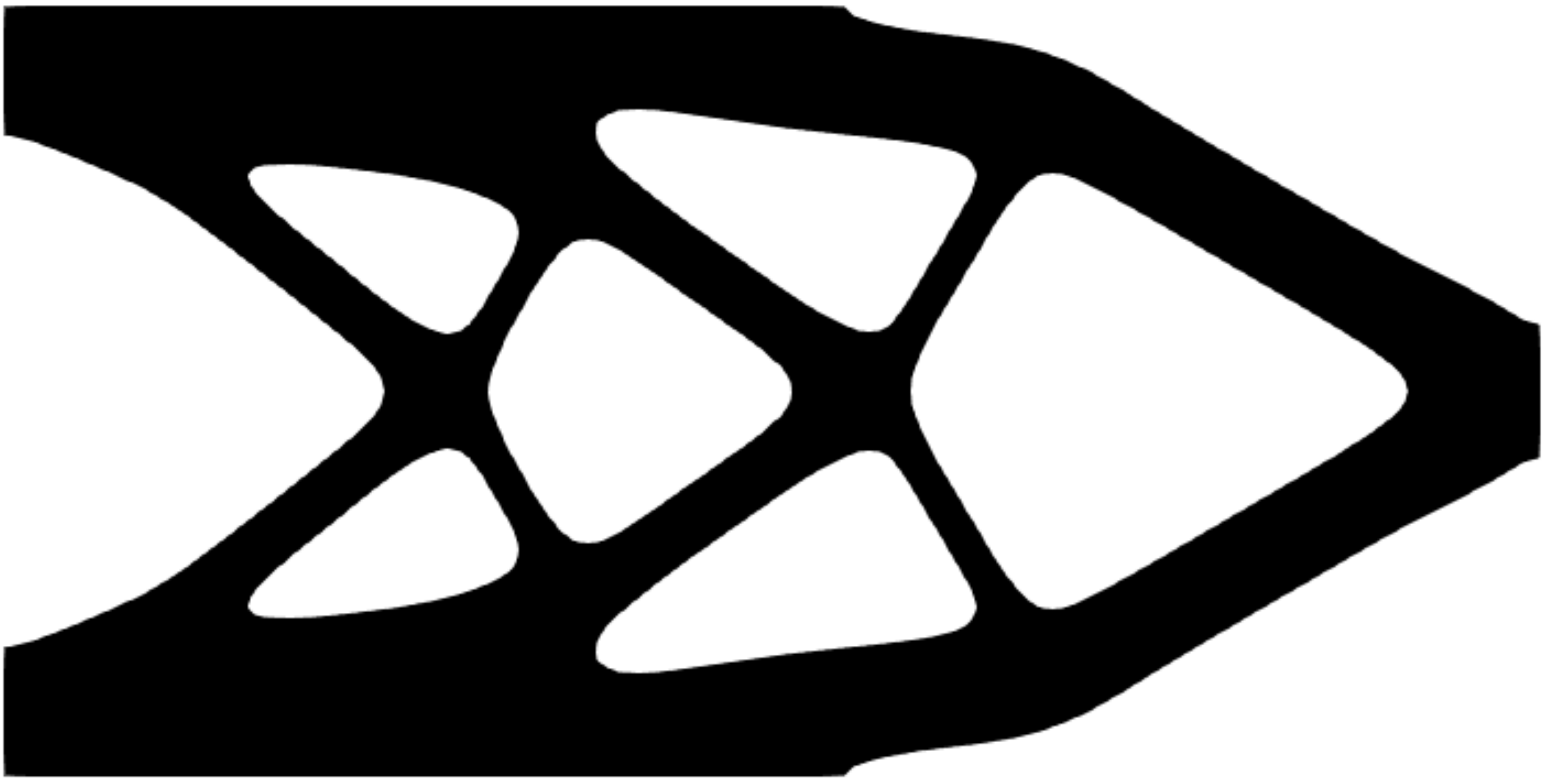}
        \subcaption{Step\,50}
        \label{1-f}
      \end{minipage} 
      \begin{minipage}[t]{0.2\hsize}
        \centering
        \includegraphics[keepaspectratio, scale=0.09]{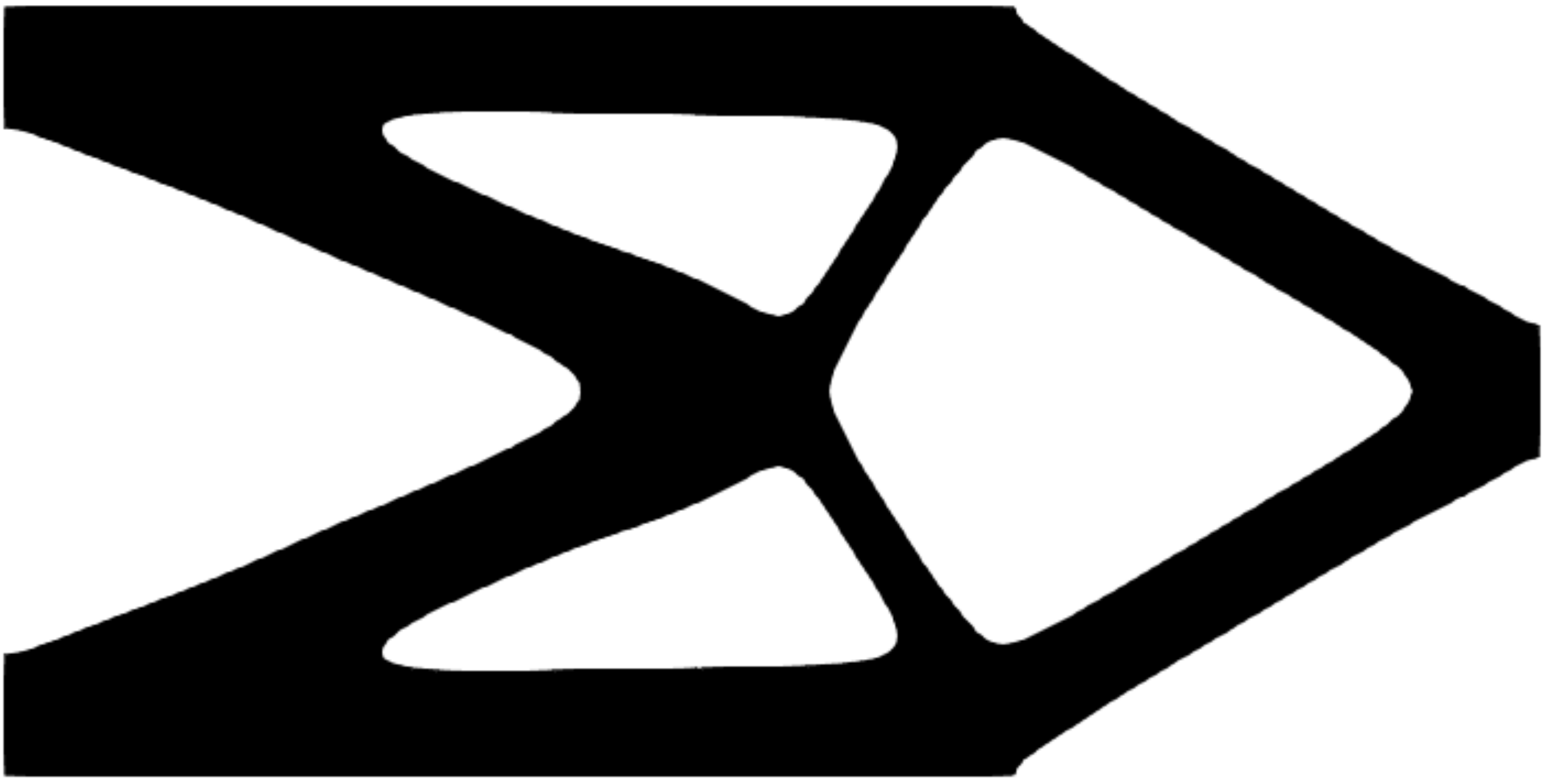}
        \subcaption{Step\,150}
        \label{1-g}
      \end{minipage} 
       \begin{minipage}[t]{0.2\hsize}
        \centering
        \includegraphics[keepaspectratio, scale=0.09]{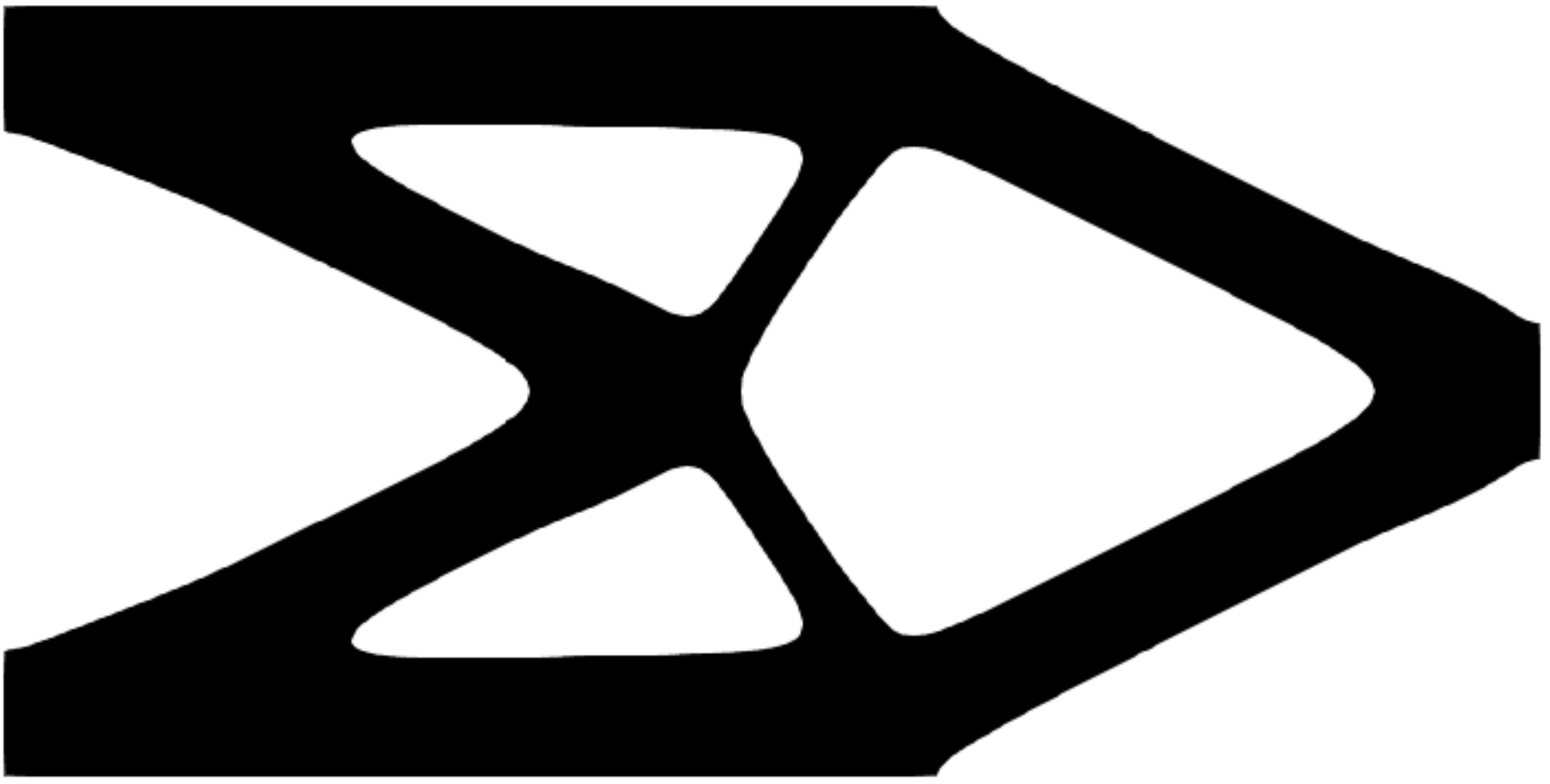}
        \subcaption{Step\,250}
        \label{1-h}
      \end{minipage} 
           \begin{minipage}[t]{0.2\hsize}
        \centering
        \includegraphics[keepaspectratio, scale=0.09]{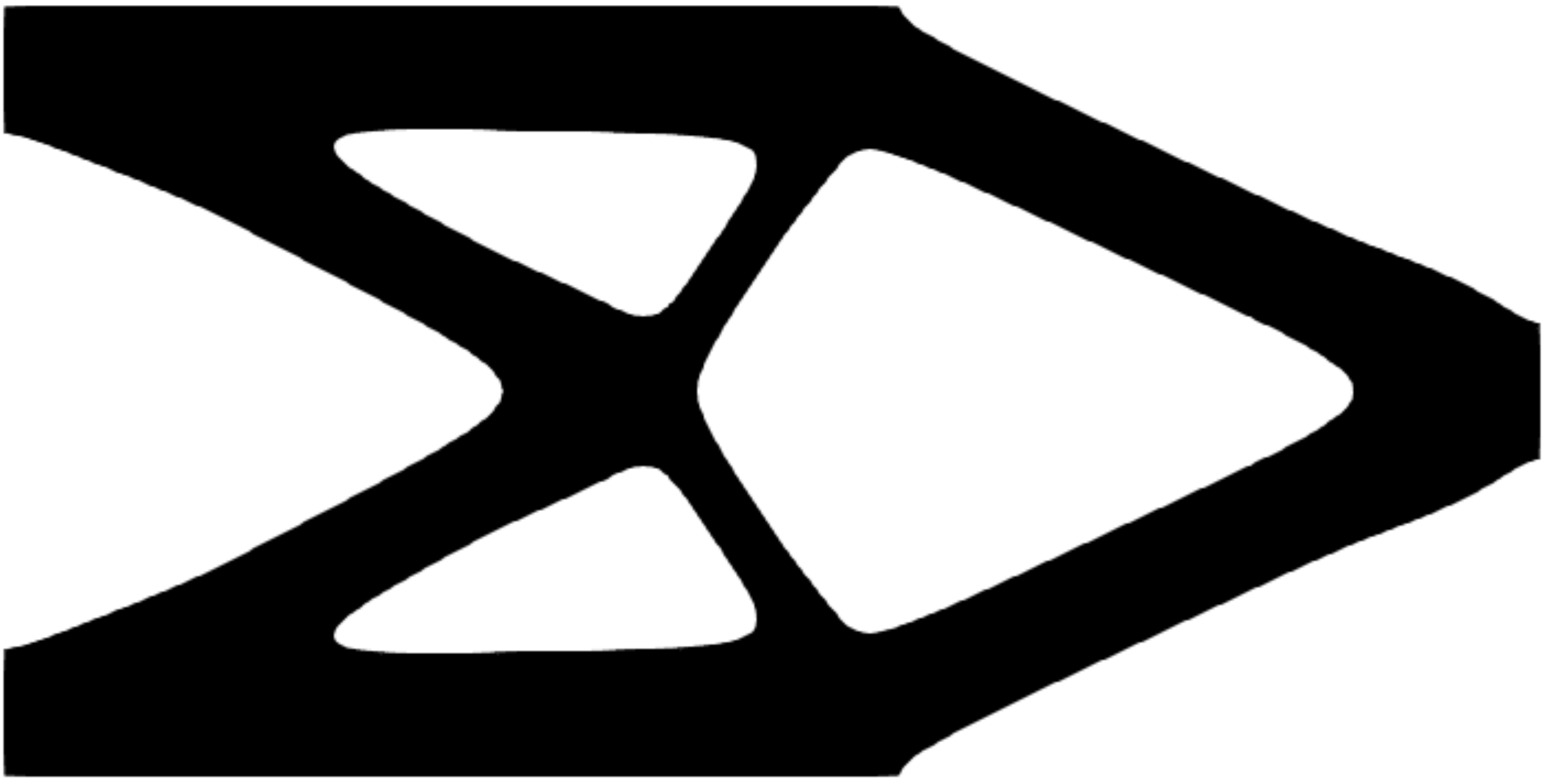}
        \subcaption{Step\,599$^{\#}$}
        \label{2-h}
      \end{minipage} 
    \end{tabular}
     \caption{Configuration $\Omega_{\phi_n}\subset D$ for the case where the initial configuration is the whole domain. 
     Figures (a)--(e) and (f)--(j) represent the results of (RD) and (NLHP), respectively.     
The symbol $^{\#}$ implies the final step.  }
     \label{case1}
  \end{figure*}

\begin{figure*}[htbp]
    \begin{tabular}{ccc}
      \hspace*{-5mm} 
      \begin{minipage}[t]{0.49\hsize}
        \centering
        \includegraphics[keepaspectratio, scale=0.33]{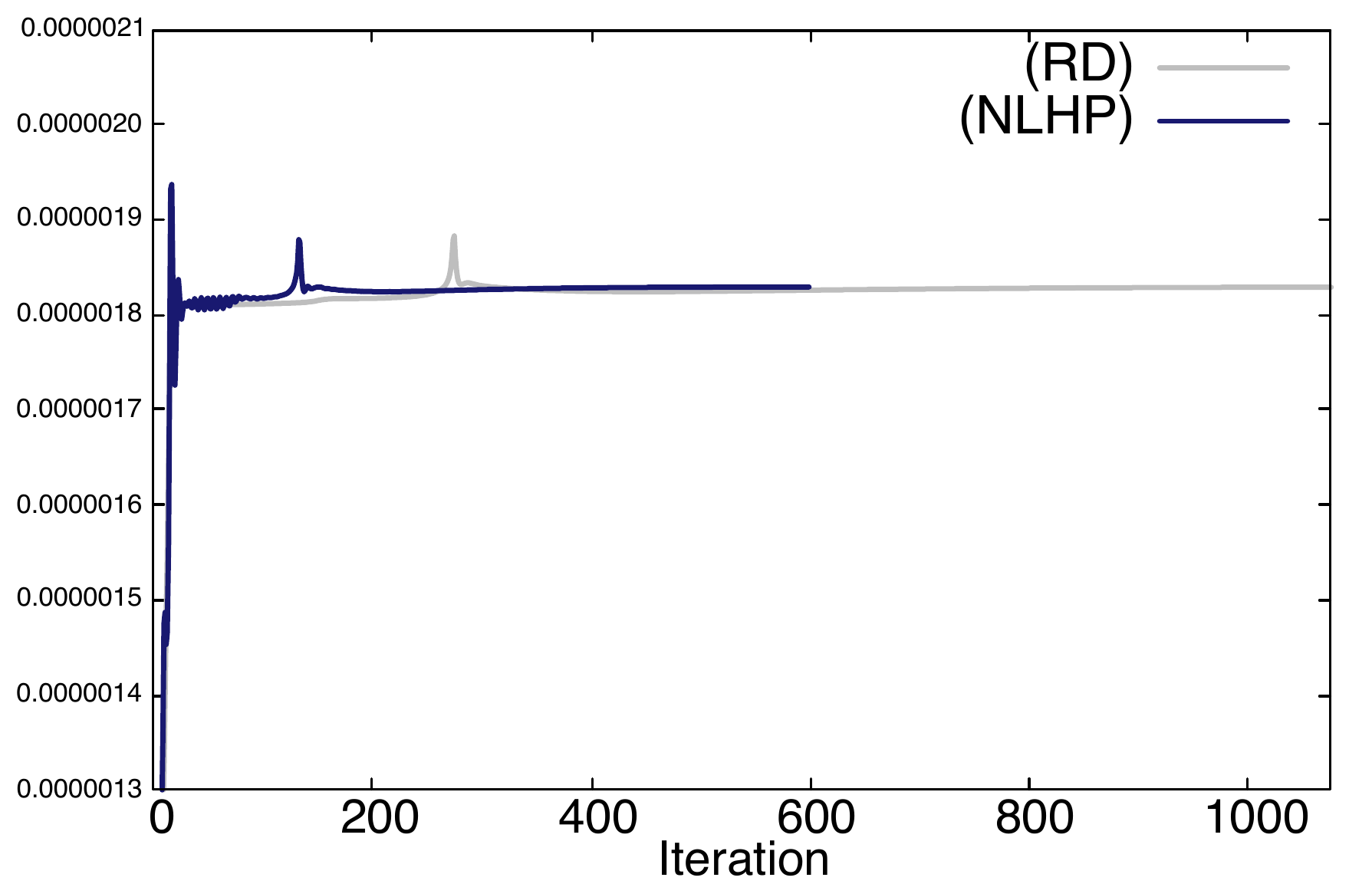}
        \subcaption{$F(\phi_n)$}
        \label{ca2-a}
      \end{minipage} 
      \begin{minipage}[t]{0.49\hsize}
        \centering
        \includegraphics[keepaspectratio, scale=0.33]{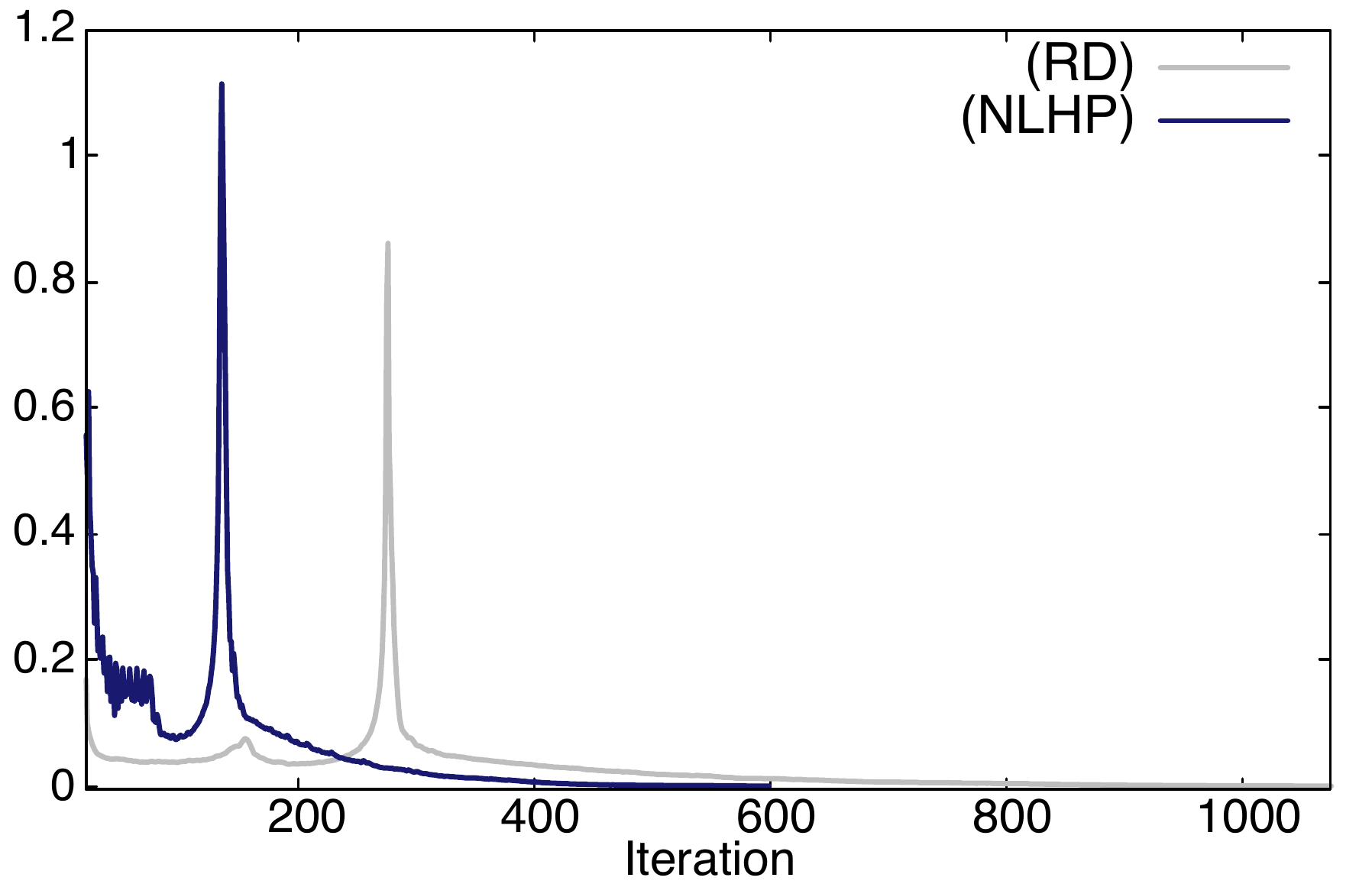}
        \subcaption{$\|\phi_{n+1}-\phi_n\|_{L^{\infty}(D)}$}
        \label{ca2-b}
      \end{minipage} &
      \end{tabular}
       \caption{ Objective functional and convergence condition for \S \ref{S:ca}-(ii).}
    \label{fig:ca2}
  \end{figure*}

\vspace{2mm}
\noindent{\bf Case (ii) (Whole domain).\,} 
We { next} consider the case where the initial configuration is the whole domain $D\subset \R^2$.
Then  Figures \ref{case1} and \ref{fig:ca2} are  obtained.  
In this case, the difference in methods obviously arises; indeed, there is 
no considerable difference up to Step\,50, but their topologies do not coincide at Step\,150. Furthermore, at Step\,250, (NLHP) is as close as possible to the optimal configuration. However, in (RD), even the topology is different from the optimal configuration.

\begin{figure*}[htbp]
\hspace*{-5mm}
    \begin{tabular}{ccccc}
      \begin{minipage}[t]{0.2\hsize}
        \centering
        \includegraphics[keepaspectratio, scale=0.09]{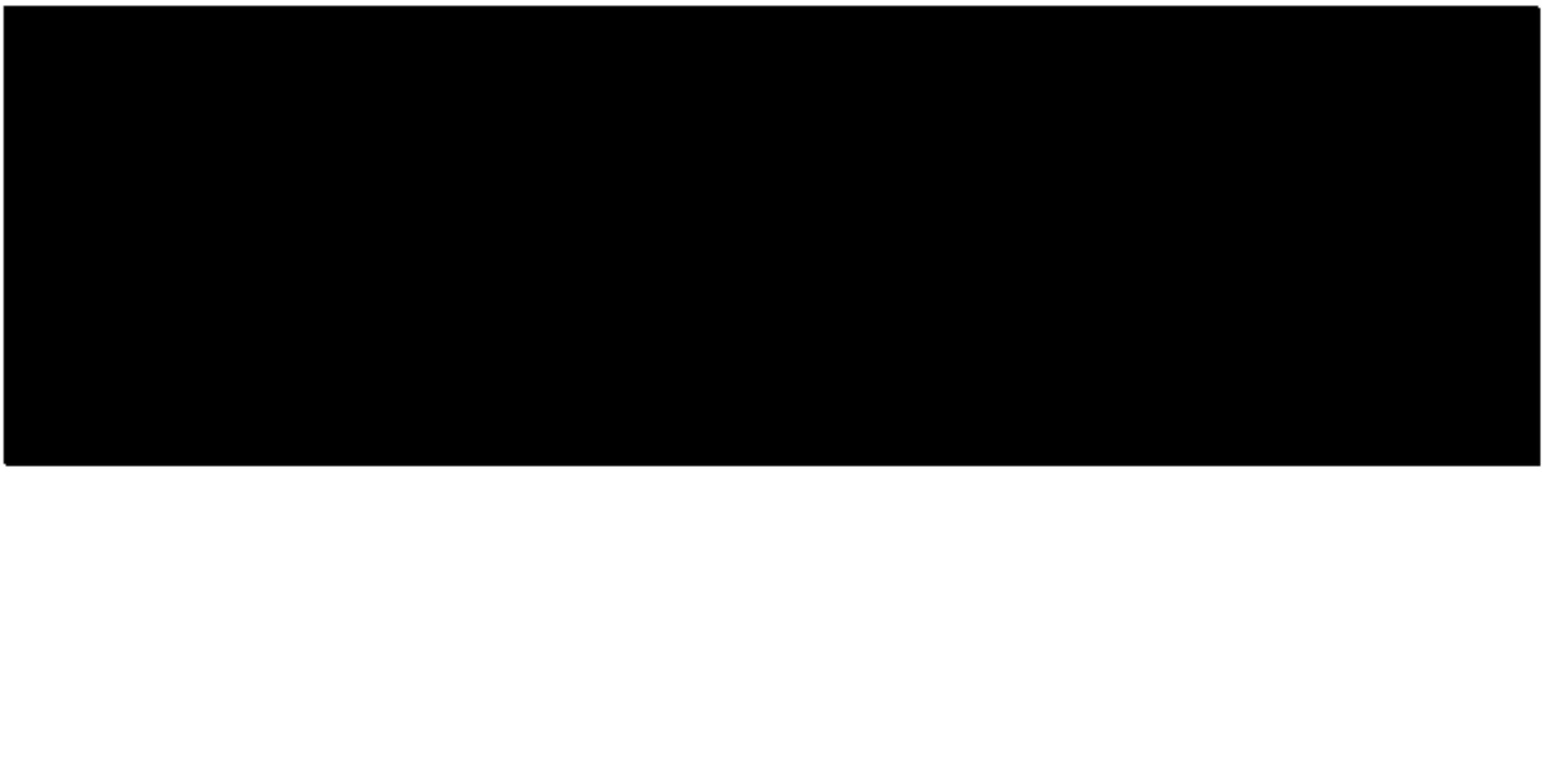}
        \subcaption{Step\,0}
        \label{3-a}
      \end{minipage} 
      \begin{minipage}[t]{0.2\hsize}
        \centering
        \includegraphics[keepaspectratio, scale=0.09]{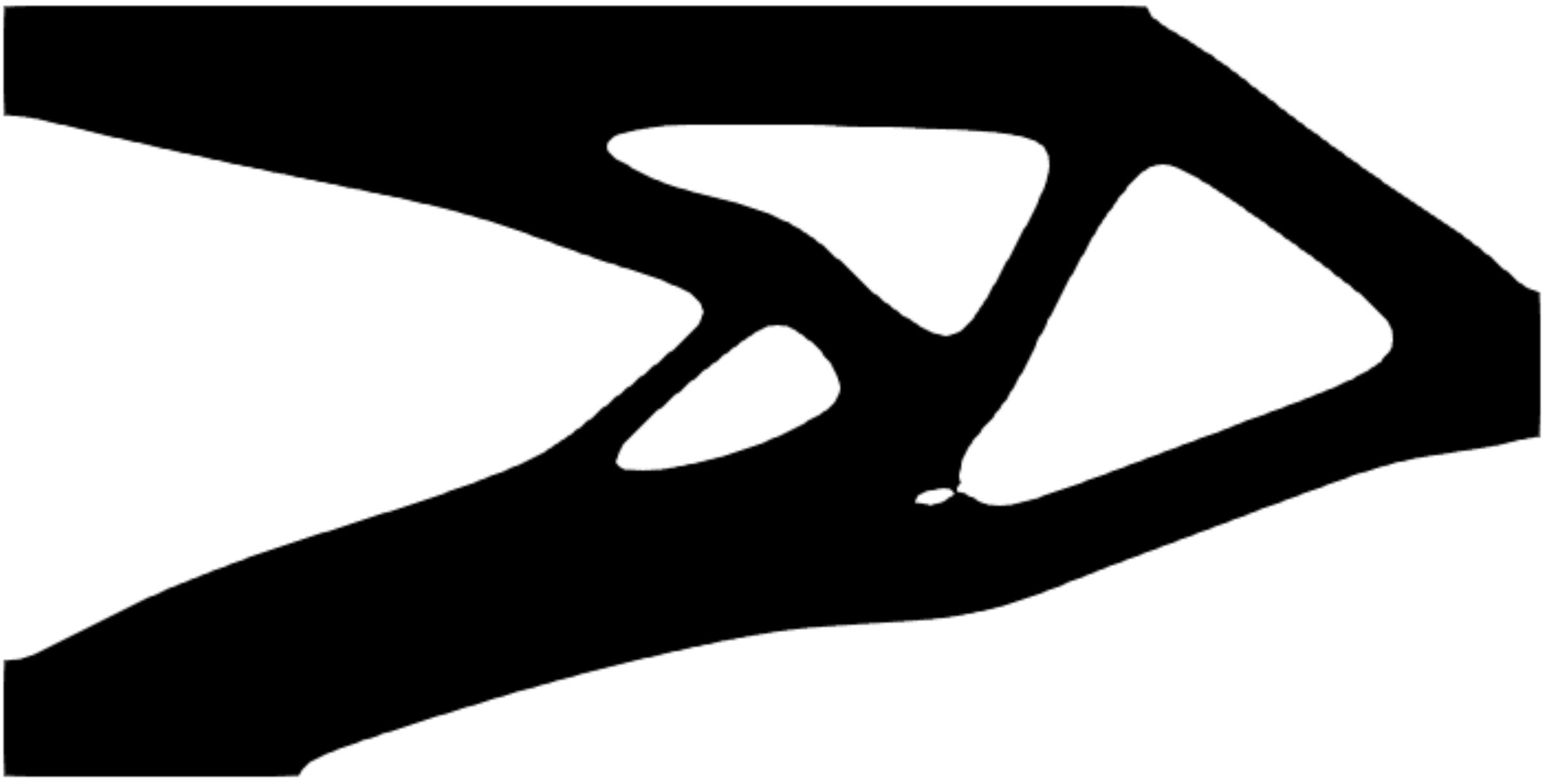}
        \subcaption{Step\,50}
        \label{3-b}
      \end{minipage} 
      \begin{minipage}[t]{0.2\hsize}
        \centering
        \includegraphics[keepaspectratio, scale=0.09]{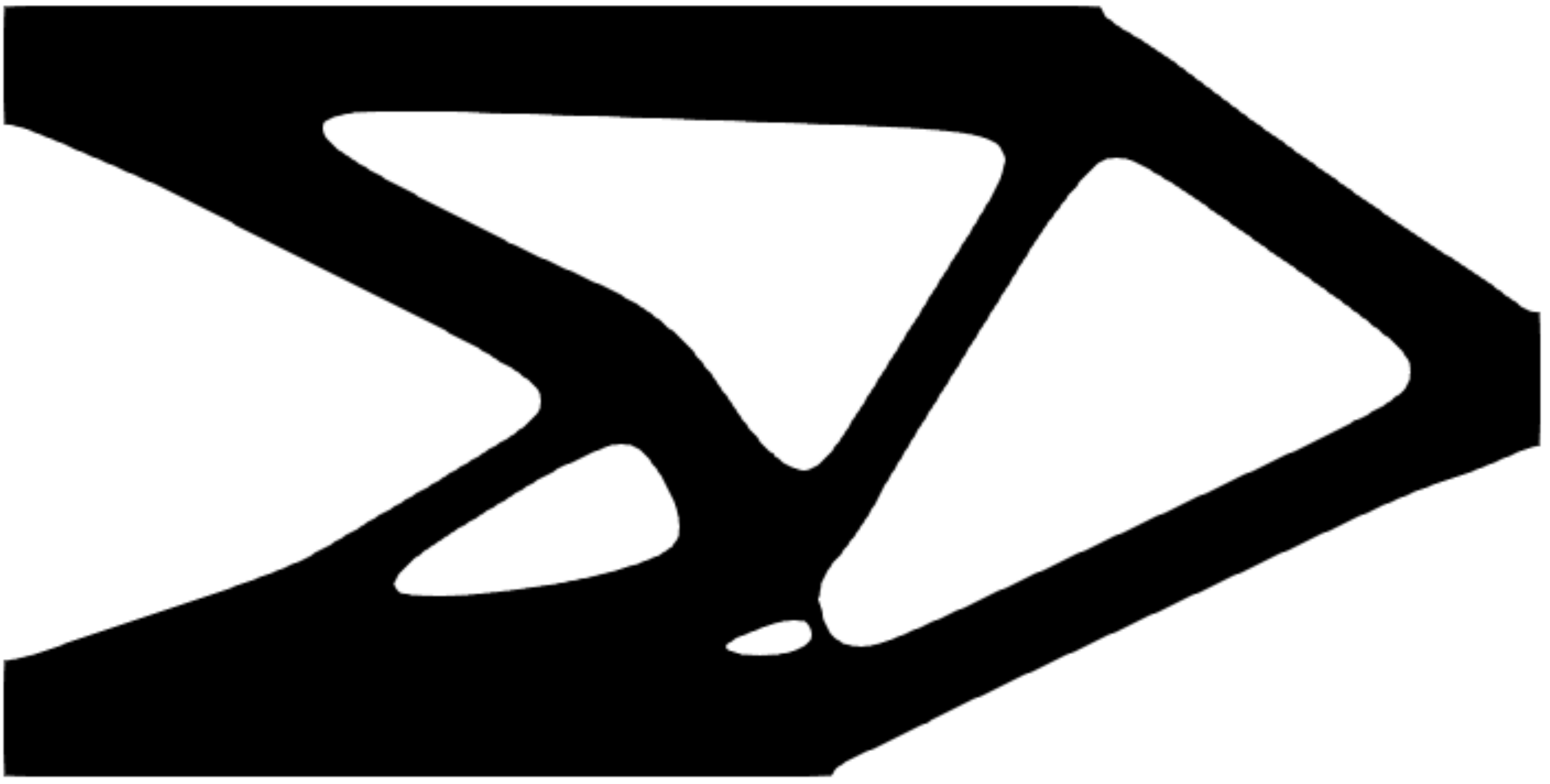}
        \subcaption{Step\,150}
        \label{3-c}
      \end{minipage} 
         \begin{minipage}[t]{0.2\hsize}
        \centering
        \includegraphics[keepaspectratio, scale=0.09]{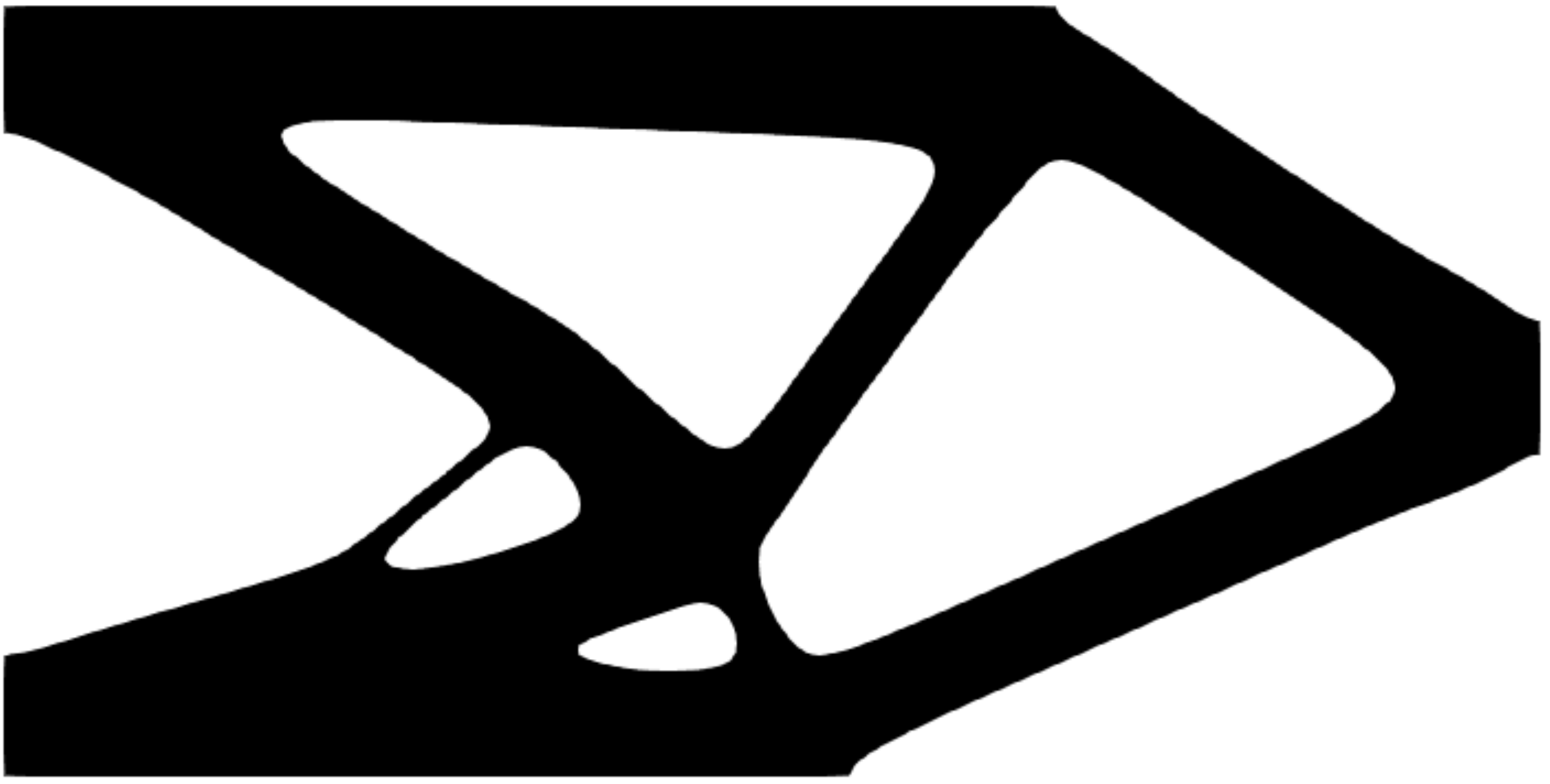}
        \subcaption{Step\,300}
        \label{3-d}
      \end{minipage} 
                 \begin{minipage}[t]{0.2\hsize}
        \centering
        \includegraphics[keepaspectratio, scale=0.09]{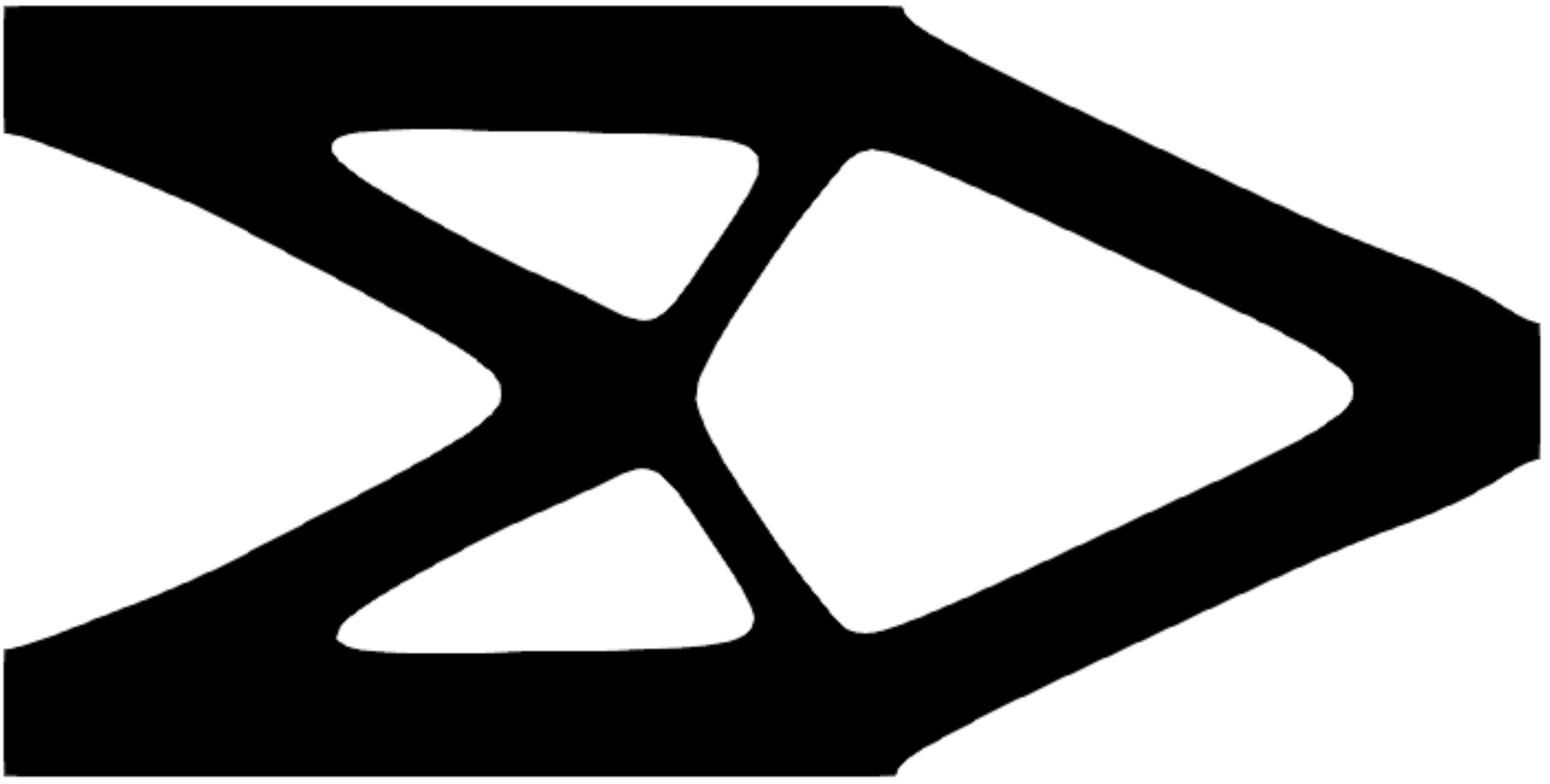}
        \subcaption{Step\,1309$^{\#}$}
      \label{3-e}
      \end{minipage} 
      \\
    \begin{minipage}[t]{0.2\hsize}
        \centering
        \includegraphics[keepaspectratio, scale=0.09]{cau0.pdf}
        \subcaption{Step\,0}
        \label{3-f}
      \end{minipage} 
      \begin{minipage}[t]{0.2\hsize}
        \centering
        \includegraphics[keepaspectratio, scale=0.09]{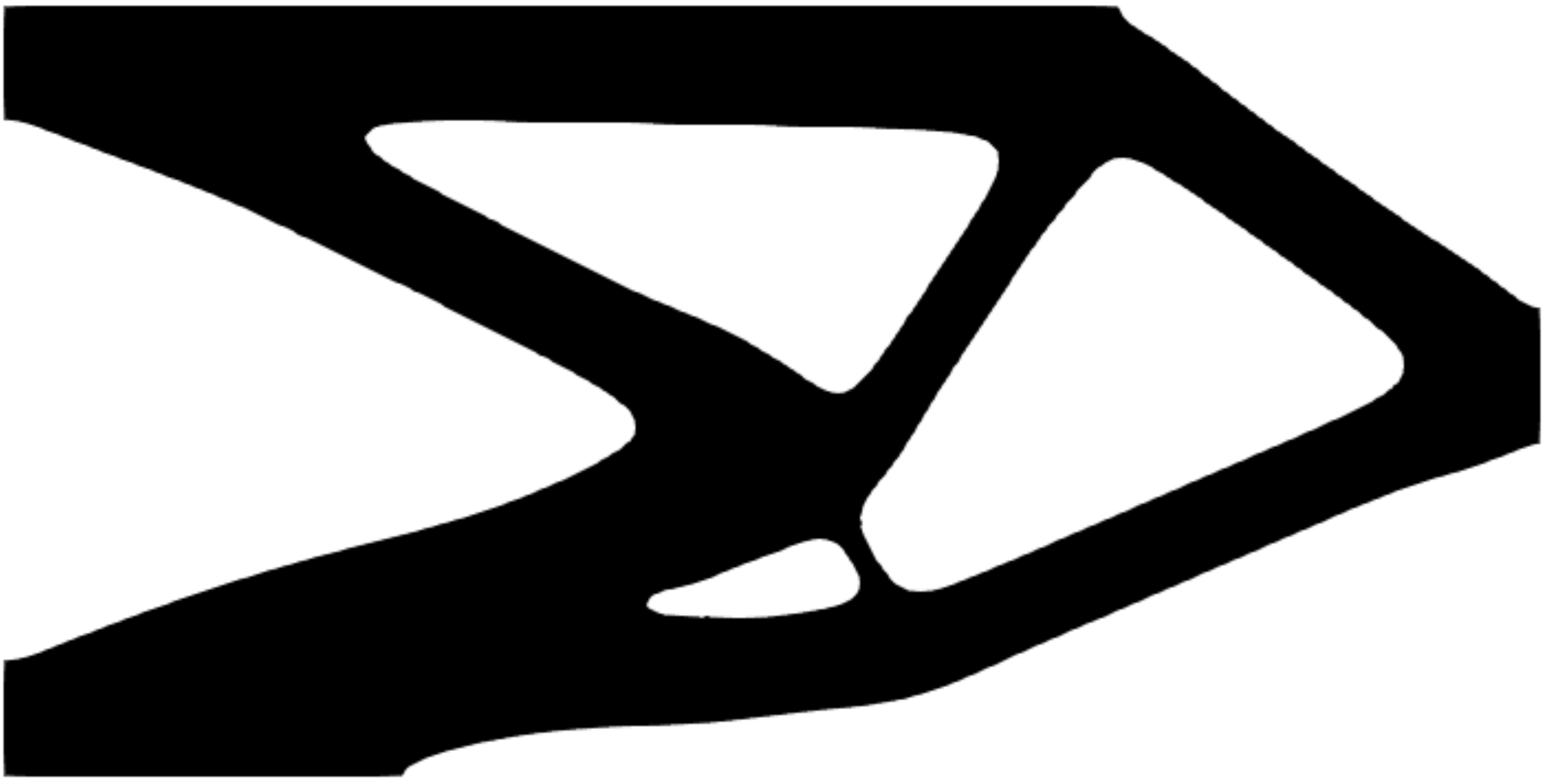}
        \subcaption{Step\,50}
        \label{3-g}
      \end{minipage} 
      \begin{minipage}[t]{0.2\hsize}
        \centering
        \includegraphics[keepaspectratio, scale=0.09]{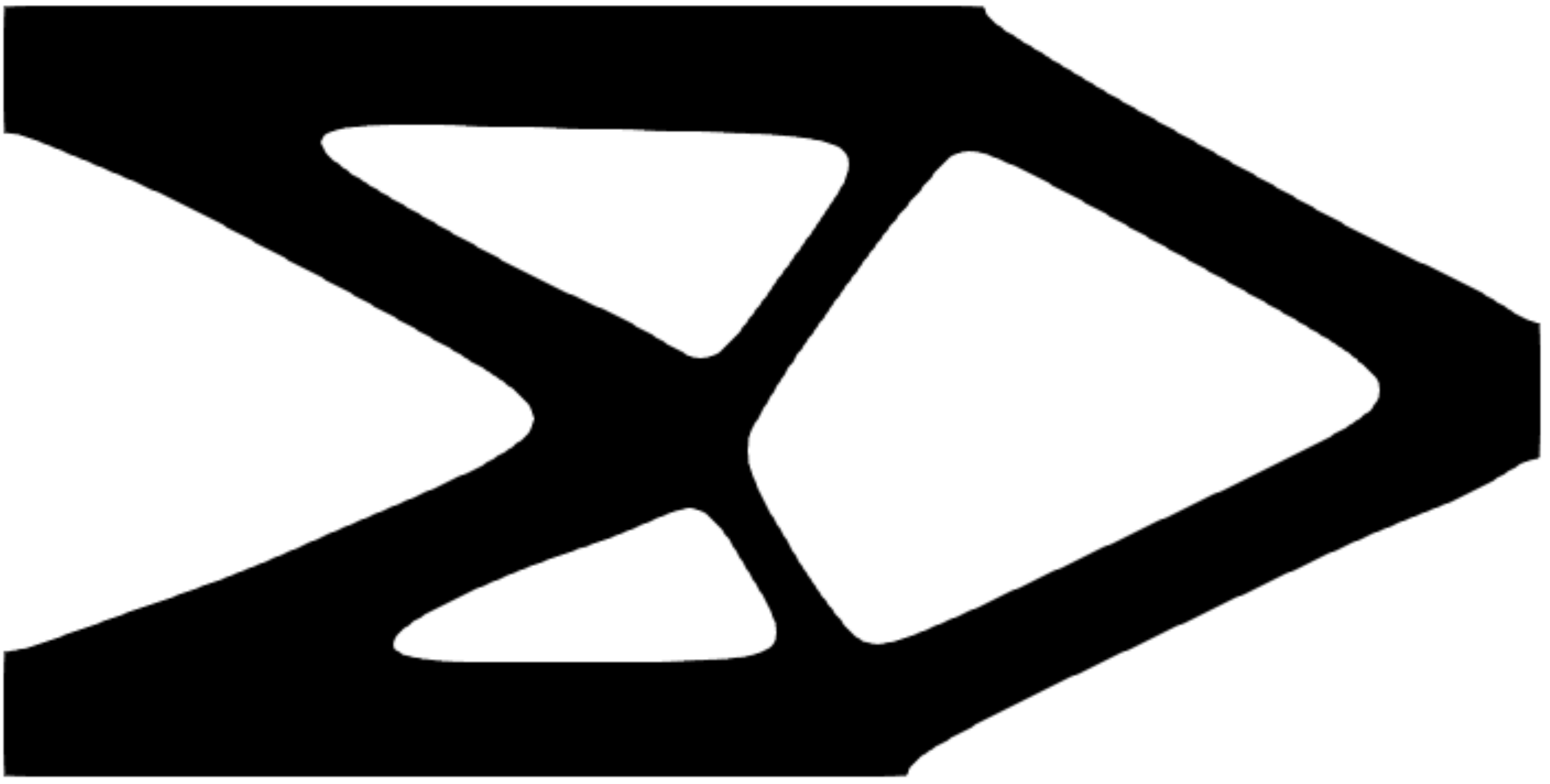}
        \subcaption{Step\,150}
        \label{3-h}
      \end{minipage} 
       \begin{minipage}[t]{0.2\hsize}
        \centering
        \includegraphics[keepaspectratio, scale=0.09]{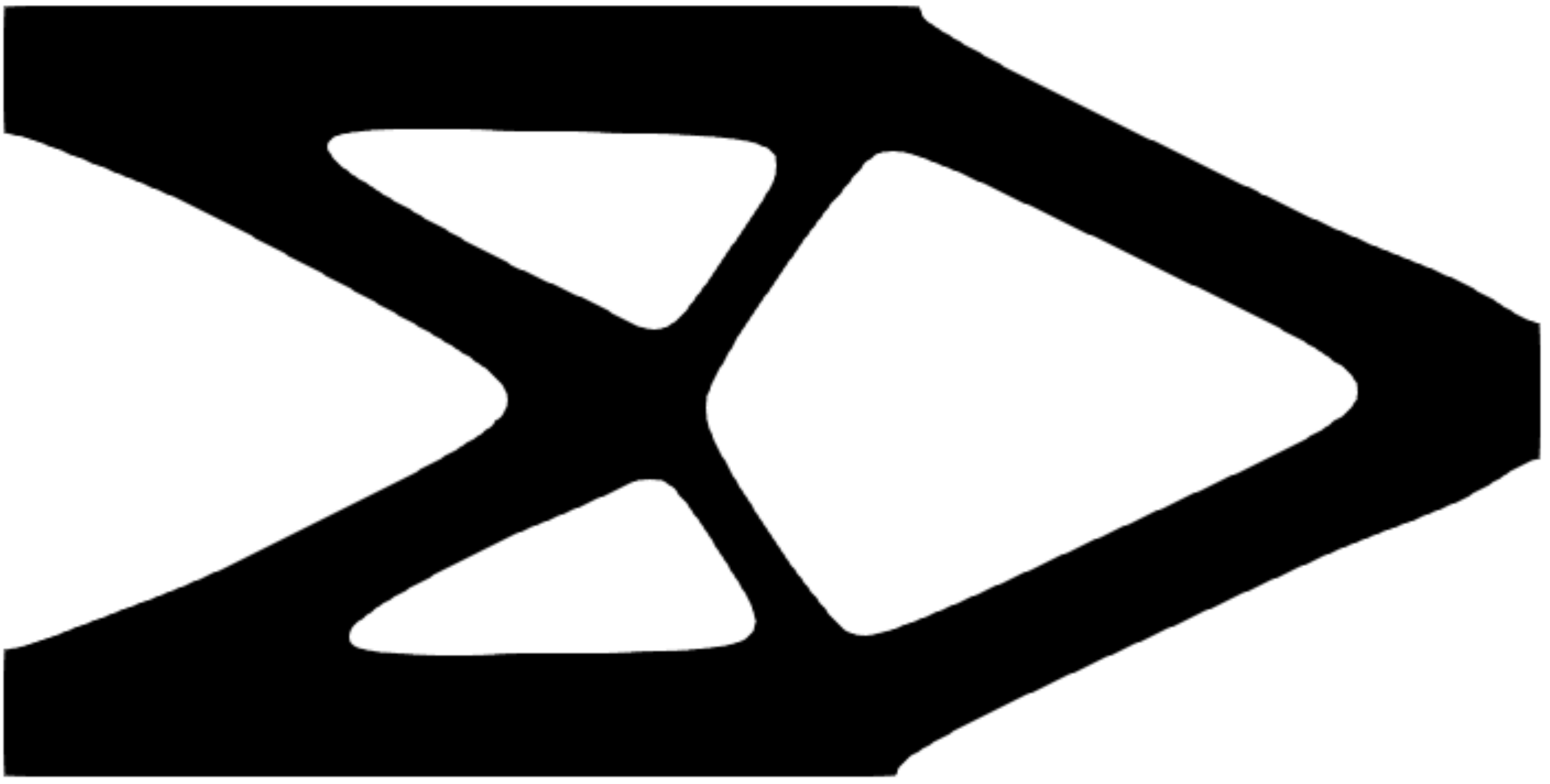}
        \subcaption{Step\,300}
        \label{3-i}
      \end{minipage}
                 \begin{minipage}[t]{0.2\hsize}
        \centering
        \includegraphics[keepaspectratio, scale=0.09]{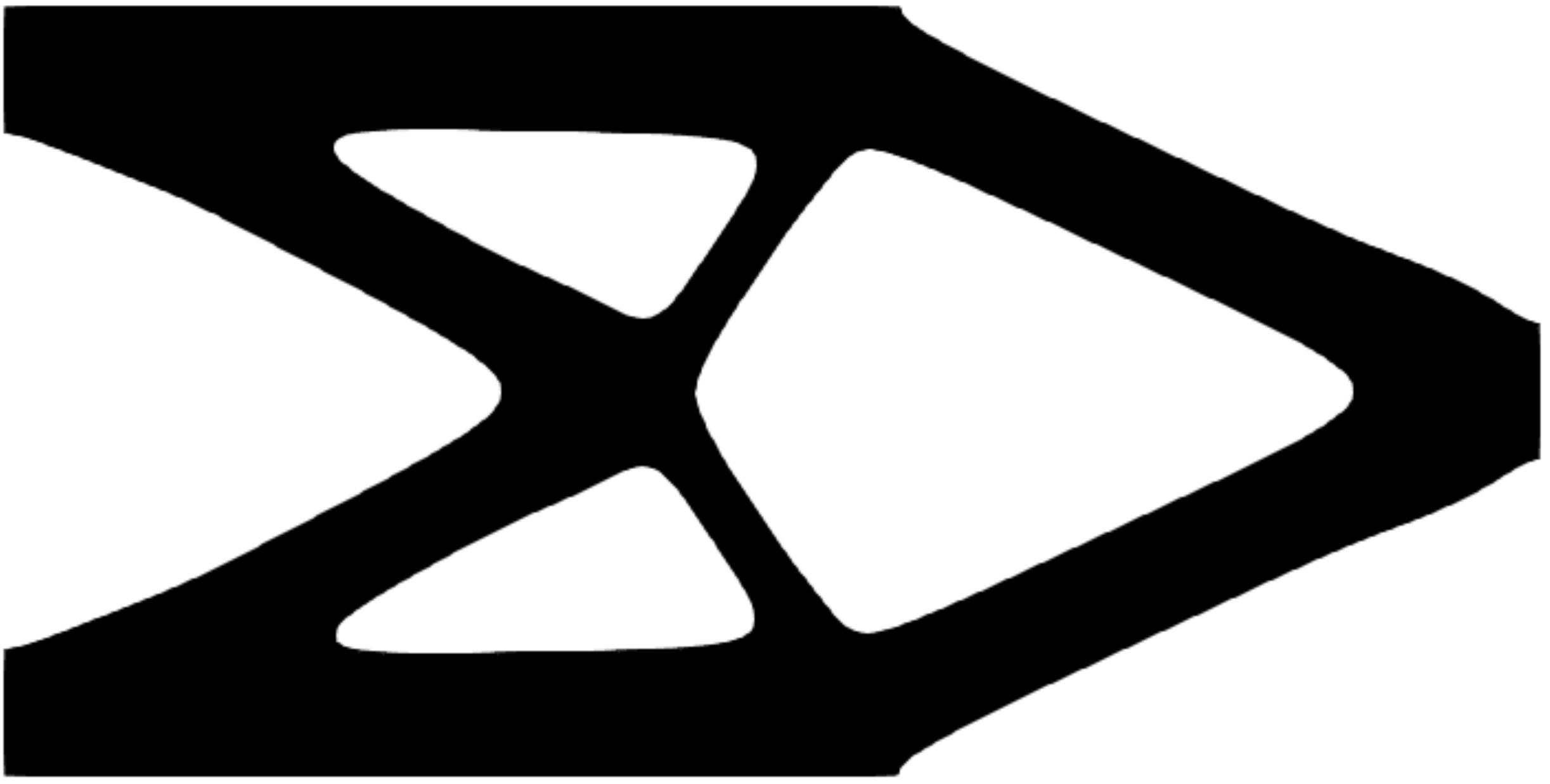}
        \subcaption{Step\,633$^{\#}$}
     \label{3-j}
      \end{minipage}  
    \end{tabular}
     \caption{Configuration $\Omega_{\phi_n}\subset D$ for the case where the initial configuration is the upper domain. 
     Figures (a)--(e) and (f)--(j) represent the results of (RD) and (NLHP), respectively.     
The symbol $^{\#}$ implies the final step. }
     \label{case3}
  \end{figure*}

\begin{figure*}[htbp]
    \begin{tabular}{ccc}
      \hspace*{-5mm} 
      \begin{minipage}[t]{0.49\hsize}
        \centering
        \includegraphics[keepaspectratio, scale=0.33]{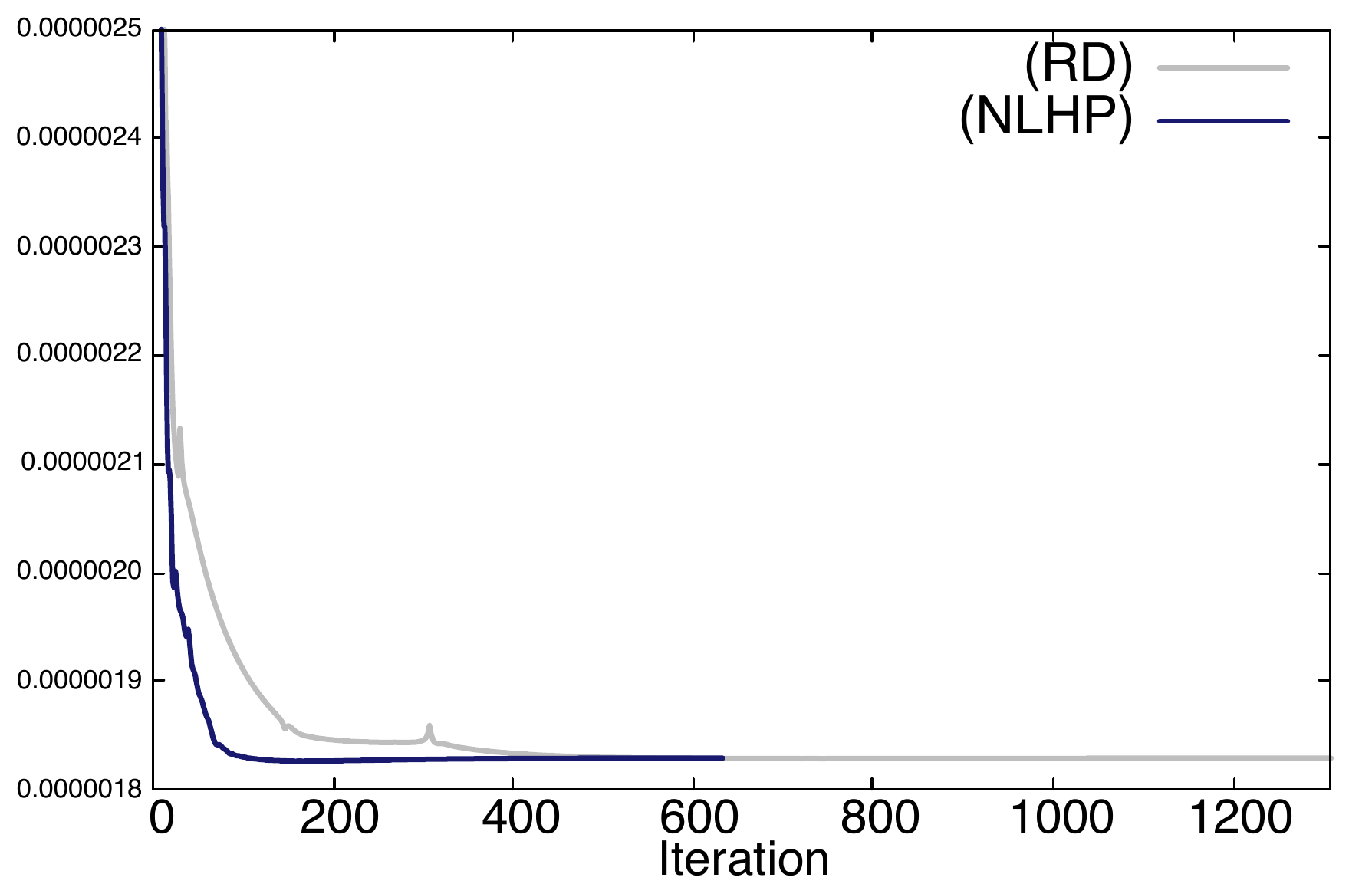}
        \subcaption{$F(\phi_n)$}
        \label{ca2-a}
      \end{minipage} 
      \begin{minipage}[t]{0.49\hsize}
        \centering
        \includegraphics[keepaspectratio, scale=0.33]{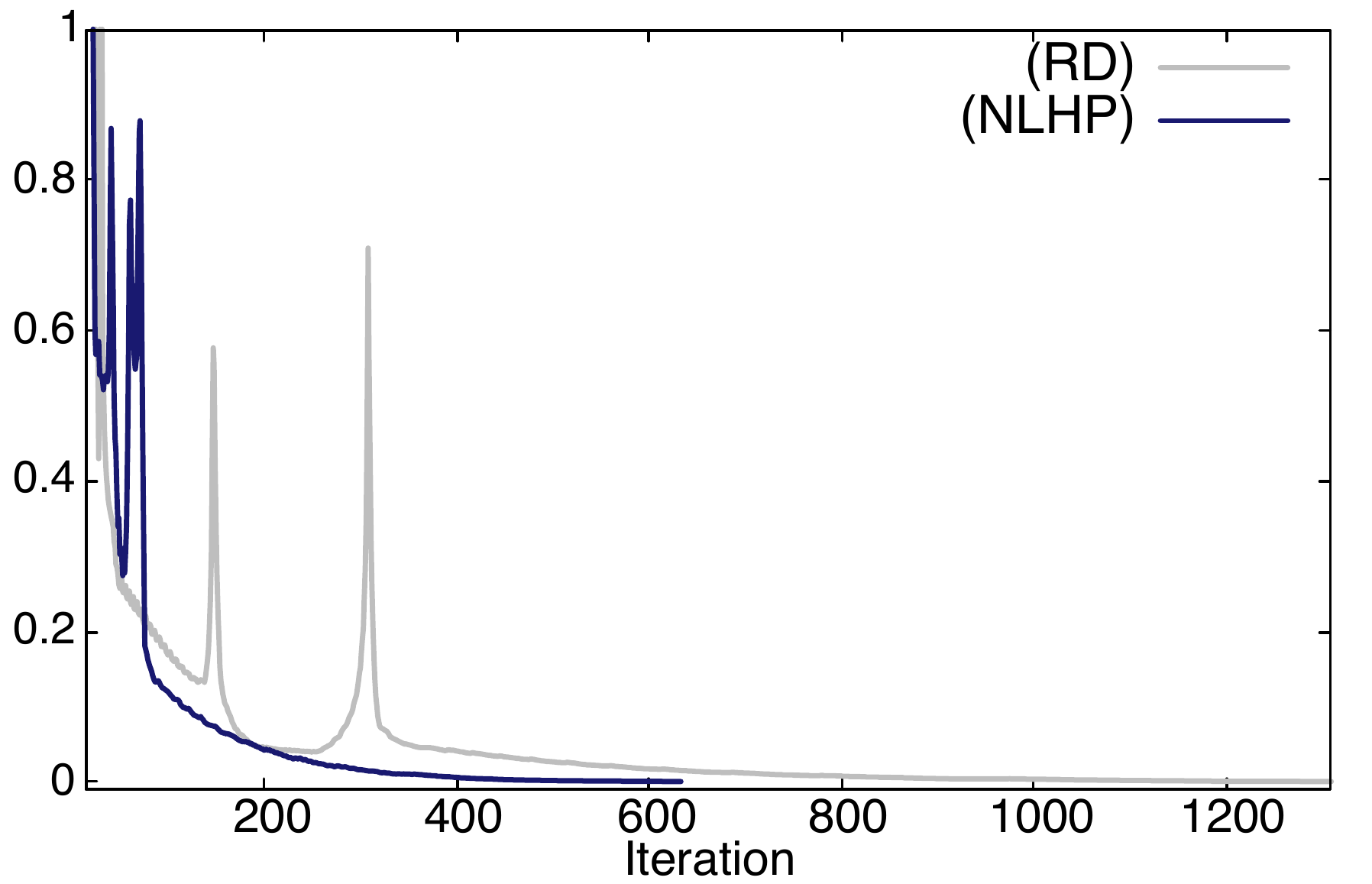}
        \subcaption{$\|\phi_{n+1}-\phi_n\|_{L^{\infty}(D)}$}
        \label{ca2-b}
      \end{minipage} &
      \end{tabular}
       \caption{ Objective functional and convergence condition for \S \ref{S:ca}-(iii).}
    \label{fig:ca3}
  \end{figure*}

\vspace{2mm}
\noindent{\bf Case (iii) (Upper domain).\,} 
We consider the case where the initial configuration is an upper domain.
Then  Figures \ref{case3} and \ref{fig:ca3} are  obtained. 
In this case, we first note that the topology and the shape must be significantly modified. 
Large differences exist in convergence among the methods (see Figure \ref{fig:ca3}). In particular, at Step\,300, Figure \ref{3-i} is considerably closer to the optimal configuration than Figure \ref{3-d}.

\begin{figure*}[htbp]
\hspace*{-5mm}
    \begin{tabular}{cccc}
      \begin{minipage}[t]{0.24\hsize}
        \centering
        \includegraphics[keepaspectratio, scale=0.105]{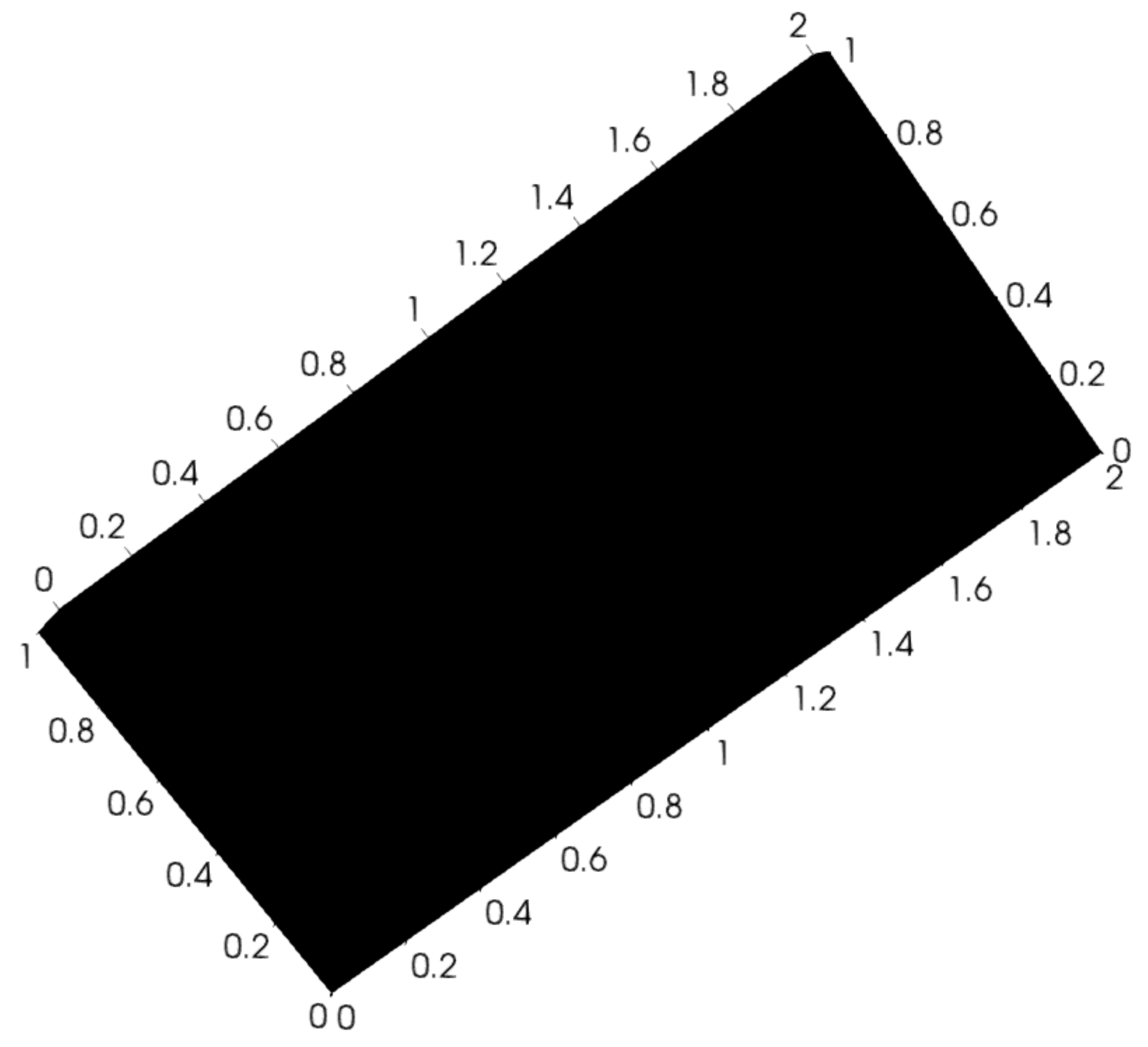}
        \subcaption{Step\,0}
        \label{3dca-a}
      \end{minipage} 
      \begin{minipage}[t]{0.24\hsize}
        \centering
        \includegraphics[keepaspectratio, scale=0.105]{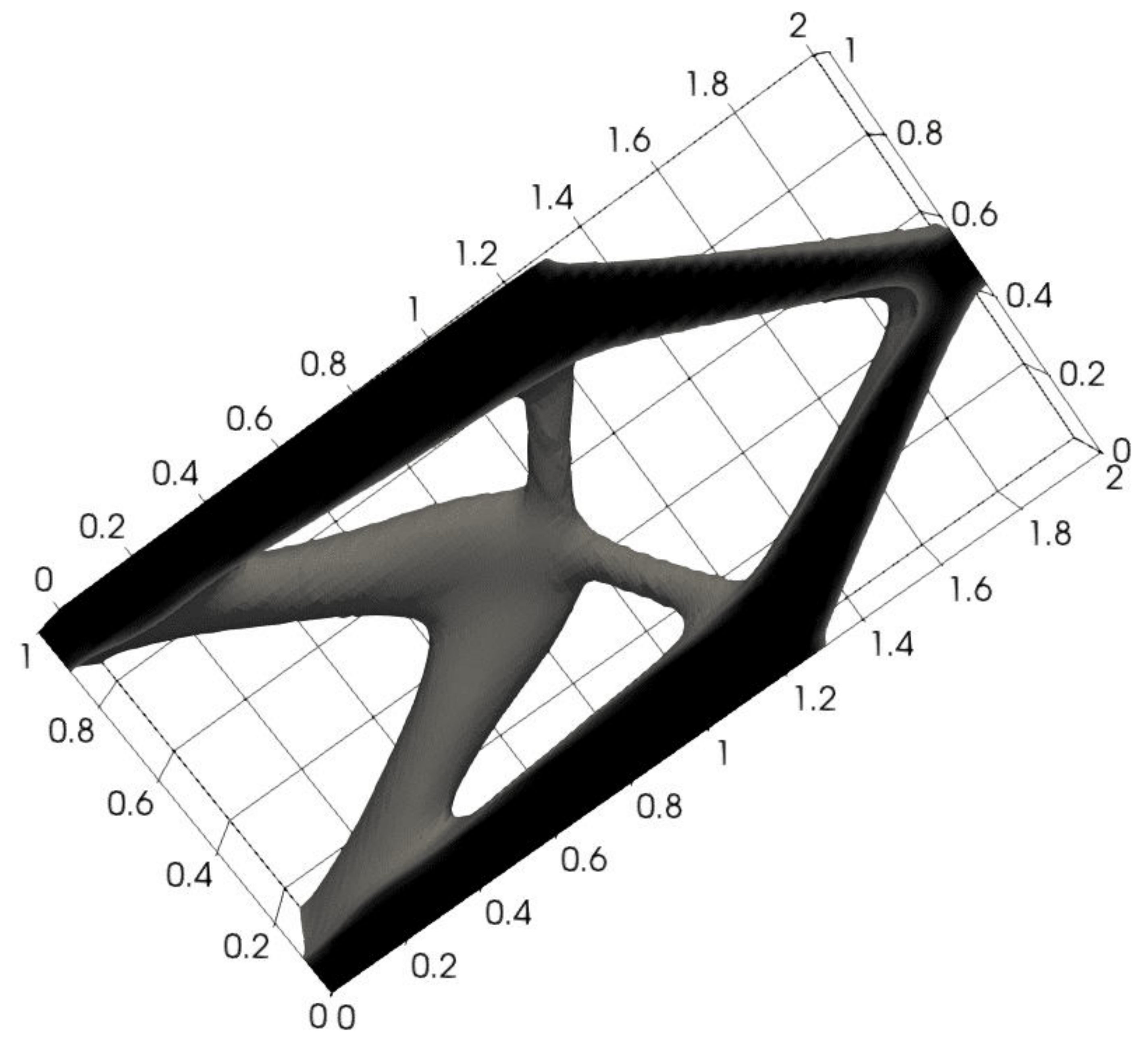}
        \subcaption{Step\,150}
        \label{3dca-b}
      \end{minipage} 
      \begin{minipage}[t]{0.24\hsize}
        \centering
        \includegraphics[keepaspectratio, scale=0.105]{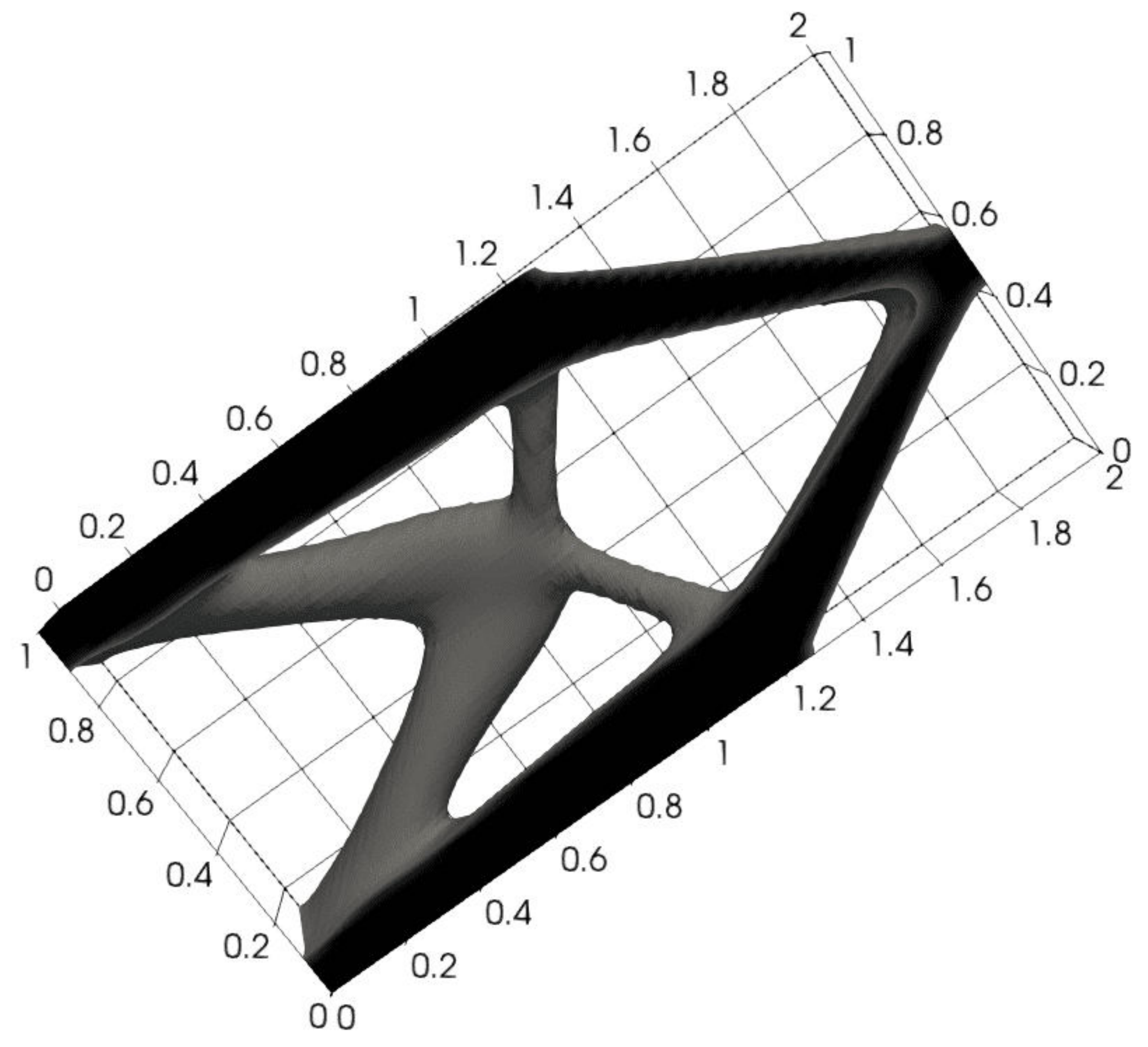}
        \subcaption{Step\,300}
        \label{3dca-c}
      \end{minipage} 
         \begin{minipage}[t]{0.24\hsize}
        \centering
        \includegraphics[keepaspectratio, scale=0.105]{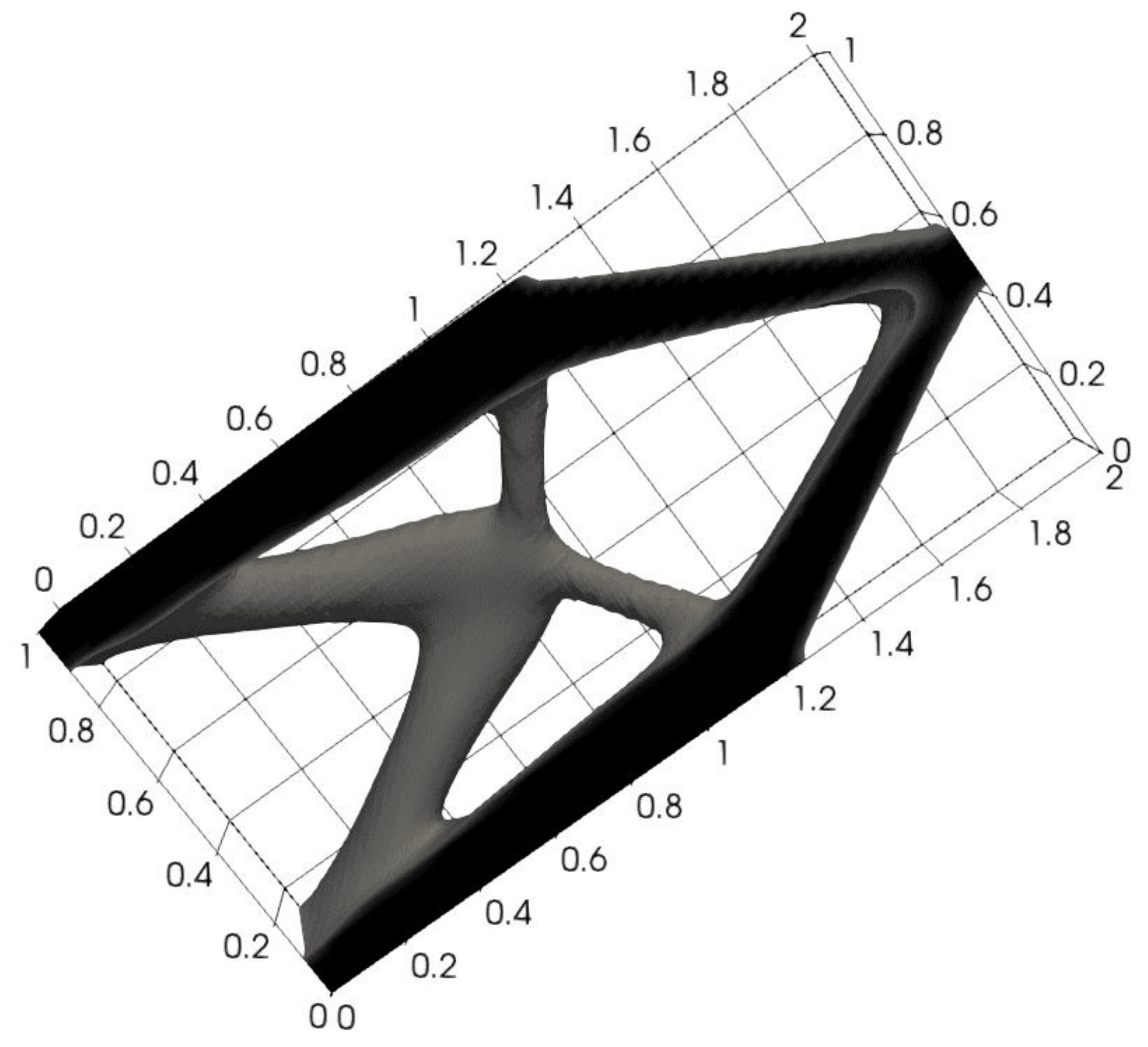}
        \subcaption{Step\,606$^\#$}
        \label{3dca-d}
      \end{minipage} 
         \\
     \begin{minipage}[t]{0.24\hsize}
        \centering
        \includegraphics[keepaspectratio, scale=0.105]{ca3d0.pdf}
        \subcaption{Step\,0}
        \label{3dca-e}
      \end{minipage} 
      \begin{minipage}[t]{0.24\hsize}
        \centering
        \includegraphics[keepaspectratio, scale=0.105]{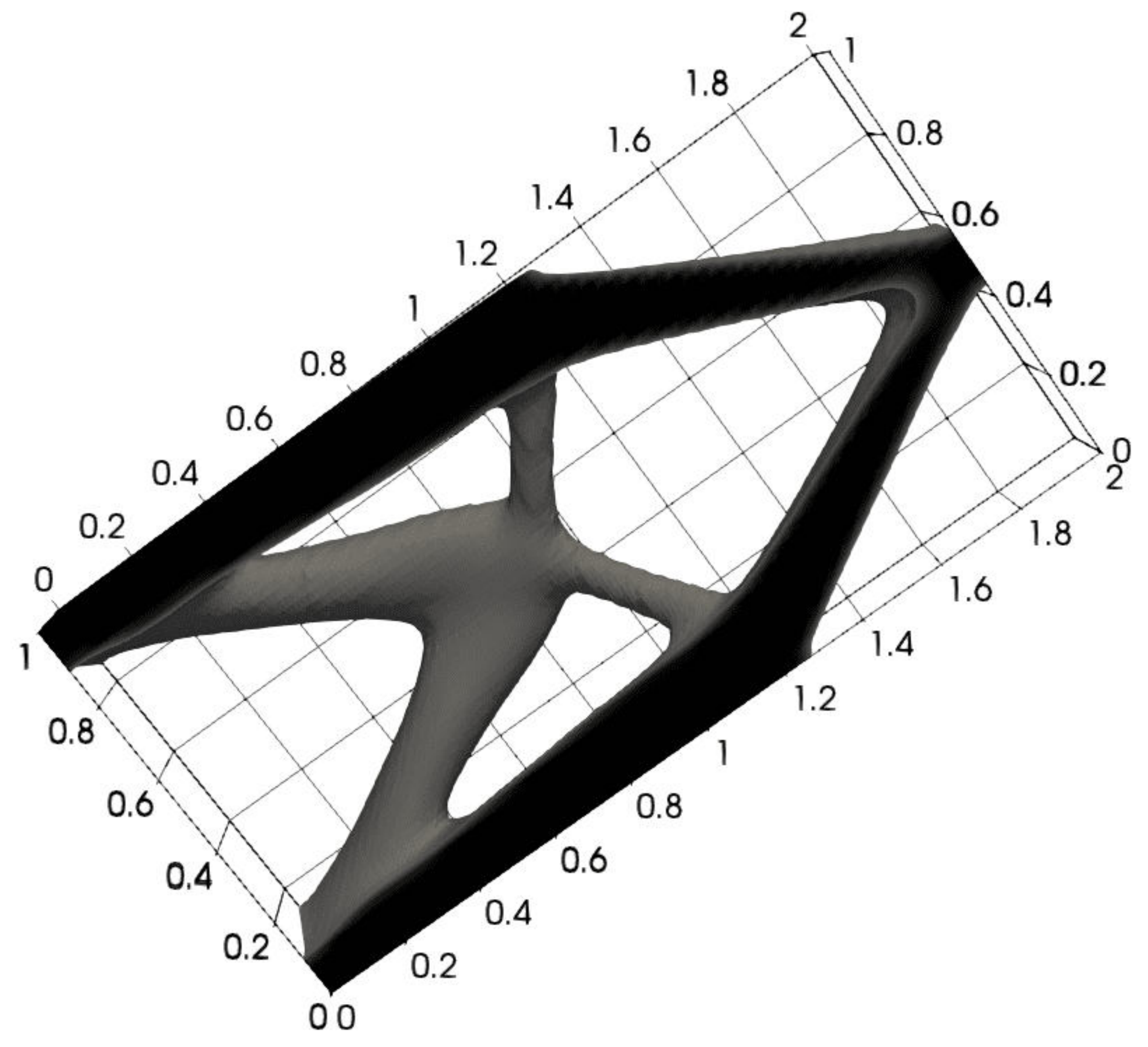}
        \subcaption{Step\,150}
        \label{3dca-f}
      \end{minipage} 
      \begin{minipage}[t]{0.24\hsize}
        \centering
        \includegraphics[keepaspectratio, scale=0.105]{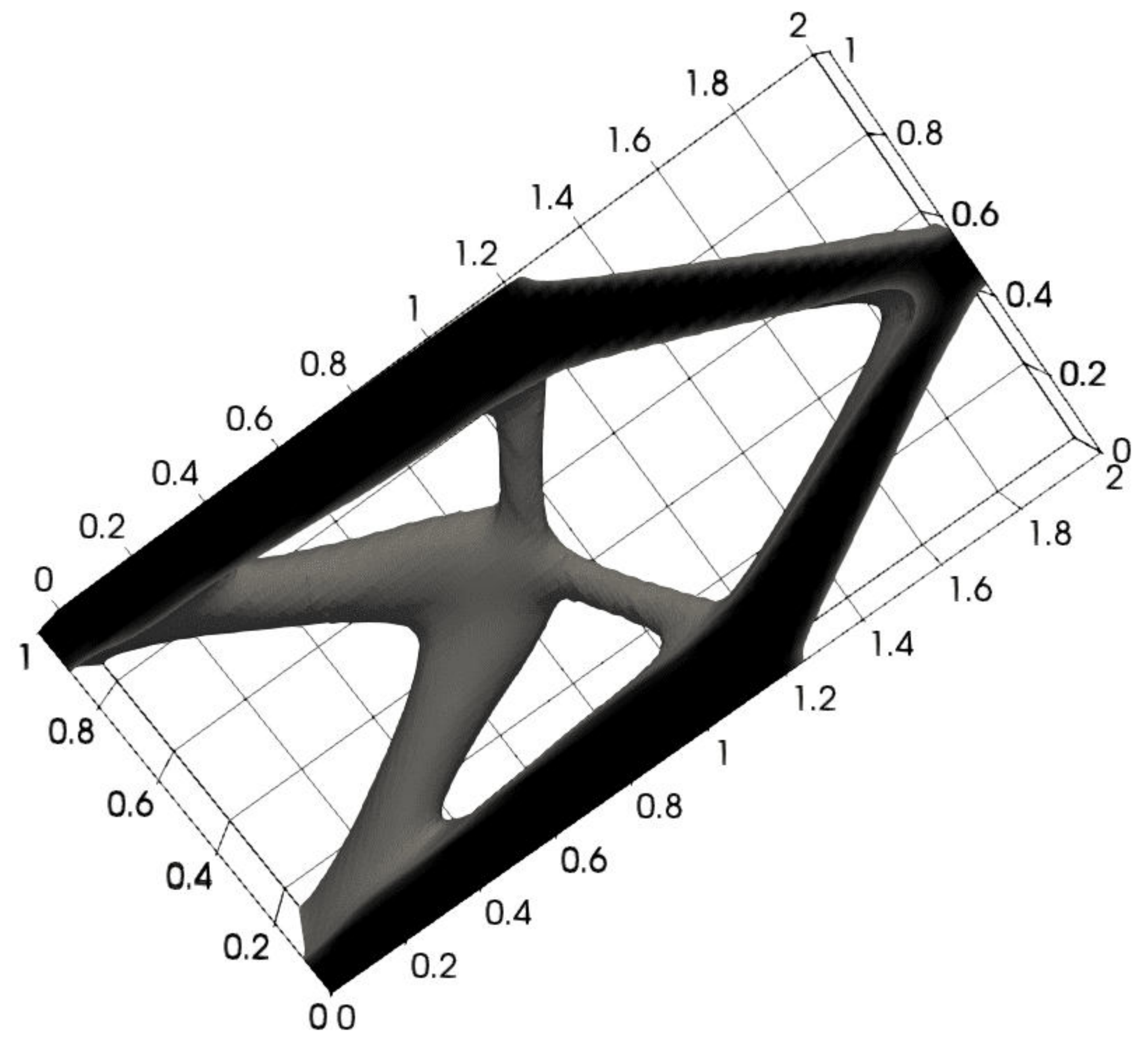}
        \subcaption{Step\,300}
        \label{3dca-g}
      \end{minipage} 
         \begin{minipage}[t]{0.24\hsize}
        \centering
        \includegraphics[keepaspectratio, scale=0.105]{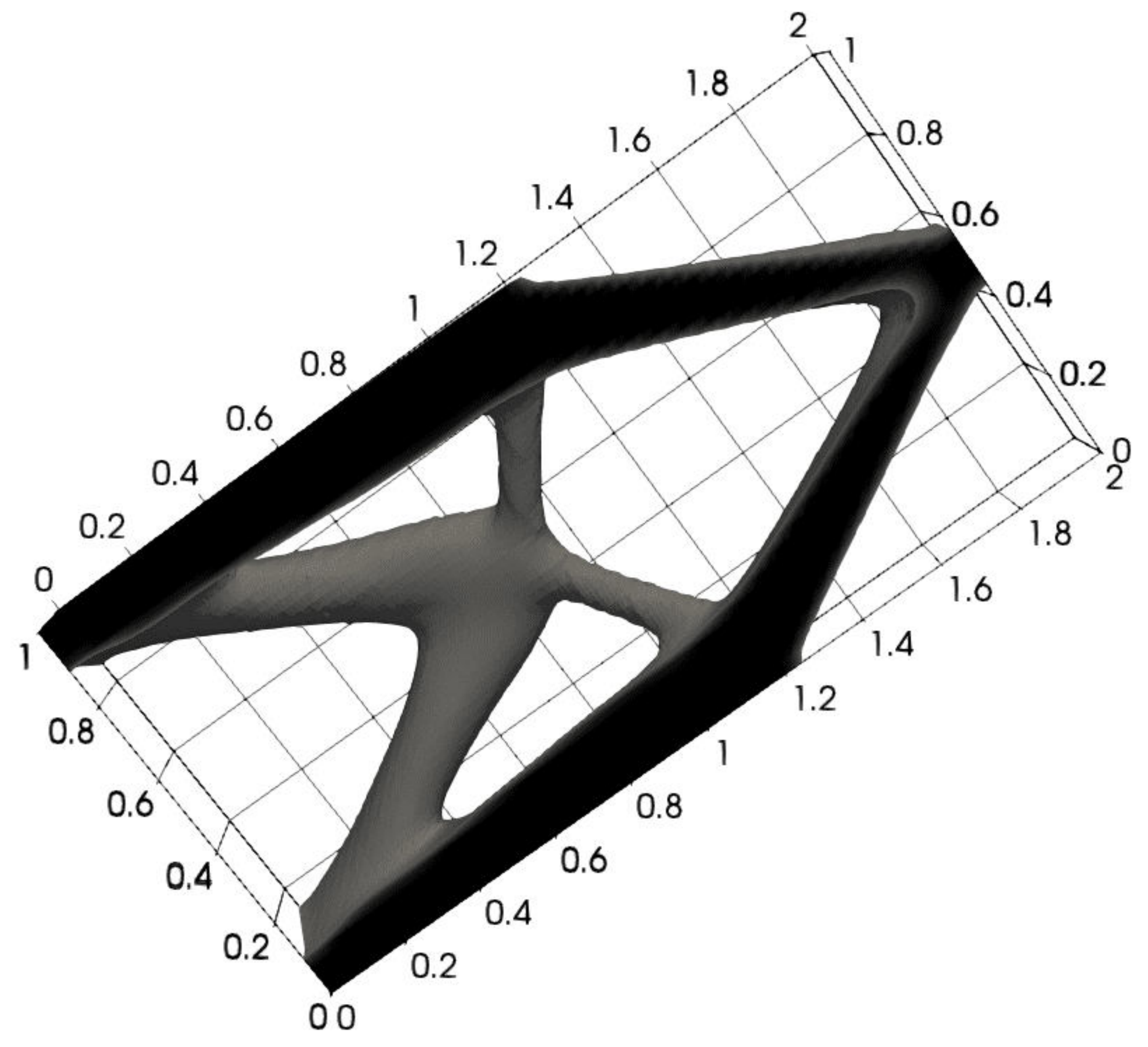}
        \subcaption{Step\,347$^\#$}
        \label{3dca-h}
      \end{minipage} 
       \end{tabular}
     \caption{ Configuration $\Omega_{\phi_n}\subset D\subset \R^3$ for the case where the initial configuration is the whole domain. 
     Figures (a)--(d) and (e)--(h) represent the results of (RD) and (NLHP), respectively.     
The symbol $^{\#}$ implies the final step. Here the depth of $D$ is set to $0.2$.}
     \label{fig:3dca}
  \end{figure*}

\begin{figure*}[htbp]
    \begin{tabular}{ccc}
      \hspace*{-5mm} 
      \begin{minipage}[t]{0.49\hsize}
        \centering
        \includegraphics[keepaspectratio, scale=0.33]{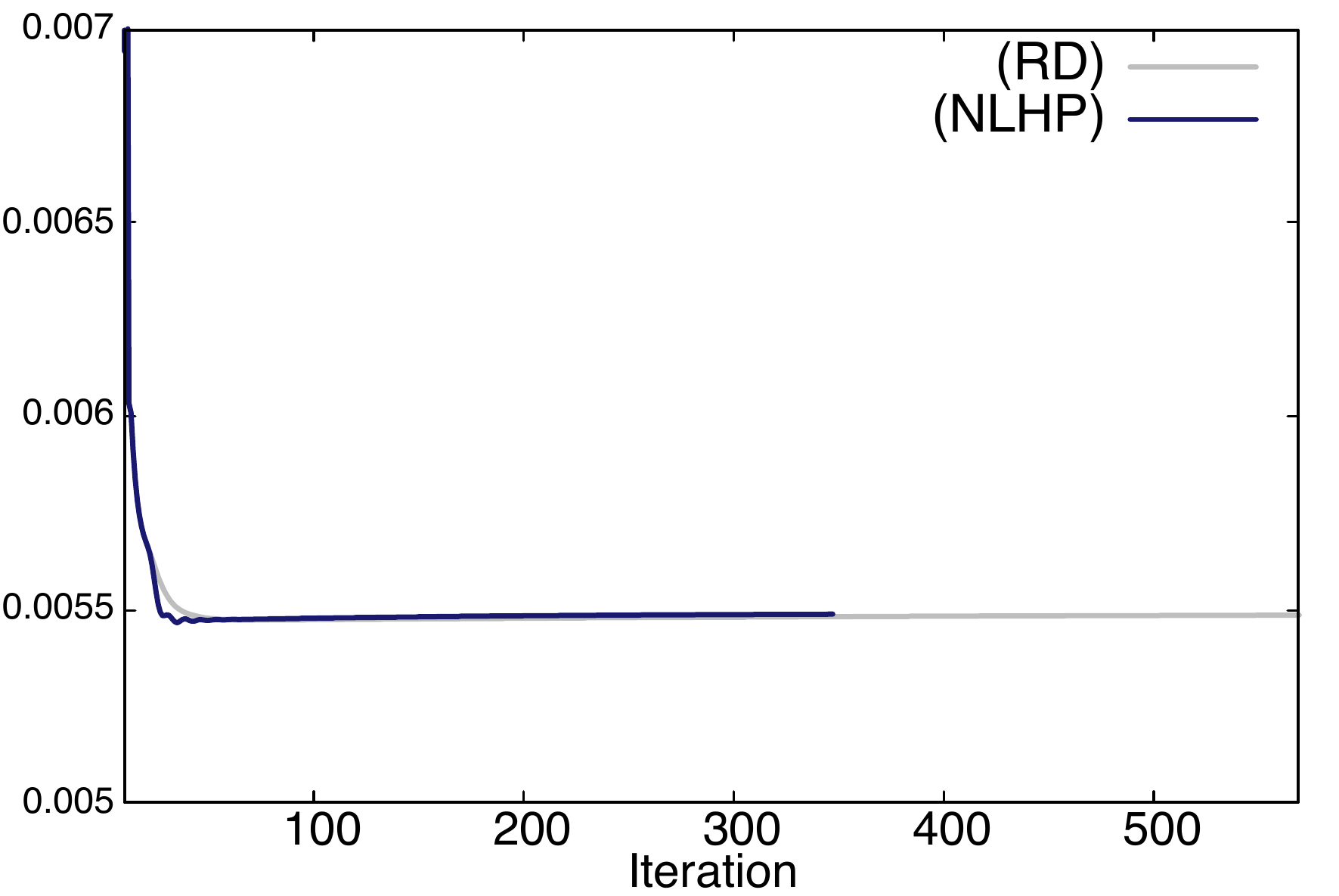}
        \subcaption{$F(\phi_n)$}
        \label{3dca-1}
      \end{minipage} 
      \begin{minipage}[t]{0.49\hsize}
        \centering
        \includegraphics[keepaspectratio, scale=0.33]{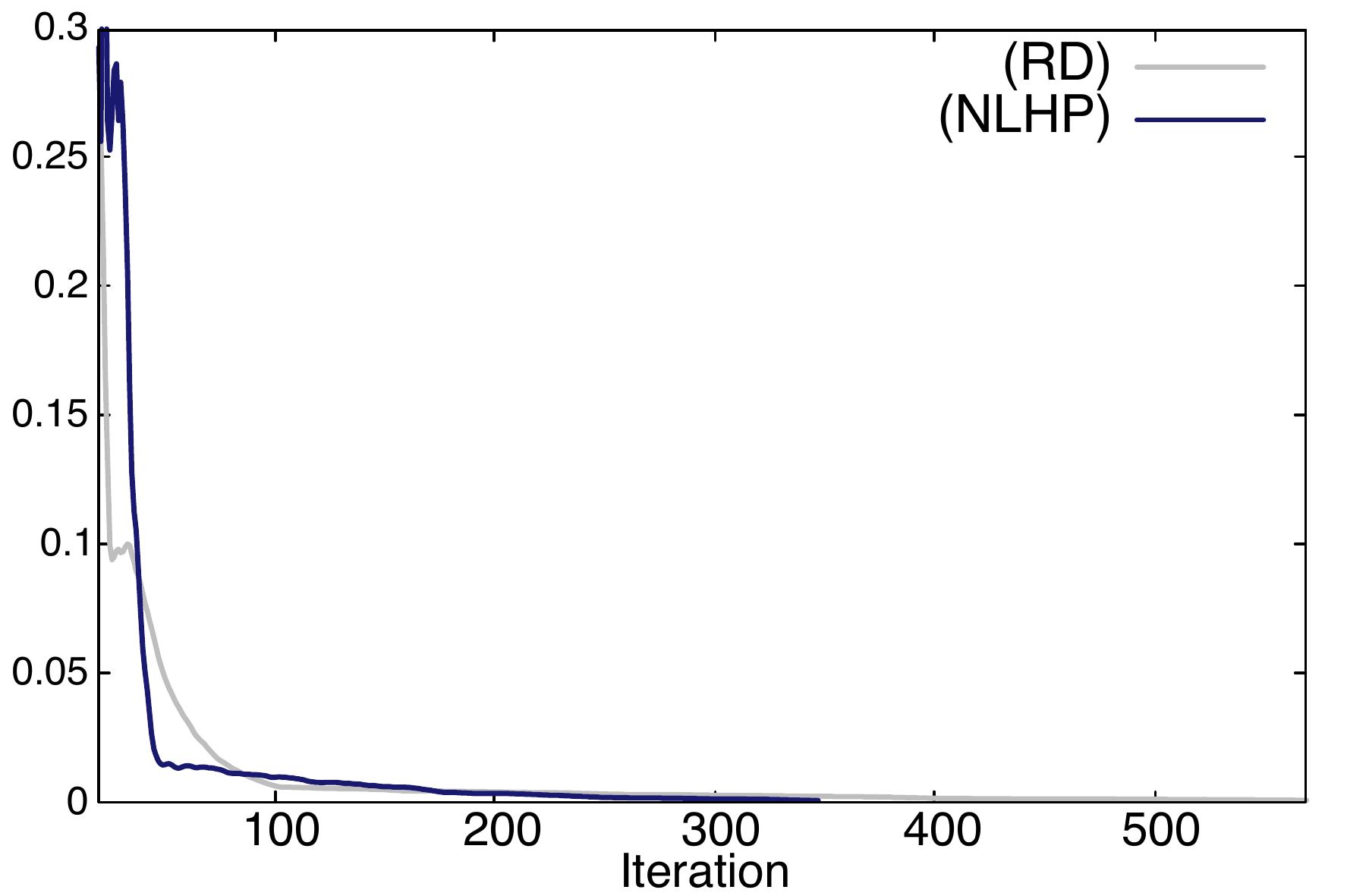}
        \subcaption{$\|\phi_{n+1}-\phi_n\|_{L^{\infty}(D)}$}
        \label{3dca-2}
      \end{minipage} &
      \end{tabular}
       \caption{ Objective functional and convergence condition for \S \ref{S:ca}-(iv).}
    \label{3dca}
  \end{figure*}

\vspace{2mm}
\noindent{\bf Case (iv) (Three-dimensional domain).\,} 
Let us finally consider the corresponding three-dimensional case. 
Here we set $(n_t,h_{\rm max})=(187280,0.0422
)$ and $(\tau, G_{\rm max})=(5.0\times 10^{-4}, 0.3)$. 
Then Figures \ref{fig:3dca} and \ref{3dca} are obtained.
Obviously, Figure \ref{3dca-2} shows that (NLHP) converges faster than (RD). In particular, at Step\,300, (NLHP) is almost identical to the final configuration, and therefore, we see that the boundary structure in (NLHP) is moving faster than that in (RD). 


\subsection{Bridge model} \label{S:br}
As another boundary condition, we next consider the so-called \emph{bridge model} (see Figure \ref{ibr}) and show numerically that the same assertion in the previous subsection is obtained.
In this subsection, we set $(n_t,h_{\rm max})=(35000, 0.0141)$ and $(\tau, G_{\rm max})=(8.0\times 10^{-5},0.35)$.
Here and henceforth, the convergence criterion is set to $\emph{eps0pt}=1.0\times 10^{-2}$ in terms of practicality.

\begin{figure*}[htbp]
\hspace*{-5mm}
    \begin{tabular}{ccccc}
      \begin{minipage}[t]{0.2\hsize}
        \centering
        \includegraphics[keepaspectratio, scale=0.09]{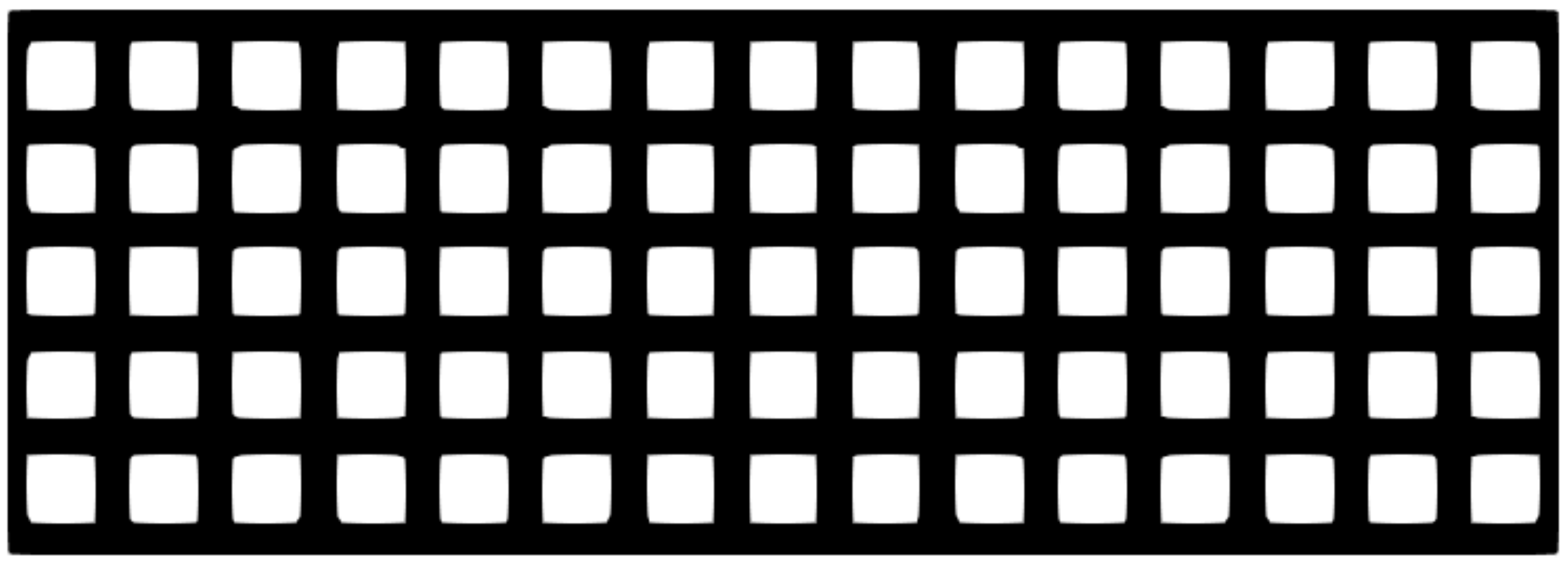}
        \subcaption{Step\,0}
        \label{b1-a}
      \end{minipage} 
      \begin{minipage}[t]{0.2\hsize}
        \centering
        \includegraphics[keepaspectratio, scale=0.09]{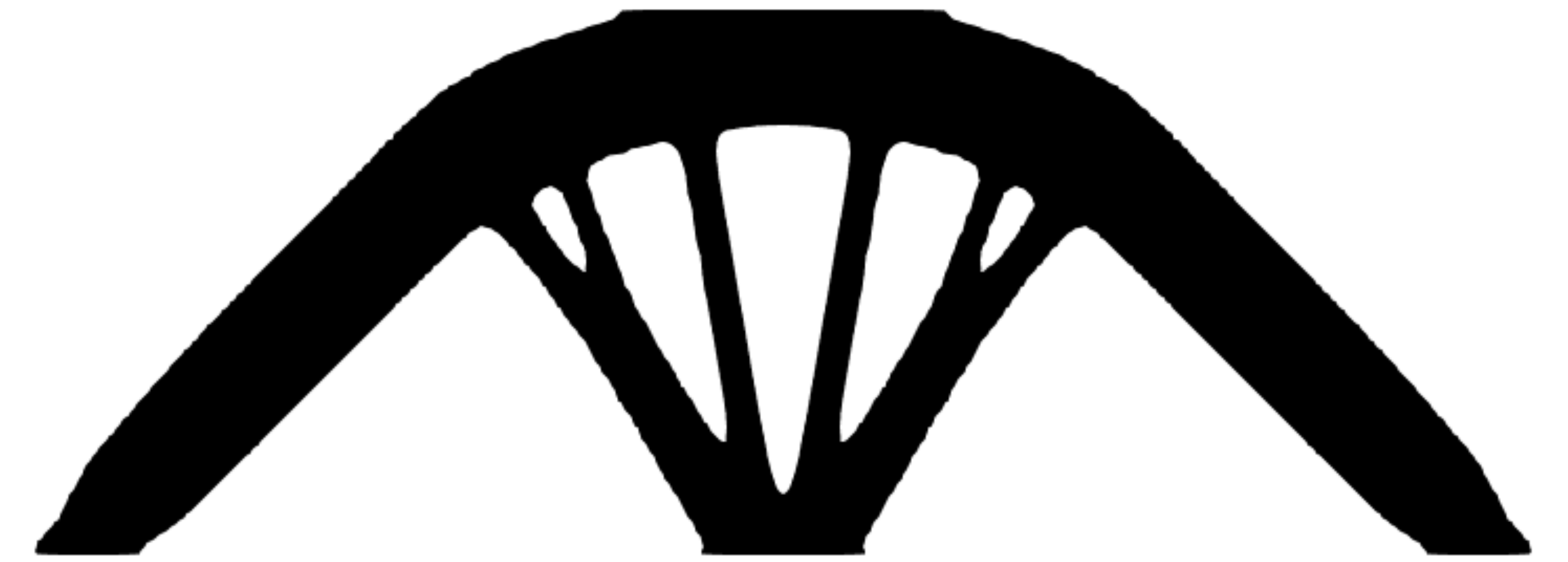}
        \subcaption{Step\,180}
        \label{b1-b}
      \end{minipage} 
      \begin{minipage}[t]{0.2\hsize}
        \centering
        \includegraphics[keepaspectratio, scale=0.09]{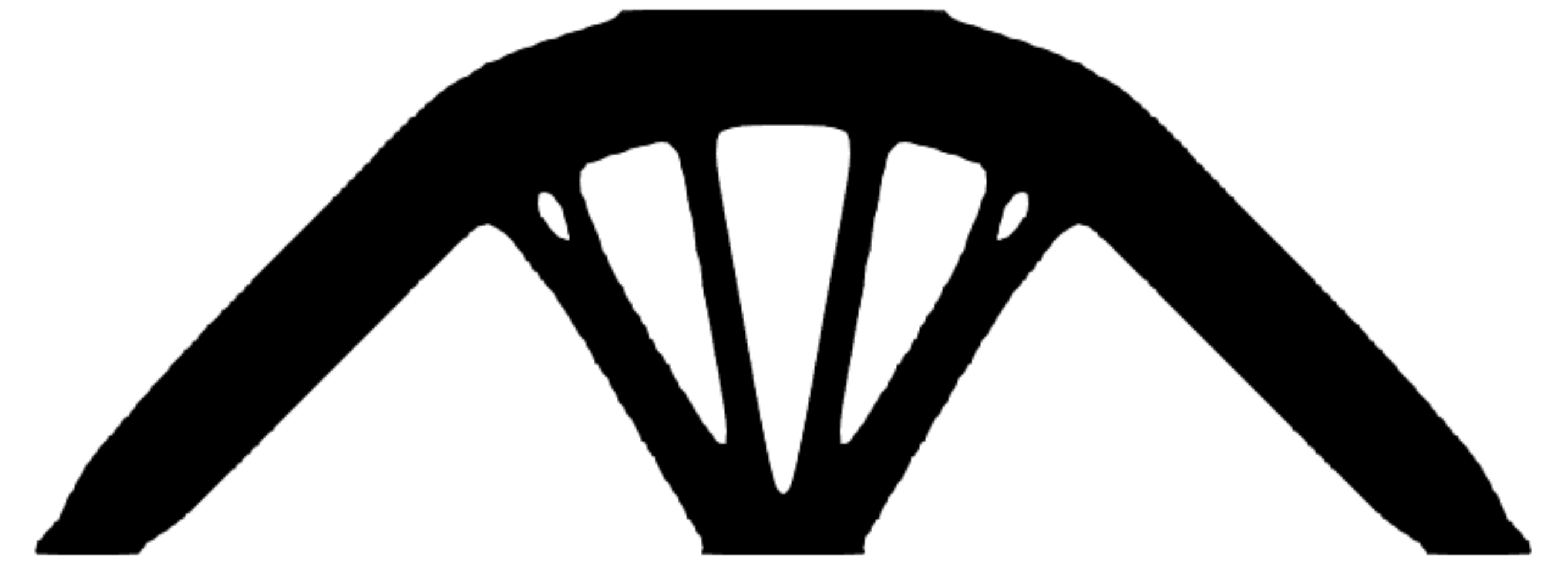}
        \subcaption{Step\,200}
        \label{b1-c}
      \end{minipage} 
         \begin{minipage}[t]{0.2\hsize}
        \centering
        \includegraphics[keepaspectratio, scale=0.09]{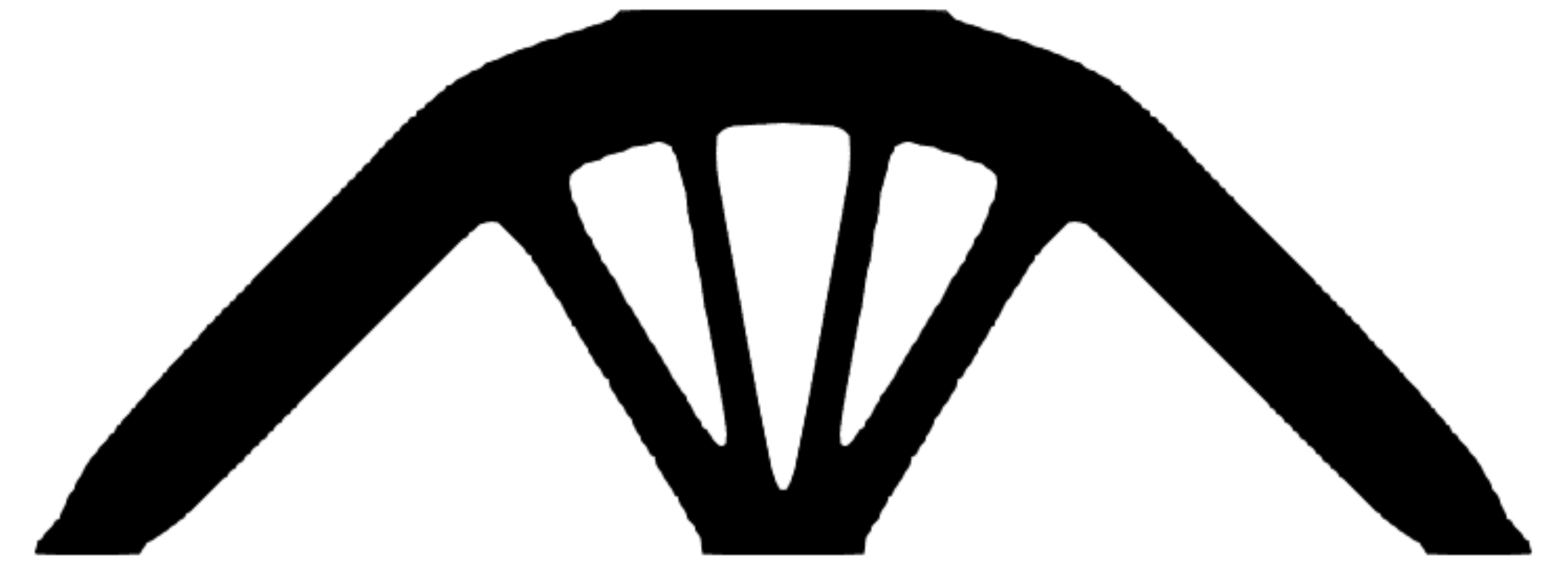}
        \subcaption{Step\,220}
        \label{b1-d}
      \end{minipage} 
                 \begin{minipage}[t]{0.2\hsize}
        \centering
        \includegraphics[keepaspectratio, scale=0.09]{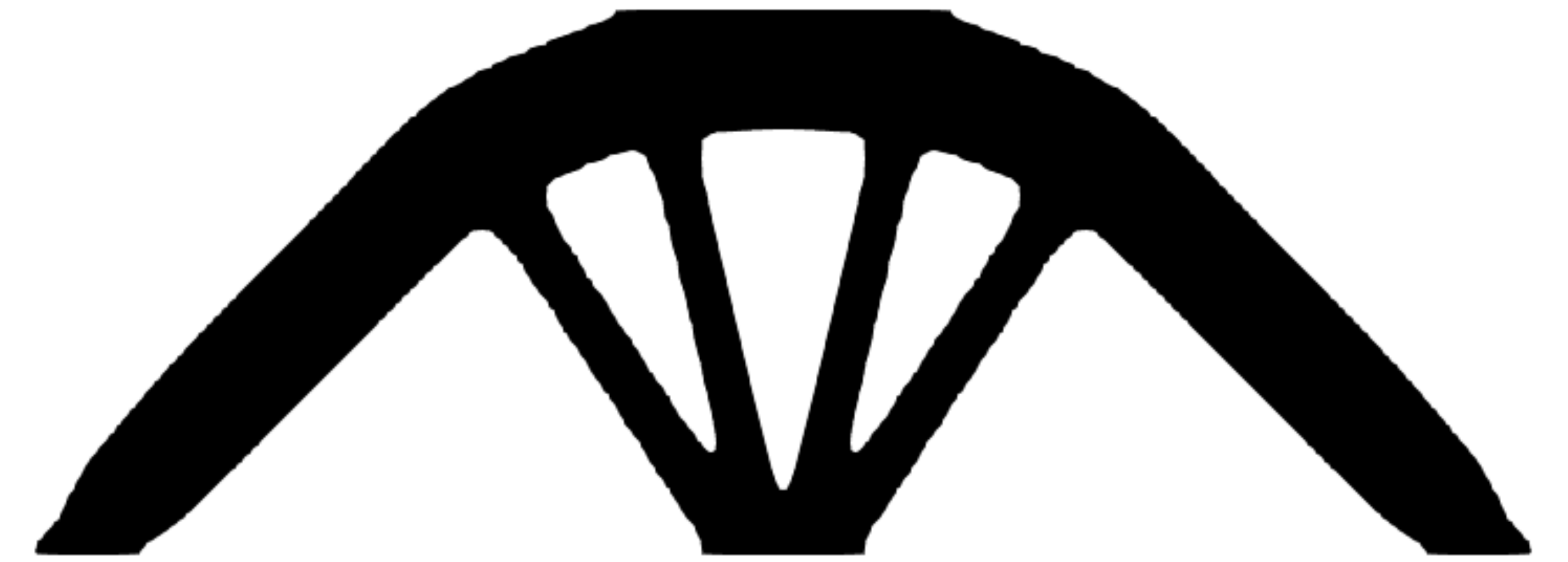} 
       \subcaption{Step\,591$^{\#}$}
       \label{b1-e}
      \end{minipage} 
      \\
    \begin{minipage}[t]{0.2\hsize}
        \centering
        \includegraphics[keepaspectratio, scale=0.09]{mbb20.pdf}
        \subcaption{Step\,0}
        \label{b1-f}
      \end{minipage} 
      \begin{minipage}[t]{0.2\hsize}
        \centering
        \includegraphics[keepaspectratio, scale=0.09]{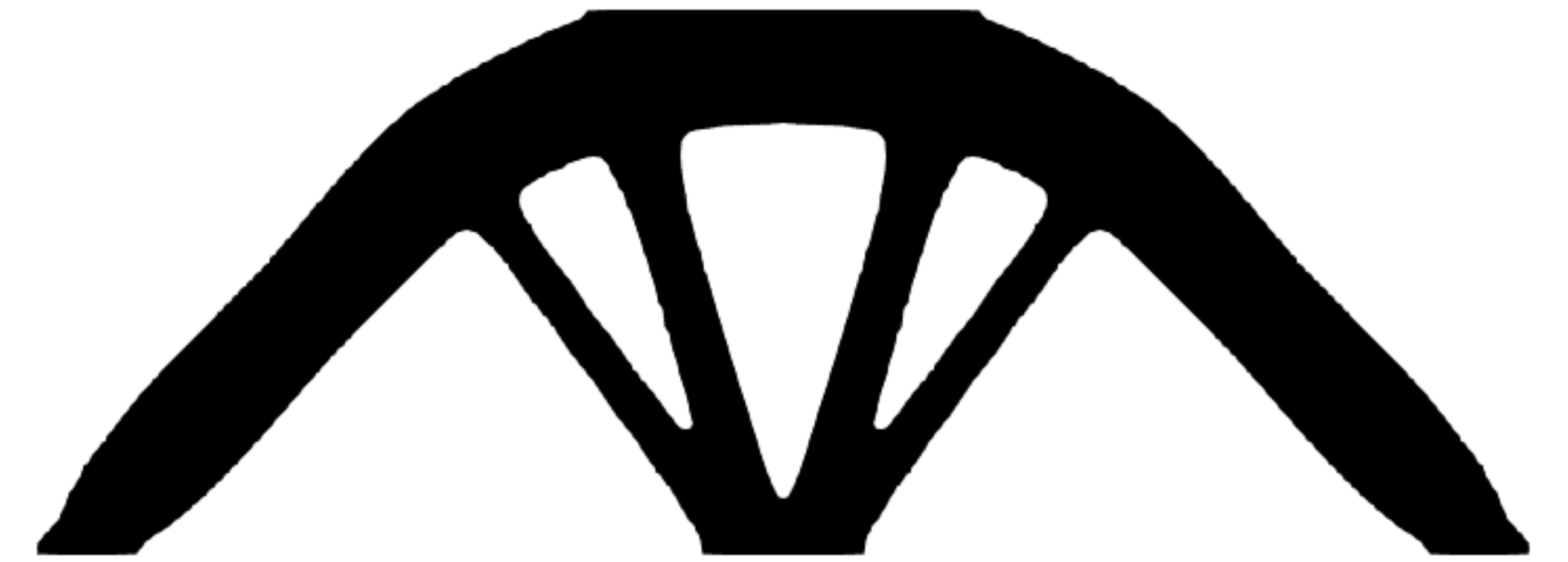}
        \subcaption{Step\,180}
        \label{b1-g}
      \end{minipage} 
      \begin{minipage}[t]{0.2\hsize}
        \centering
        \includegraphics[keepaspectratio, scale=0.09]{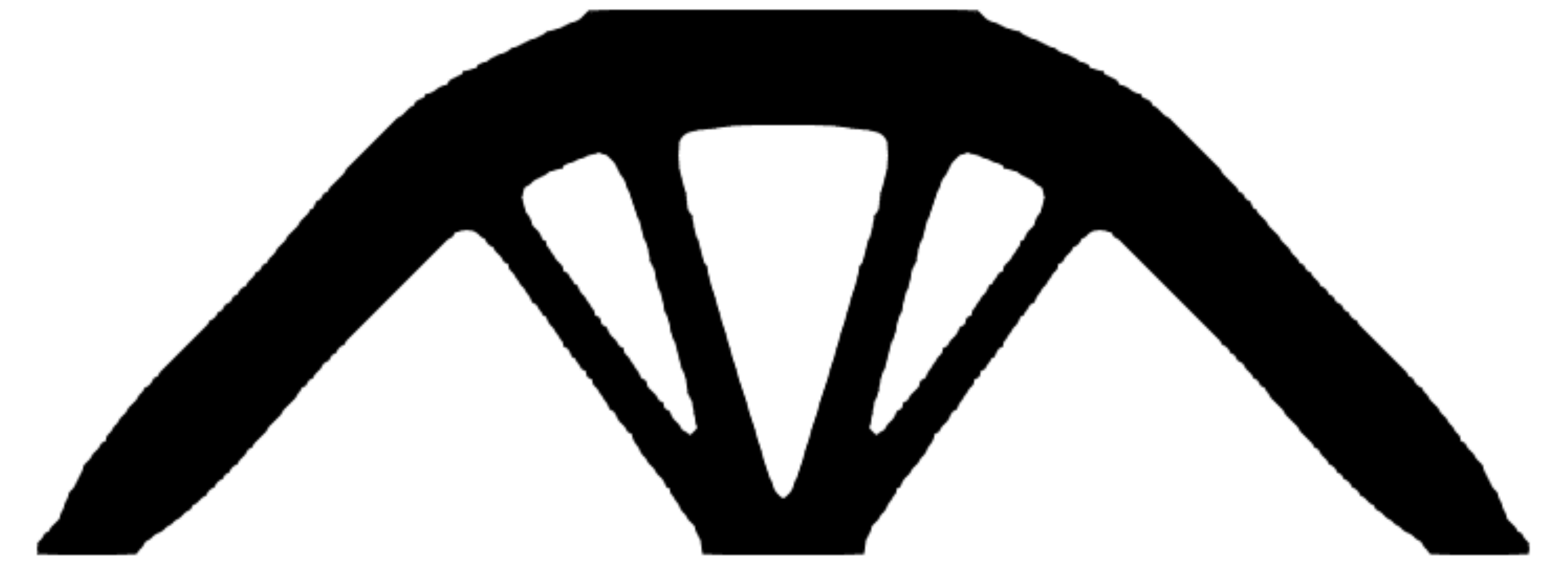}
        \subcaption{Step\,200}
        \label{b1-h}
      \end{minipage} 
       \begin{minipage}[t]{0.2\hsize}
        \centering
        \includegraphics[keepaspectratio, scale=0.09]{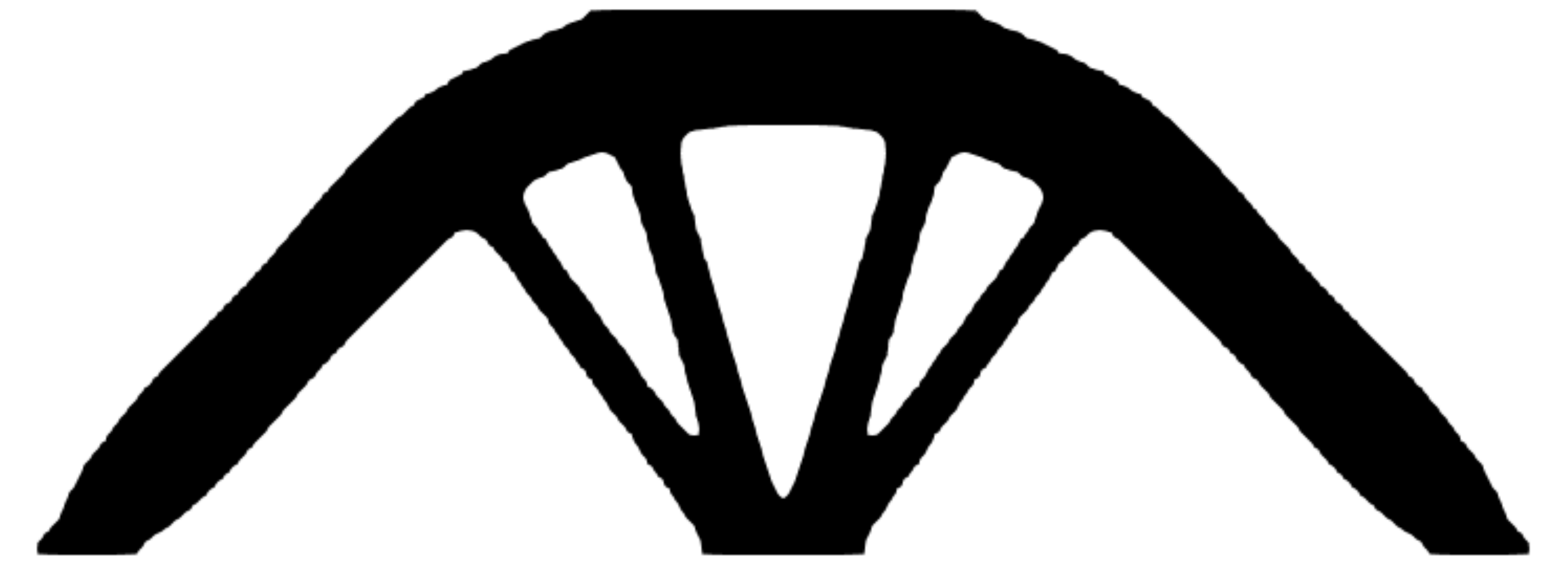}
        \subcaption{Step\,220}
        \label{b1-i}
      \end{minipage}
                 \begin{minipage}[t]{0.2\hsize}
        \centering
        \includegraphics[keepaspectratio, scale=0.09]{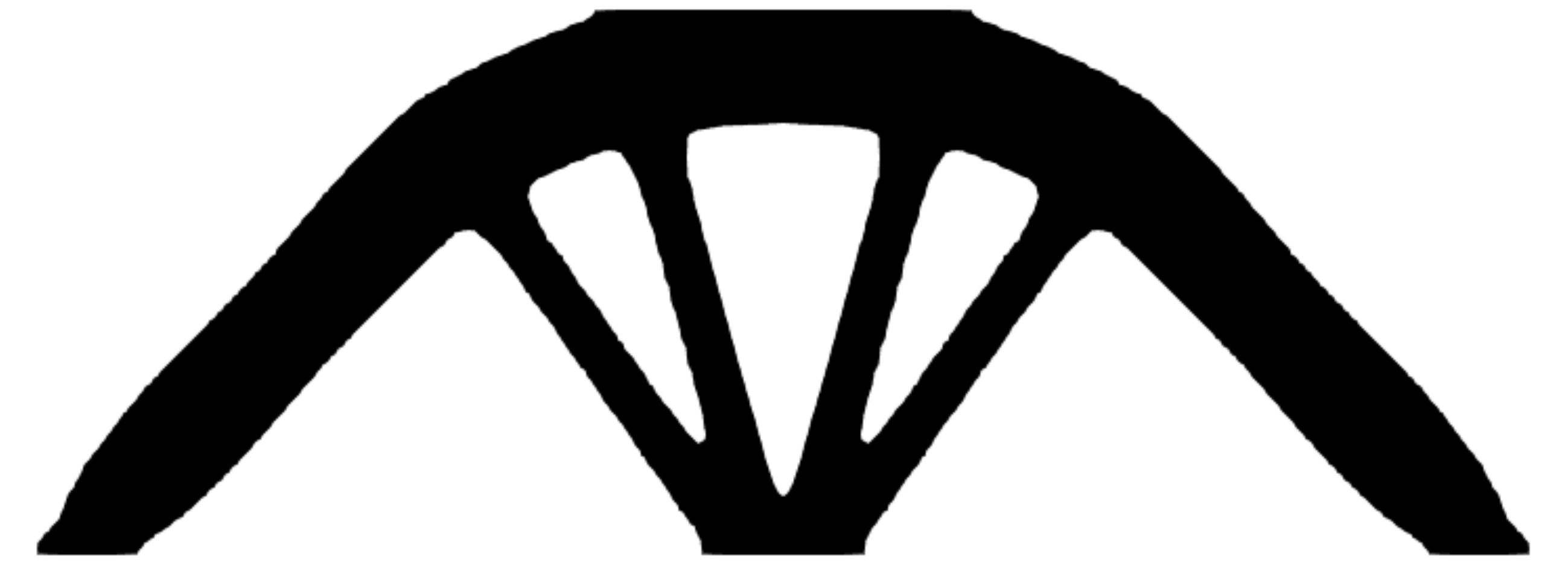}
        \subcaption{Step\,291$^{\#}$}
        \label{b1-j}
      \end{minipage}  
    \end{tabular}
     \caption{ Configuration $\Omega_{\phi_n}\subset D$ for the case where the initial configuration is the periodically perforated domain. Figures (a)--(e) and (f)--(j) represent the results of (RD) and (NLHP), respectively.     
The symbol $^{\#}$ implies the final step.}
     \label{fig:b1}
  \end{figure*}

\begin{figure*}[htbp]
    \begin{tabular}{ccc}
      \hspace*{-5mm} 
      \begin{minipage}[t]{0.49\hsize}
        \centering
        \includegraphics[keepaspectratio, scale=0.33]{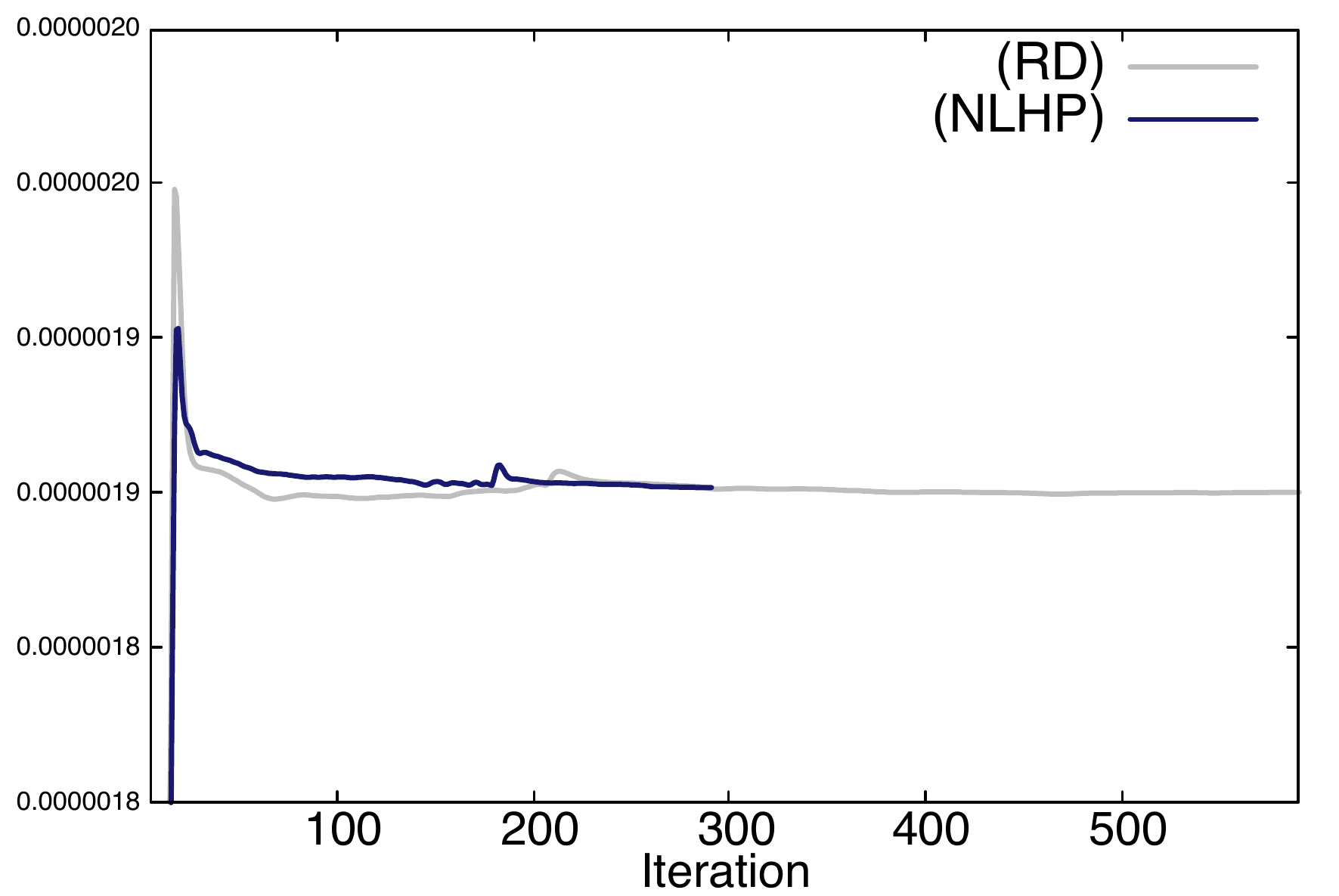}
        \subcaption{$F(\phi_n)$}
        \label{b1-1}
      \end{minipage} 
      \begin{minipage}[t]{0.49\hsize}
        \centering
        \includegraphics[keepaspectratio, scale=0.33]{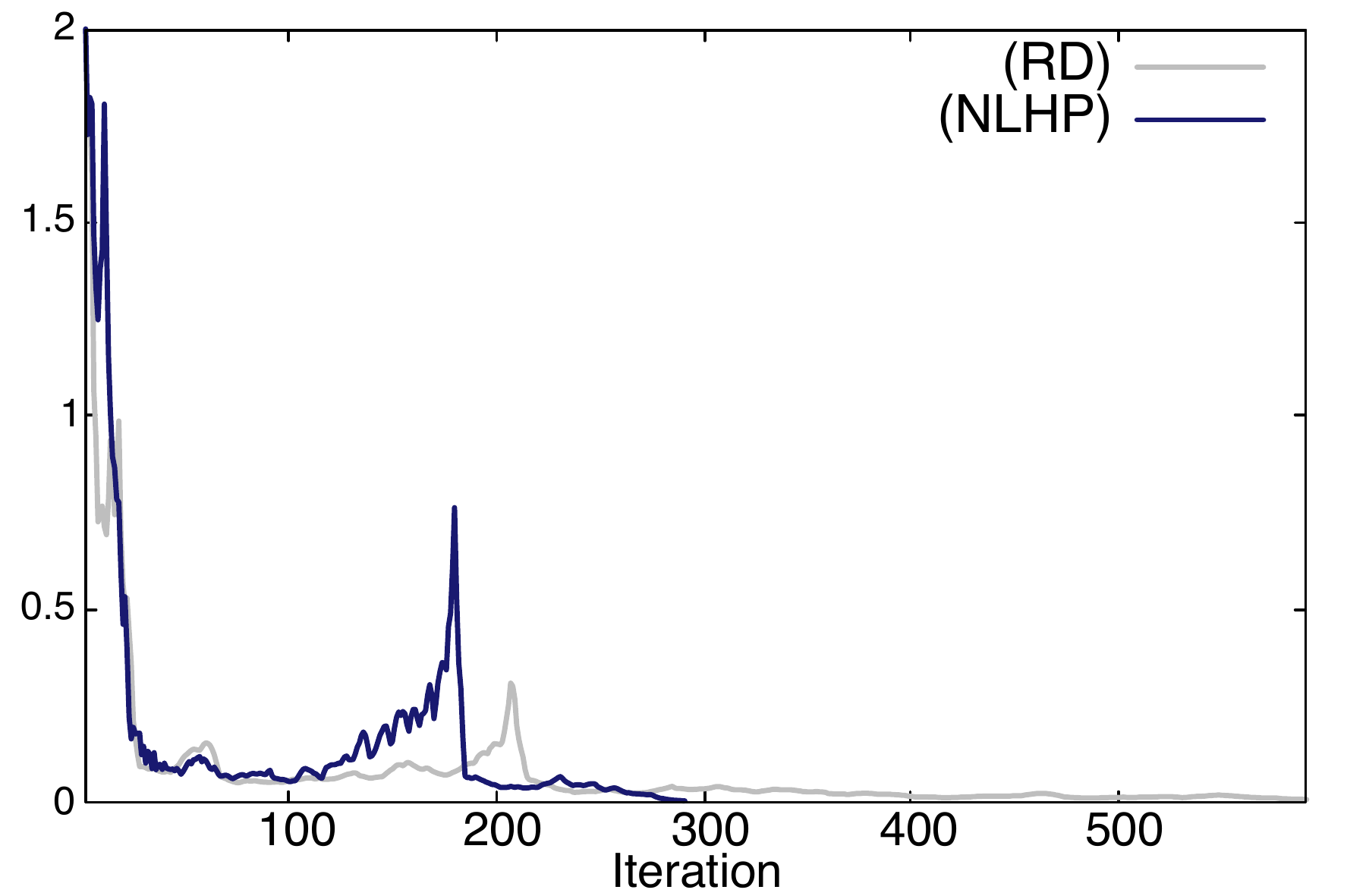}
        \subcaption{$\|\phi_{n+1}-\phi_n\|_{L^{\infty}(D)}$}
        \label{b1-2}
      \end{minipage} &
      \end{tabular}
       \caption{ Objective functional and convergence condition for \S \ref{S:br}-(i).}
    \label{b1}
  \end{figure*}

\vspace{2mm}
\noindent{\bf Case (i) (Periodically perforated domain).\,} 
As in \S \ref{S:ca}-(i), one takes the initial configuration as the periodically perforated domain. Then Figures \ref{fig:b1} and \ref{b1} ensure the assertion in this study; indeed, at Step\,180 (see Figures \ref{b1-b} and \ref{b1-g}), the topology of $\Omega_{\phi_n}\subset D$ in (NLHP) can be optimized, and moreover, (NLHP) satisfies the convergence condition in at least half the number of iterations for (RD).

\begin{figure*}[htbp]
\hspace*{-5mm}
    \begin{tabular}{ccccc}
      \begin{minipage}[t]{0.2\hsize}
        \centering
        \includegraphics[keepaspectratio, scale=0.09]{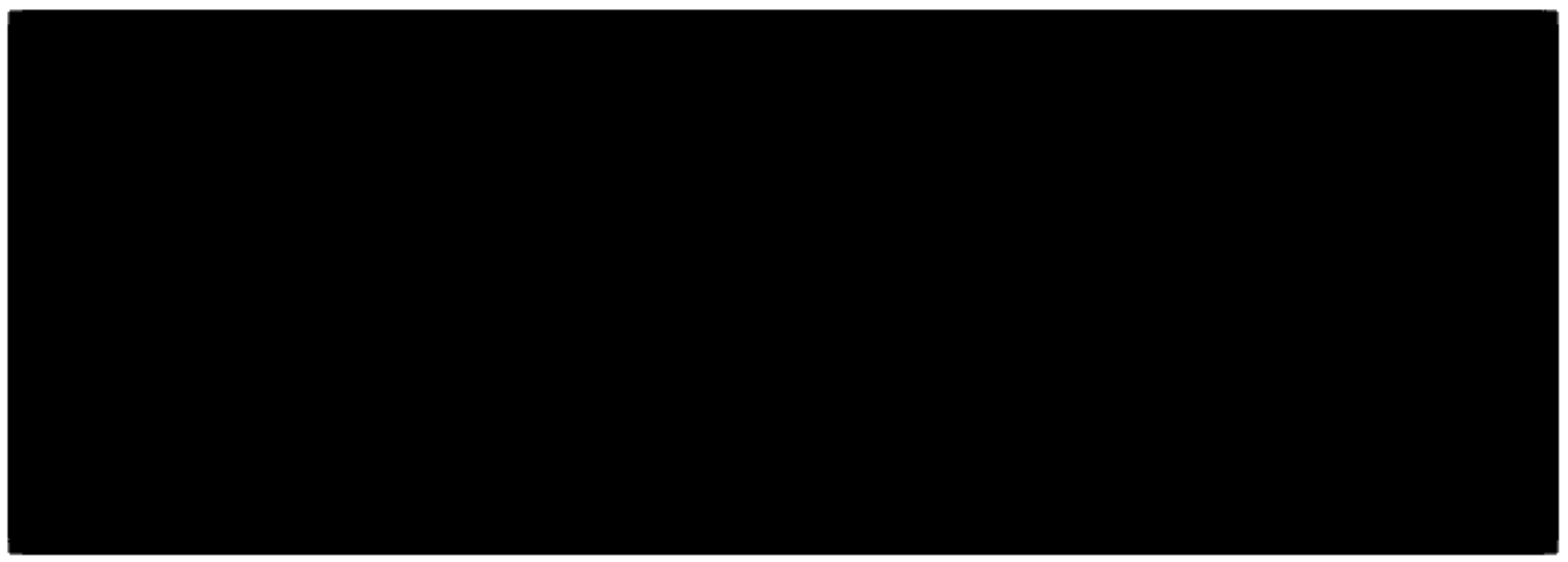}
        \subcaption{Step\,0}
        \label{b2-a}
      \end{minipage} 
      \begin{minipage}[t]{0.2\hsize}
        \centering
        \includegraphics[keepaspectratio, scale=0.09]{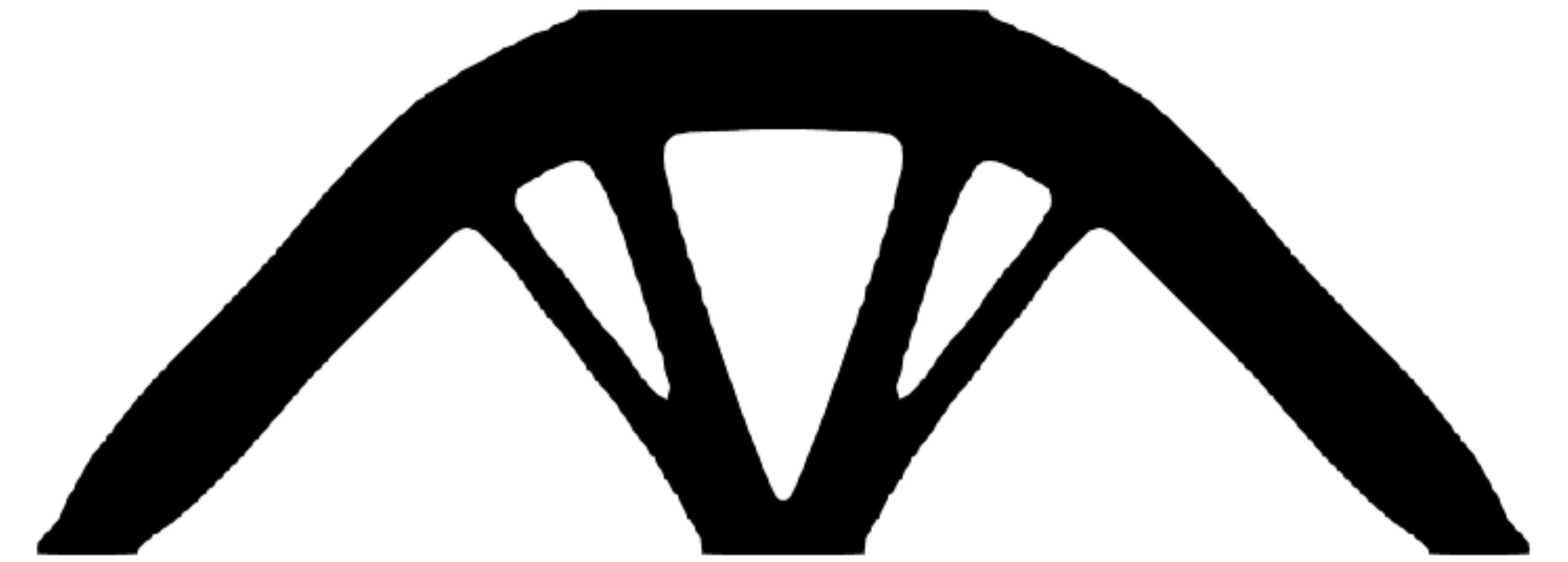}
        \subcaption{Step\,80}
        \label{b2-b}
      \end{minipage} 
      \begin{minipage}[t]{0.2\hsize}
        \centering
        \includegraphics[keepaspectratio, scale=0.09]{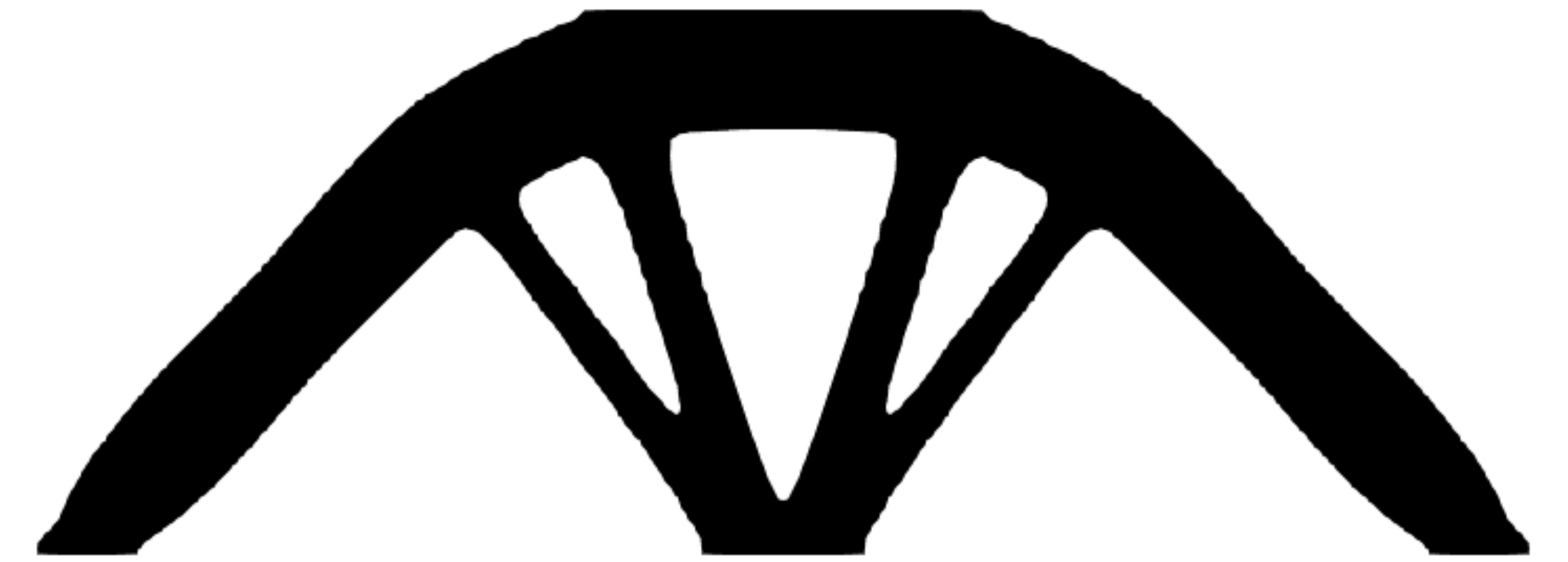}
        \subcaption{Step\,120}
        \label{b2-c}
      \end{minipage} 
         \begin{minipage}[t]{0.2\hsize}
        \centering
        \includegraphics[keepaspectratio, scale=0.09]{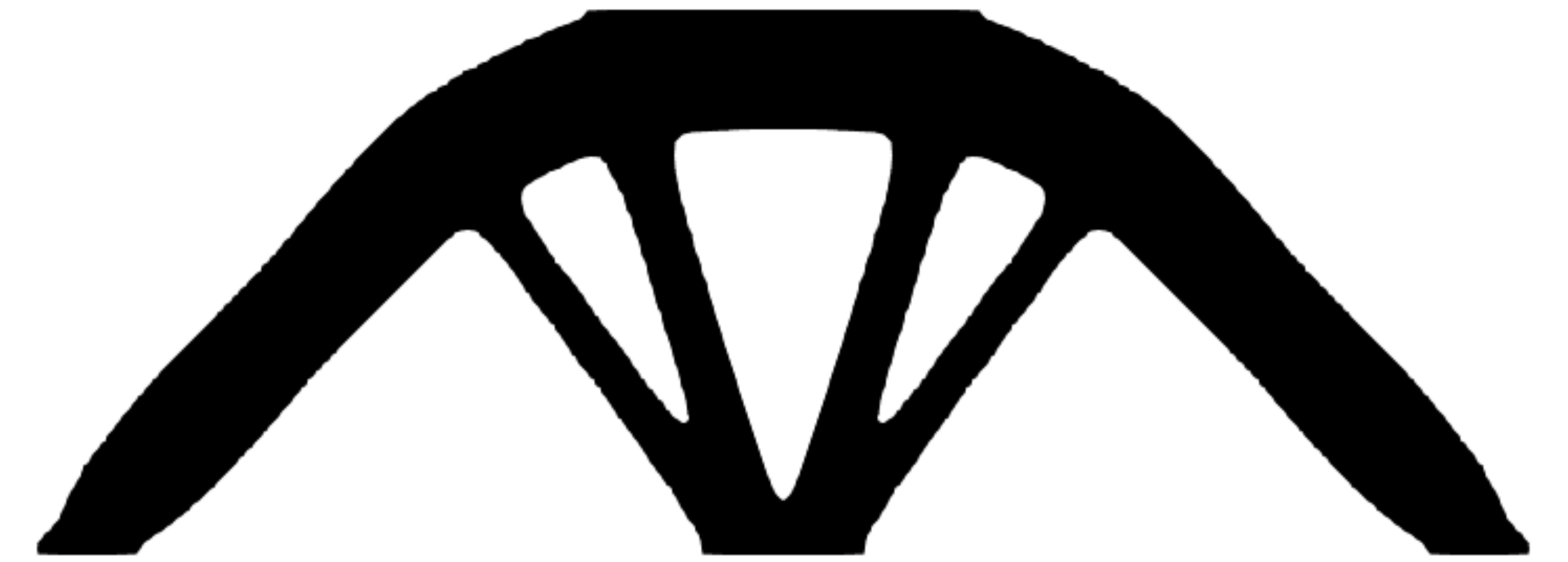}
        \subcaption{Step\,160}
        \label{b2-d}
      \end{minipage} 
                 \begin{minipage}[t]{0.2\hsize}
        \centering
        \includegraphics[keepaspectratio, scale=0.09]{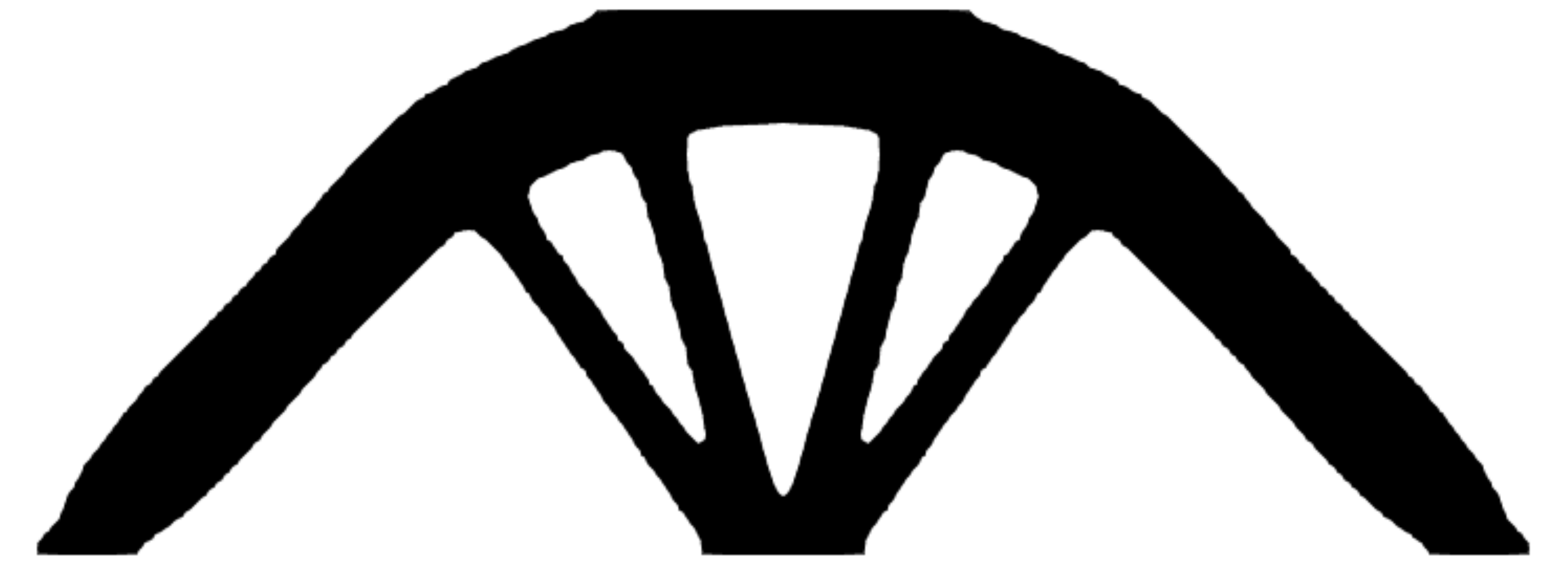}
        \subcaption{Step\,416$^{\#}$}
        \label{b2-e}
      \end{minipage} 
      \\
    \begin{minipage}[t]{0.2\hsize}
        \centering
        \includegraphics[keepaspectratio, scale=0.09]{mbb2w0.pdf}
        \subcaption{Step\,0}
        \label{b2-f}
      \end{minipage} 
      \begin{minipage}[t]{0.2\hsize}
        \centering
        \includegraphics[keepaspectratio, scale=0.09]{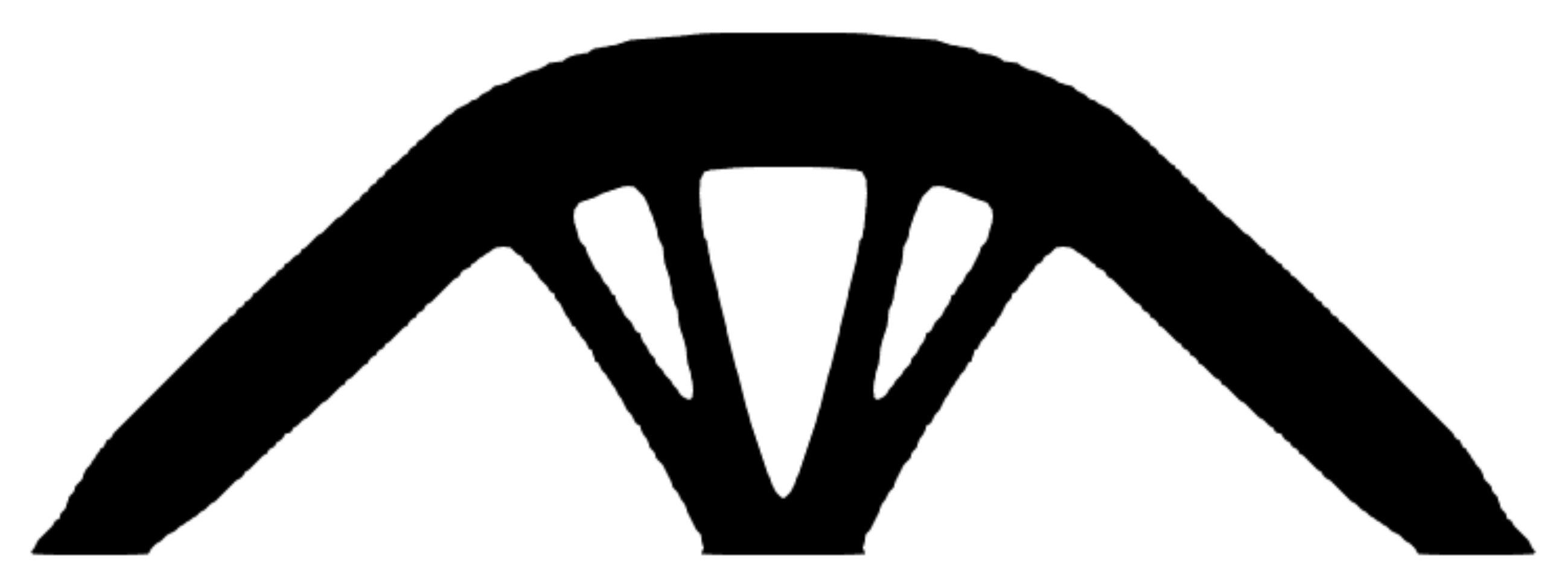}
        \subcaption{Step\,80}
        \label{b2-g}
      \end{minipage} 
      \begin{minipage}[t]{0.2\hsize}
        \centering
        \includegraphics[keepaspectratio, scale=0.09]{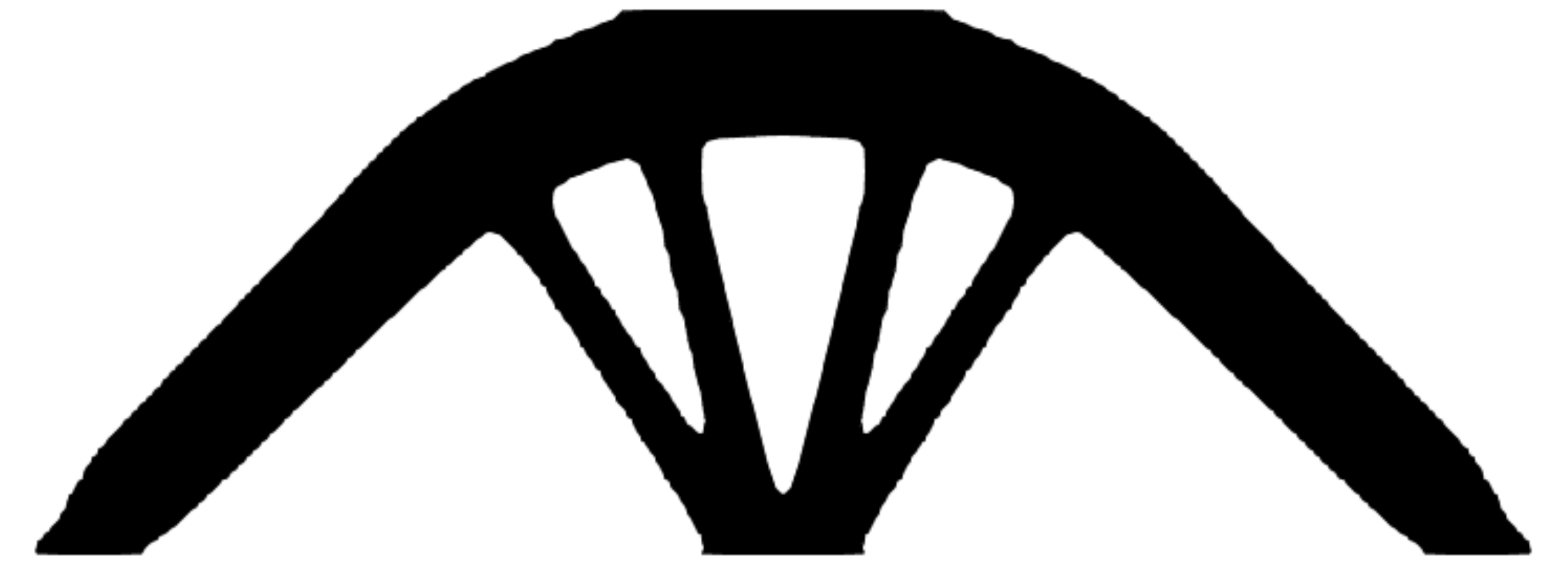}
        \subcaption{Step\,120}
        \label{3-h}
      \end{minipage} 
       \begin{minipage}[t]{0.2\hsize}
        \centering
        \includegraphics[keepaspectratio, scale=0.09]{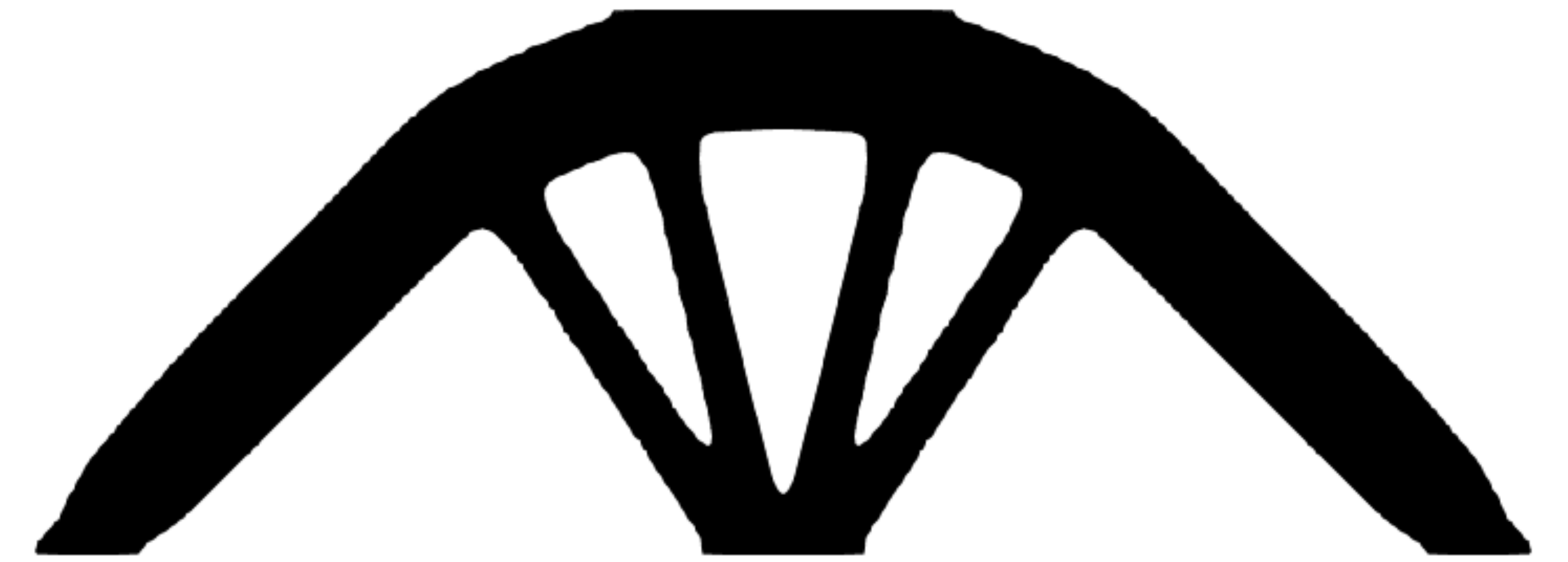}
        \subcaption{Step\,160}
        \label{b2-i}
      \end{minipage}
                 \begin{minipage}[t]{0.2\hsize}
        \centering
        \includegraphics[keepaspectratio, scale=0.09]{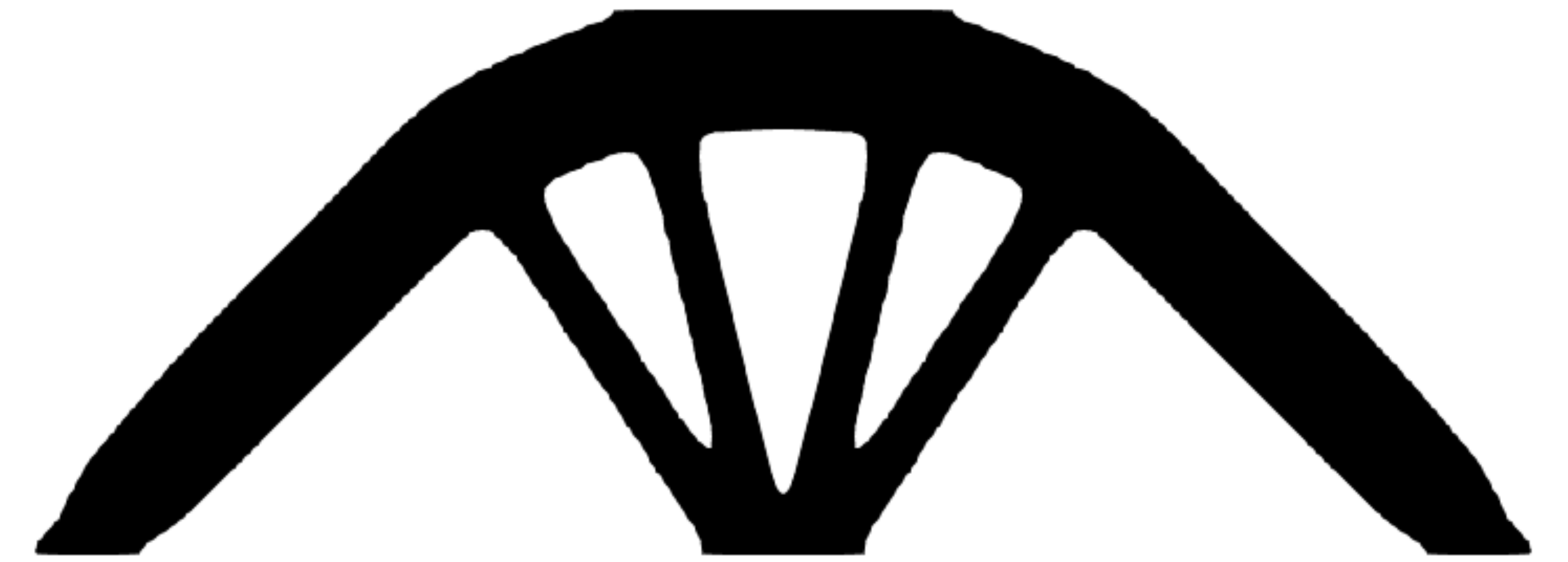}
        \subcaption{Step\,180$^{\#}$}
        \label{b2-j}
      \end{minipage}  
    \end{tabular}
     \caption{ Configuration $\Omega_{\phi_n}\subset D$ for the case where the initial configuration is the whole domain. 
     Figures (a)--(e) and (f)--(j) represent the results of (RD) and (NLHP), respectively.     
The symbol $^{\#}$ implies the final step. }
     \label{fig:b2}
  \end{figure*}

\begin{figure*}[htbp]
    \begin{tabular}{ccc}
      \hspace*{-5mm} 
      \begin{minipage}[t]{0.49\hsize}
        \centering
        \includegraphics[keepaspectratio, scale=0.33]{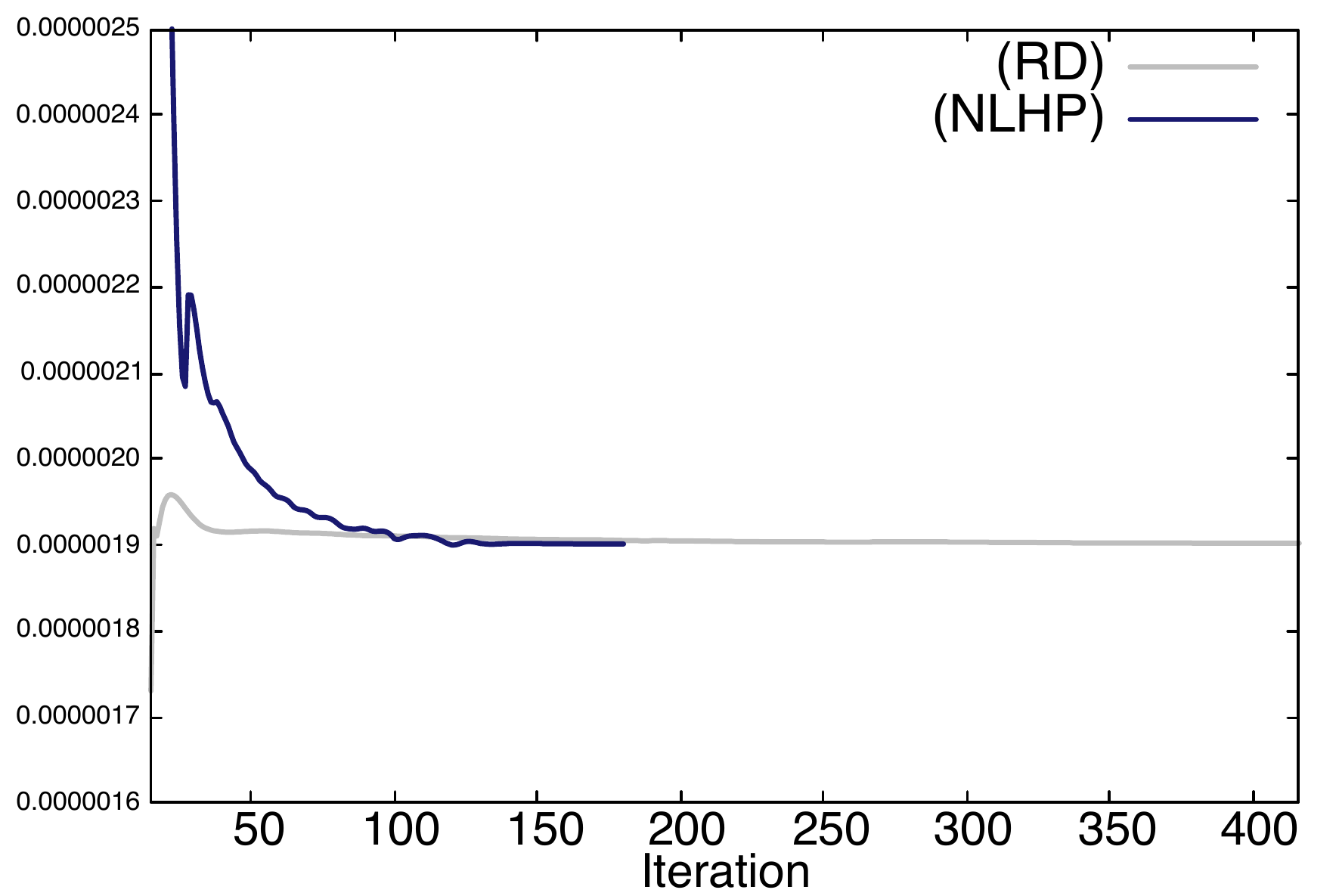}
        \subcaption{$F(\phi_n)$}
        \label{b2-1}
      \end{minipage} 
      \begin{minipage}[t]{0.49\hsize}
        \centering
        \includegraphics[keepaspectratio, scale=0.33]{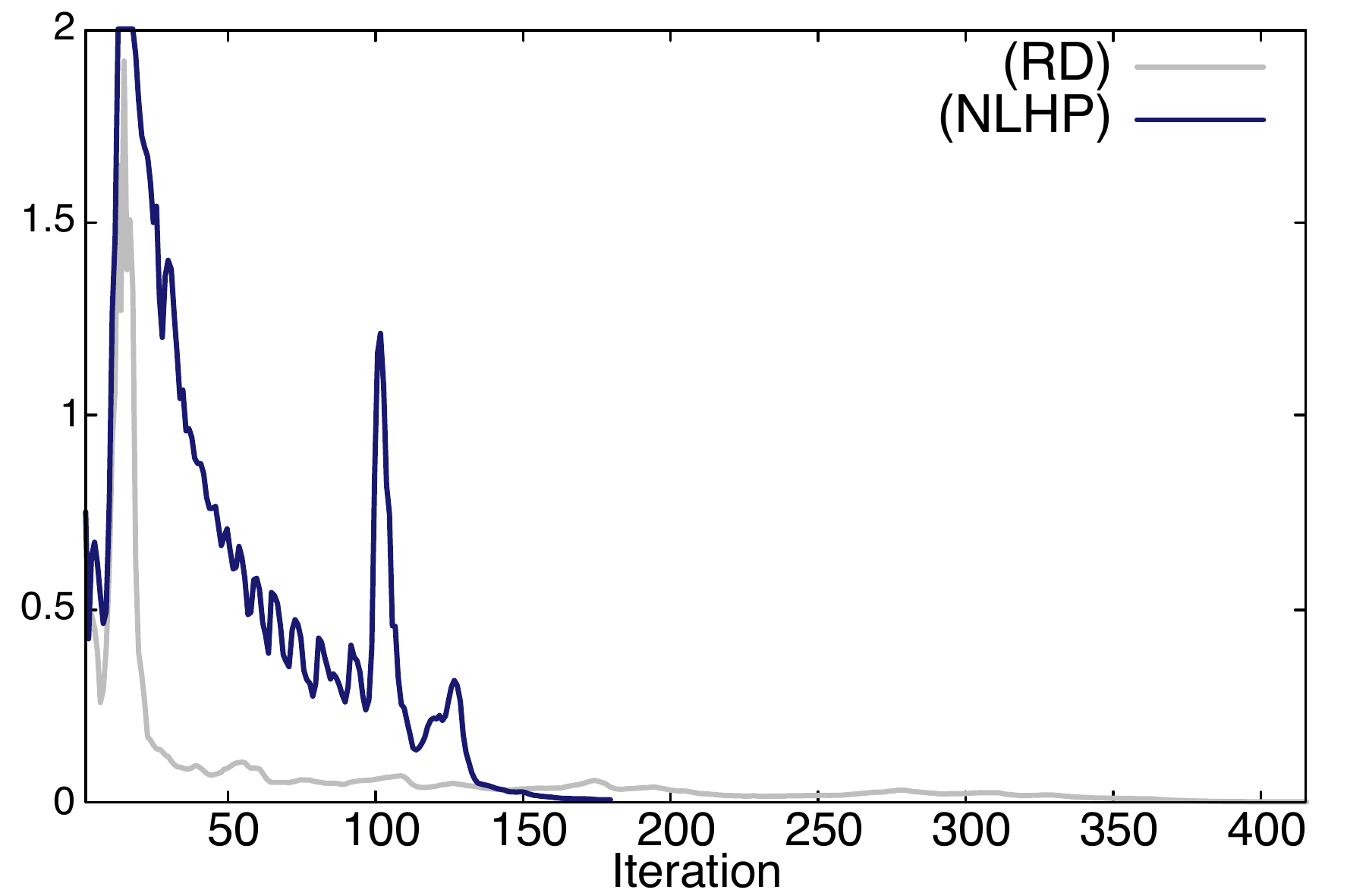}
        \subcaption{$\|\phi_{n+1}-\phi_n\|_{L^{\infty}(D)}$}
        \label{b2-2}
      \end{minipage} &
      \end{tabular}
       \caption{ Objective functional and convergence condition for \S \ref{S:br}-(ii).}
    \label{b2}
  \end{figure*}

\vspace{2mm}
\noindent{\bf Case (ii) (Whole domain).\,} 
We next choose the whole domain as the initial configuration. 
By comparing Figure \ref{b2-b} with \ref{b2-g}, the topology is optimized at the same speed, and there seems to be no difference in the speed for convergence; however, from Figure \ref{b2}, one can see that the boundary structure moves and converges faster in (NLHP) than in (RD). 
As in the previous case, it is noteworthy that more than twice the number of iterations can be improved.

\begin{figure*}[htbp]
\hspace*{-5mm}
    \begin{tabular}{ccccc}
      \begin{minipage}[t]{0.2\hsize}
        \centering
        \includegraphics[keepaspectratio, scale=0.09]{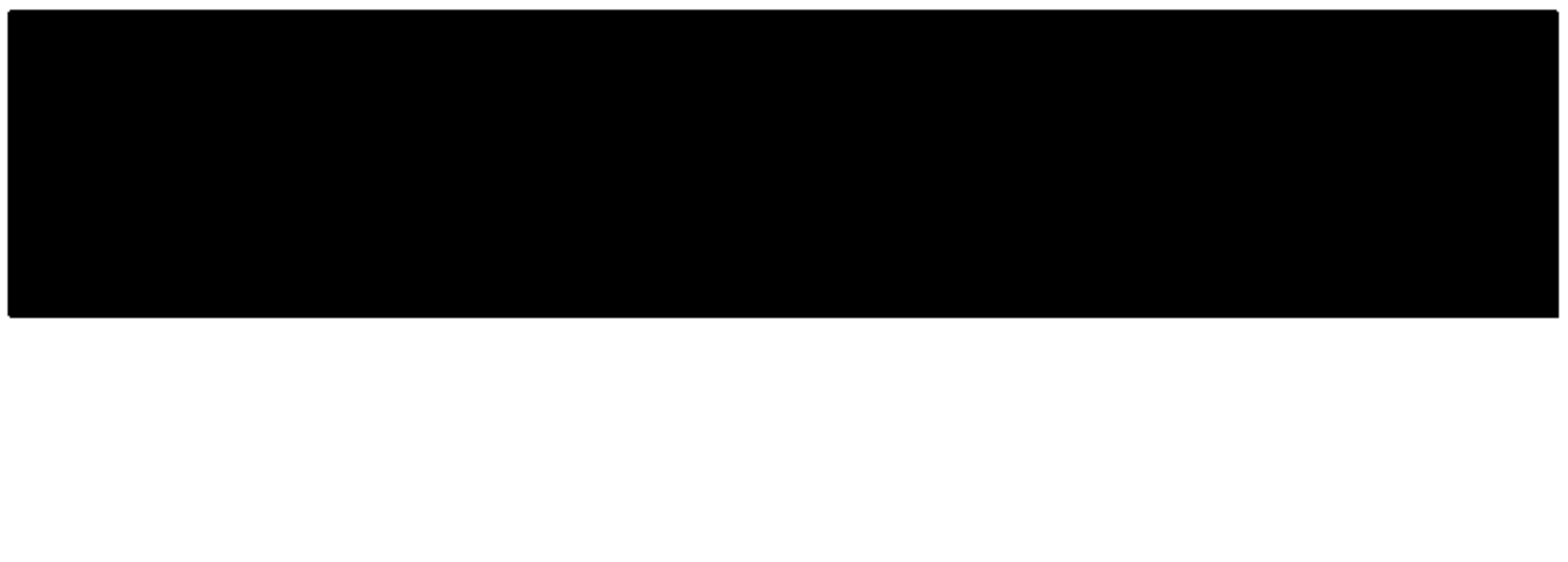}
        \subcaption{Step\,0}
        \label{b3-a}
      \end{minipage} 
      \begin{minipage}[t]{0.2\hsize}
        \centering
        \includegraphics[keepaspectratio, scale=0.09]{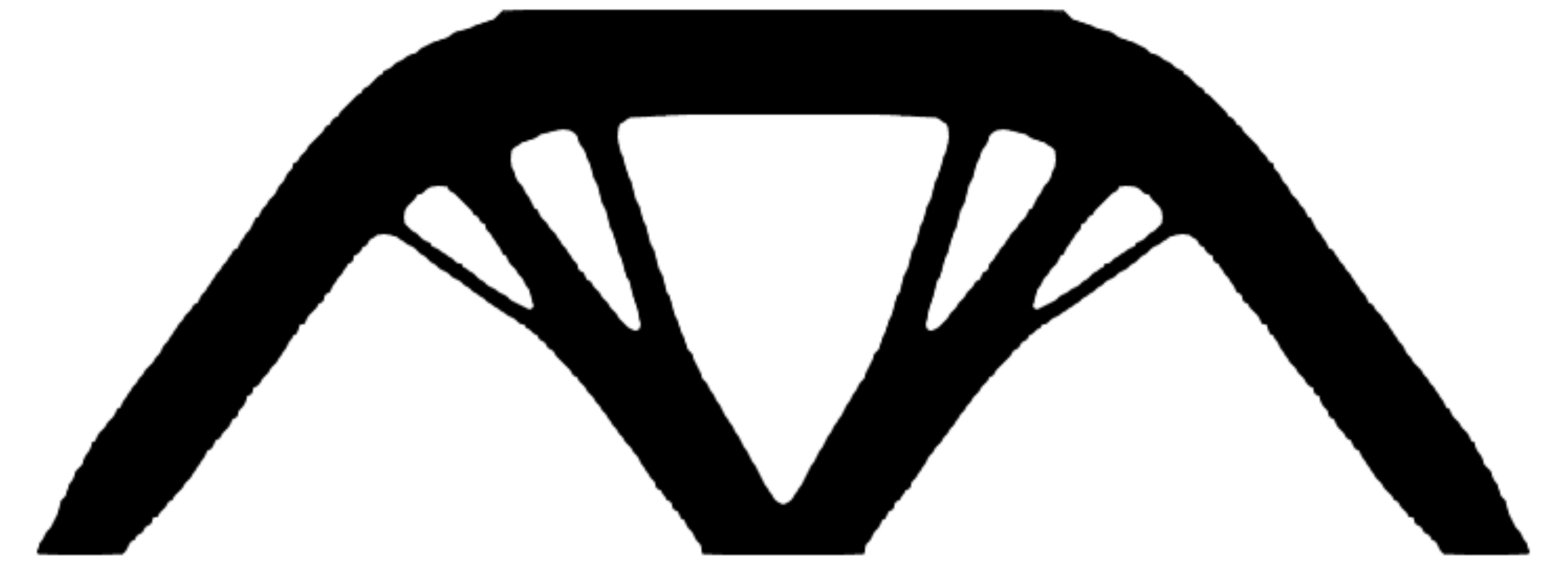}
        \subcaption{Step\,30}
        \label{b3-b}
      \end{minipage} 
      \begin{minipage}[t]{0.2\hsize}
        \centering
        \includegraphics[keepaspectratio, scale=0.09]{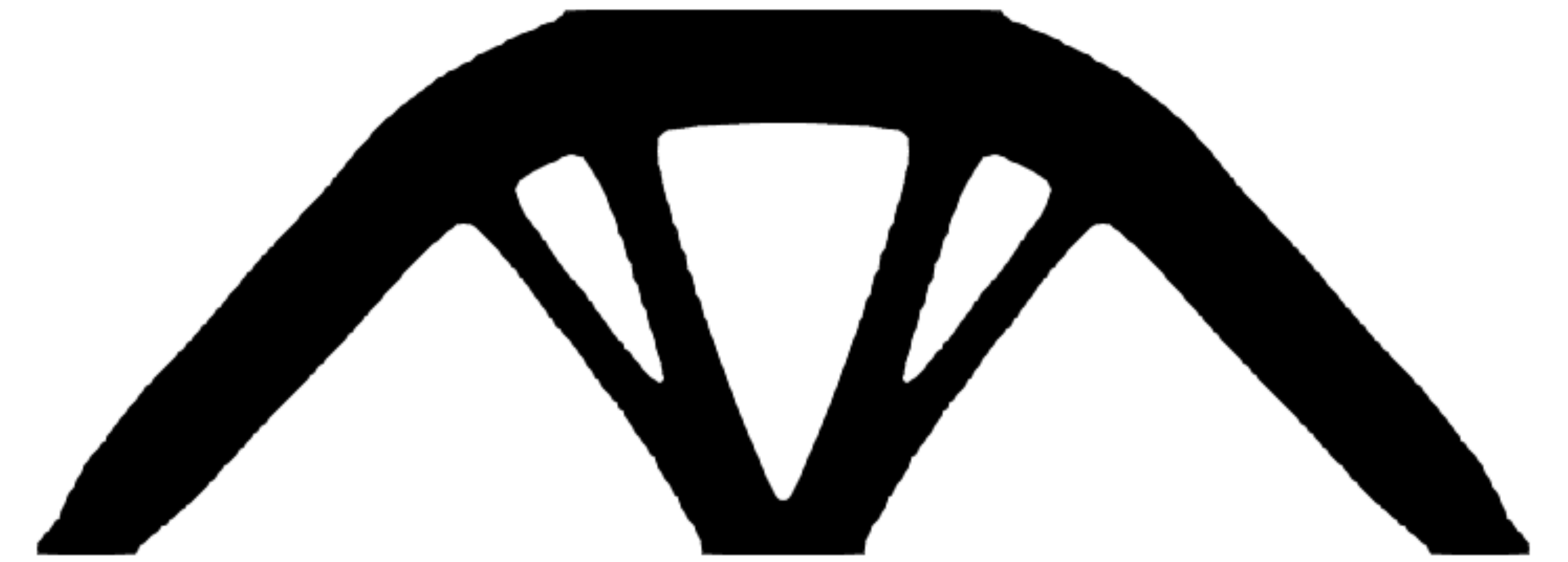}
        \subcaption{Step\,130}
        \label{b3-c}
      \end{minipage} 
         \begin{minipage}[t]{0.2\hsize}
        \centering
        \includegraphics[keepaspectratio, scale=0.09]{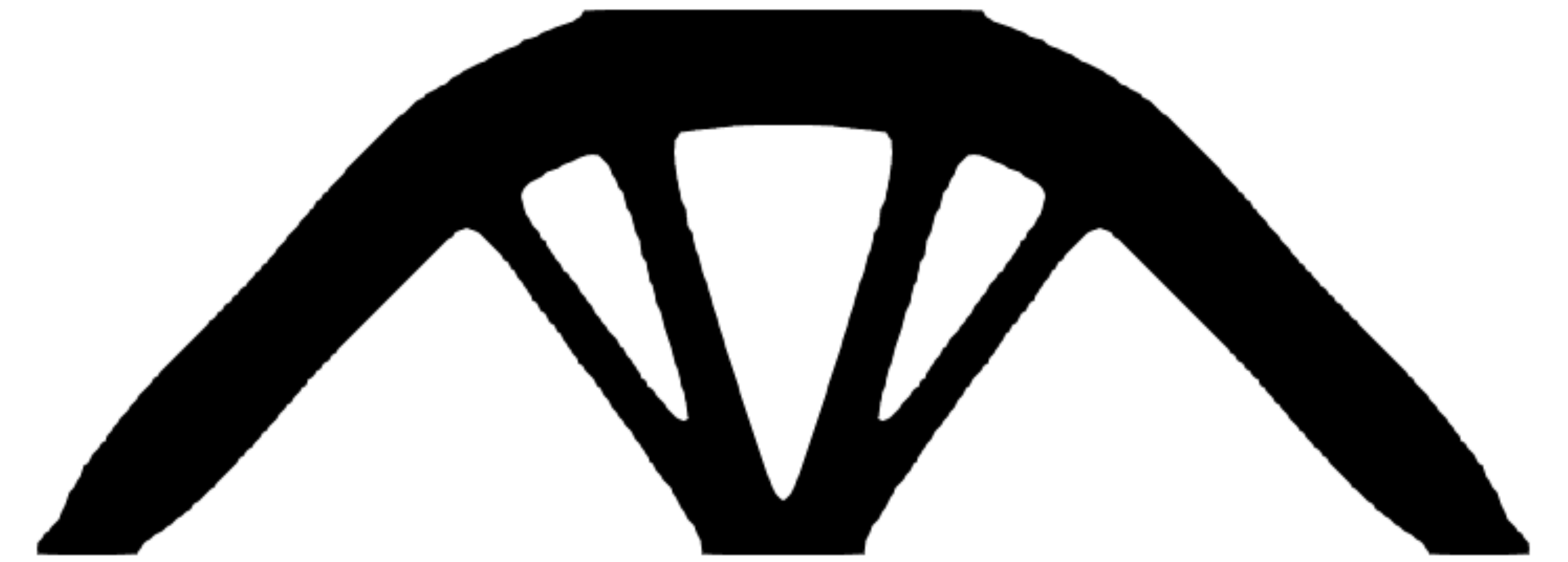}
        \subcaption{Step\,230}
        \label{b3-d}
      \end{minipage} 
                 \begin{minipage}[t]{0.2\hsize}
        \centering
        \includegraphics[keepaspectratio, scale=0.09]{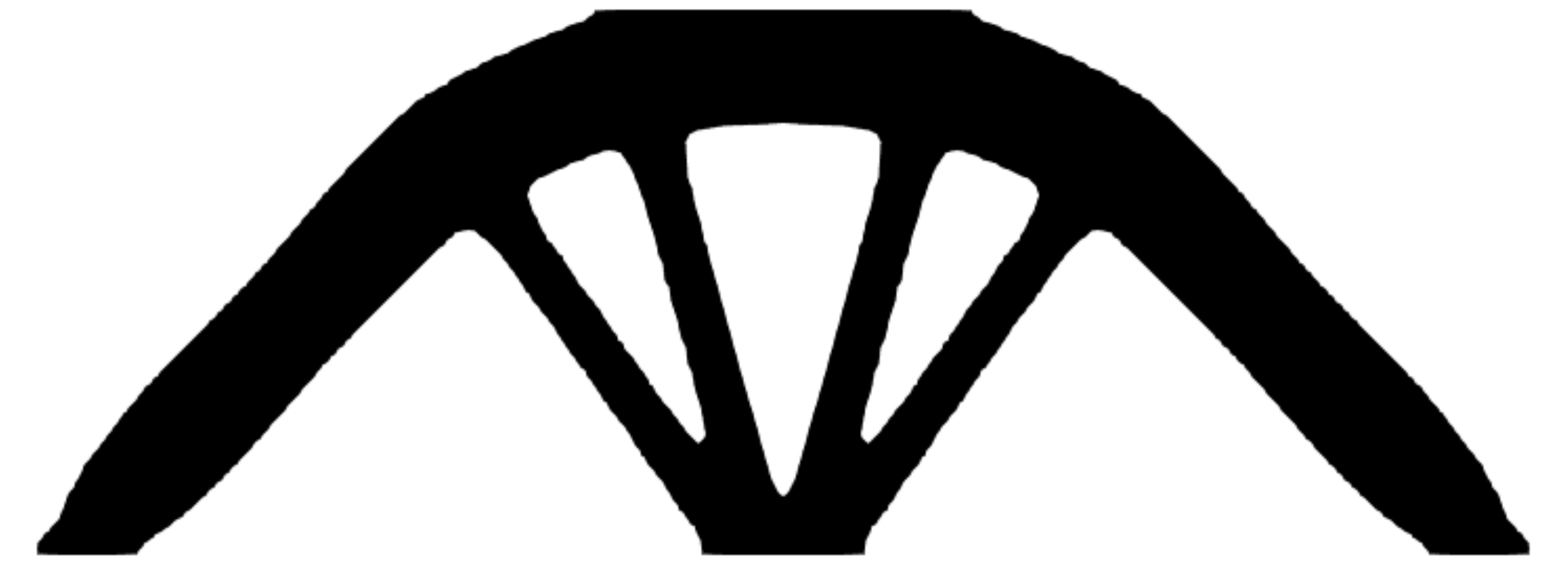}
        \subcaption{Step\,433$^{\#}$}
        \label{b3-e}
      \end{minipage} 
      \\
    \begin{minipage}[t]{0.2\hsize}
        \centering
        \includegraphics[keepaspectratio, scale=0.09]{mbb2u0.pdf}
        \subcaption{Step\,0}
        \label{b3-f}
      \end{minipage} 
      \begin{minipage}[t]{0.2\hsize}
        \centering
        \includegraphics[keepaspectratio, scale=0.09]{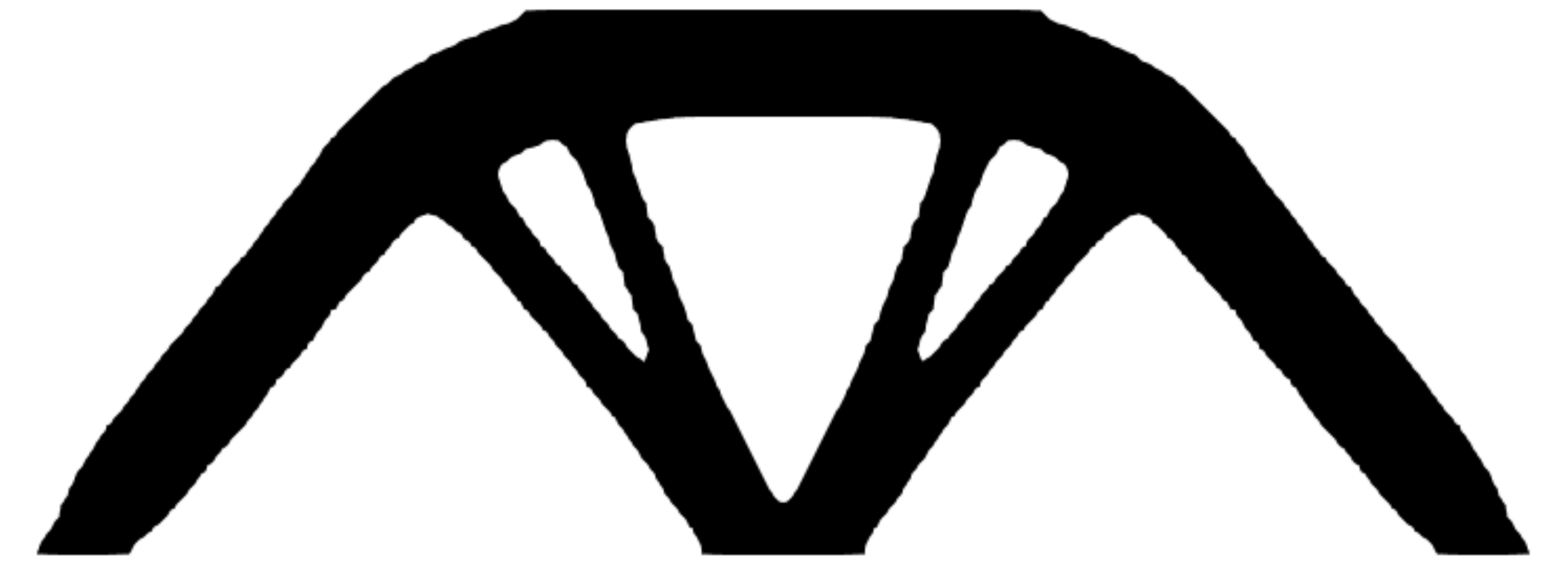}
        \subcaption{Step\,30}
        \label{b3-g}
      \end{minipage} 
      \begin{minipage}[t]{0.2\hsize}
        \centering
        \includegraphics[keepaspectratio, scale=0.09]{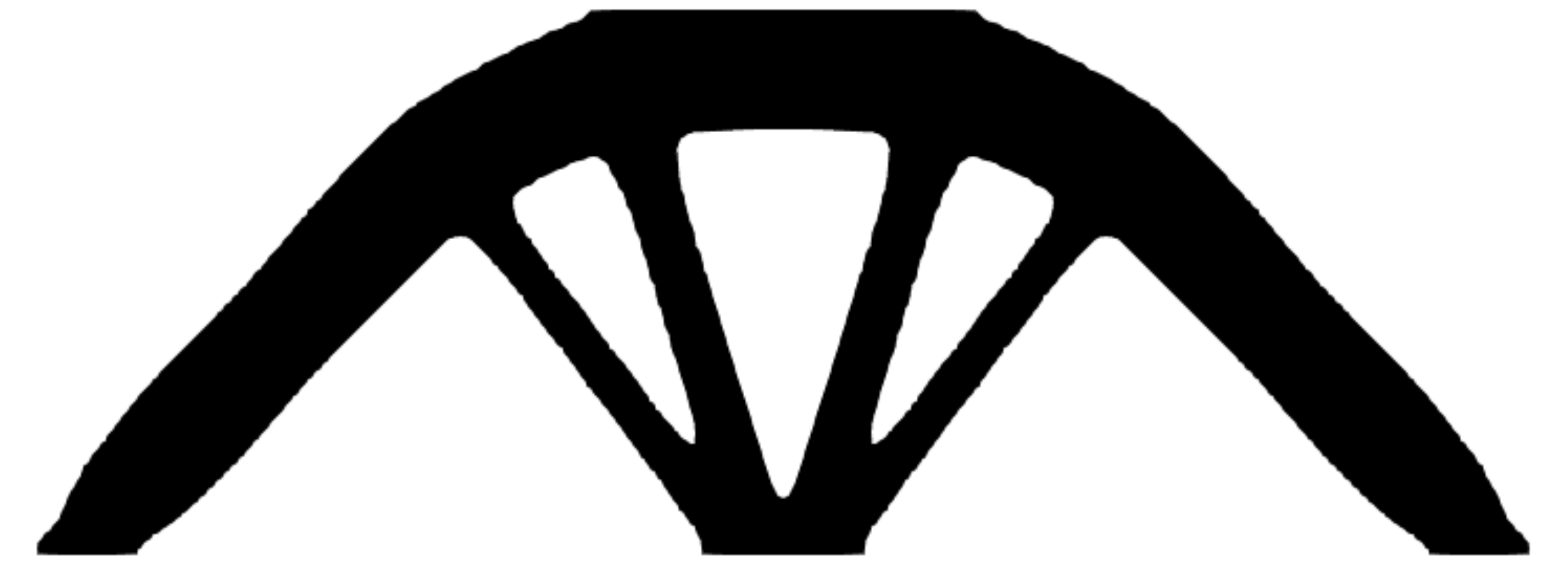}
        \subcaption{Step\,130}
        \label{b3-h}
      \end{minipage} 
       \begin{minipage}[t]{0.2\hsize}
        \centering
        \includegraphics[keepaspectratio, scale=0.09]{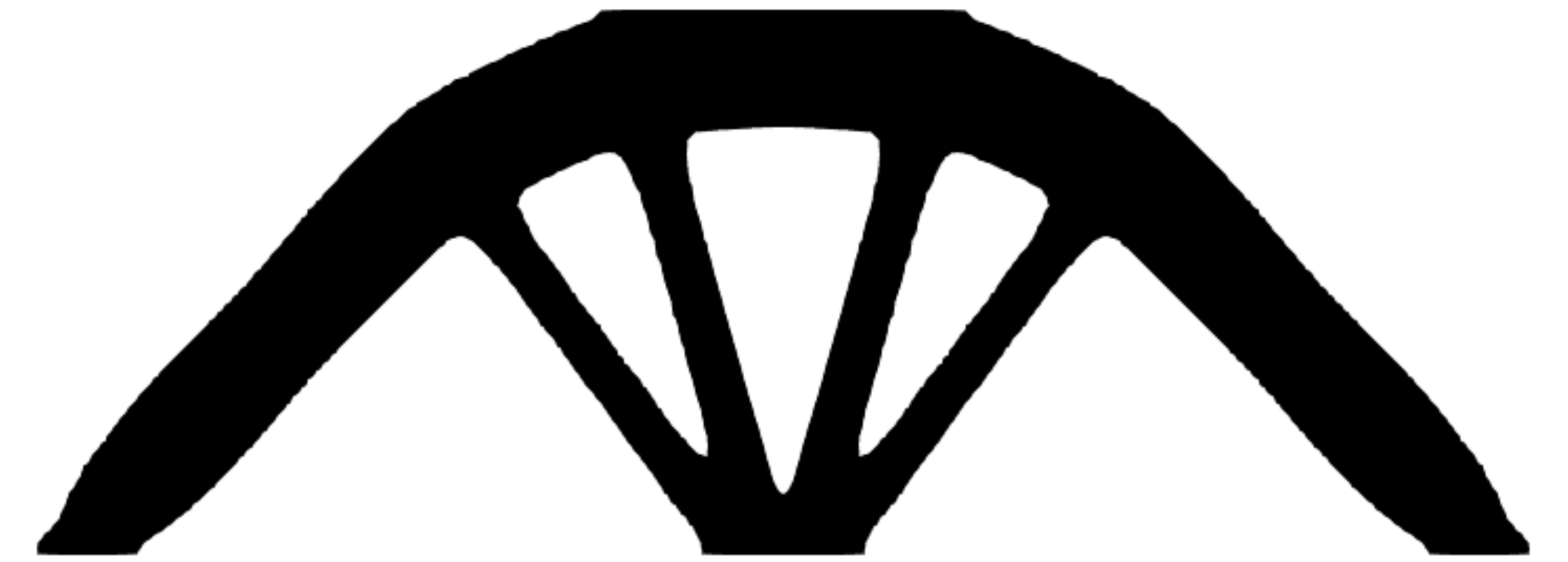}
        \subcaption{Step\,230}
        \label{b3-i}
      \end{minipage}
                 \begin{minipage}[t]{0.2\hsize}
        \centering
        \includegraphics[keepaspectratio, scale=0.09]{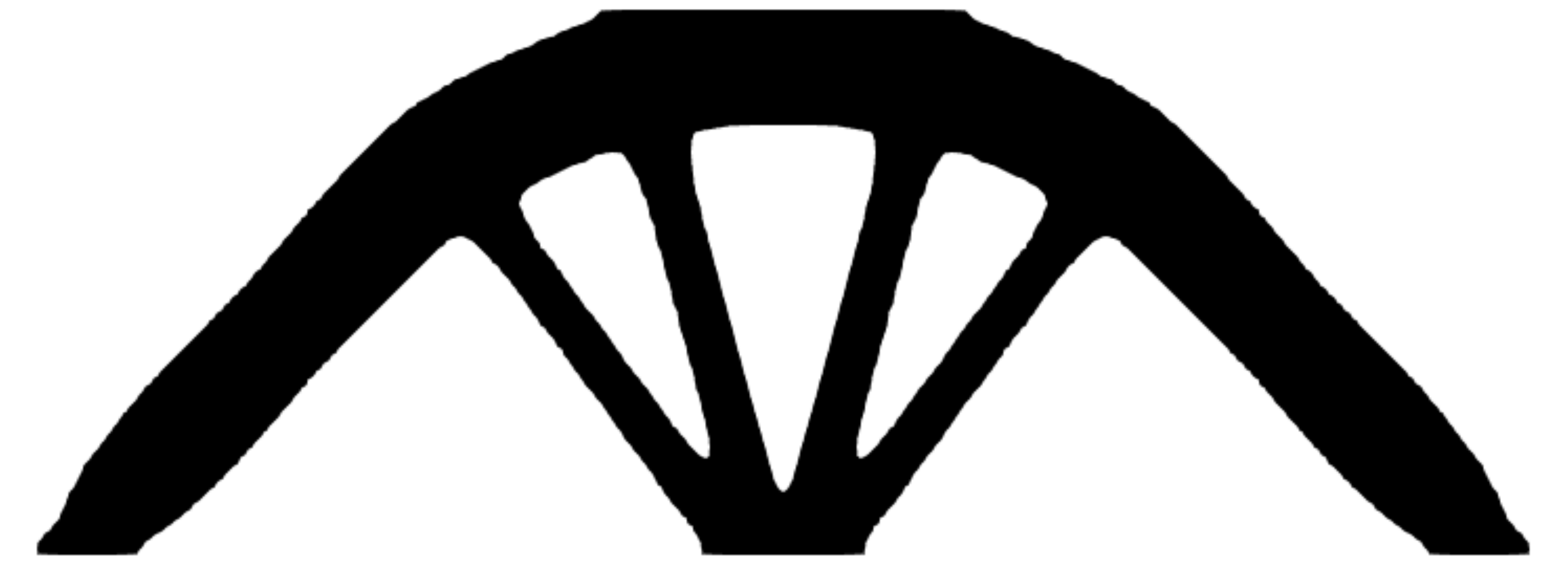}
        \subcaption{Step\,309$^{\#}$}
      \end{minipage}  
    \end{tabular}
     \caption{ Configuration $\Omega_{\phi_n}\subset D$ for the case where the initial configuration is the upper domain. 
     Figures (a)--(e) and (f)--(j) represent the results of (RD) and (NLHP), respectively.     
The symbol $^{\#}$ implies the final step. }
     \label{fig:b3}
  \end{figure*}

\begin{figure*}[htbp]
    \begin{tabular}{ccc}
      \hspace*{-5mm} 
      \begin{minipage}[t]{0.49\hsize}
        \centering
        \includegraphics[keepaspectratio, scale=0.33]{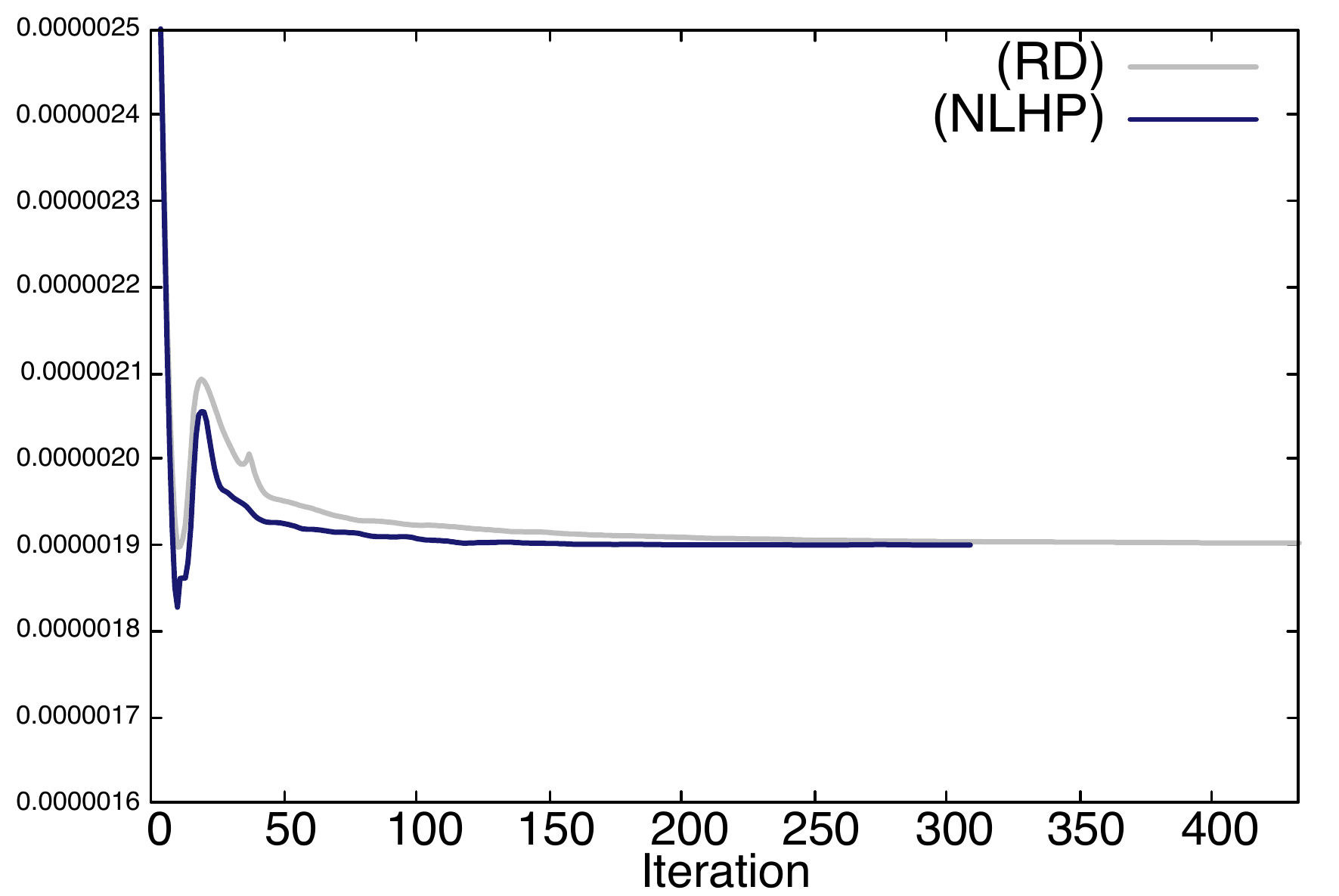}
        \subcaption{$F(\phi_n)$}
        \label{b3-1}
      \end{minipage} 
      \begin{minipage}[t]{0.49\hsize}
        \centering
        \includegraphics[keepaspectratio, scale=0.33]{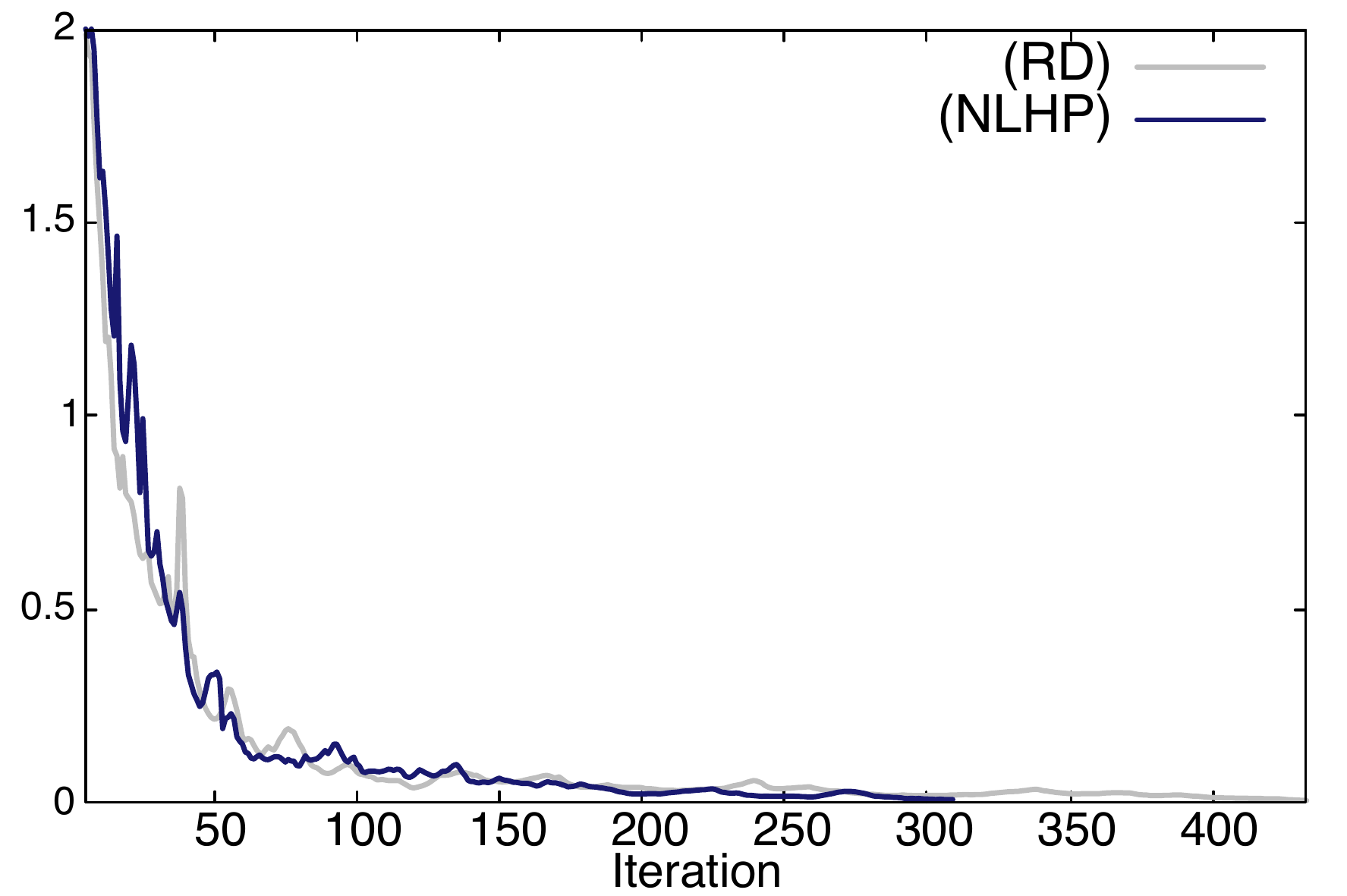}
        \subcaption{$\|\phi_{n+1}-\phi_n\|_{L^{\infty}(D)}$}
        \label{b3-2}
      \end{minipage} &
      \end{tabular}
       \caption{ Objective functional and convergence condition for \S \ref{S:br}-(iii).}
    \label{b3}
  \end{figure*}

\vspace{2mm}
\noindent{\bf Case (iii) (Upper domain).\,} 
As for the case where the initial configuration is the upper domain, 
Figures \ref{fig:b3} and \ref{b3} yield the assertion in this study; indeed, it is confirmed from Figures \ref{b3-b} and \ref{b3-g} that (NLHP) optimizes the topology faster than (RD).   
Moreover, we can also see that the boundary structure in (NLHP) is moving faster than that of (RD) in terms of satisfying the convergence condition more quickly.

\begin{figure*}[htbp]
\hspace*{-5mm}
    \begin{tabular}{cccc}
      \begin{minipage}[t]{0.24\hsize}
        \centering
        \includegraphics[keepaspectratio, scale=0.13]{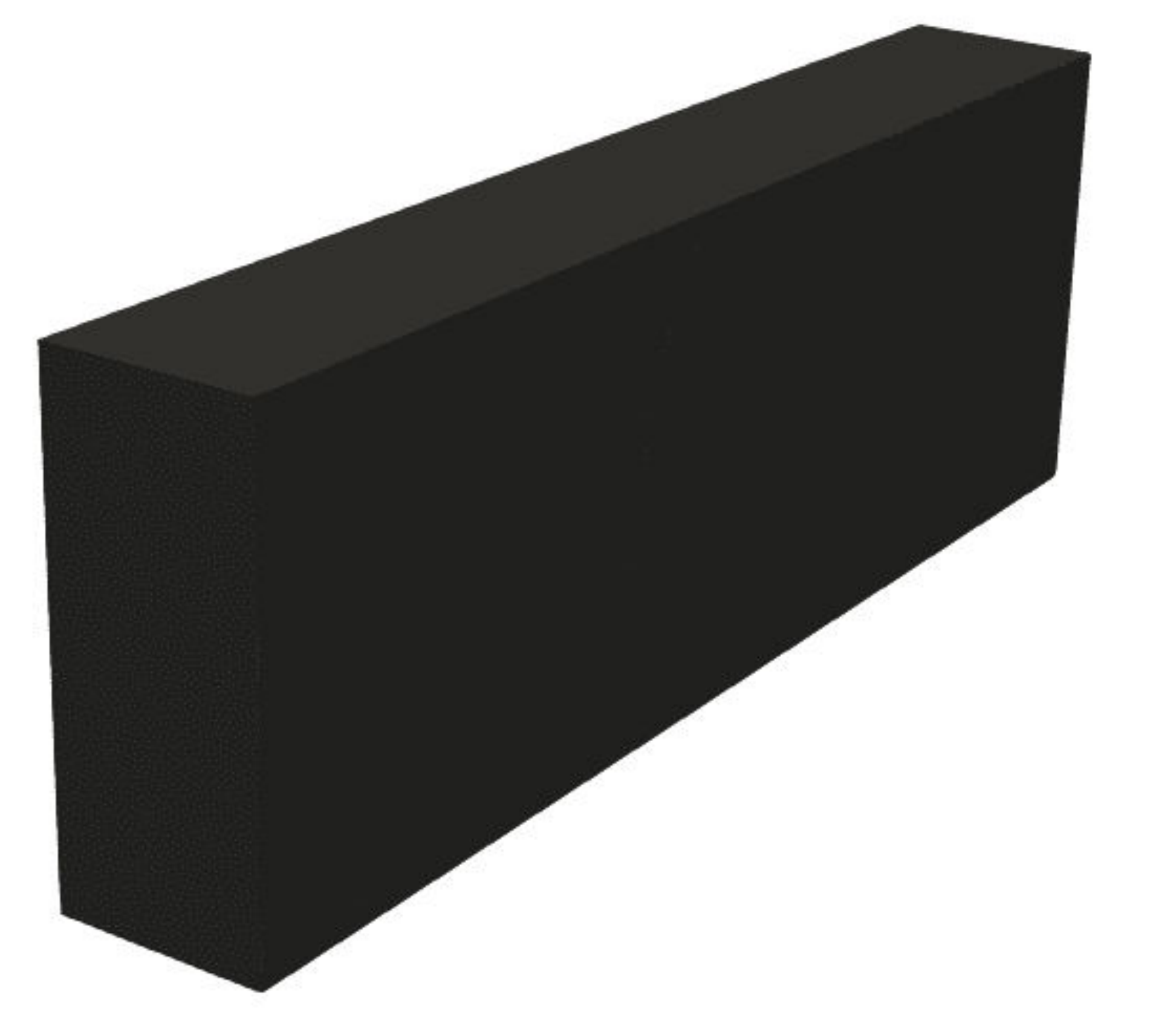}
        \subcaption{Step\,0}
        \label{b3d-a}
      \end{minipage} 
      \begin{minipage}[t]{0.24\hsize}
        \centering
        \includegraphics[keepaspectratio, scale=0.13]{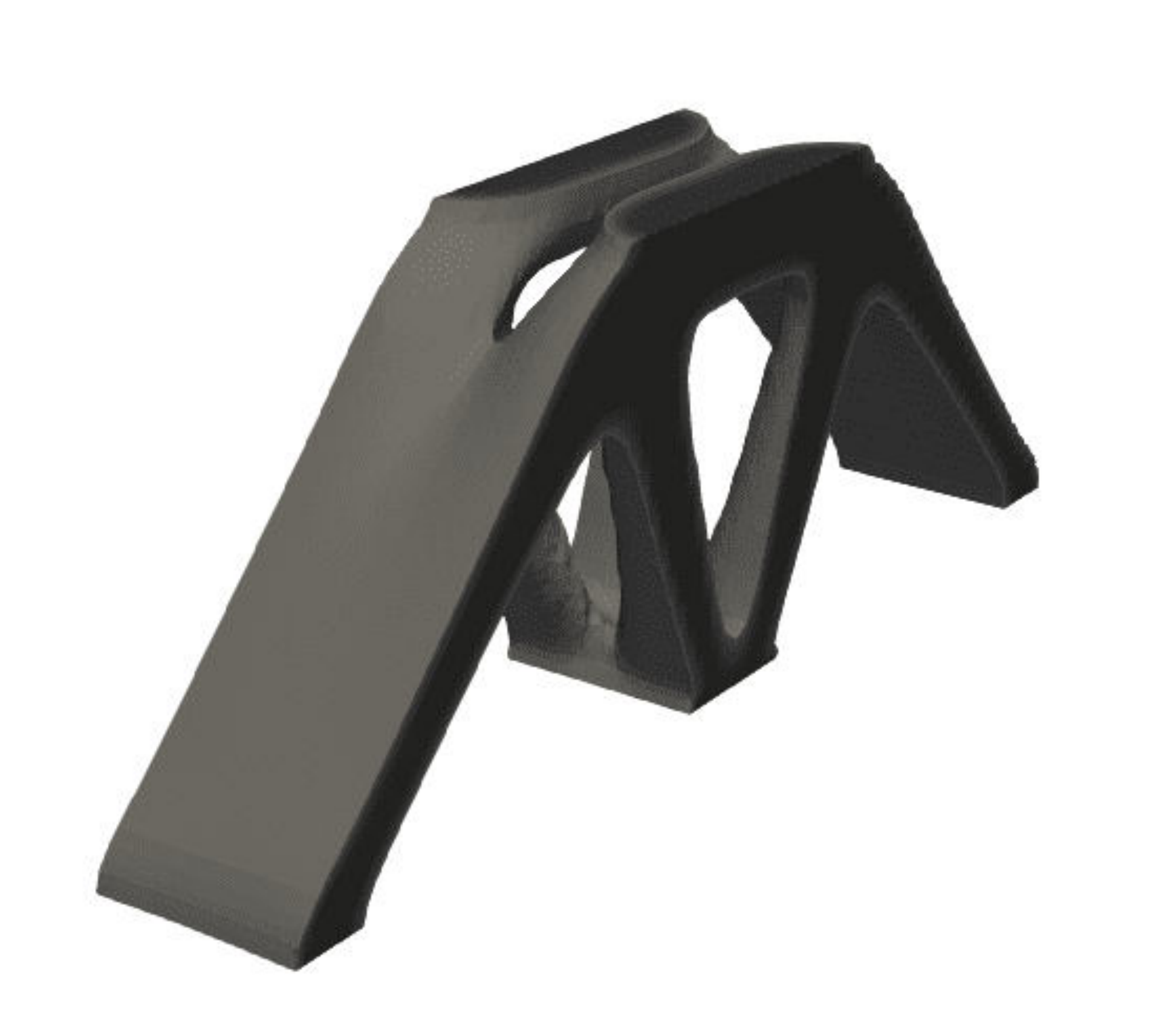}
        \subcaption{Step\,20}
        \label{b3d-b}
      \end{minipage} 
      \begin{minipage}[t]{0.24\hsize}
        \centering
        \includegraphics[keepaspectratio, scale=0.13]{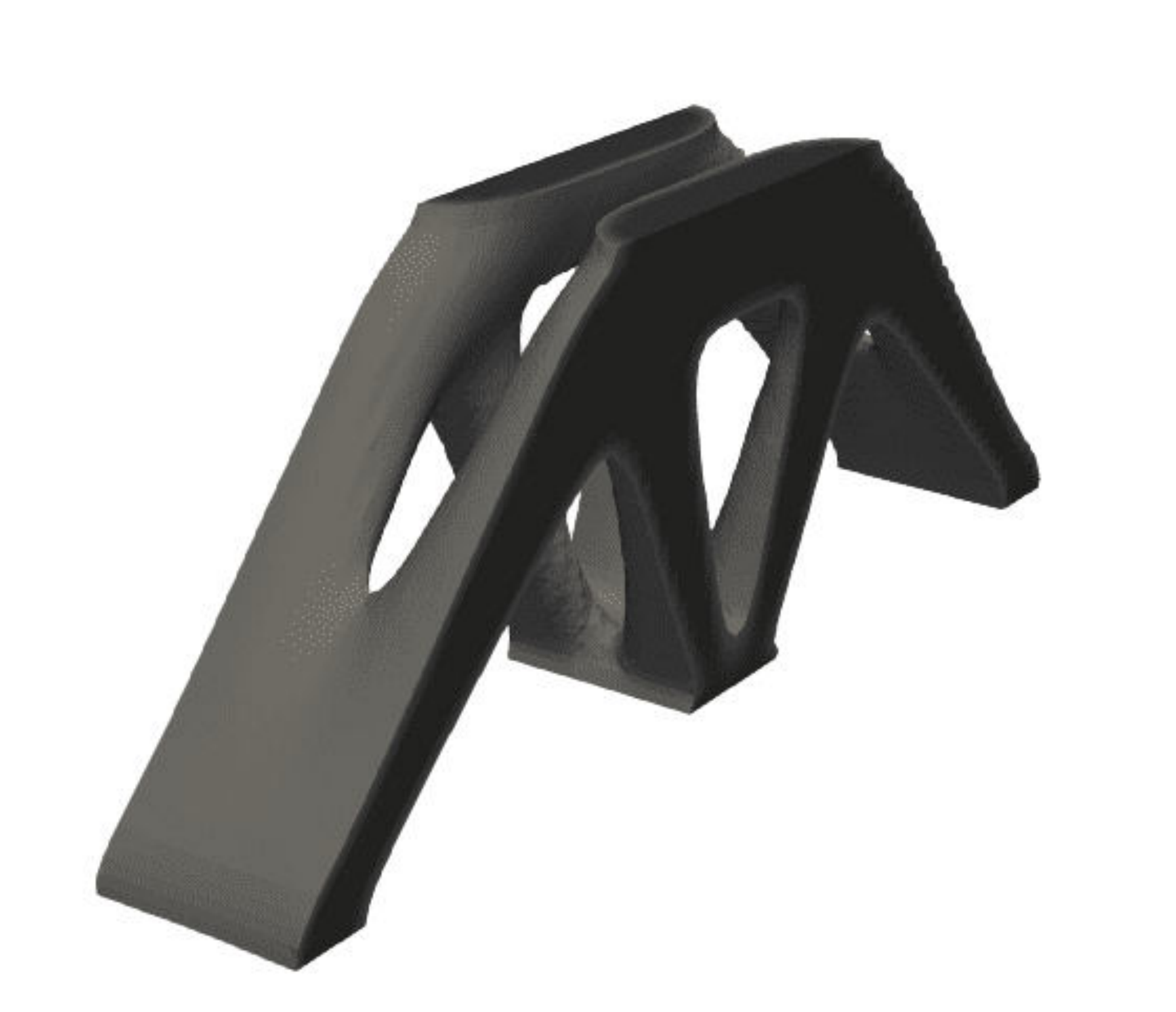}
        \subcaption{Step\,40}
        \label{b3d-c}
      \end{minipage} 
         \begin{minipage}[t]{0.24\hsize}
        \centering
        \includegraphics[keepaspectratio, scale=0.13]{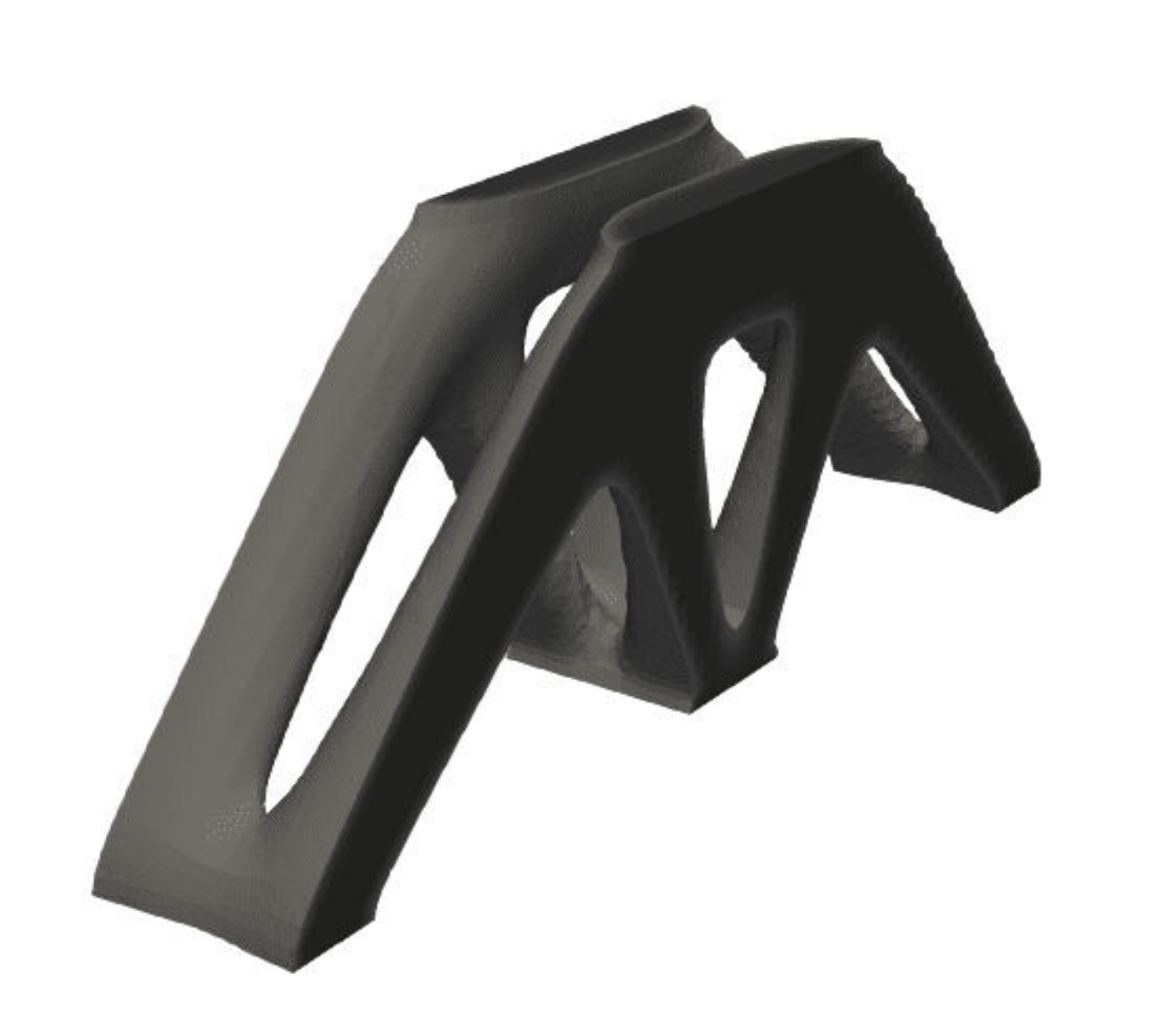}
        \subcaption{Step\,104$^\#$}
        \label{b3d-d}
      \end{minipage} 
         \\
     \begin{minipage}[t]{0.24\hsize}
        \centering
        \includegraphics[keepaspectratio, scale=0.13]{mbb3d0.pdf}
        \subcaption{Step\,0}
        \label{b3d-e}
      \end{minipage} 
      \begin{minipage}[t]{0.24\hsize}
        \centering
        \includegraphics[keepaspectratio, scale=0.13]{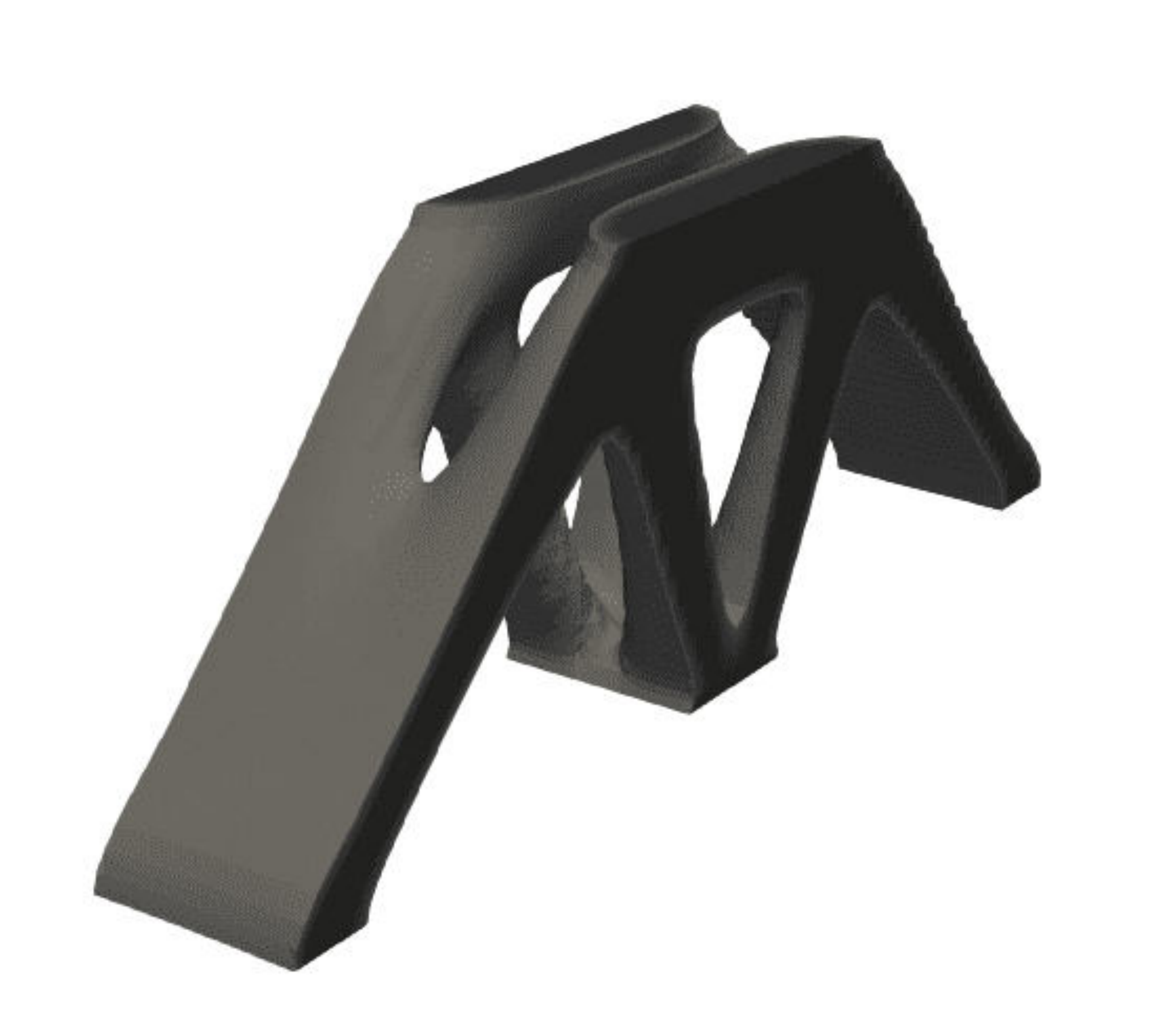}
        \subcaption{Step\,20}
        \label{b3d-f}
      \end{minipage} 
      \begin{minipage}[t]{0.24\hsize}
        \centering
        \includegraphics[keepaspectratio, scale=0.13]{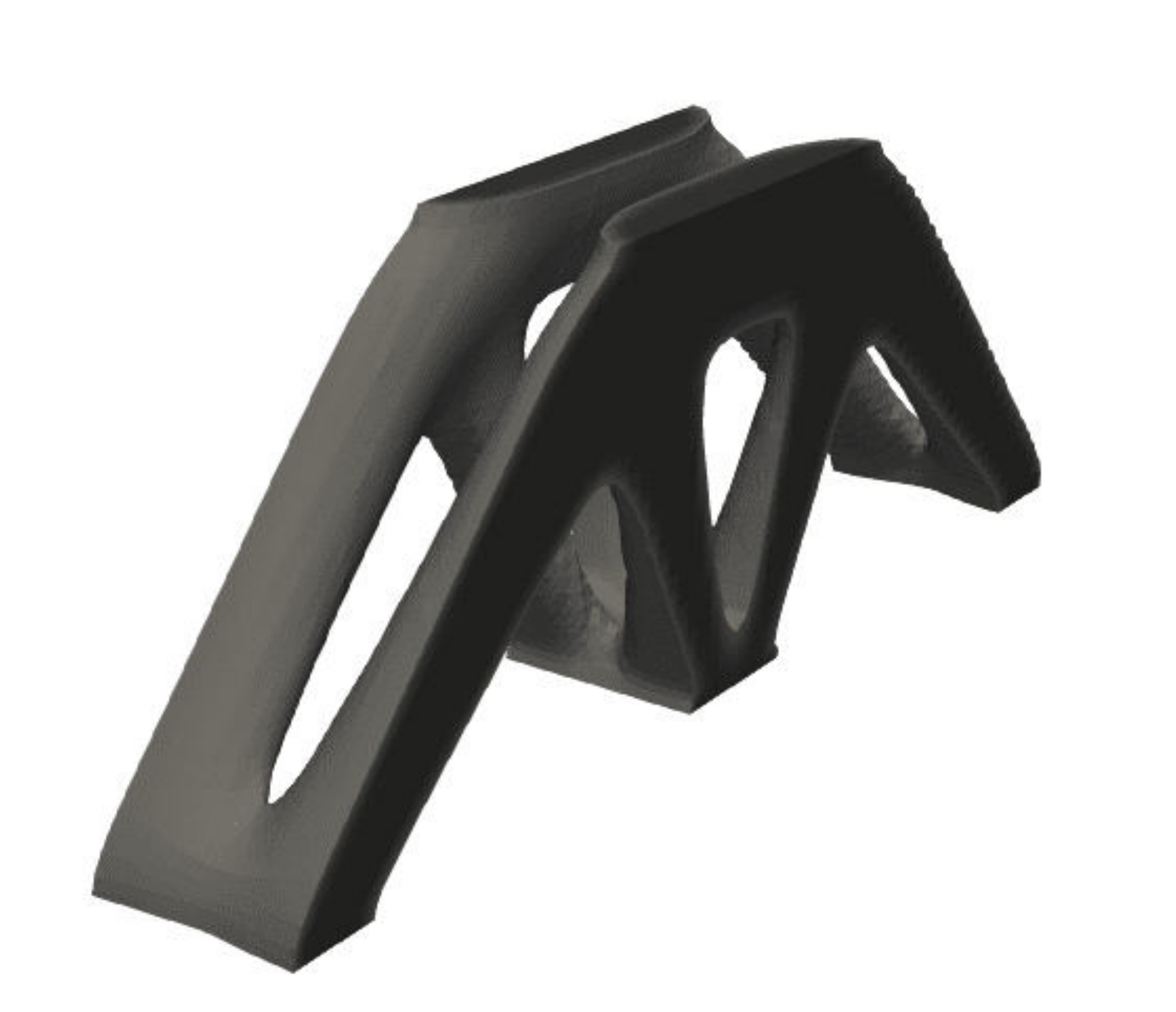}
        \subcaption{Step\,40}
        \label{b3d-g}
      \end{minipage} 
         \begin{minipage}[t]{0.24\hsize}
        \centering
        \includegraphics[keepaspectratio, scale=0.13]{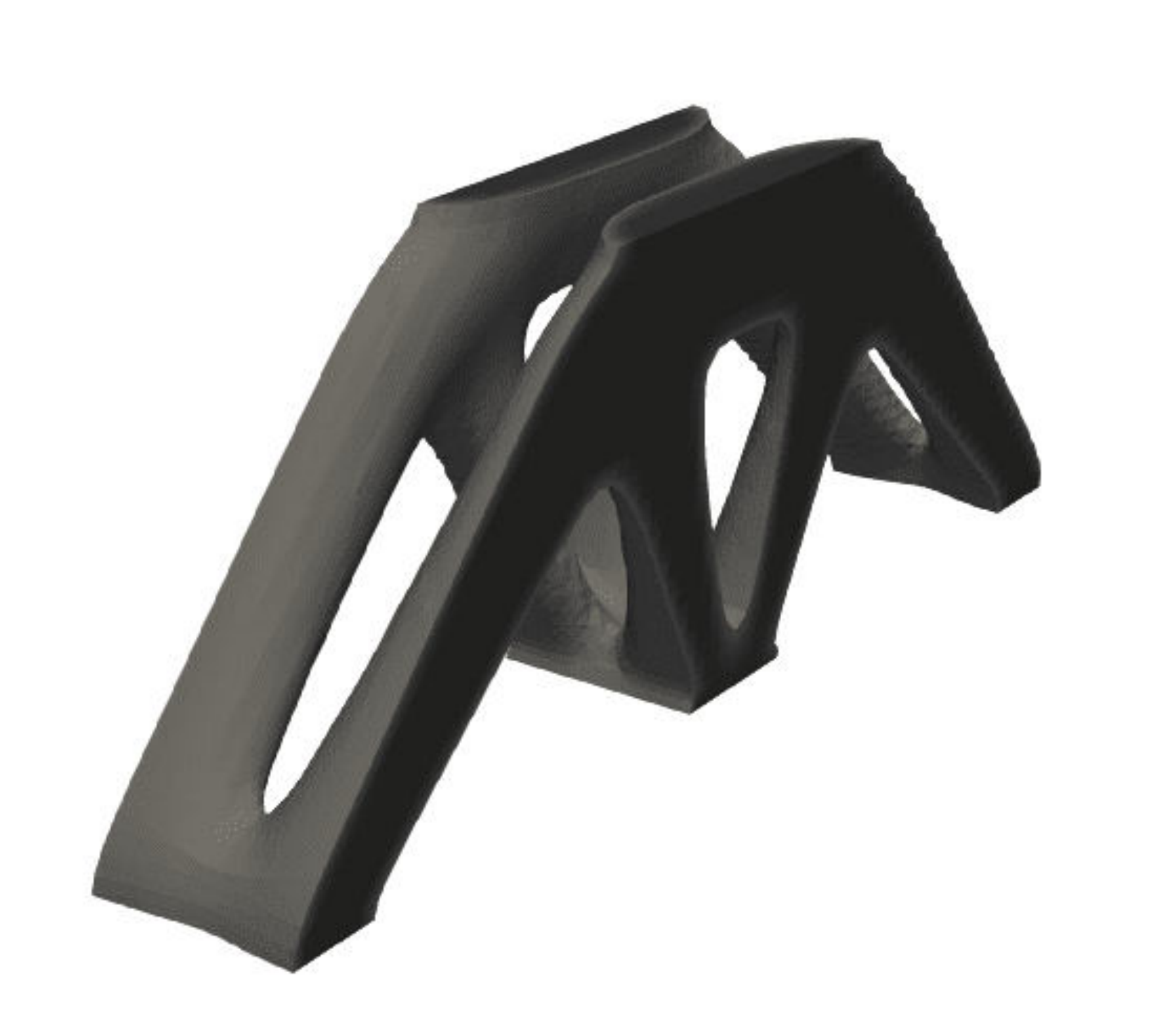}
        \subcaption{Step\,59$^\#$}
        \label{b3d-h}
      \end{minipage} 
       \end{tabular}
     \caption{ Configuration $\Omega_{\phi_n}\subset D\subset \R^3$ for the case where the initial configuration is the whole domain. 
     Figures (a)--(d) and (e)--(h) represent the results of (RD) and (NLHP), respectively.     
The symbol $^{\#}$ implies the final step. Here the depth of $D$ is set to $0.3$.}
     \label{fig:b3d}
  \end{figure*}

\begin{figure*}[htbp]
    \begin{tabular}{ccc}
      \hspace*{-5mm} 
      \begin{minipage}[t]{0.49\hsize}
        \centering
        \includegraphics[keepaspectratio, scale=0.33]{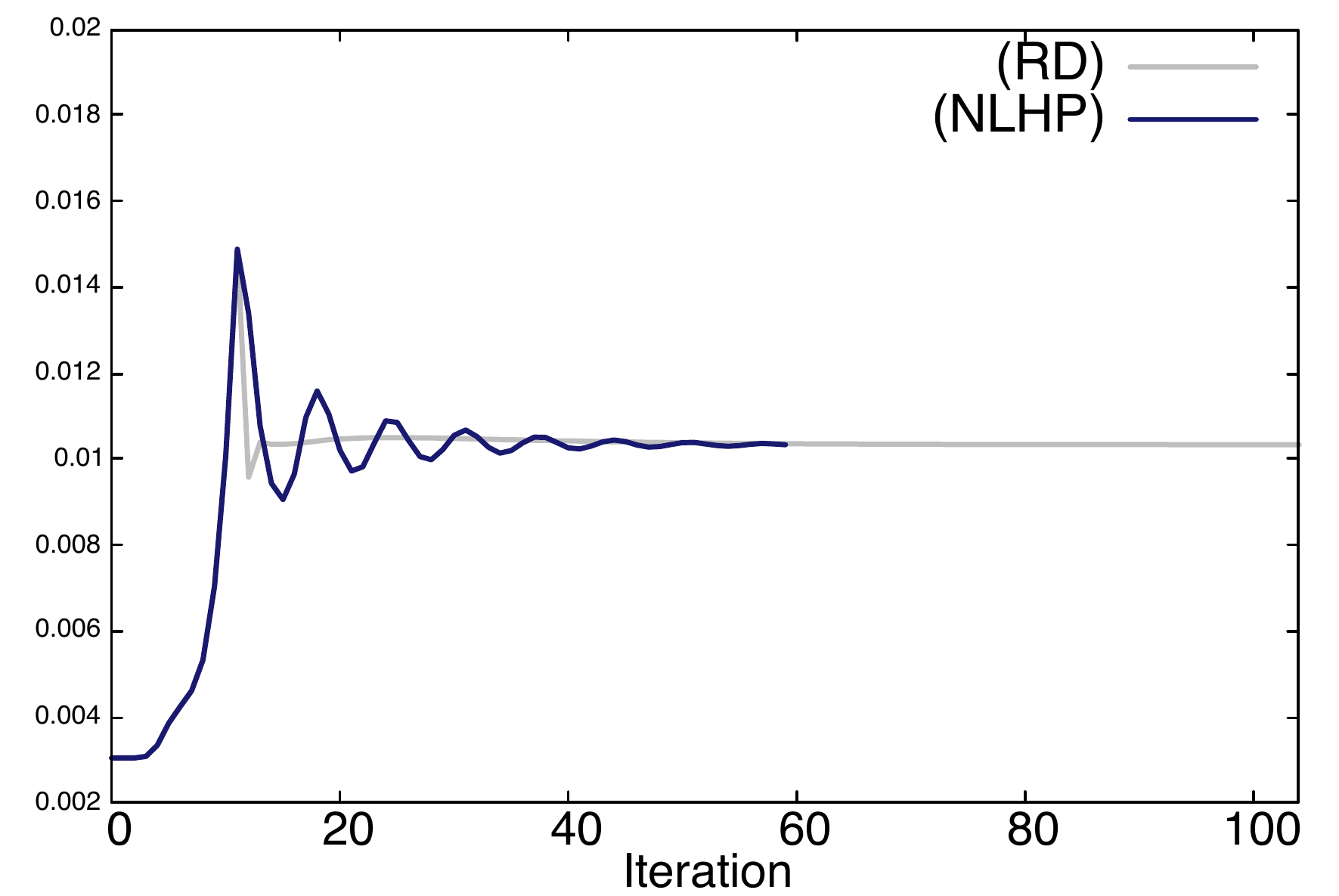}
        \subcaption{$F(\phi_n)$}
        \label{b3d-a}
      \end{minipage} 
      \begin{minipage}[t]{0.49\hsize}
        \centering
        \includegraphics[keepaspectratio, scale=0.33]{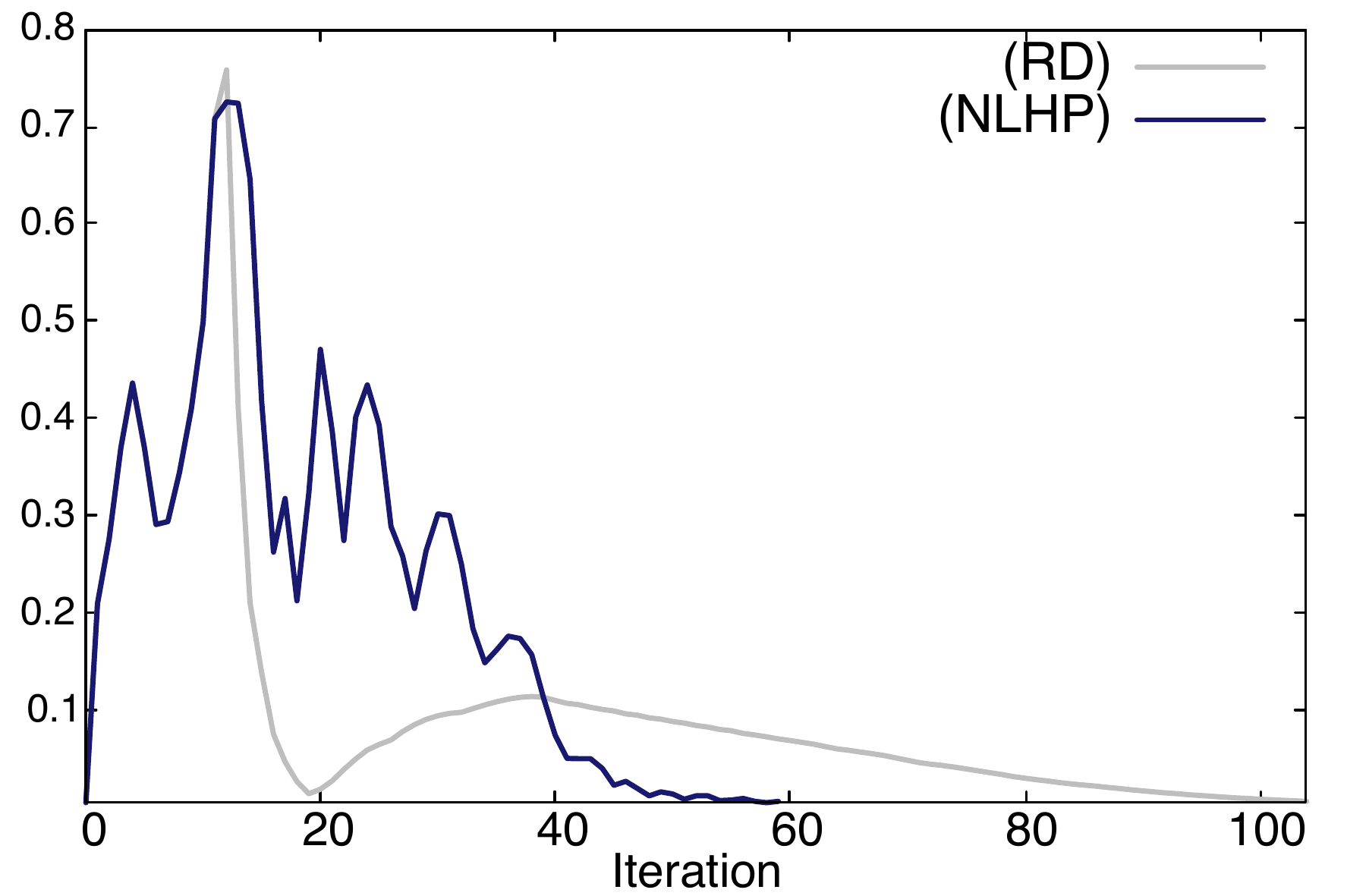}
        \subcaption{$\|\phi_{n+1}-\phi_n\|_{L^{\infty}(D)}$}
        \label{b3d-b}
      \end{minipage} &
      \end{tabular}
       \caption{ Objective functional and convergence condition for \S \ref{S:br}-(iv).}
    \label{b3d}
  \end{figure*}

\vspace{2mm}
\noindent{\bf Case (iv) (Three-dimensional domain).\,} 
The same results as in the two-dimensional cases can also be obtained for the three-dimensional case.
Actually, we readily deduce from Figure \ref{fig:b3d} that 
the boundary structure in (NLHP) is moving faster than that of (RD).
Moreover, Figure \ref{b3d} ensures the validity of the proposed method.
Here we set $(n_t,h_{\rm max})=(176629,0.0433)$ and $(\tau, G_{\rm max})=(1.0\times 10^{-3}, 0.2)$.


\subsection{Radiator model} \label{S:ra}
We finally consider the so-called \emph{radiator model} as a more complex geometry model (see, e.g.,~\cite[Example 5.9]{ACMOY19} and \cite{Y22}).
Here we set $(n_t,h_{\rm max})=(25600, 0.0088)$ and $(\tau, G_{\rm max})=(2.0\times 10^{-5},0.5)$. The setting for the fixed domain $D\subset \R^2$ and boundary conditions are as in Figure \ref{ira}.

\begin{figure*}[htbp]
\hspace*{-5mm}
    \begin{tabular}{ccccc}
      \begin{minipage}[t]{0.2\hsize}
        \centering
        \includegraphics[keepaspectratio, scale=0.09]{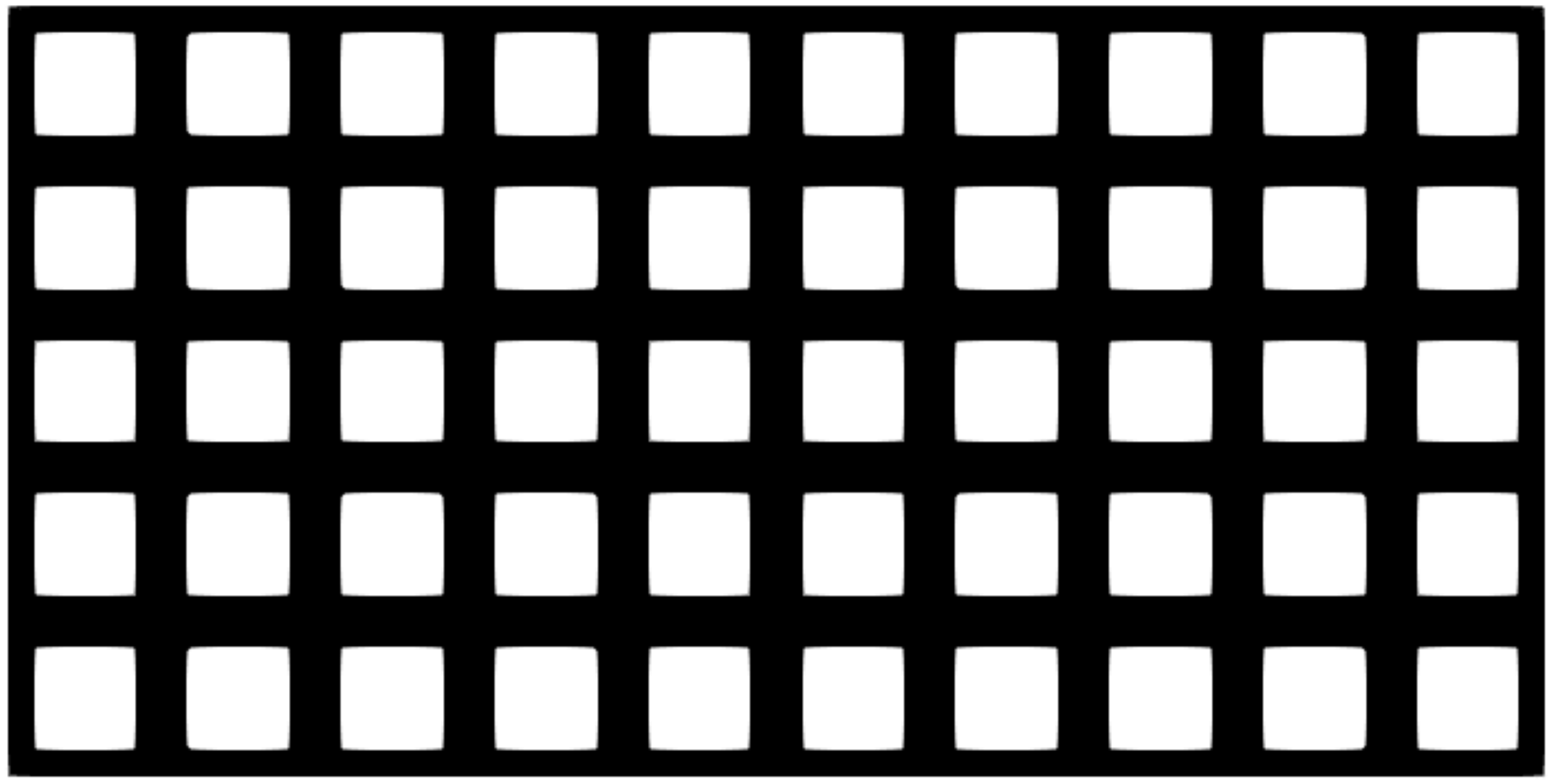}
        \subcaption{Step\,0}
        \label{ra1-a}
      \end{minipage} 
      \begin{minipage}[t]{0.2\hsize}
        \centering
        \includegraphics[keepaspectratio, scale=0.09]{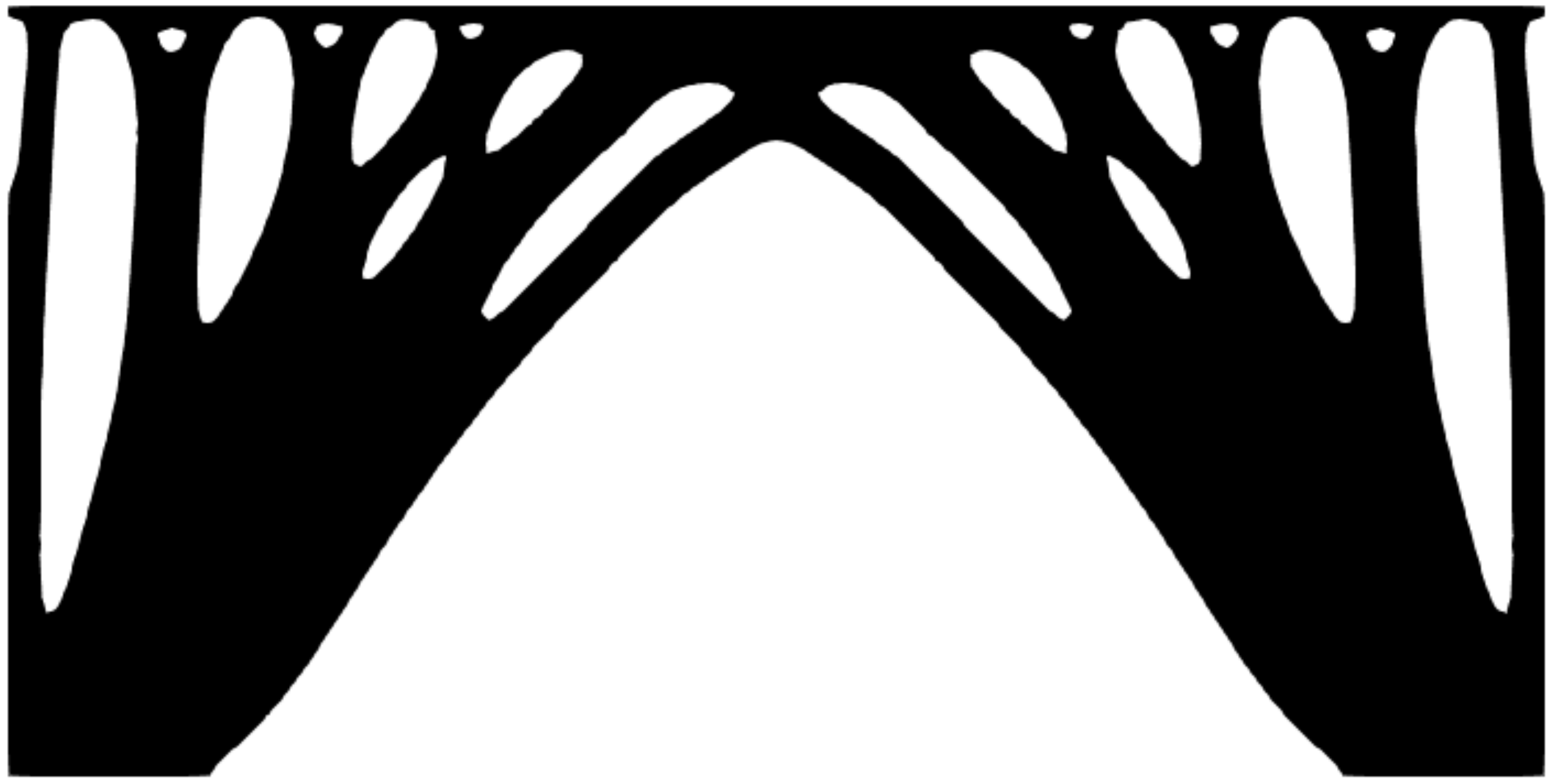}
        \subcaption{Step\,50}
        \label{ra1-b}
      \end{minipage} 
      \begin{minipage}[t]{0.2\hsize}
        \centering
        \includegraphics[keepaspectratio, scale=0.09]{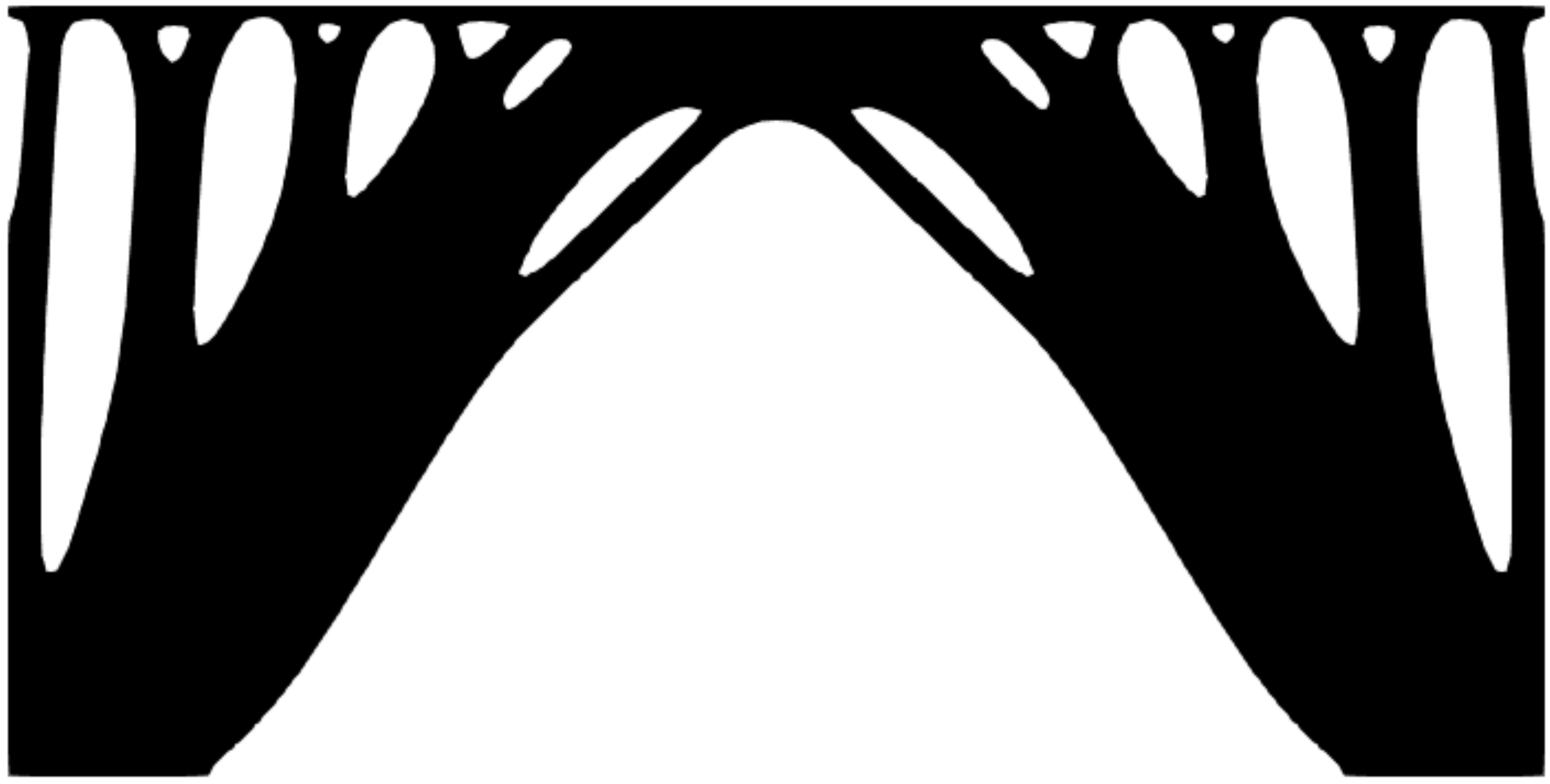}
        \subcaption{Step\,100}
        \label{ra1-c}
      \end{minipage} 
         \begin{minipage}[t]{0.2\hsize}
        \centering
        \includegraphics[keepaspectratio, scale=0.09]{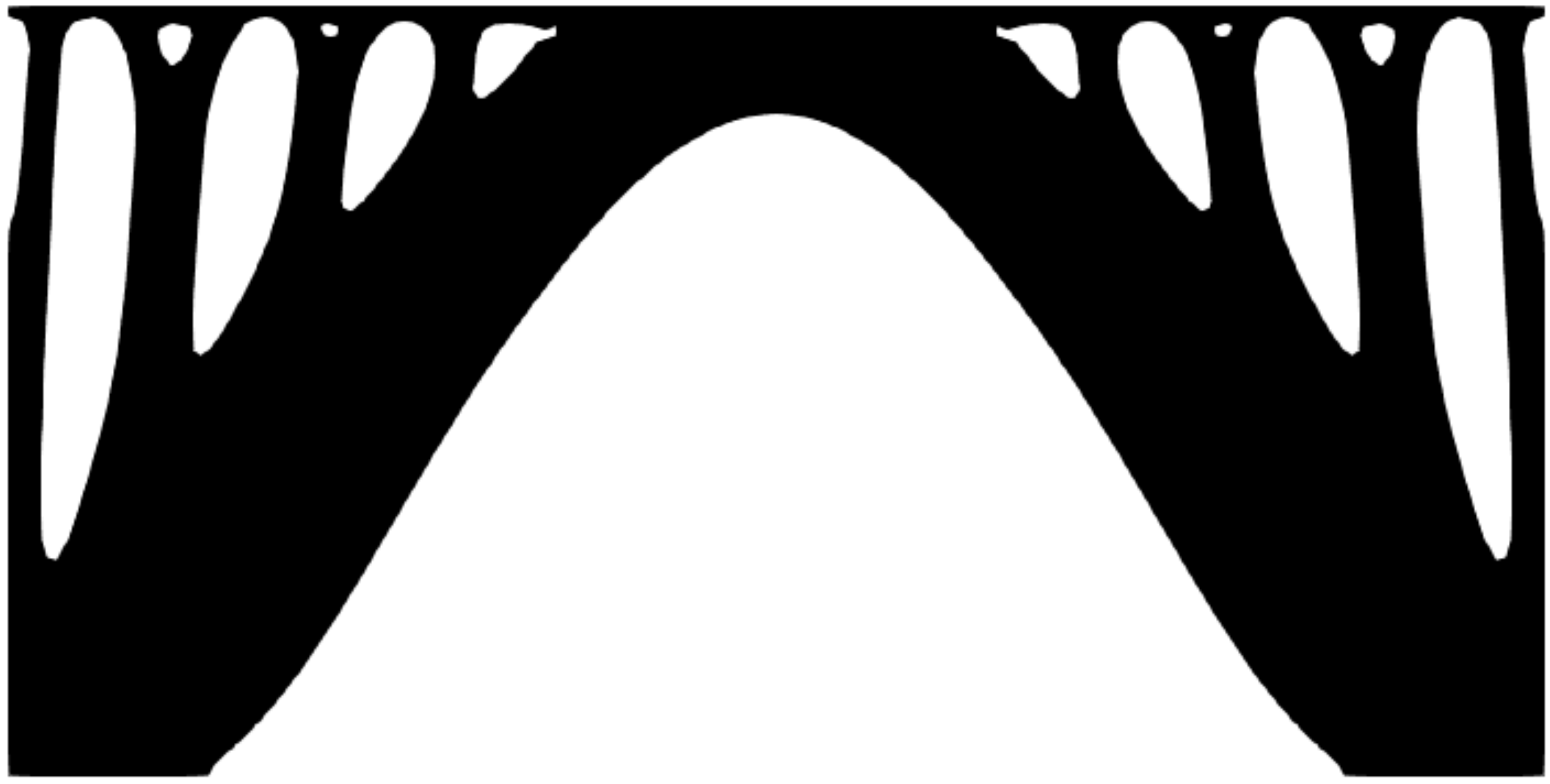}
        \subcaption{Step\,150}
        \label{ra1-d}
      \end{minipage} 
                 \begin{minipage}[t]{0.2\hsize}
        \centering
        \includegraphics[keepaspectratio, scale=0.09]{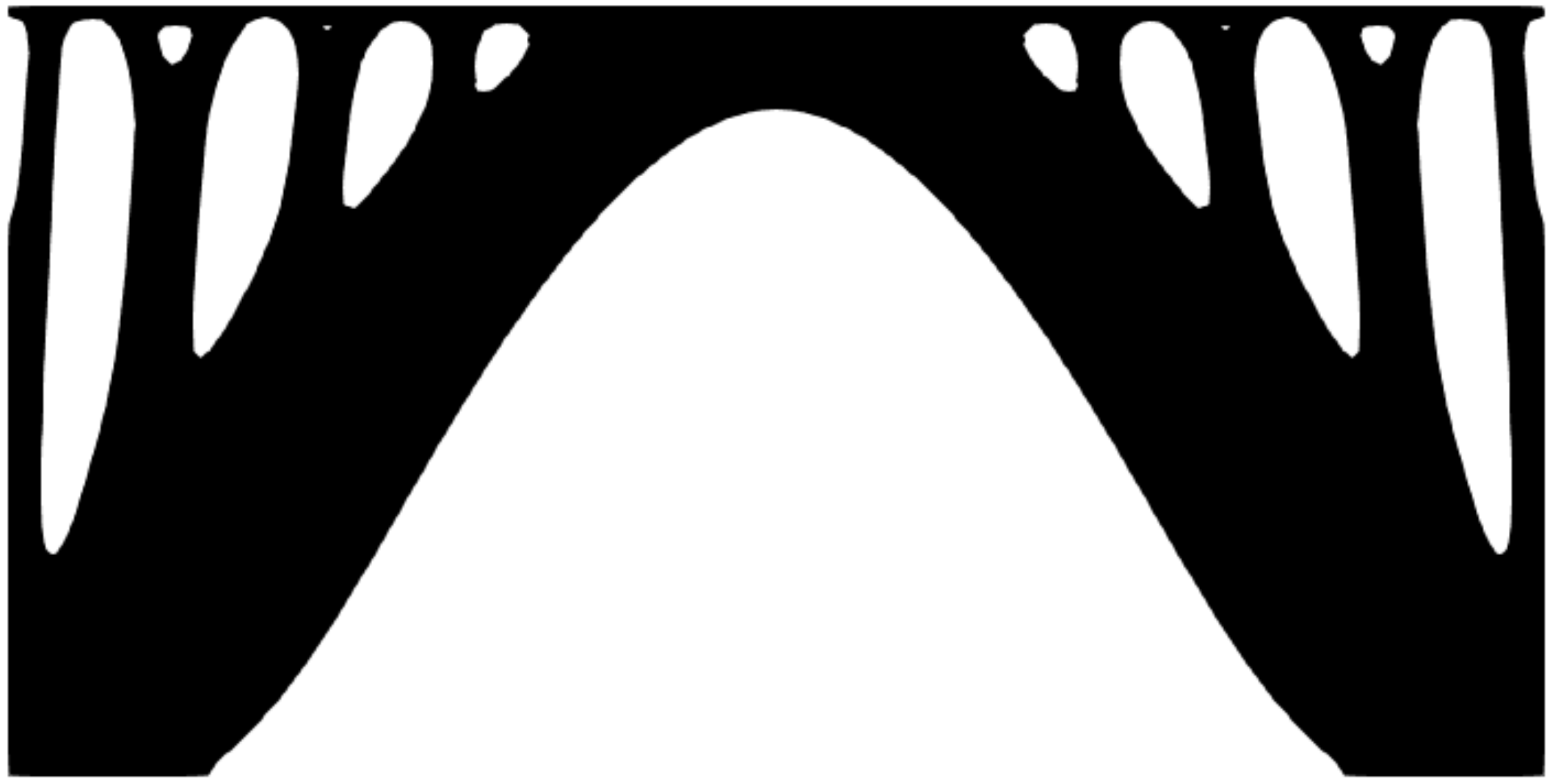}
        \subcaption{Step\,189$^{\#}$}
      \end{minipage} 
      \\
    \begin{minipage}[t]{0.2\hsize}
        \centering
        \includegraphics[keepaspectratio, scale=0.09]{rap0.pdf}
        \subcaption{Step\,0}
        \label{ra1-f}
      \end{minipage} 
      \begin{minipage}[t]{0.2\hsize}
        \centering
        \includegraphics[keepaspectratio, scale=0.09]{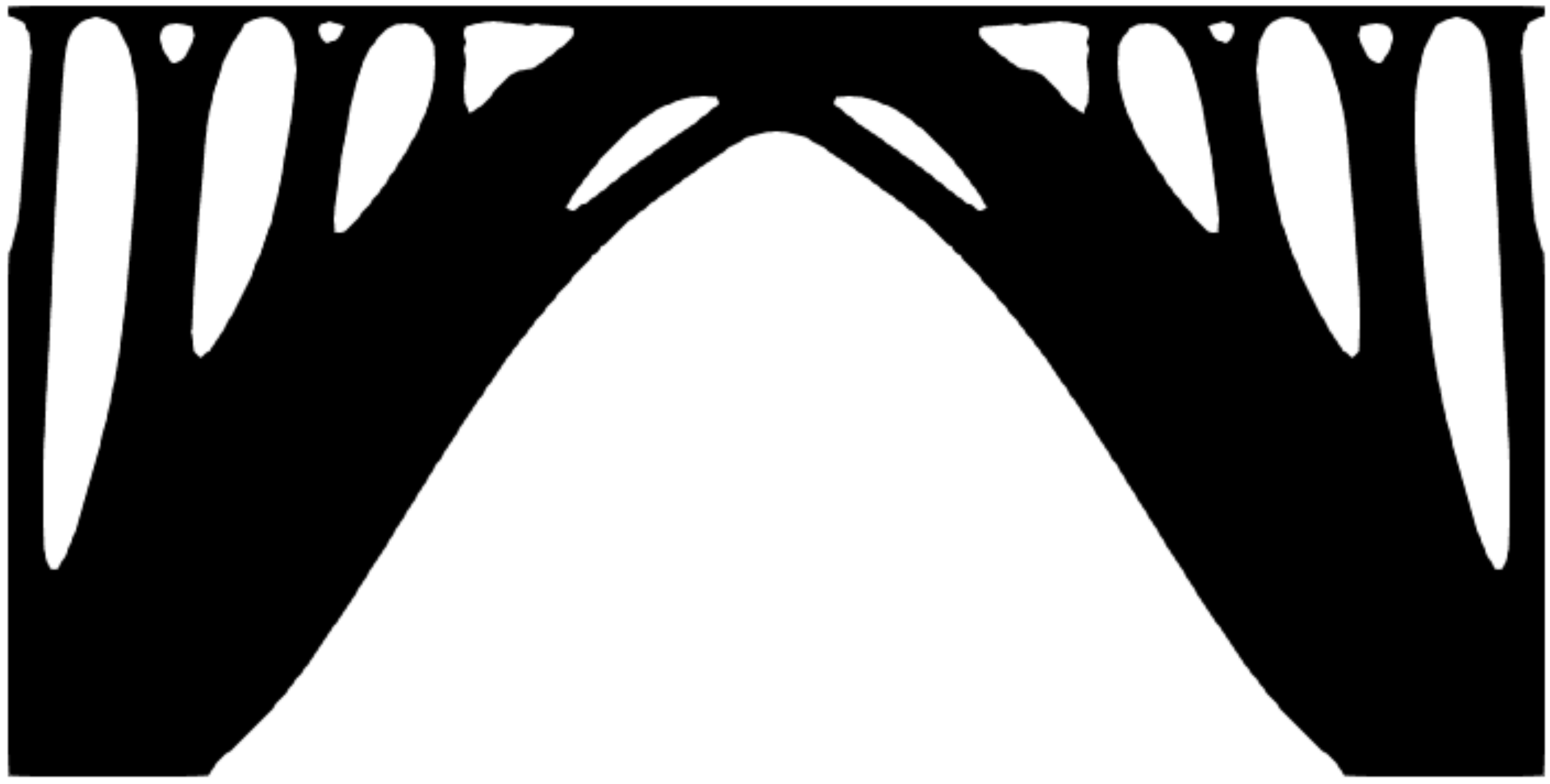}
        \subcaption{Step\,50}
        \label{ra1-g}
      \end{minipage} 
      \begin{minipage}[t]{0.2\hsize}
        \centering
        \includegraphics[keepaspectratio, scale=0.09]{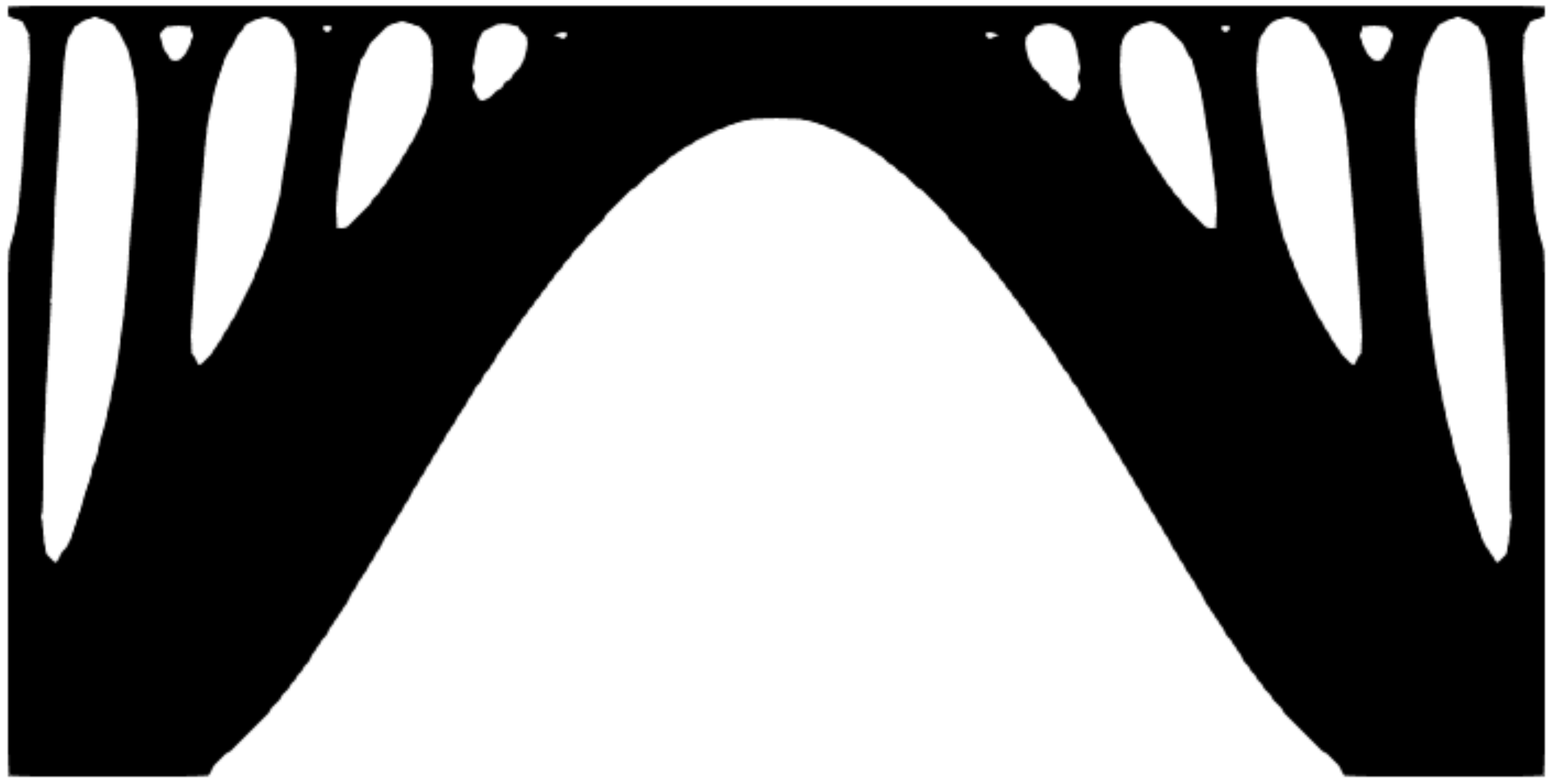}
        \subcaption{Step\,100}
        \label{ra1-h}
      \end{minipage} 
       \begin{minipage}[t]{0.2\hsize}
        \centering
        \includegraphics[keepaspectratio, scale=0.09]{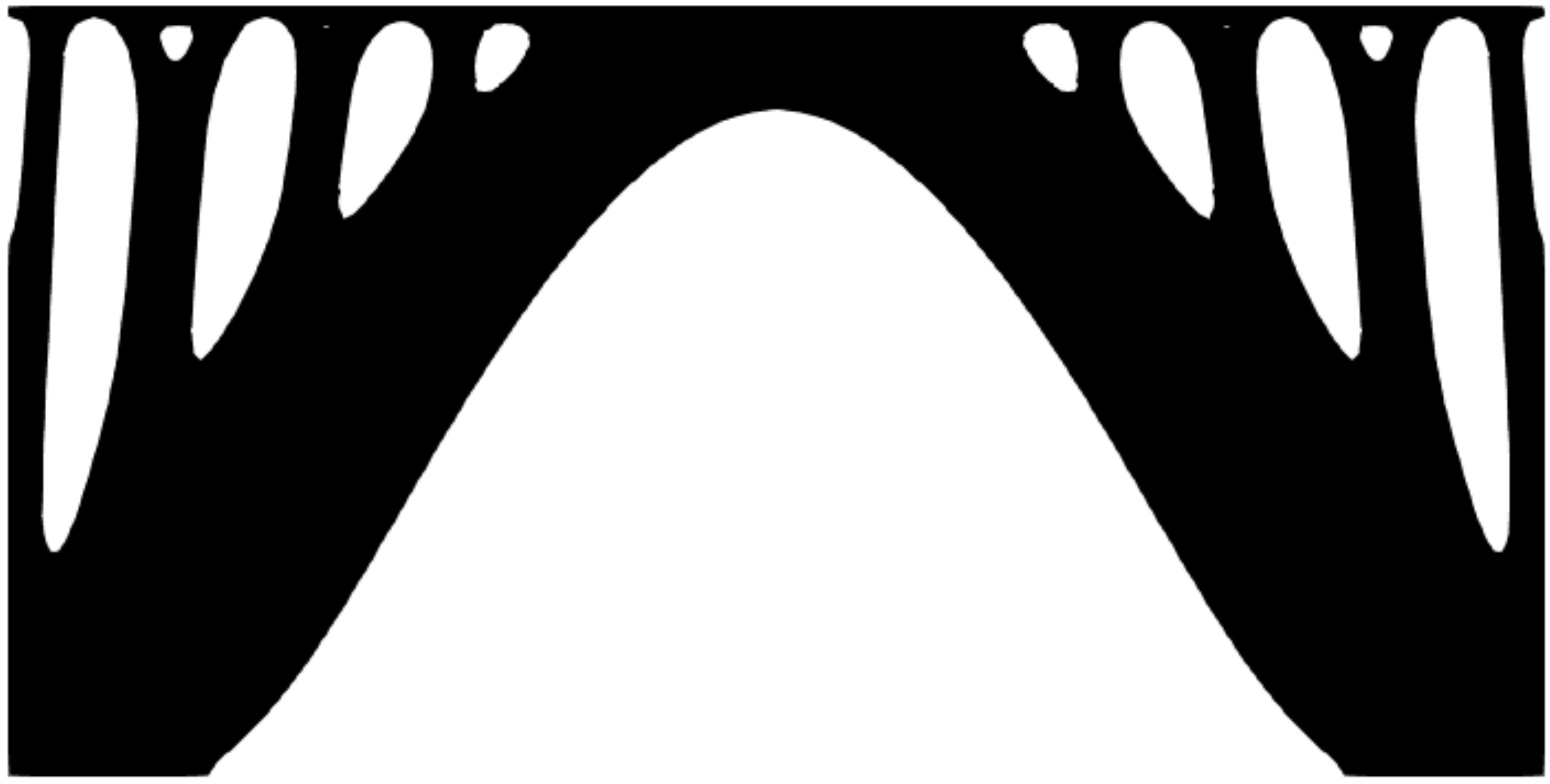}
        \subcaption{Step\,150}
        \label{ra1-i}
      \end{minipage}
                 \begin{minipage}[t]{0.2\hsize}
        \centering
        \includegraphics[keepaspectratio, scale=0.09]{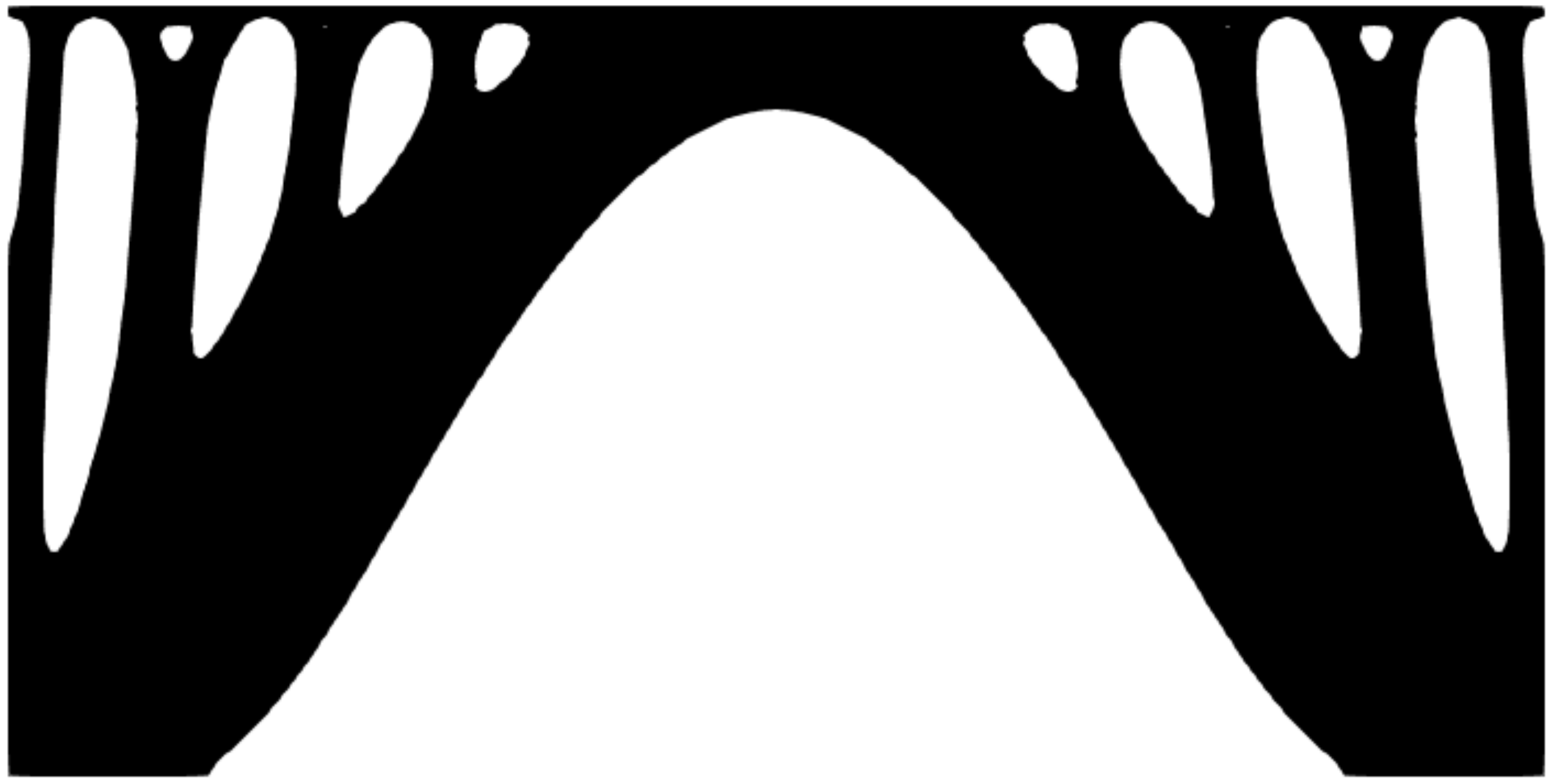}
        \subcaption{Step\,157$^{\#}$}
      \end{minipage}  
    \end{tabular}
     \caption{ Configuration $\Omega_{\phi_n}\subset D$ for the case where the initial configuration is the periodically perforated domain. Figures (a)--(e) and (f)--(j) represent the results of (RD) and (NLHP), respectively.     
The symbol $^{\#}$ implies the final step.}
     \label{fig:ra1}
  \end{figure*}

\begin{figure*}[htbp]
    \begin{tabular}{ccc}
      \hspace*{-5mm} 
      \begin{minipage}[t]{0.49\hsize}
        \centering
        \includegraphics[keepaspectratio, scale=0.33]{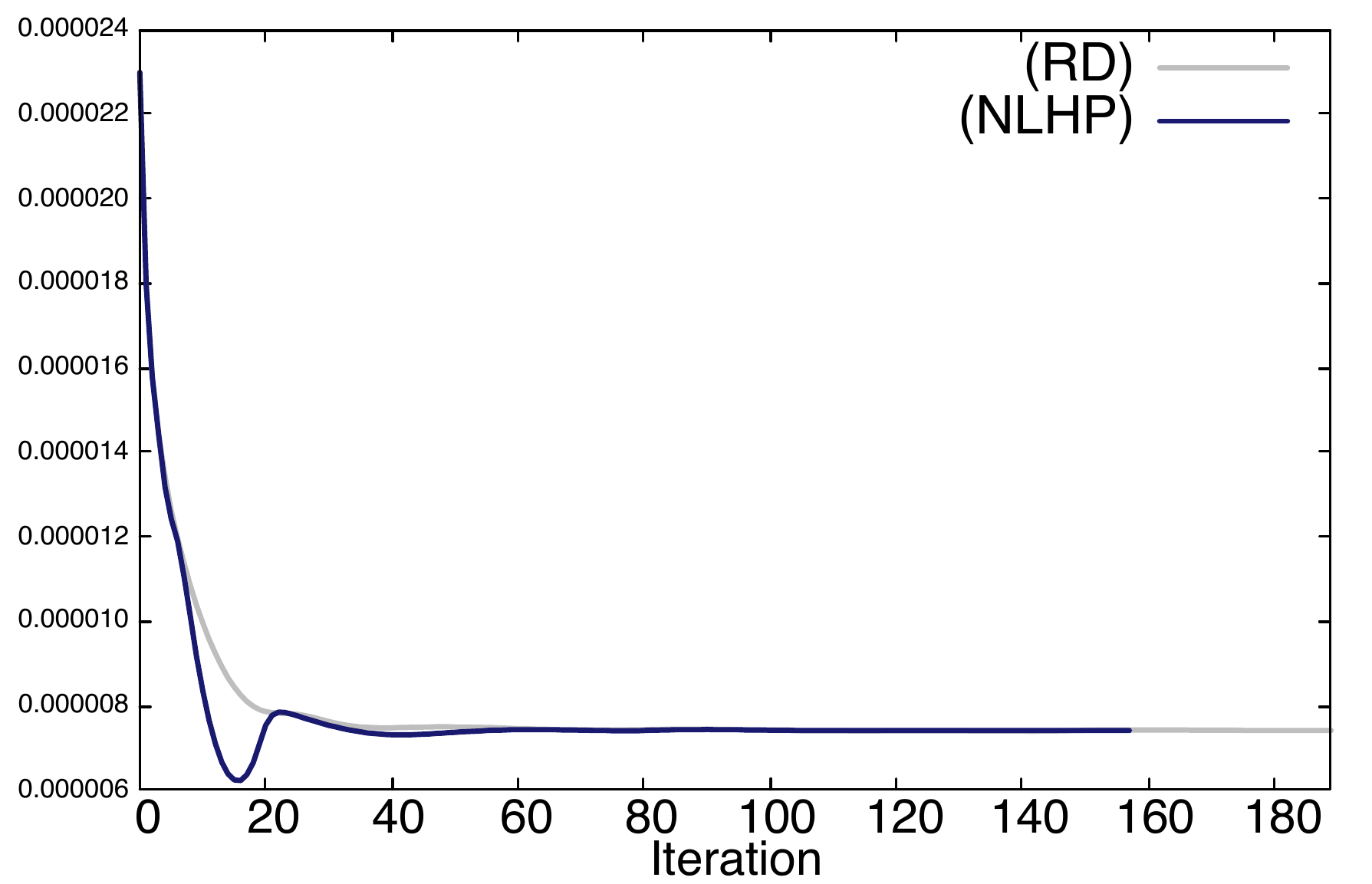}
        \subcaption{$F(\phi_n)$}
        \label{ra1-1}
      \end{minipage} 
      \begin{minipage}[t]{0.49\hsize}
        \centering
        \includegraphics[keepaspectratio, scale=0.33]{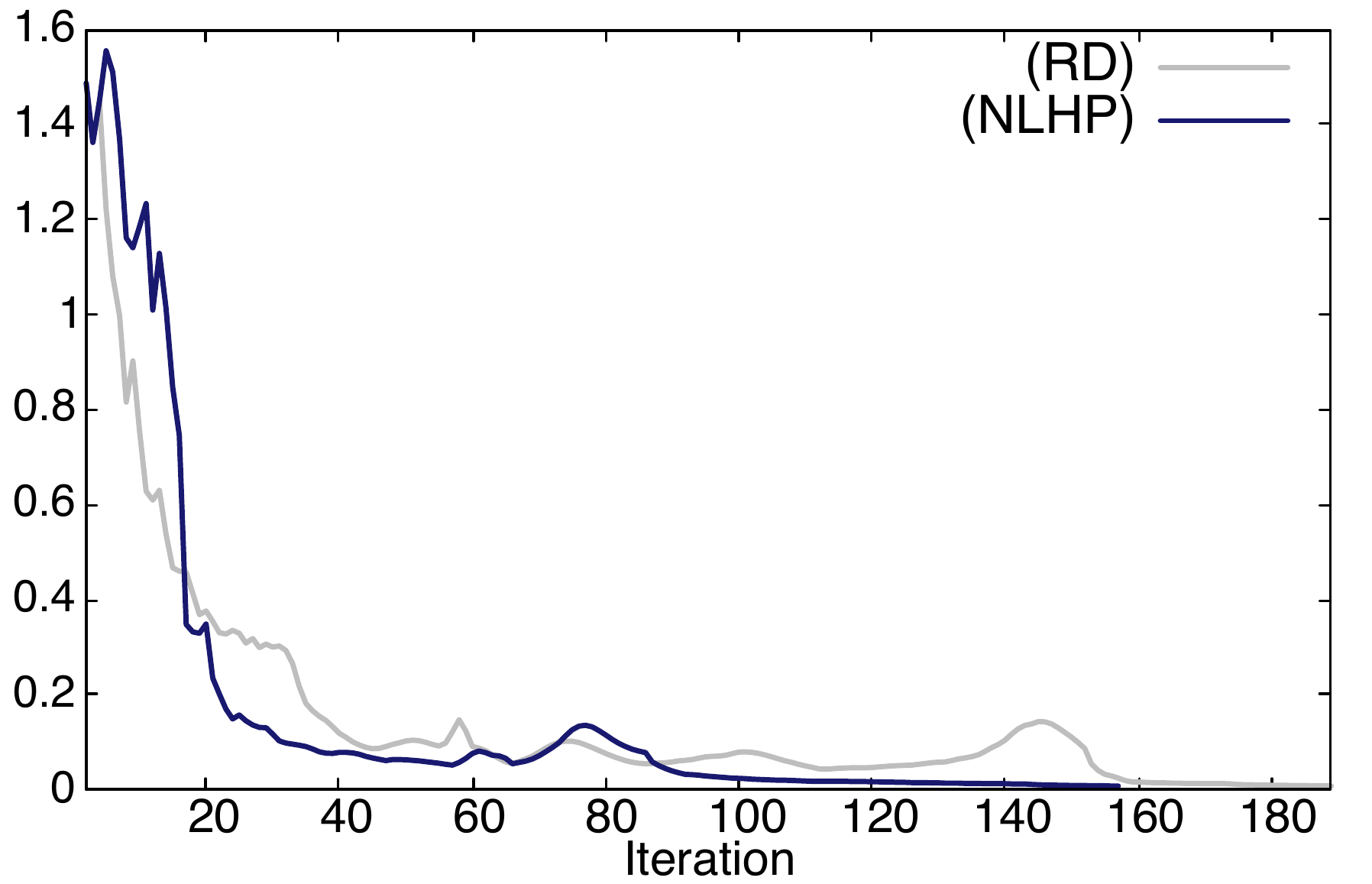}
        \subcaption{$\|\phi_{n+1}-\phi_n\|_{L^{\infty}(D)}$}
        \label{ra1-2}
      \end{minipage} &
      \end{tabular}
       \caption{ Objective functional and convergence condition for \S \ref{S:ra}-(i).}
    \label{ra1}
  \end{figure*}

\vspace{2mm}
\noindent{\bf Case (i) (Periodically perforated domain).\,} 
Let us consider the case where the initial configuration is the periodically perforated domain.
From Figures \ref{ra1-c} and \ref{ra1-h}, we infer that the topology of $\Omega_{\phi_n}\subset D$ in (NLHP) can be optimized faster than that of (RD), and then Figure \ref{ra1} yields the assertion in this study. Here we switched to the reaction-diffusion in (NLHP) after $15$ steps to avoid oscillation. 

\begin{figure*}[htbp]
\hspace*{-5mm}
    \begin{tabular}{ccccc}
      \begin{minipage}[t]{0.2\hsize}
        \centering
        \includegraphics[keepaspectratio, scale=0.09]{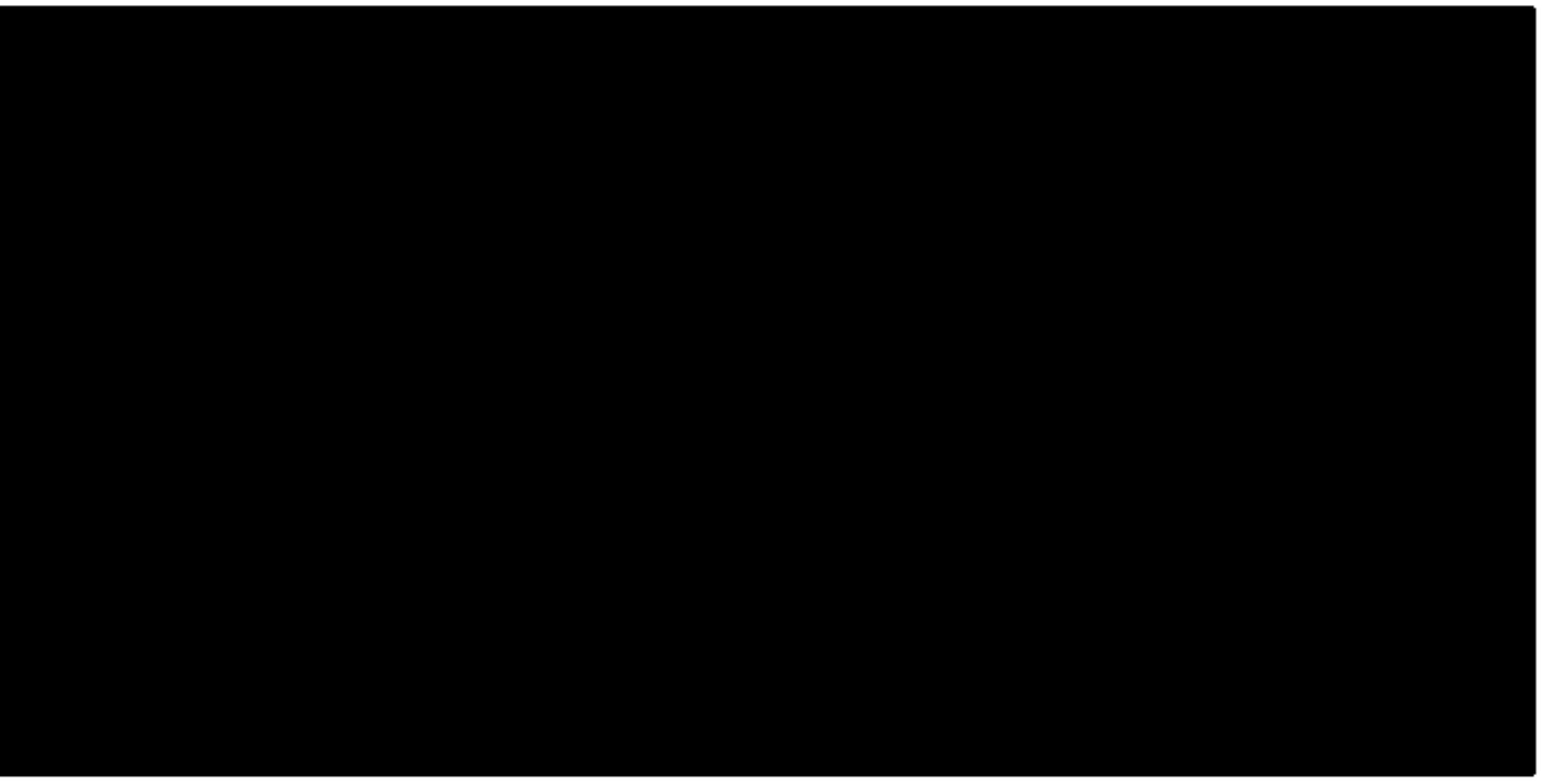}
        \subcaption{Step\,0}
        \label{ra2-a}
      \end{minipage} 
      \begin{minipage}[t]{0.2\hsize}
        \centering
        \includegraphics[keepaspectratio, scale=0.09]{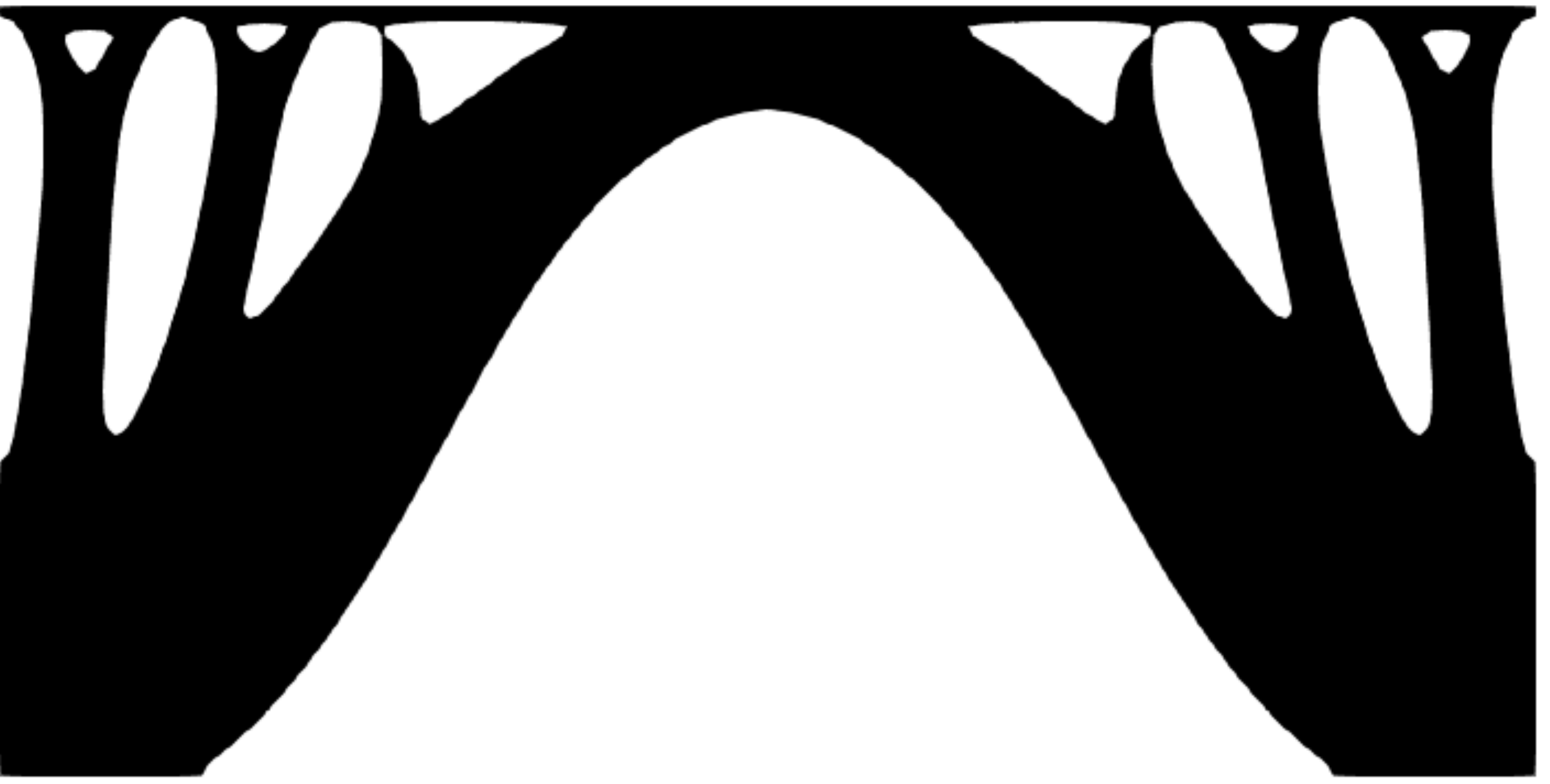}
        \subcaption{Step\,40}
        \label{ra2-b}
      \end{minipage} 
      \begin{minipage}[t]{0.2\hsize}
        \centering
        \includegraphics[keepaspectratio, scale=0.09]{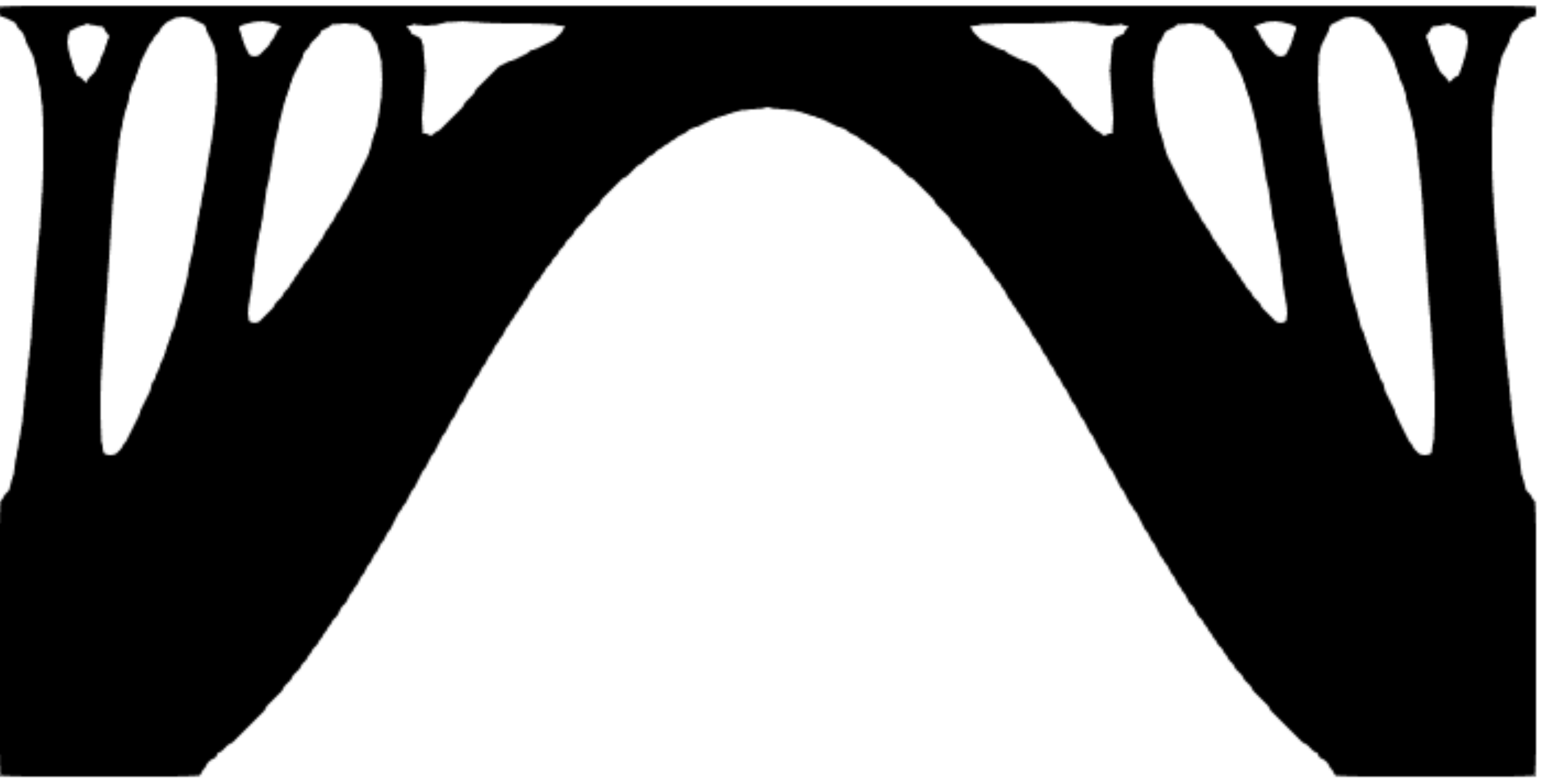}
        \subcaption{Step\,55}
        \label{ra2-c}
      \end{minipage} 
         \begin{minipage}[t]{0.2\hsize}
        \centering
        \includegraphics[keepaspectratio, scale=0.09]{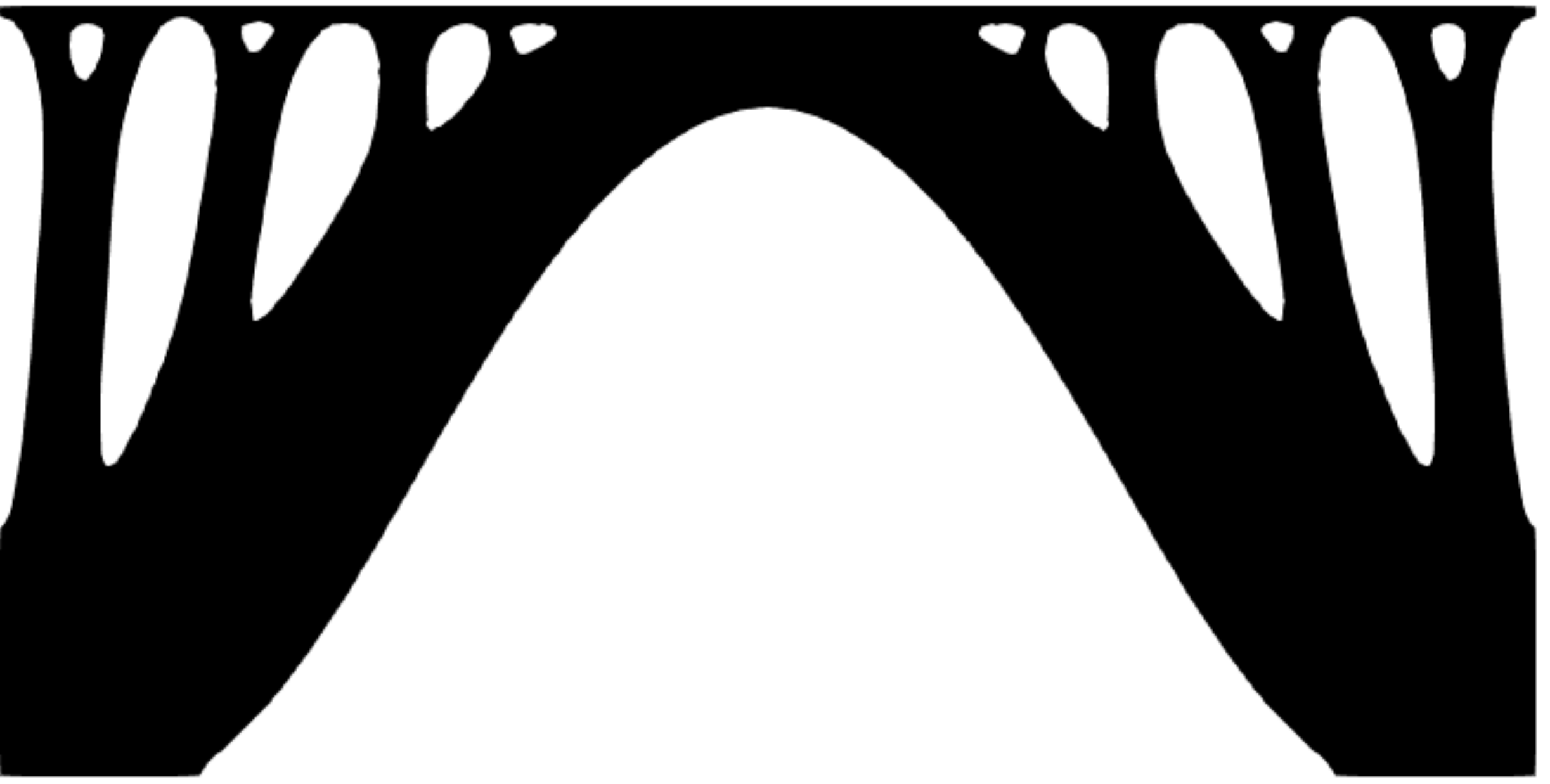}
        \subcaption{Step\,70}
        \label{ra2-d}
      \end{minipage} 
                 \begin{minipage}[t]{0.2\hsize}
        \centering
        \includegraphics[keepaspectratio, scale=0.09]{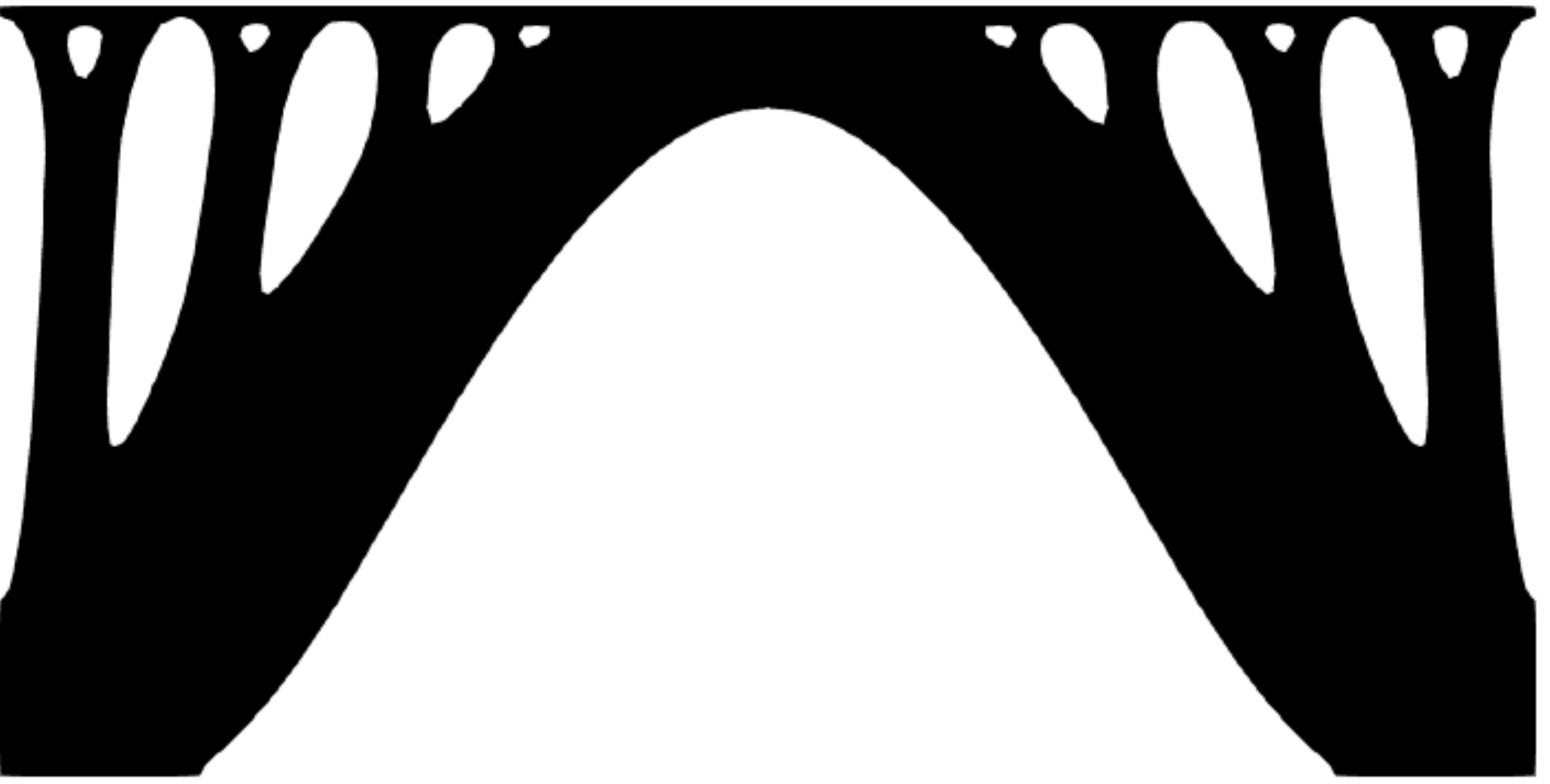}
        \subcaption{Step\,185$^{\#}$}
         \label{ra2-e}
      \end{minipage} 
      \\
    \begin{minipage}[t]{0.2\hsize}
        \centering
        \includegraphics[keepaspectratio, scale=0.09]{raw0.pdf}
        \subcaption{Step\,0}
        \label{ra2-f}
      \end{minipage} 
      \begin{minipage}[t]{0.2\hsize}
        \centering
        \includegraphics[keepaspectratio, scale=0.09]{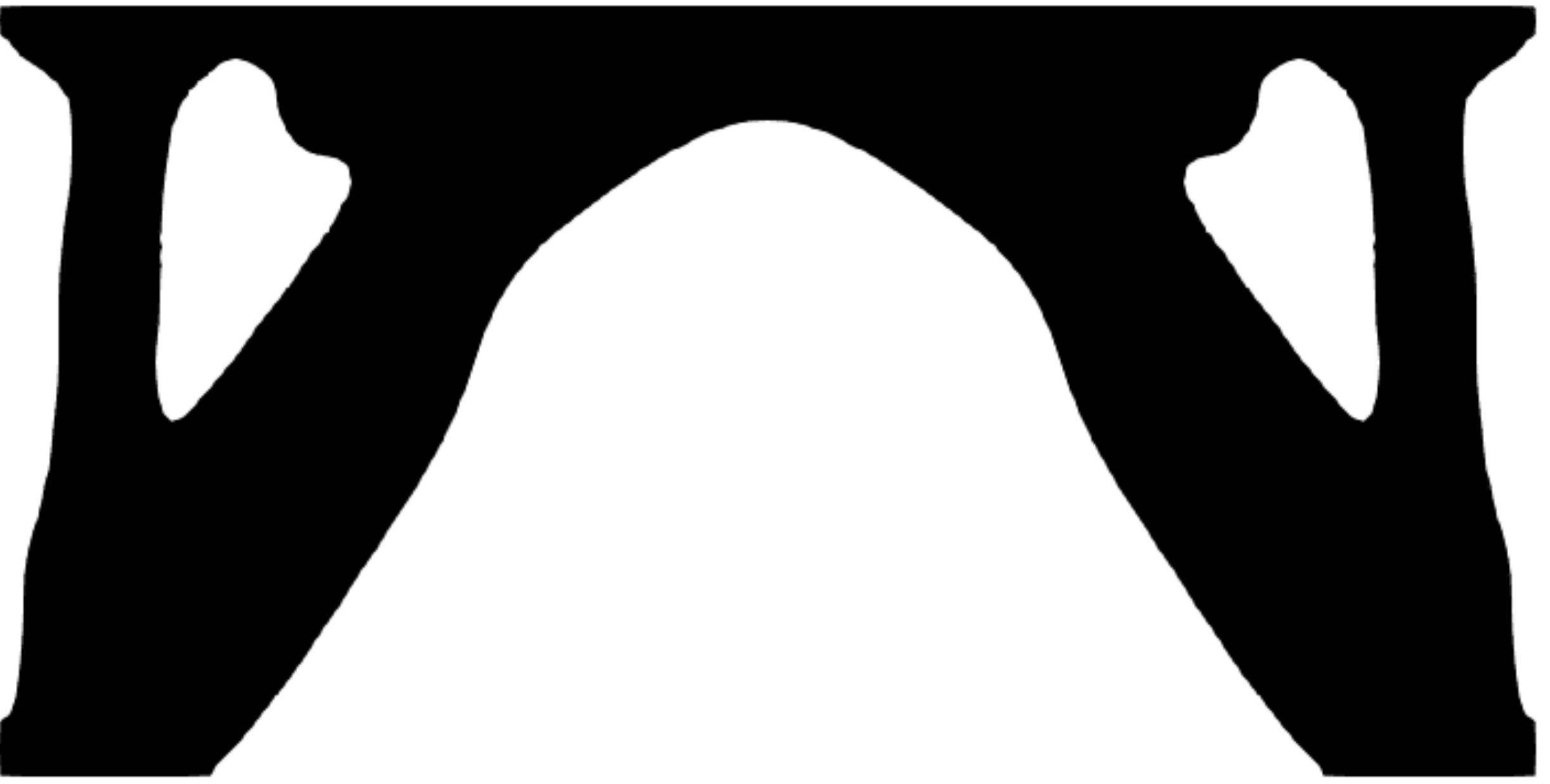}
        \subcaption{Step\,40}
        \label{ra2-g}
      \end{minipage} 
      \begin{minipage}[t]{0.2\hsize}
        \centering
        \includegraphics[keepaspectratio, scale=0.09]{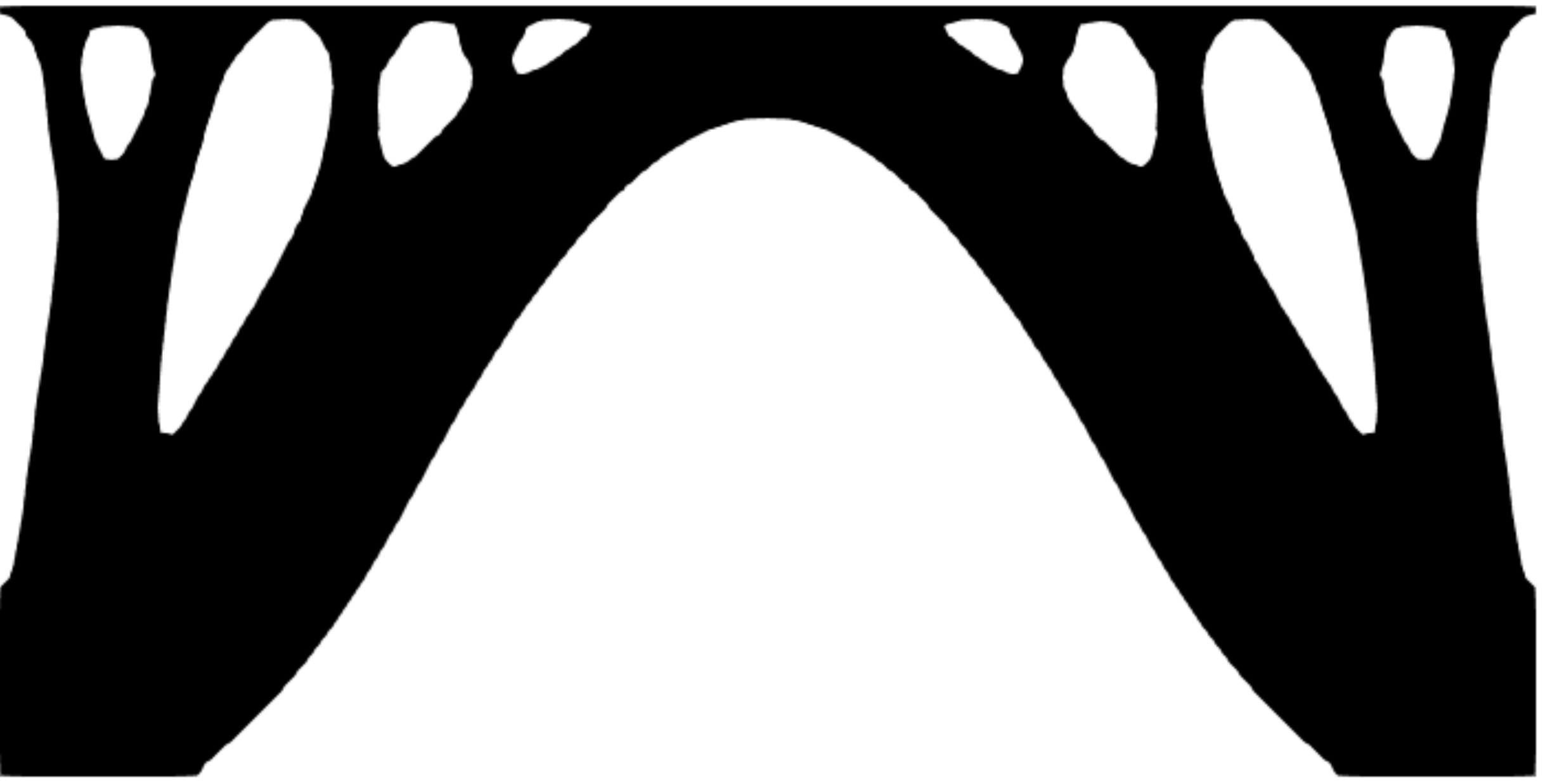}
        \subcaption{Step\,55}
        \label{ra2-h}
      \end{minipage} 
       \begin{minipage}[t]{0.2\hsize}
        \centering
        \includegraphics[keepaspectratio, scale=0.09]{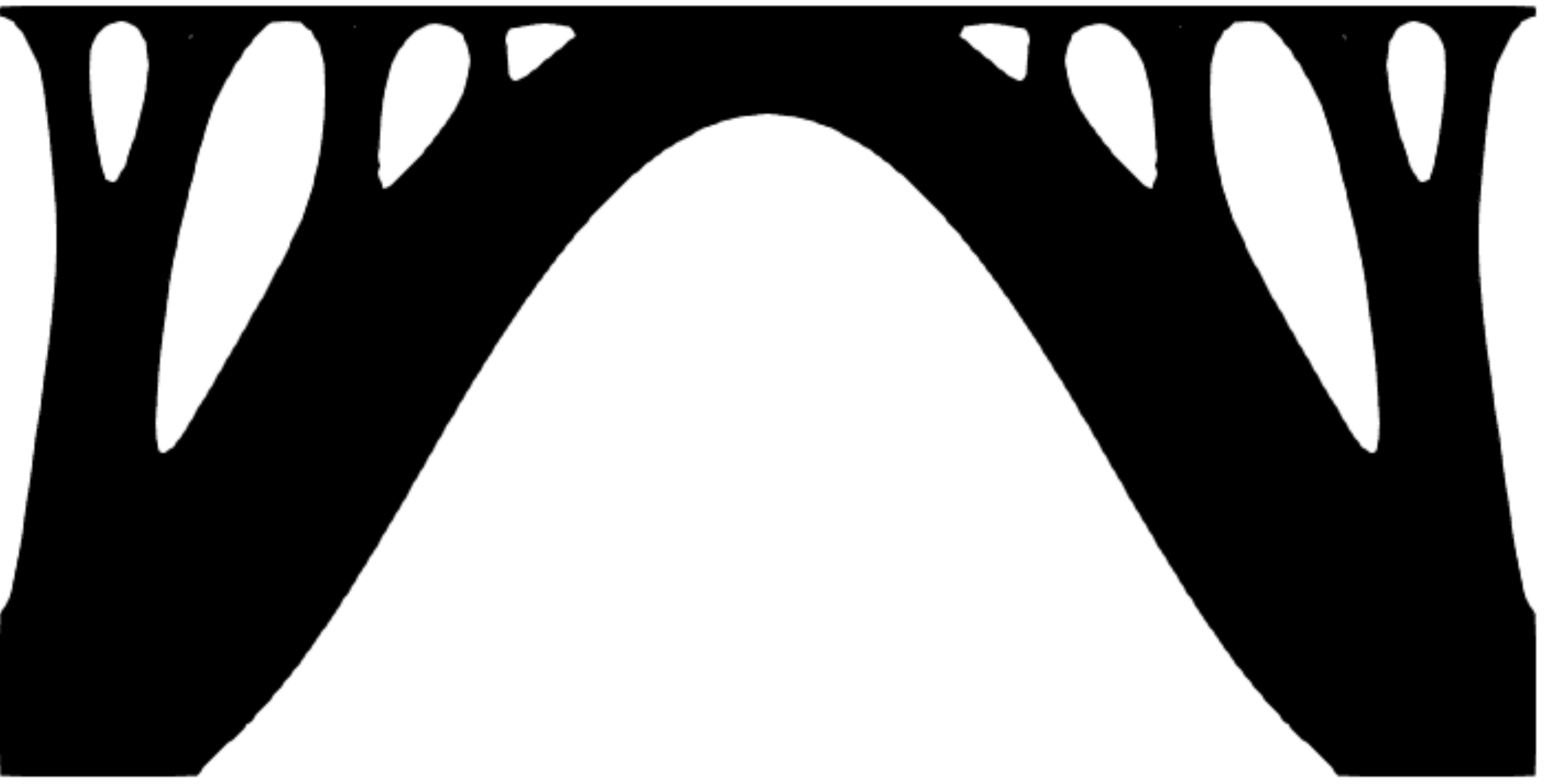}
        \subcaption{Step\,70}
        \label{ra2-i}
      \end{minipage}
                 \begin{minipage}[t]{0.2\hsize}
        \centering
        \includegraphics[keepaspectratio, scale=0.09]{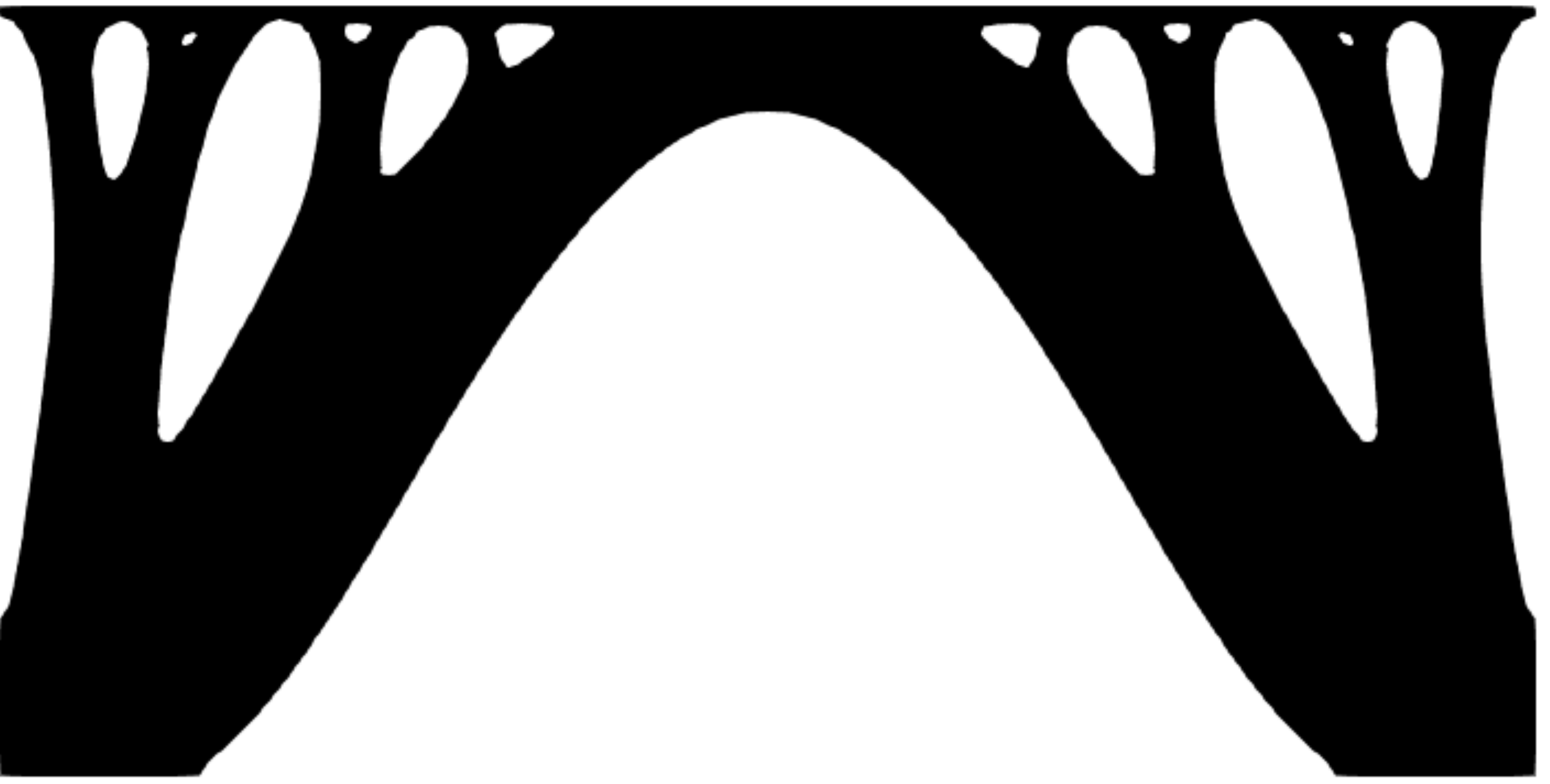}
        \subcaption{Step\,114$^{\#}$}
        \label{ra2-j}
      \end{minipage}  
    \end{tabular}
     \caption{ Configuration $\Omega_{\phi_n}\subset D$ for the case where the initial configuration is the whole domain. 
     Figures (a)--(e) and (f)--(j) represent the results of (RD) and (NLHP), respectively.     
The symbol $^{\#}$ implies the final step.}
     \label{fig:ra2}
  \end{figure*}

\begin{figure*}[htbp]
    \begin{tabular}{ccc}
      \hspace*{-5mm} 
      \begin{minipage}[t]{0.49\hsize}
        \centering
        \includegraphics[keepaspectratio, scale=0.33]{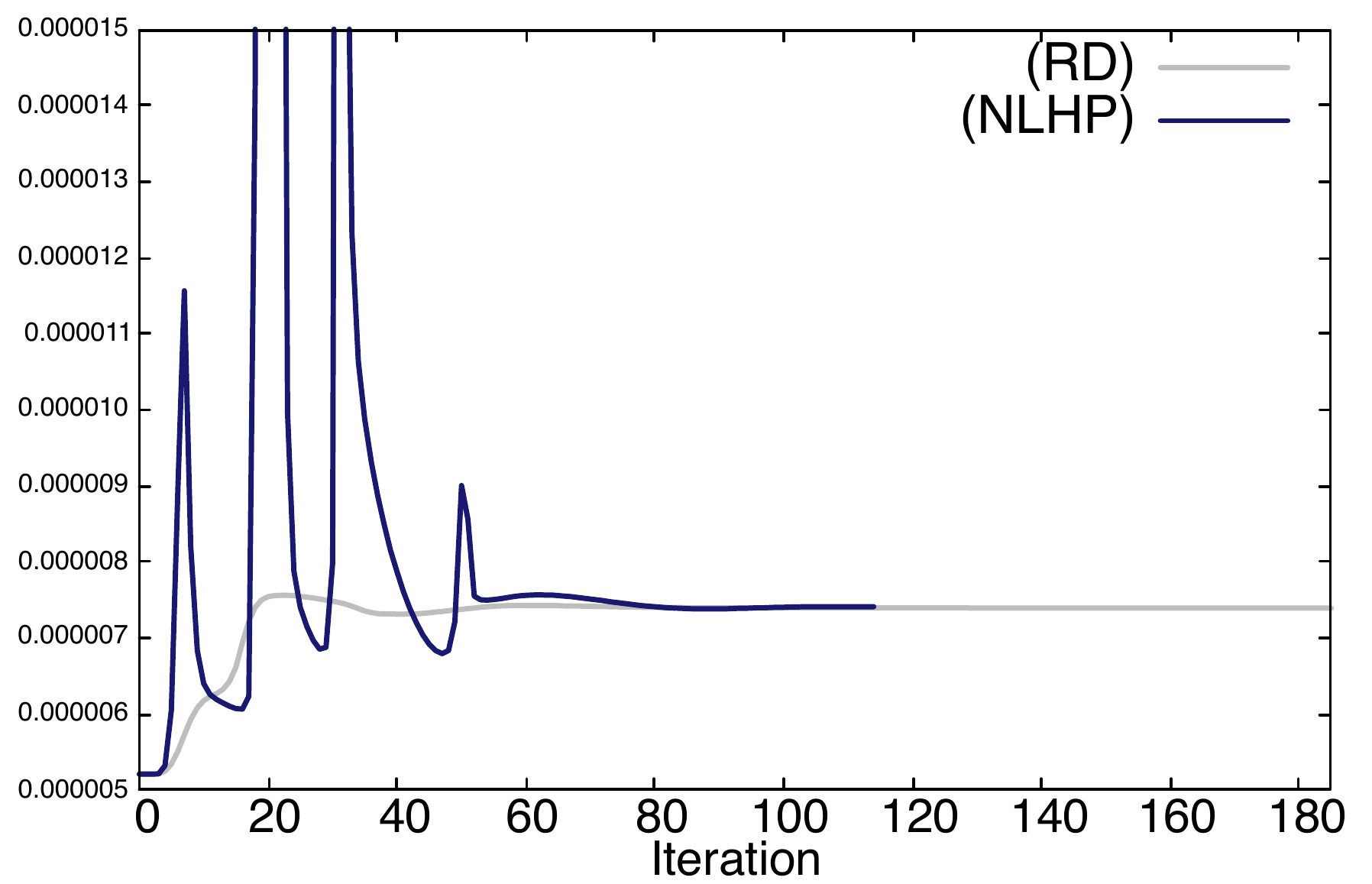}
        \subcaption{$F(\phi_n)$}
        \label{ra2-1}
      \end{minipage} 
      \begin{minipage}[t]{0.49\hsize}
        \centering
        \includegraphics[keepaspectratio, scale=0.33]{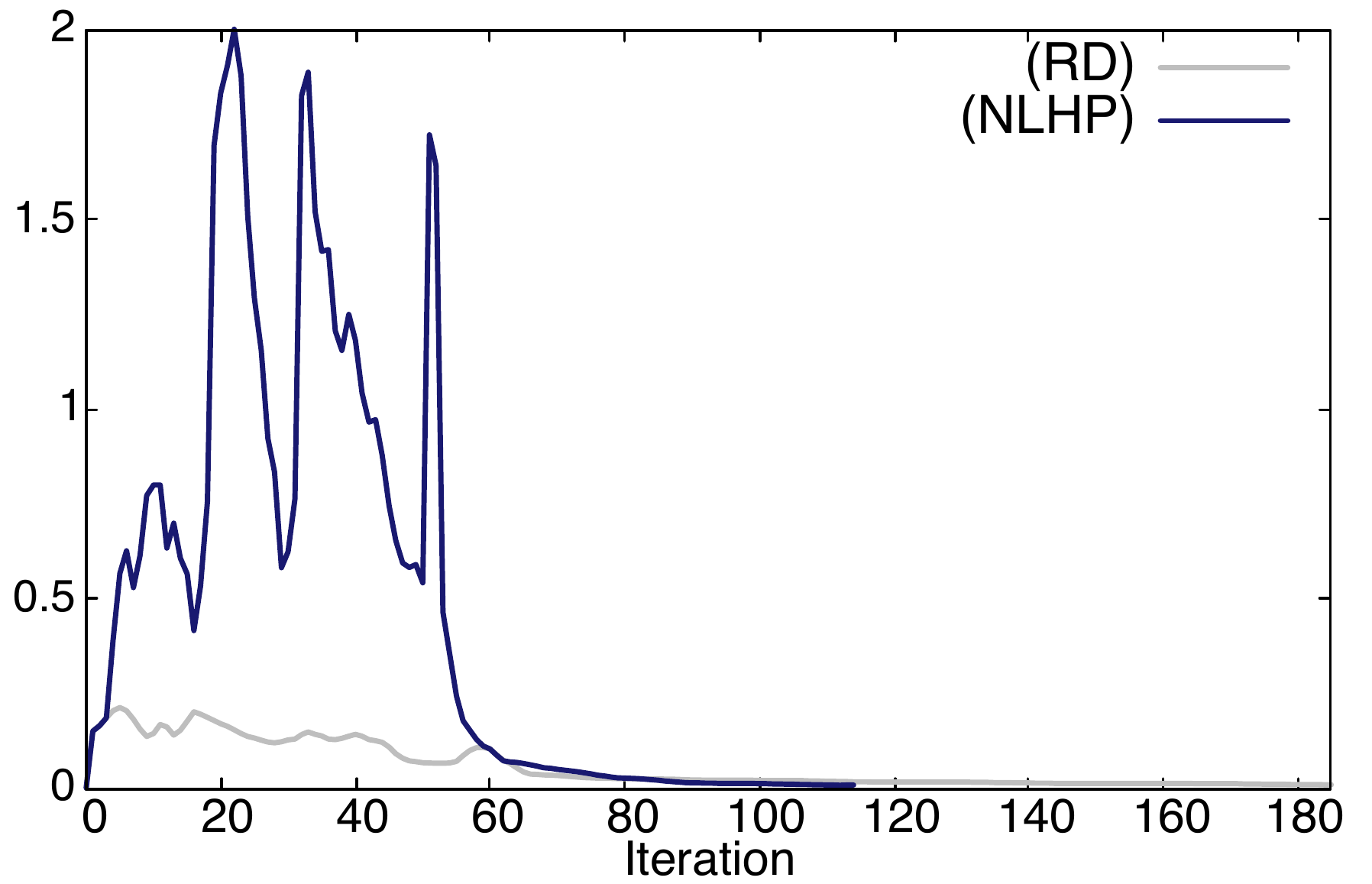}
        \subcaption{$\|\phi_{n+1}-\phi_n\|_{L^{\infty}(D)}$}
        \label{ra2-2}
      \end{minipage} &
      \end{tabular}
       \caption{ Objective functional and convergence condition \S \ref{S:ra}-(ii).}
    \label{ra2}
  \end{figure*}

\vspace{2mm}
\noindent{\bf Case (ii) (Whole domain).\,} 
As for the case where the initial configuration is the whole domain, 
we first note that optimal configurations in Figure \ref{ra2-j} are not the same as that in Figure \ref{ra2-e}; however, from Figure \ref{ra2-1}, we see that the convergence values are the same, and (NLHP) is optimized to the same degree as (RD).     
On the other hand, by noting that objective functional $F(\phi_n)$ is oscillating in Figure \ref{ra2-a}  and switching to the reaction-diffusion in 50 steps, the convergence condition can be satisfied quickly, which implies that the use of a slightly nonlinear damped wave equation improves convergence. 

\begin{figure*}[htbp]
\hspace*{-5mm}
    \begin{tabular}{ccccc}
      \begin{minipage}[t]{0.2\hsize}
        \centering
        \includegraphics[keepaspectratio, scale=0.09]{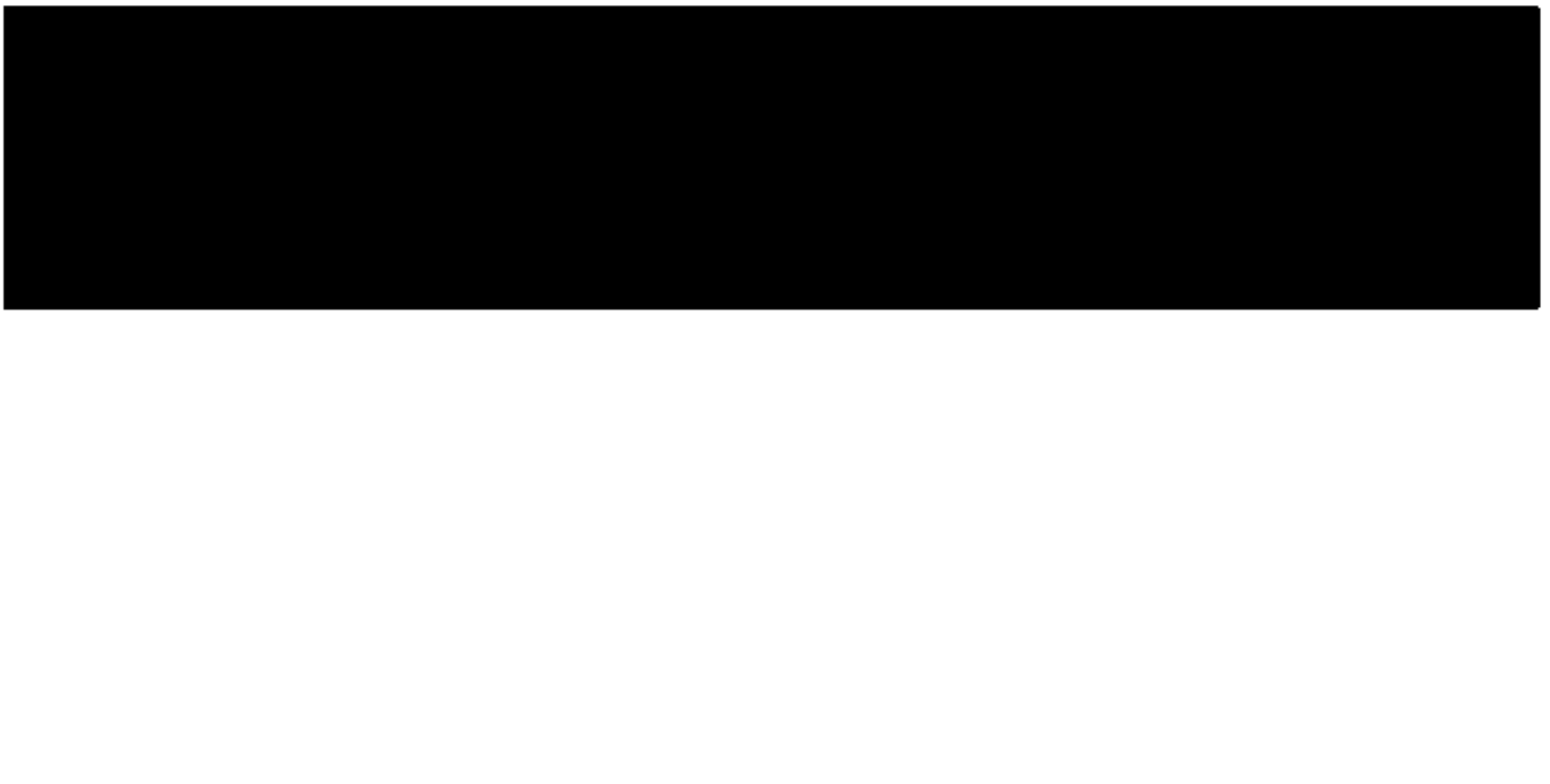}
        \subcaption{Step\,0}
        \label{ra3-a}
      \end{minipage} 
      \begin{minipage}[t]{0.2\hsize}
        \centering
        \includegraphics[keepaspectratio, scale=0.09]{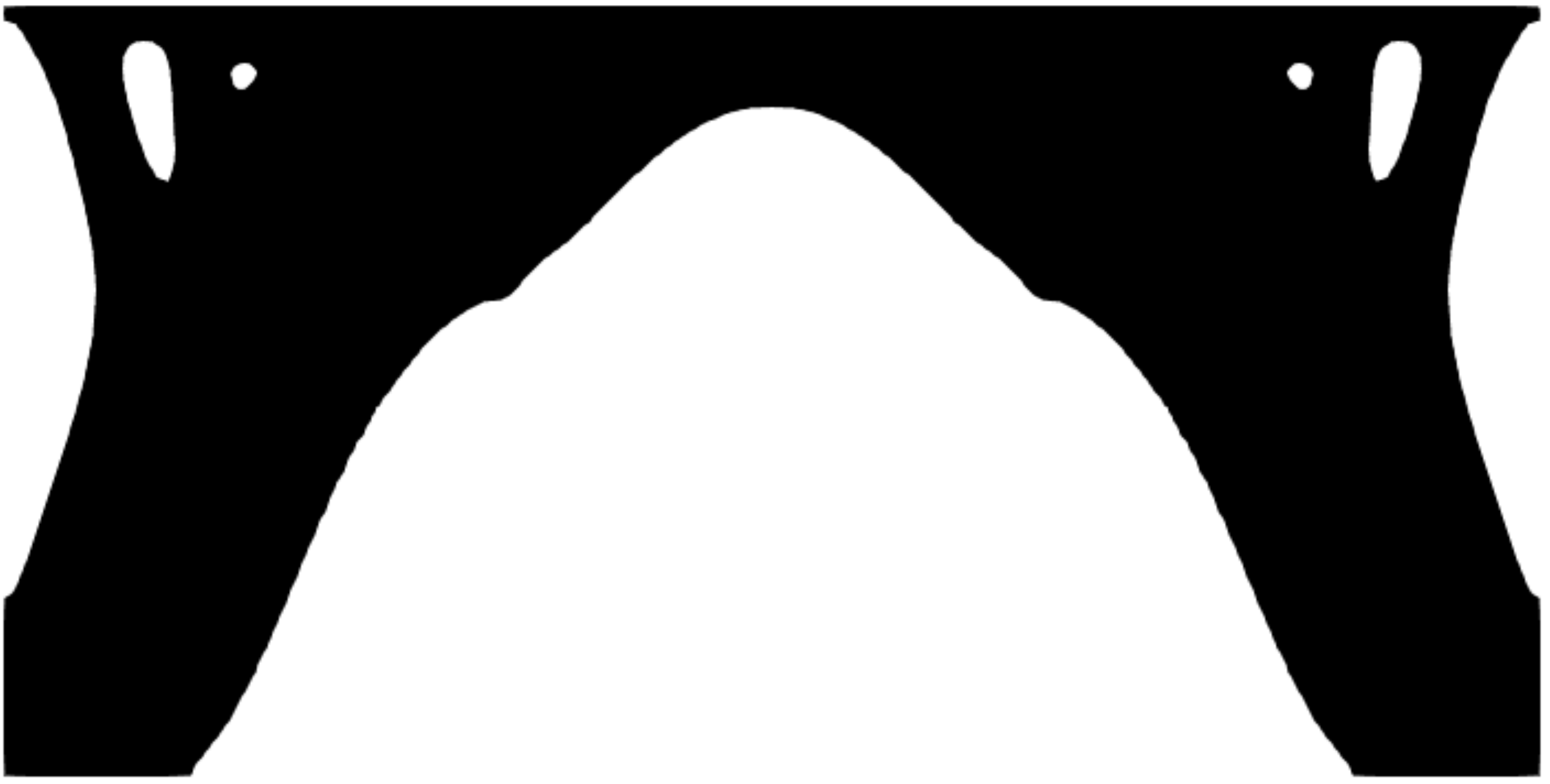}
        \subcaption{Step\,70}
        \label{ra3-b}
      \end{minipage} 
      \begin{minipage}[t]{0.2\hsize}
        \centering
        \includegraphics[keepaspectratio, scale=0.09]{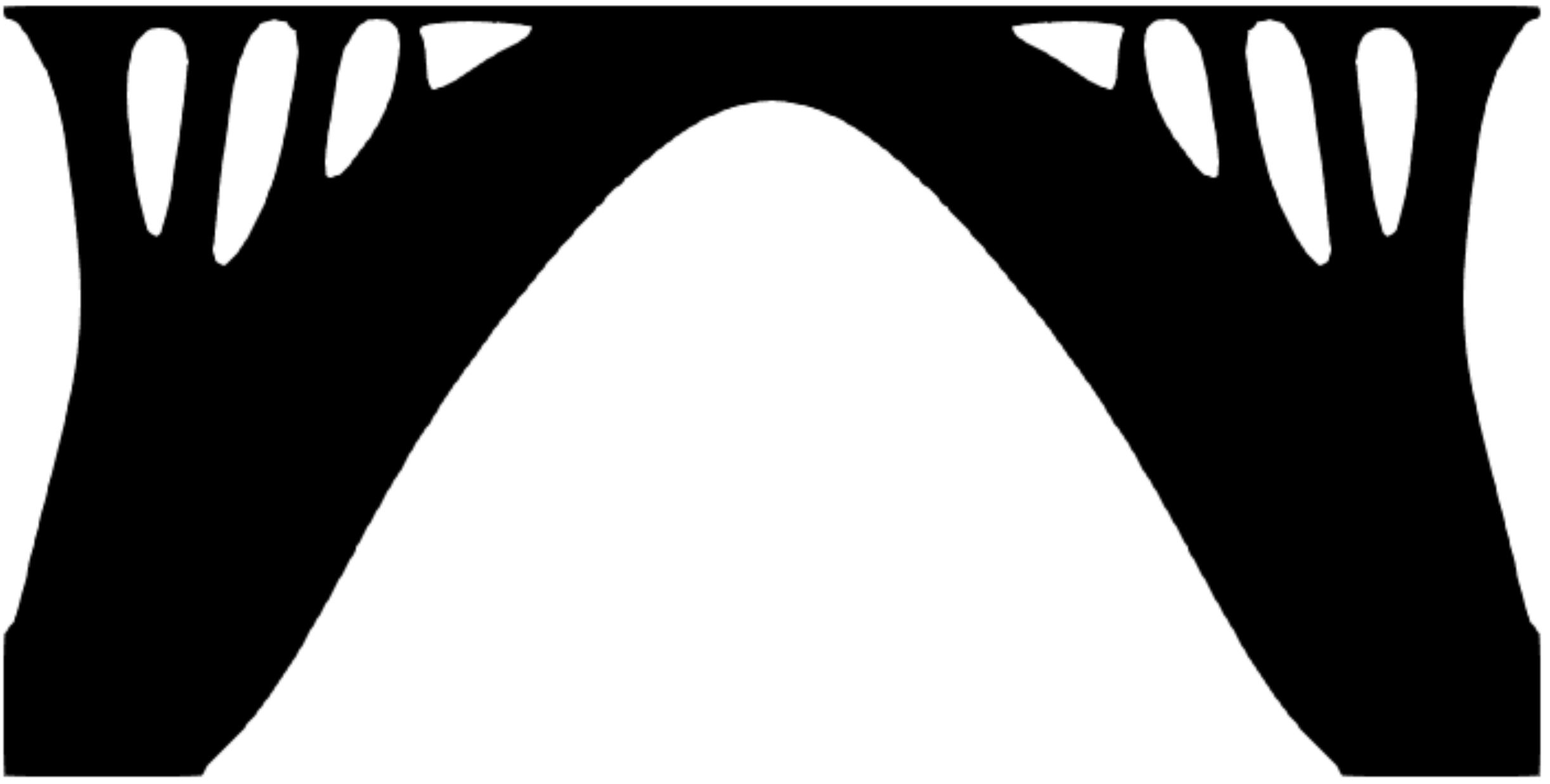}
        \subcaption{Step\,120}
        \label{ra3-c}
      \end{minipage} 
         \begin{minipage}[t]{0.2\hsize}
        \centering
        \includegraphics[keepaspectratio, scale=0.09]{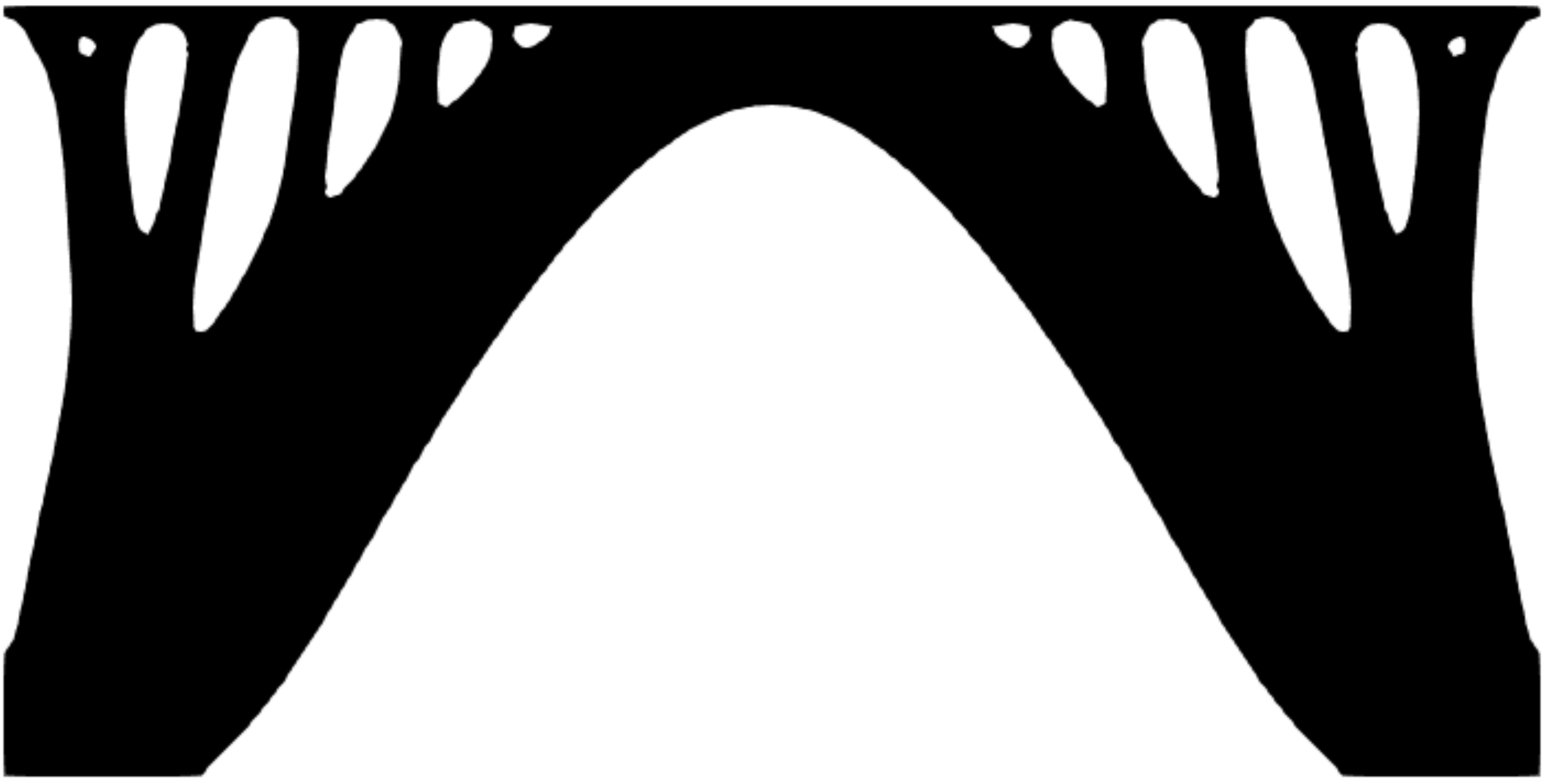}
        \subcaption{Step\,170}
        \label{ra3-d}
      \end{minipage} 
                 \begin{minipage}[t]{0.2\hsize}
        \centering
        \includegraphics[keepaspectratio, scale=0.09]{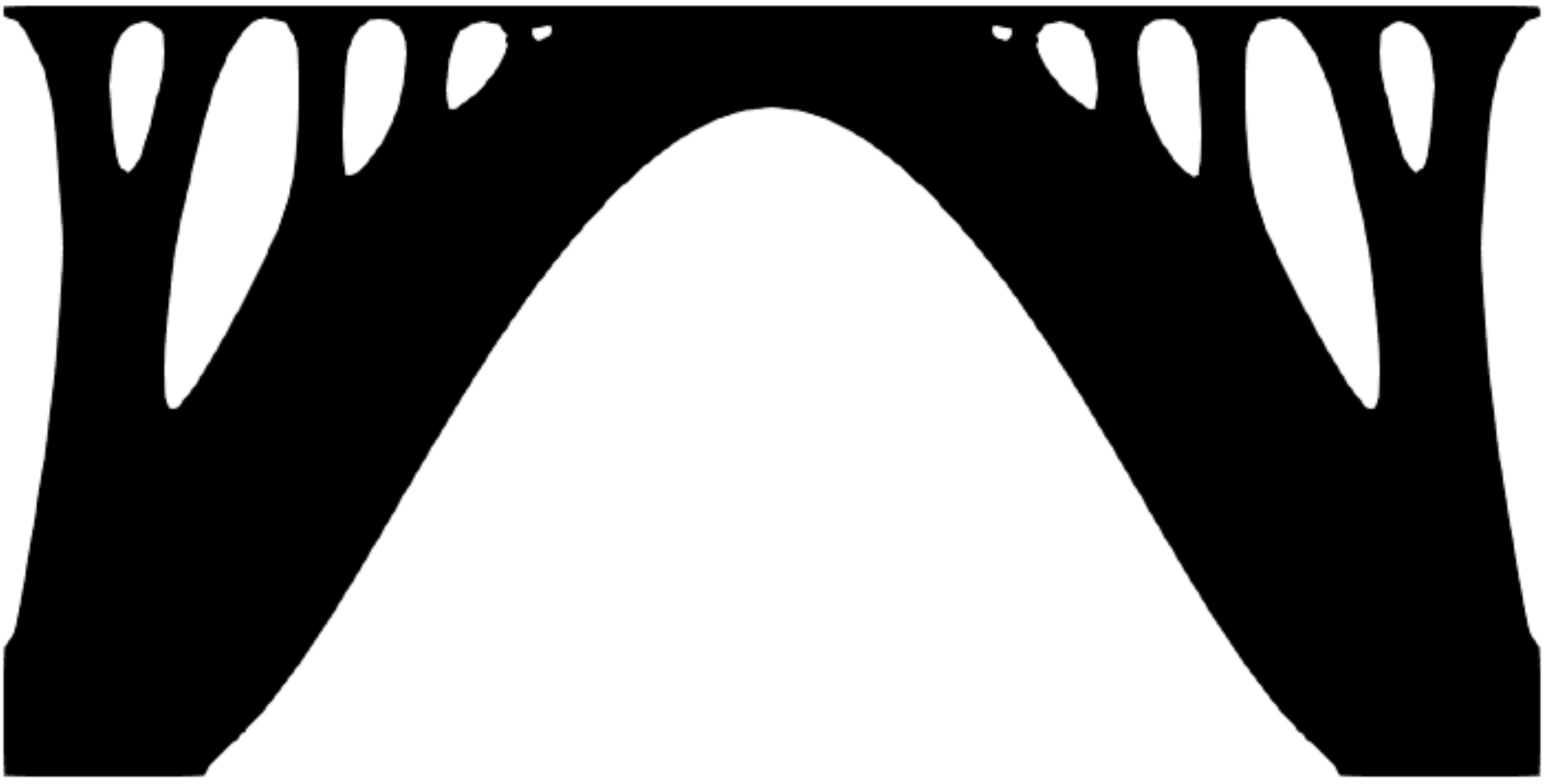}
        \subcaption{Step\,412$^{\#}$}
        \label{ra3-e}
      \end{minipage} 
      \\
    \begin{minipage}[t]{0.2\hsize}
        \centering
        \includegraphics[keepaspectratio, scale=0.09]{rau0.pdf}
        \subcaption{Step\,0}
        \label{ra3-f}
      \end{minipage} 
      \begin{minipage}[t]{0.2\hsize}
        \centering
        \includegraphics[keepaspectratio, scale=0.09]{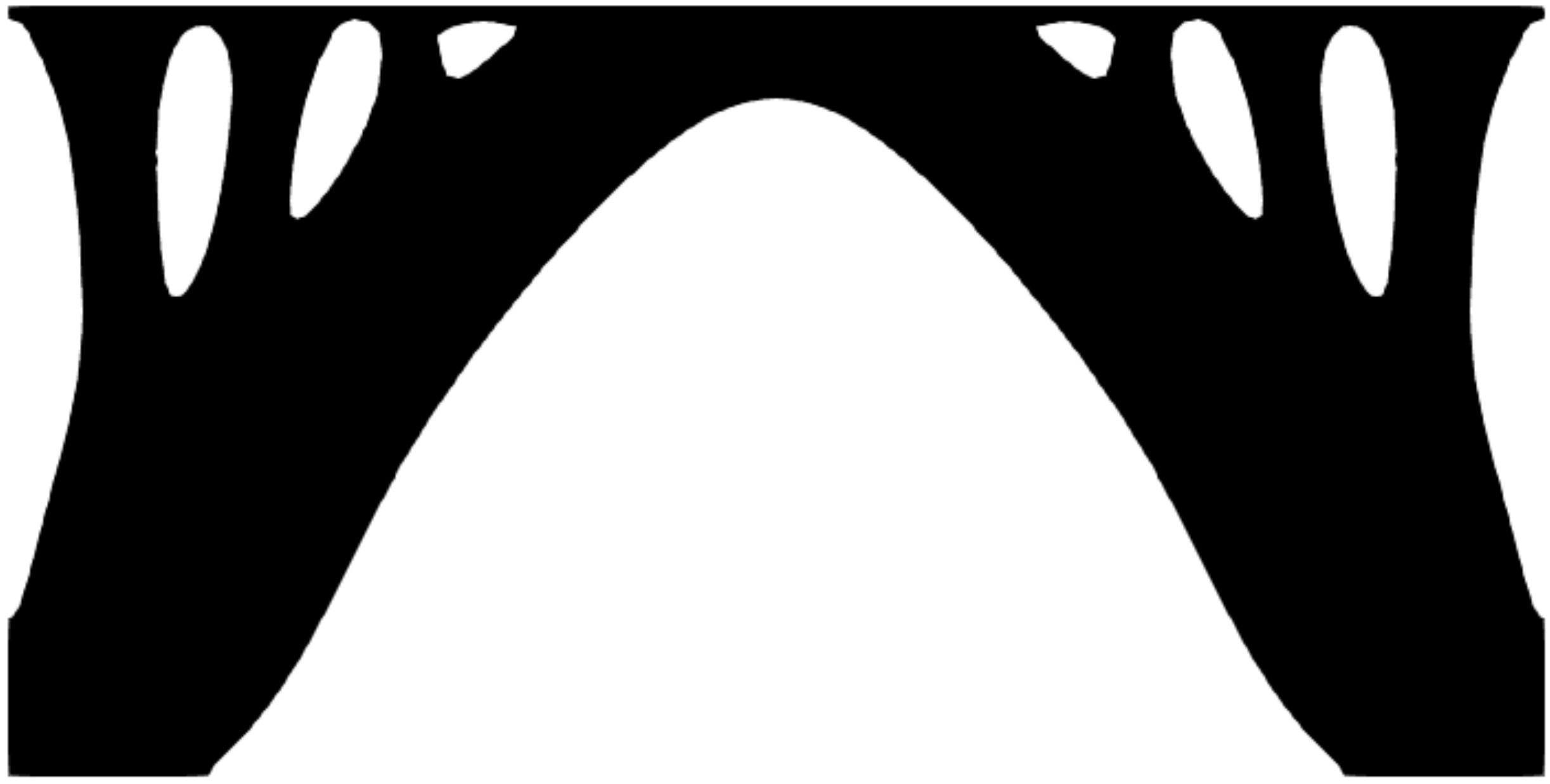}
        \subcaption{Step\,70}
        \label{ra3-g}
      \end{minipage} 
      \begin{minipage}[t]{0.2\hsize}
        \centering
        \includegraphics[keepaspectratio, scale=0.09]{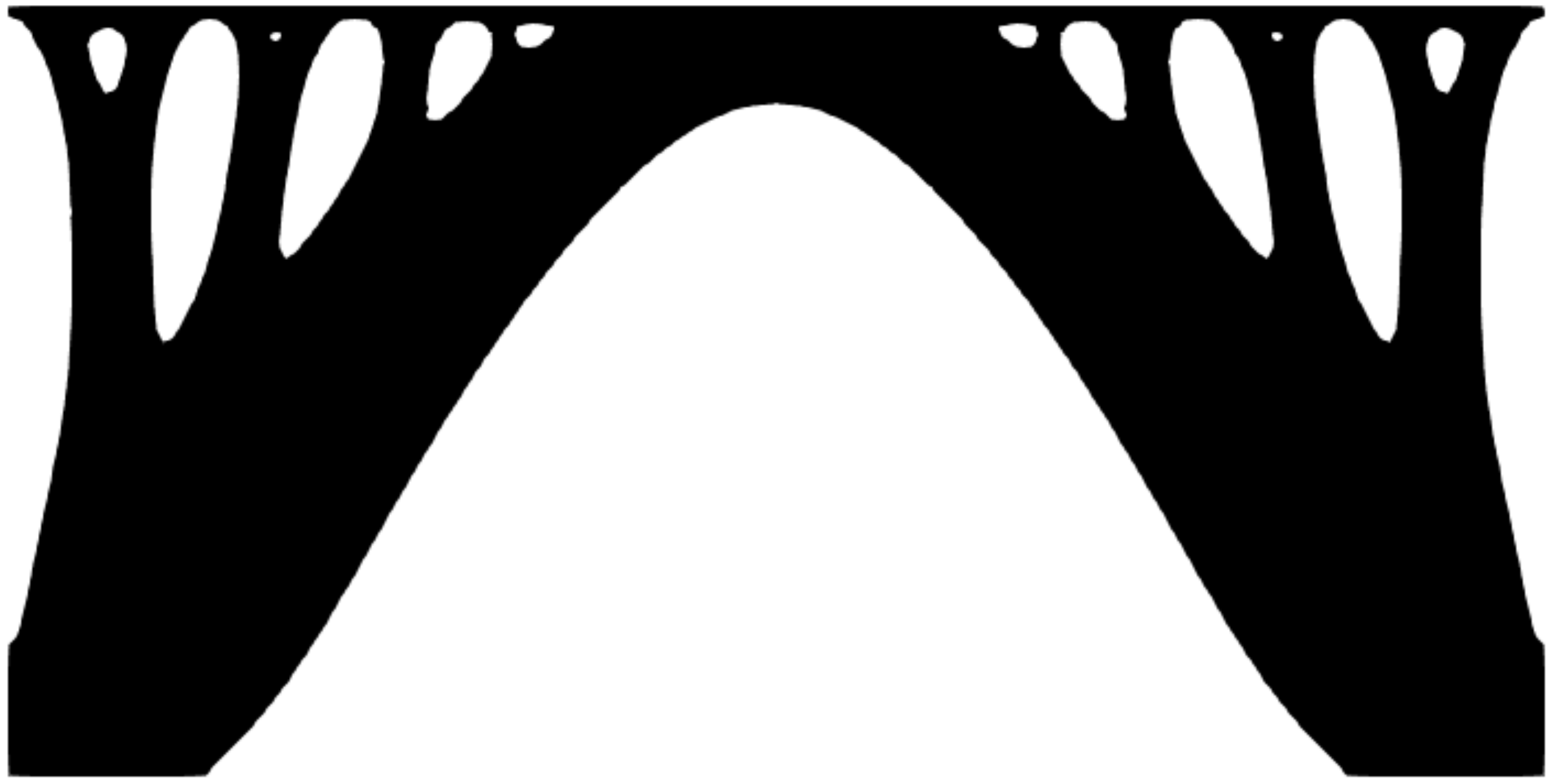}
        \subcaption{Step\,120}
        \label{ra3-h}
      \end{minipage} 
       \begin{minipage}[t]{0.2\hsize}
        \centering
        \includegraphics[keepaspectratio, scale=0.09]{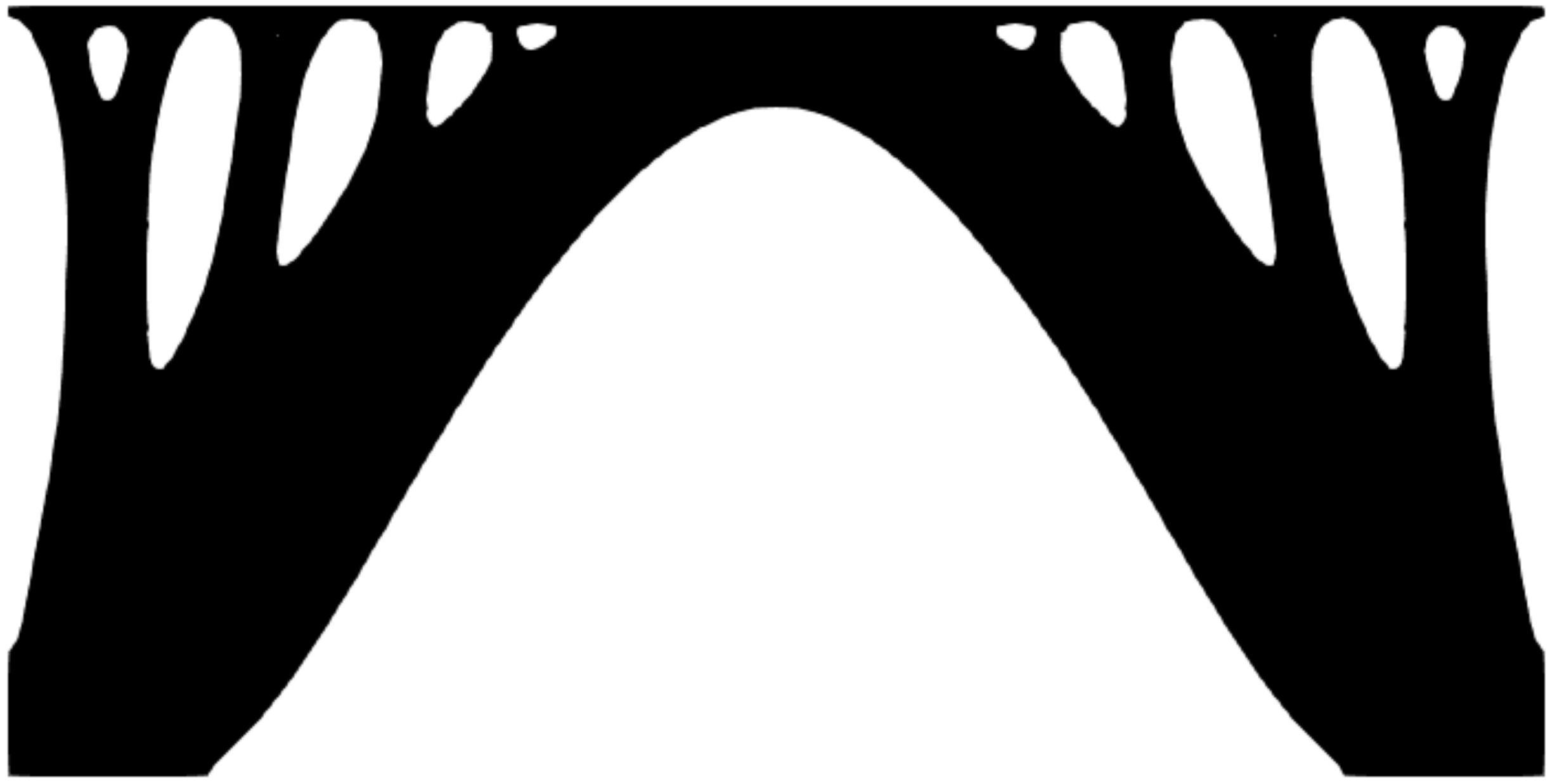}
        \subcaption{Step\,170}
        \label{ra3-i}
      \end{minipage}
                 \begin{minipage}[t]{0.2\hsize}
        \centering
        \includegraphics[keepaspectratio, scale=0.09]{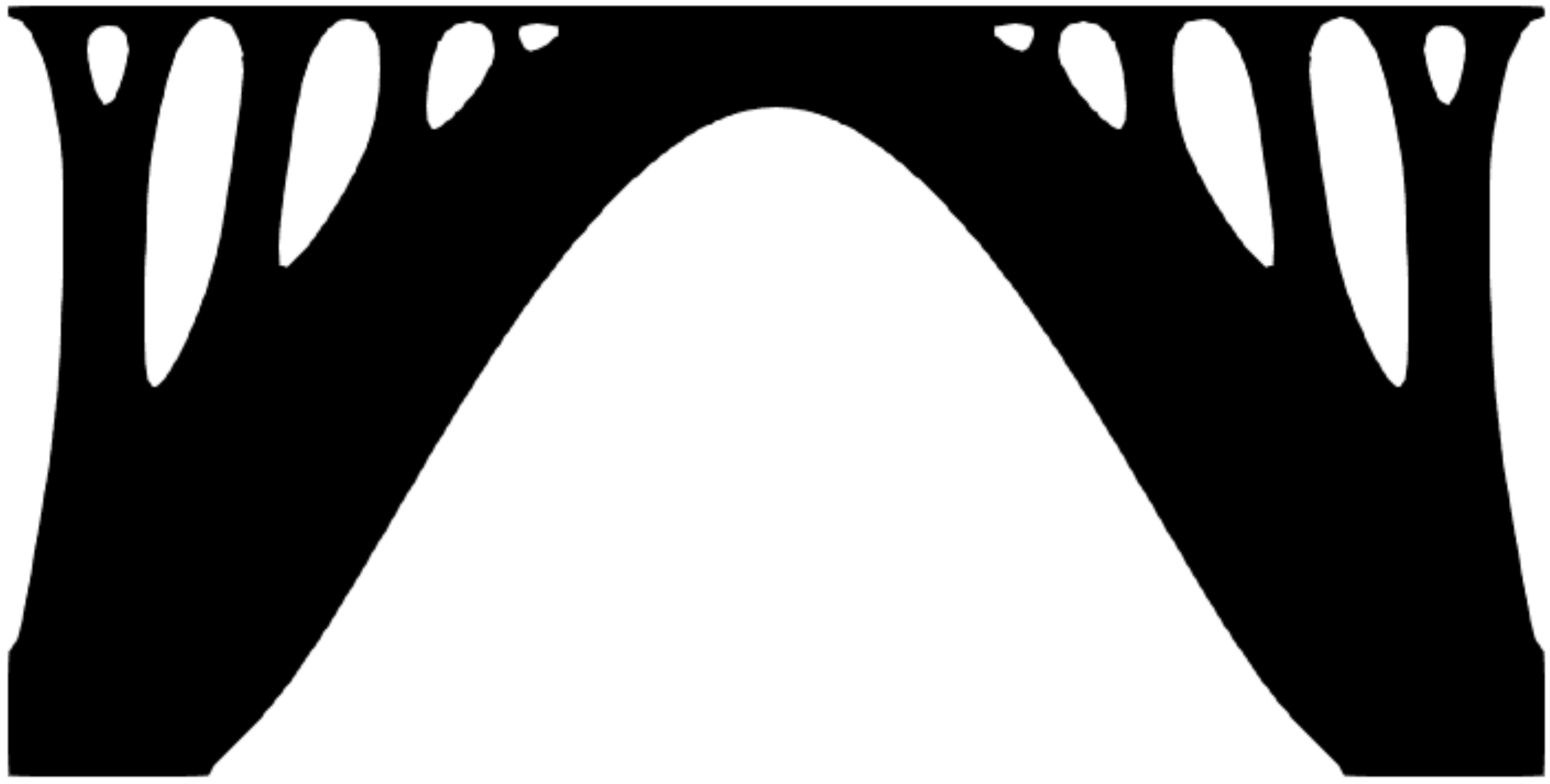}
        \subcaption{Step\,219$^{\#}$}
        \label{ra3-j}
      \end{minipage}  
    \end{tabular}
     \caption{ Configuration $\Omega_{\phi_n}\subset D$ for the case where the initial configuration is the upper domain. 
     Figures (a)--(e) and (f)--(j) represent the results of (RD) and (NLHP), respectively.     
The symbol $^{\#}$ implies the final step.}
     \label{fig:ra3}
  \end{figure*}

\begin{figure*}[htbp]
    \begin{tabular}{ccc}
      \hspace*{-5mm} 
      \begin{minipage}[t]{0.49\hsize}
        \centering
        \includegraphics[keepaspectratio, scale=0.33]{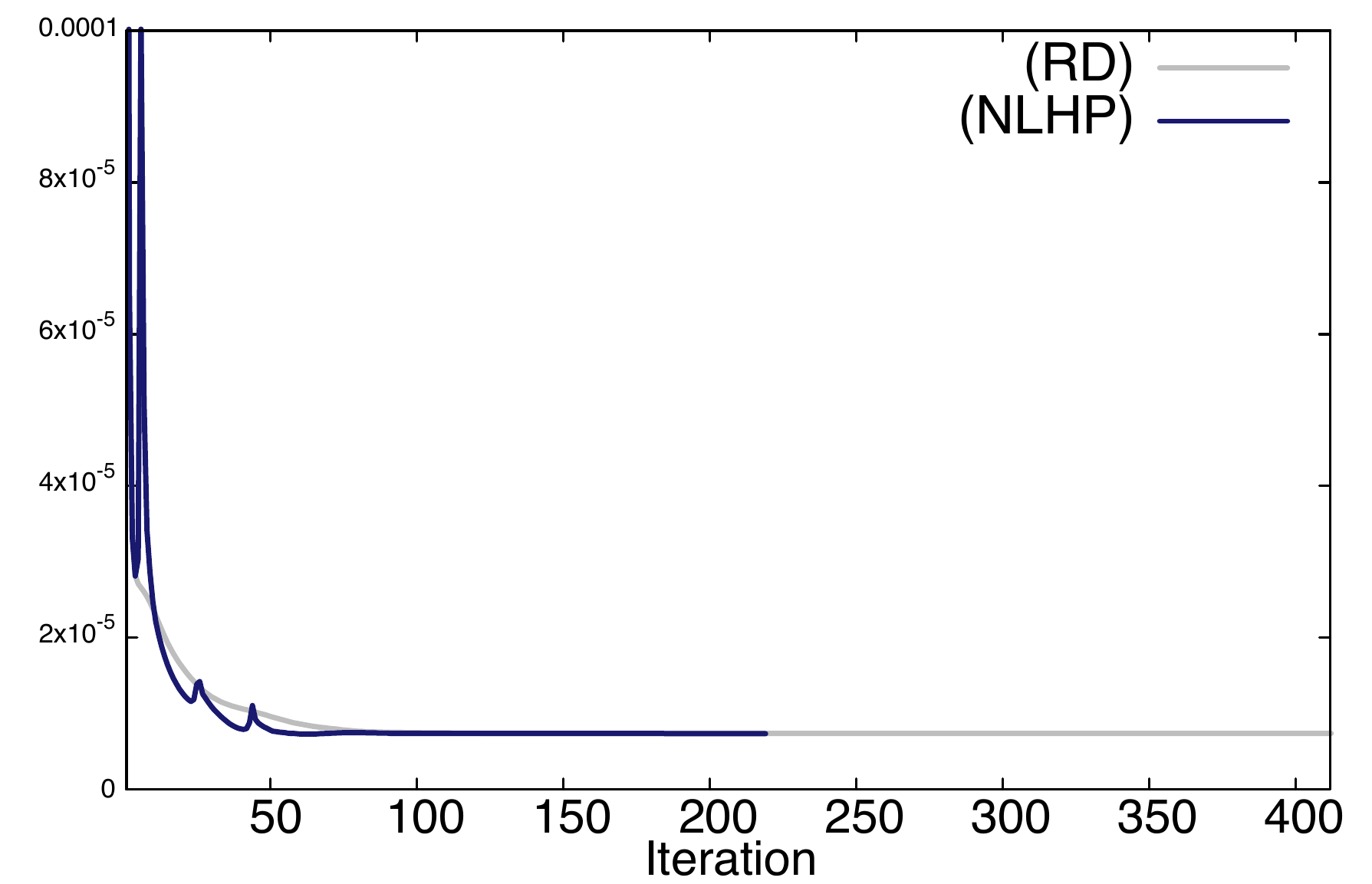}
        \subcaption{$F(\phi_n)$}
        \label{ra3-1}
      \end{minipage} 
      \begin{minipage}[t]{0.49\hsize}
        \centering
        \includegraphics[keepaspectratio, scale=0.33]{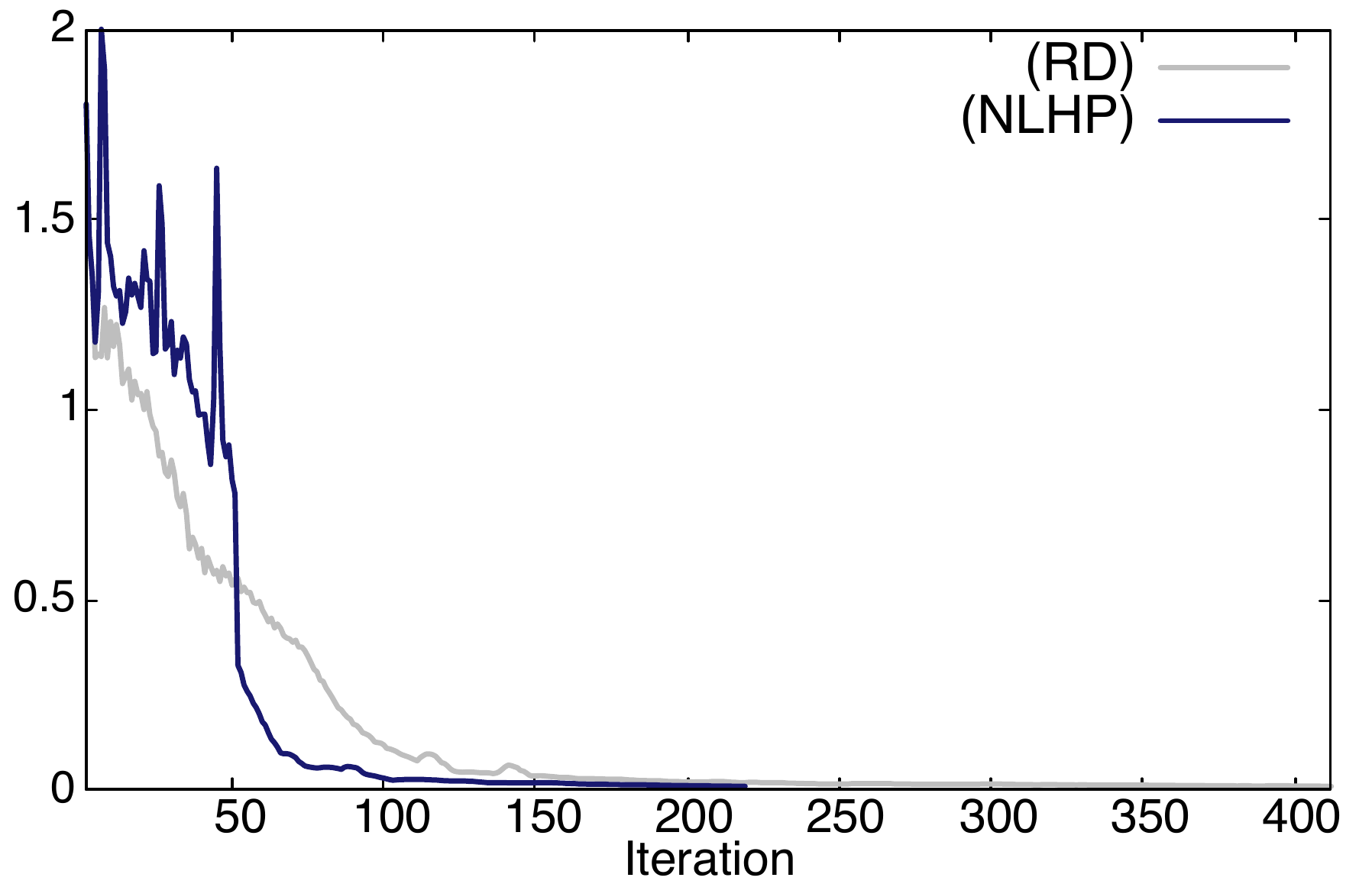}
        \subcaption{$\|\phi_{n+1}-\phi_n\|_{L^{\infty}(D)}$}
        \label{ra3-2}
      \end{minipage} &
      \end{tabular}
       \caption{ Objective functional and convergence condition \S \ref{S:ra}-(iii).}
    \label{ra3}
  \end{figure*}

\vspace{2mm}
\noindent{\bf Case (iii) (Upper domain).\,} 
We take the initial configuration as the upper domain. 
Then Figures \ref{fig:ra3} and \ref{ra3} are obtained.
As in the previous case, since the objective functional $F(\phi_n)$ in (NLHP) oscillates, we also switched to the reaction-diffusion at 50 steps.
By comparing Figures \ref{ra3-c} with \ref{ra3-h}, one can see that the topology in (NLHP) is optimized faster than that in (RD), and moreover, it is noteworthy that convergence is greatly improved.

\begin{figure*}[htbp]
\hspace*{-5mm}
    \begin{tabular}{cccc}
      \begin{minipage}[t]{0.24\hsize}
        \centering
        \includegraphics[keepaspectratio, scale=0.11]{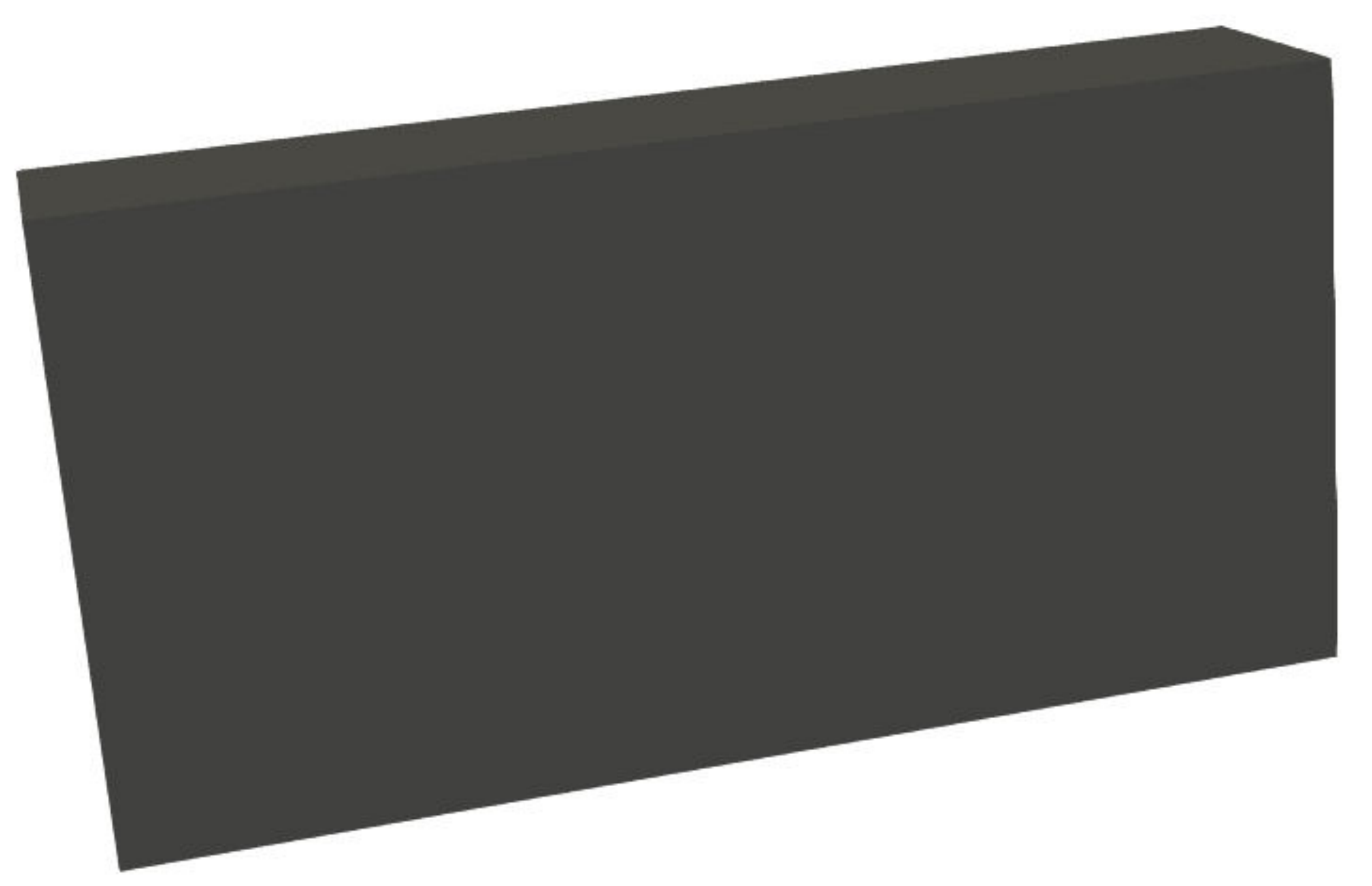}
        \subcaption{Step\,0}
        \label{ra3d-a}
      \end{minipage} 
      \begin{minipage}[t]{0.24\hsize}
        \centering
        \includegraphics[keepaspectratio, scale=0.11]{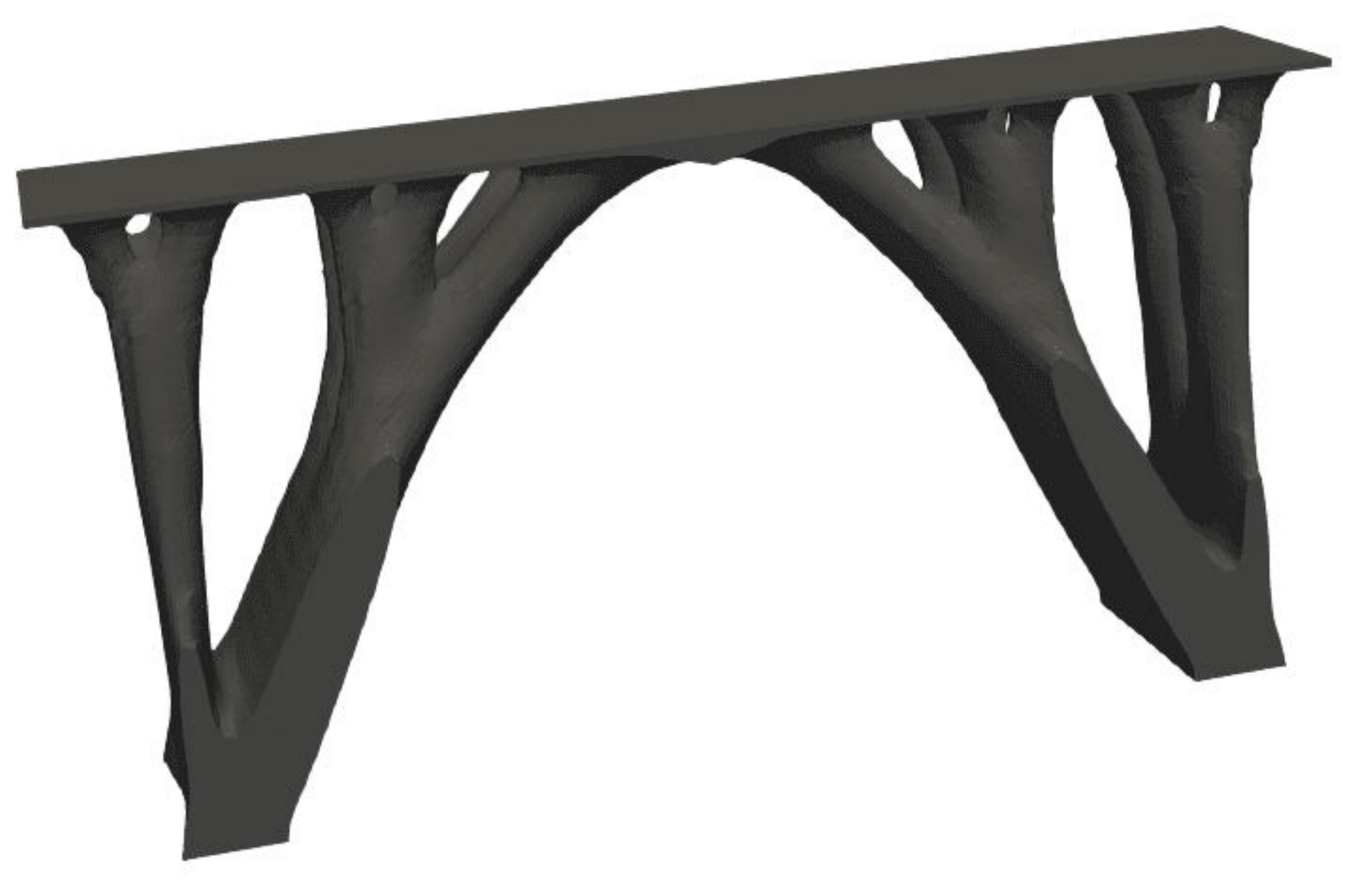}
        \subcaption{Step\,30}
        \label{ra3d-b}
      \end{minipage} 
      \begin{minipage}[t]{0.24\hsize}
        \centering
        \includegraphics[keepaspectratio, scale=0.11]{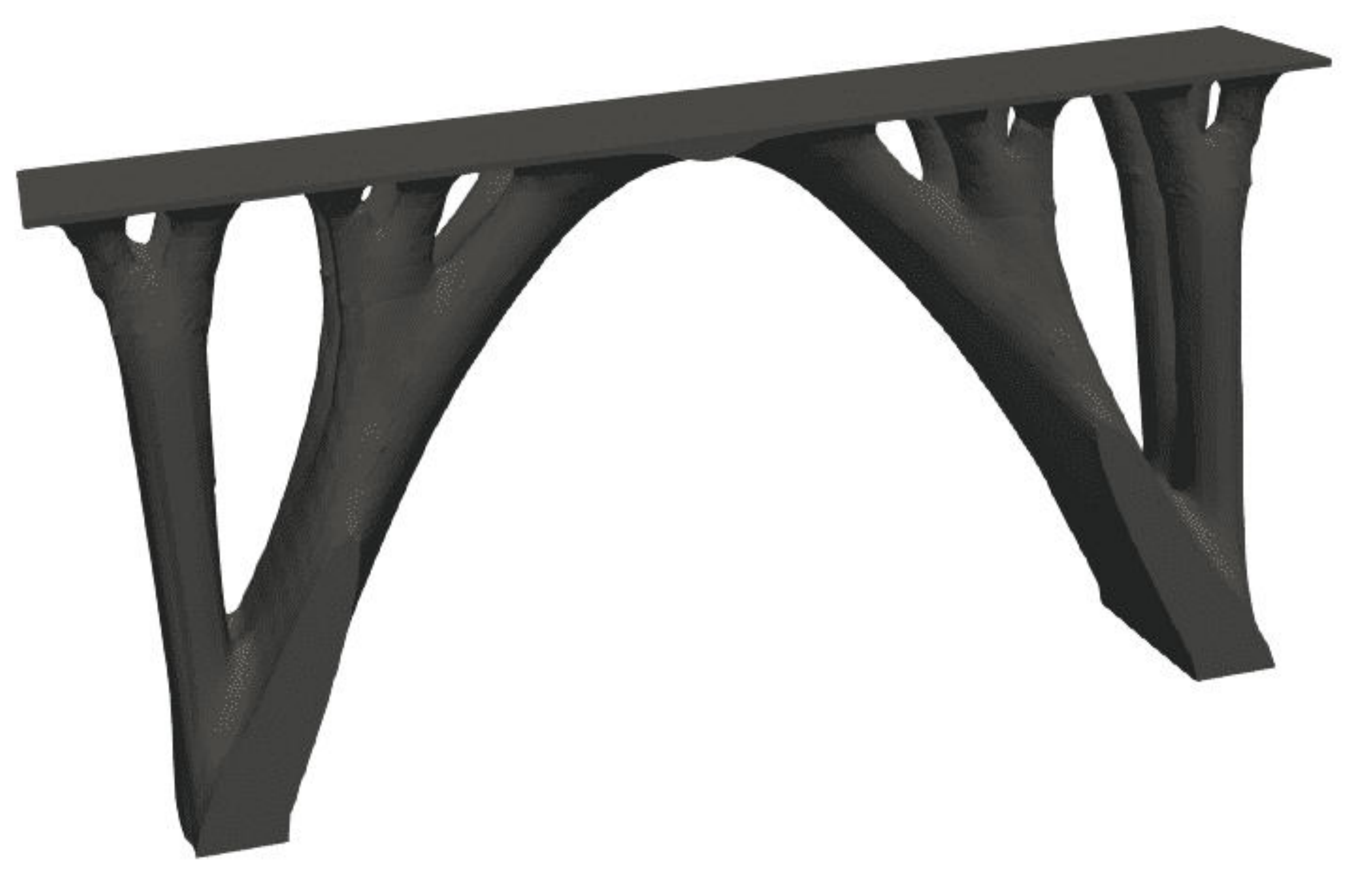}
        \subcaption{Step\,60}
        \label{ra3d-c}
      \end{minipage} 
         \begin{minipage}[t]{0.24\hsize}
        \centering
        \includegraphics[keepaspectratio, scale=0.11]{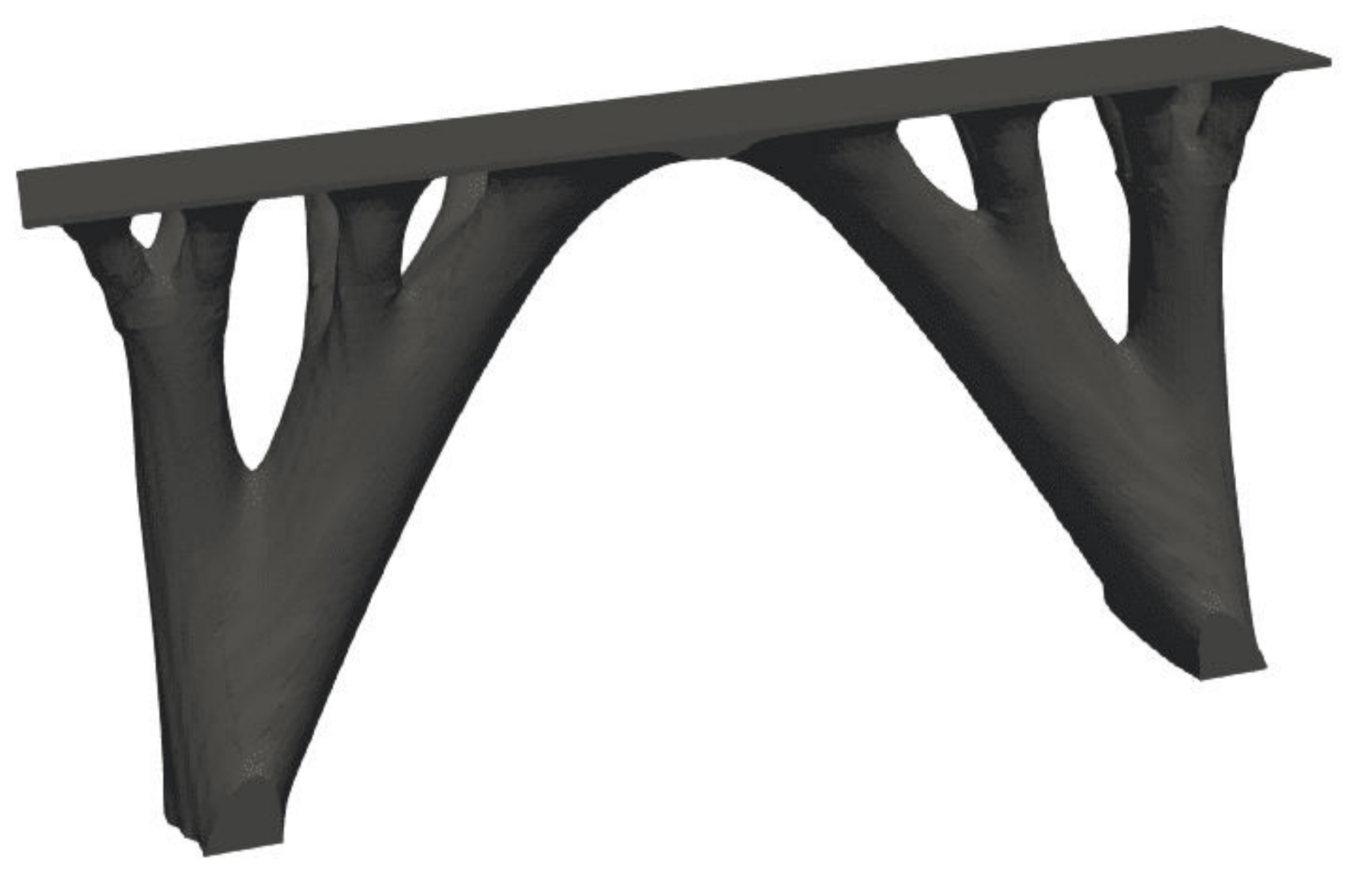}
        \subcaption{Step\,261$^\#$}
        \label{ra3d-d}
      \end{minipage} 
         \\
     \begin{minipage}[t]{0.24\hsize}
        \centering
        \includegraphics[keepaspectratio, scale=0.11]{d3d0.pdf}
        \subcaption{Step\,0}
        \label{ra3d-e}
      \end{minipage} 
      \begin{minipage}[t]{0.24\hsize}
        \centering
        \includegraphics[keepaspectratio, scale=0.11]{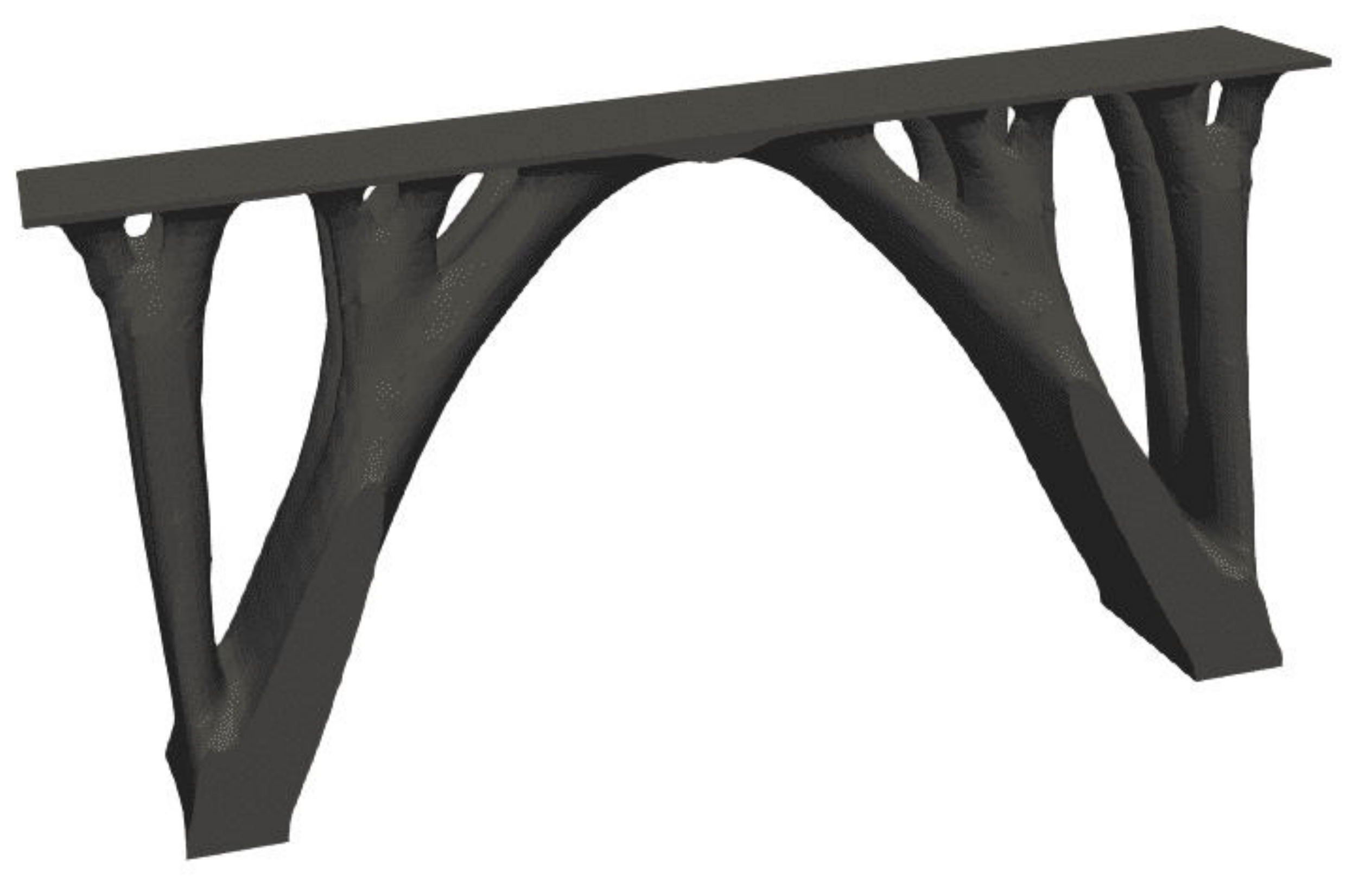}
        \subcaption{Step\,30}
        \label{ra3d-f}
      \end{minipage} 
      \begin{minipage}[t]{0.24\hsize}
        \centering
        \includegraphics[keepaspectratio, scale=0.11]{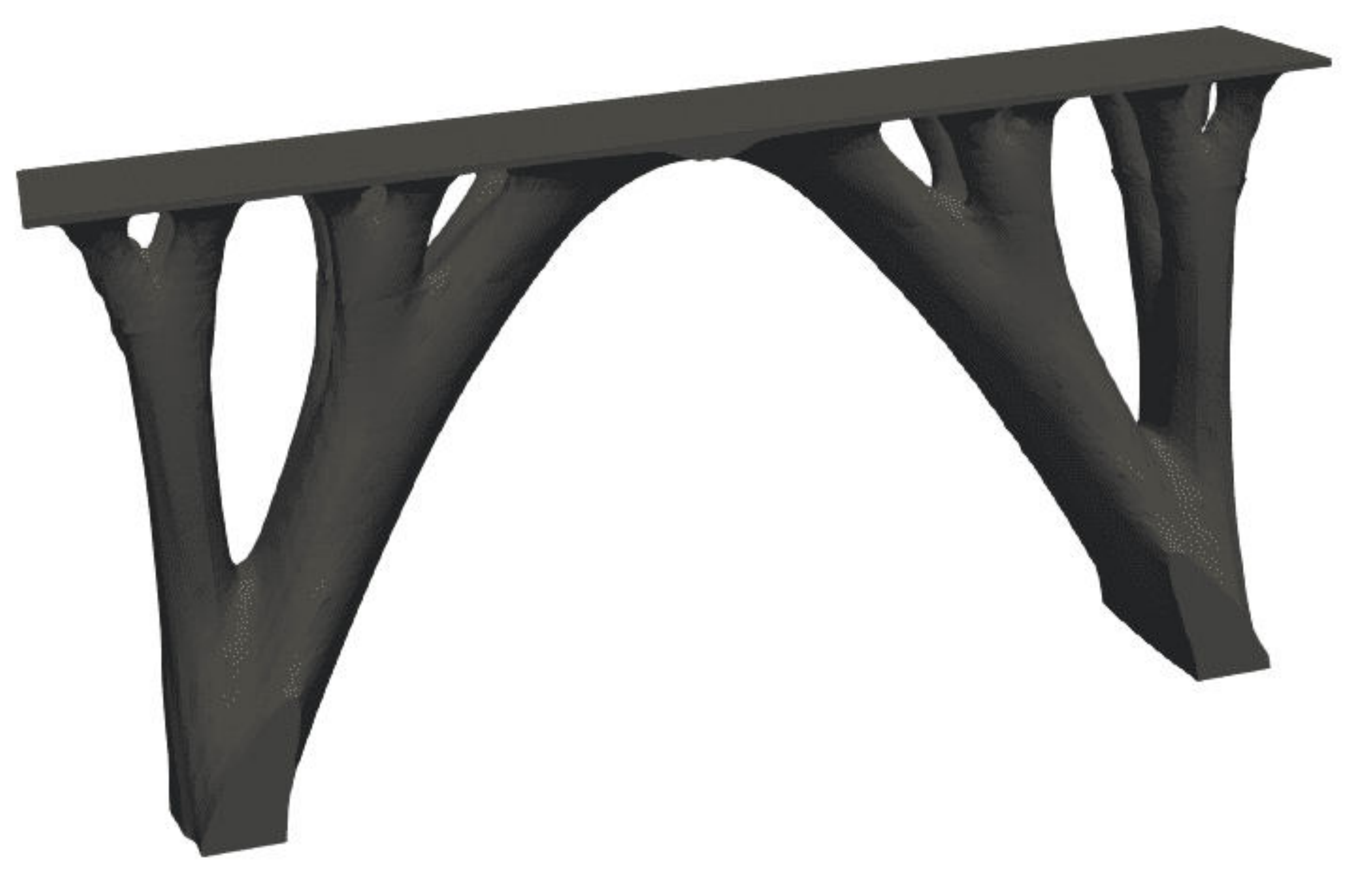}
        \subcaption{Step\,60}
        \label{ra3d-g}
      \end{minipage} 
         \begin{minipage}[t]{0.24\hsize}
        \centering
        \includegraphics[keepaspectratio, scale=0.11]{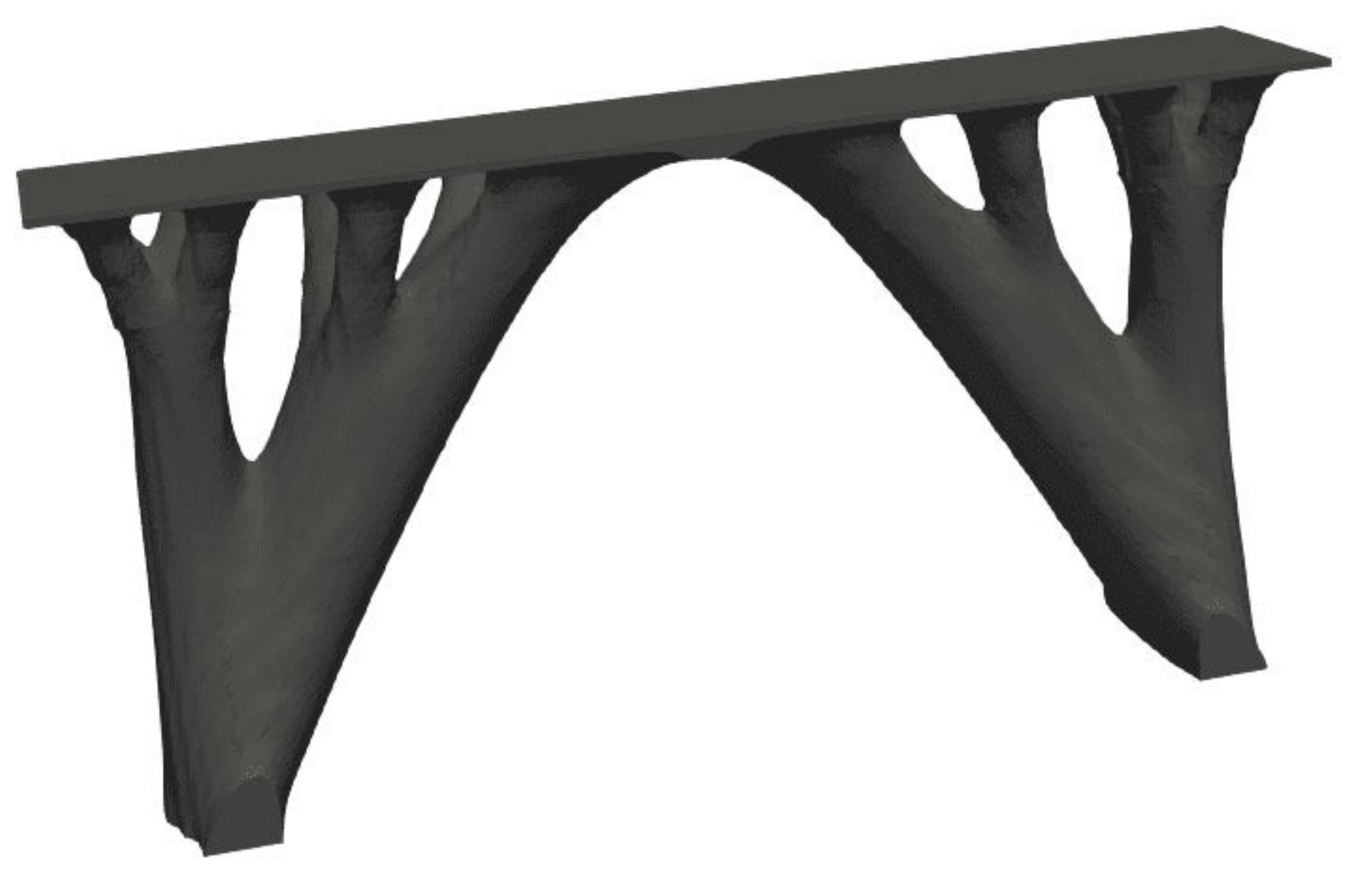}
        \subcaption{Step\,208$^\#$}
        \label{ra3d-h}
      \end{minipage} 
       \end{tabular}
     \caption{ Configuration $\Omega_{\phi_n}\subset D\subset \R^3$ for the case where the initial configuration is the whole domain. 
     Figures (a)--(d) and (e)--(h) represent the results of (RD) and (NLHP), respectively.     
The symbol $^{\#}$ implies the final step. Here the depth of $D$ is set to $0.2$.}
     \label{fig:ra3d}
  \end{figure*}

\begin{figure*}[htbp]
    \begin{tabular}{ccc}
      \hspace*{-5mm} 
      \begin{minipage}[t]{0.49\hsize}
        \centering
        \includegraphics[keepaspectratio, scale=0.33]{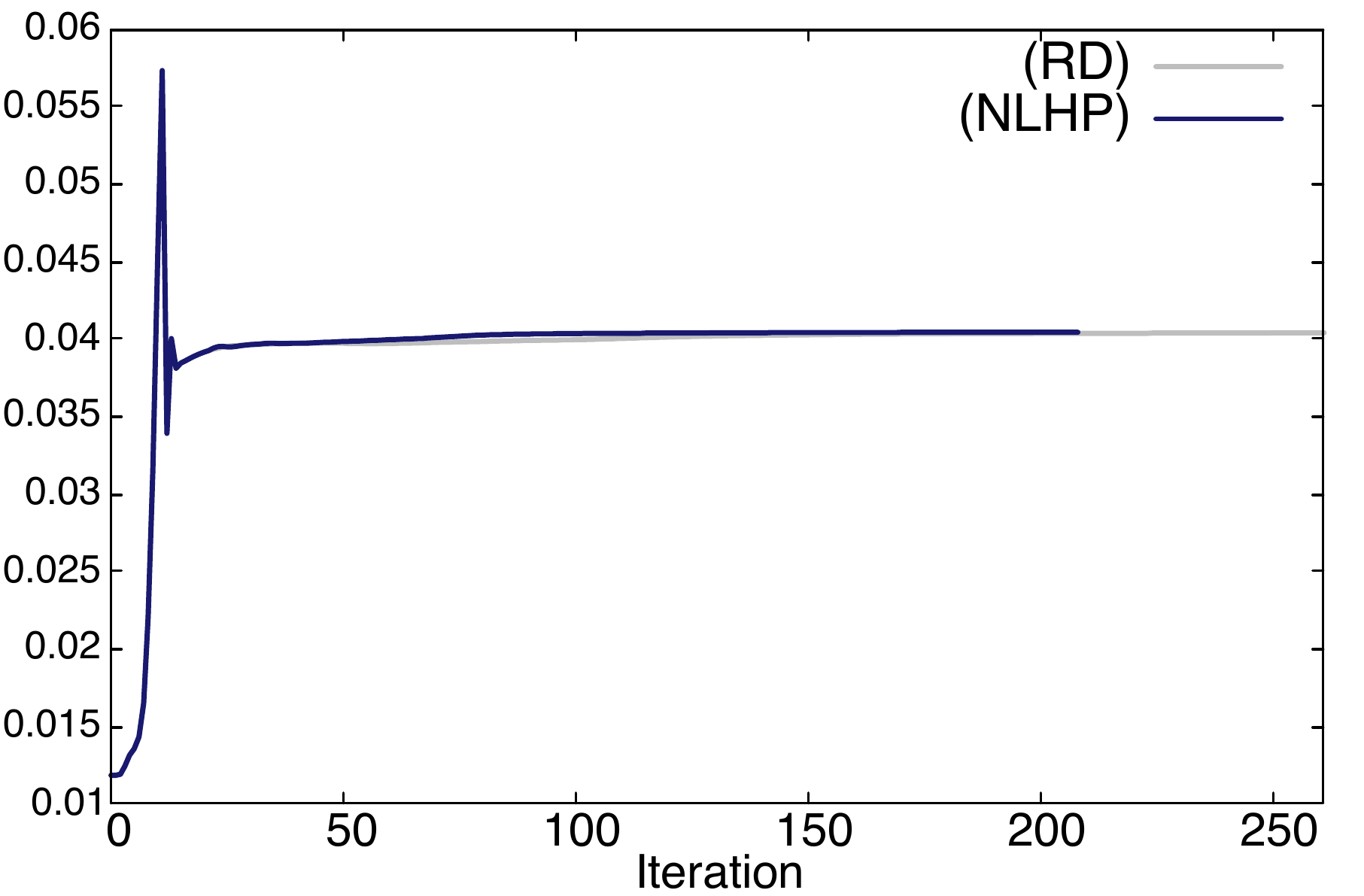}
        \subcaption{$F(\phi_n)$}
        \label{ra3d-1}
      \end{minipage} 
      \begin{minipage}[t]{0.49\hsize}
        \centering
        \includegraphics[keepaspectratio, scale=0.33]{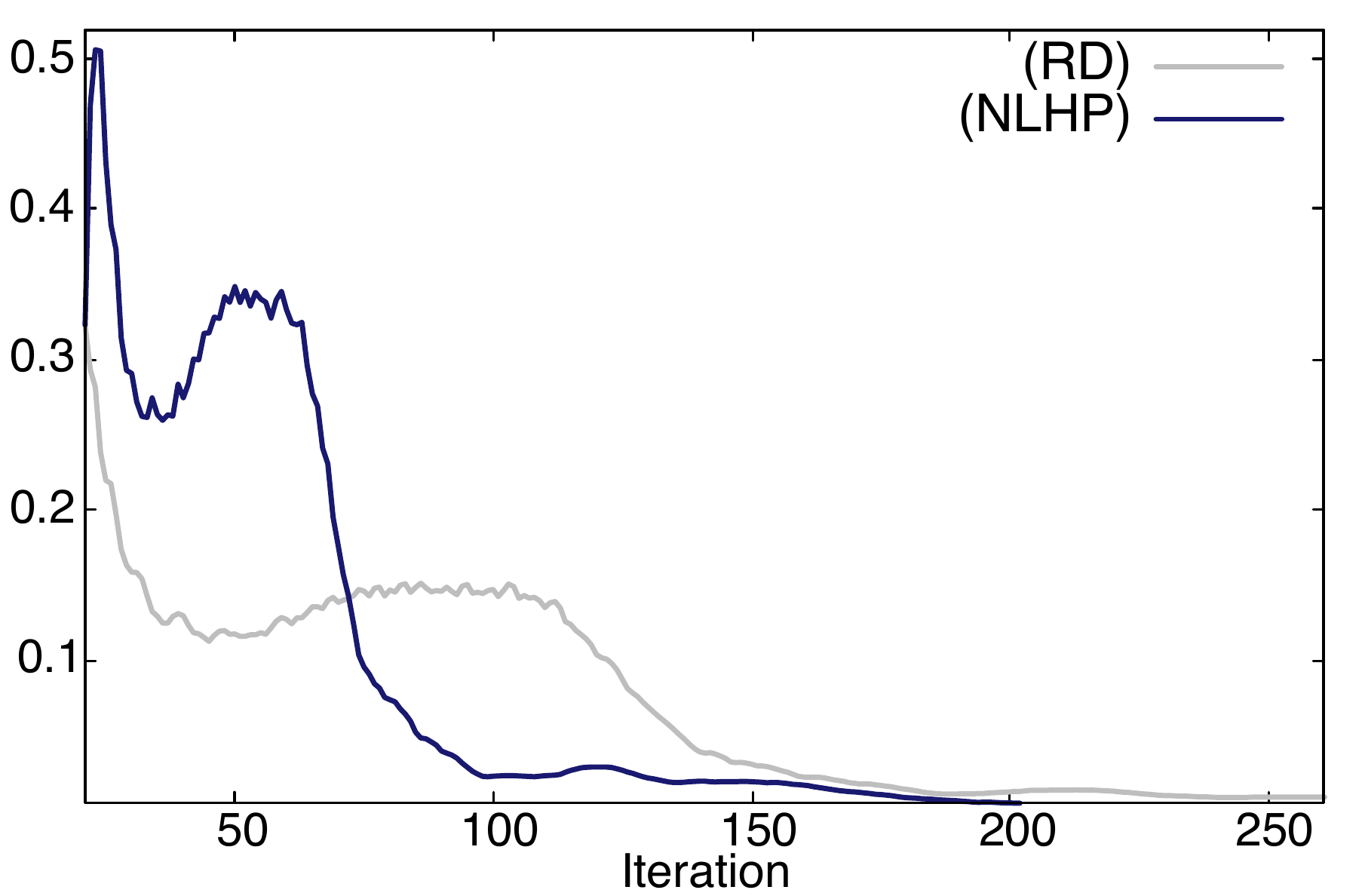}
        \subcaption{$\|\phi_{n+1}-\phi_n\|_{L^{\infty}(D)}$}
        \label{ra3d-2}
      \end{minipage} &
      \end{tabular}
       \caption{ Objective functional and convergence condition for \S \ref{S:ra}-(iv).}
    \label{ra3d}
  \end{figure*}

\vspace{2mm}
\noindent{\bf Case (iv) (Three-dimensional domain).\,} 
We finally show numerically that the assertion in this study can also be obtained for the corresponding three-dimensional case.  
Here we set  $(n_t,h_{\rm max})=(227955,0.0192)$ and $(\tau, G_{\rm max})=(1.0\times 10^{-4}, 0.2)$.
Then Figures \ref{fig:ra3d} and \ref{ra3d} are obtained, and then we see that the  boundary structure in (NLHP) is moving faster than that in (RD).
In particular, Figure \ref{ra3d-2} ensures the validity of the proposed method.
This completes the confirmation for the assertion.


\section{Conclusion}\label{S:conc}
Based on Nesterov's accelerated gradient descent method, this study devised a level set-based topology optimization, which employed reaction-diffusion and nonlinear (damped) wave equations. Furthermore, it was applied to minimum mean compliance problems and yielded numerically faster convergence to an optimized configuration than a reported result \cite{YINT10}. 
Besides, we presented our FreeFEM++ code (for the cantilever model) in the appendix.

\appendix
\section{Free FEM++ code}\label{Ap}
The code for the optimization problem of the cantilever model is as follows\/{\rm:} 
\scriptsize
\begin{lstlisting}%[caption=FEM++ code,label=CODE]
%\begin{verbatim} for arxiv

//	_/_/_/_/_/_/_/_/_/_/_/_/_/_/_/_/_/_/_/
//	_/ (c) 2021 Takayuki YAMADA        _/
//	_/ All right reserved             _/
//	_/ THE UNIVERSITY OF TOKYO       _/
//	_/_/_/_/_/_/_/_/_/_/_/_/_/_/_/_/_/
verbosity=0; 	//  Set message level
//-----------------------------------------------------------
//  parameters
//-----------------------------------------------------------
// Solid Mechanics
real E;			E = 210e9;		
// Young's modulus
real nu; 		nu = 0.3;			
// Poisson's ratio
real lambda;	lambda = E*nu*1./(1+nu)/(1-2.*nu);	
// Lame coefficient
real mu; 		mu = E/2./(1+nu);				
// Lame coefficient
real A1;		A1= E*(-3.*(1-nu)/(2.*(1+nu)*(7-5.*nu)))*
((1-14.*nu+15.*nu^2)*1./(1-2.*nu)^2);
real A2;		A2 = E*(15.*(1-nu)*1./(2.*(1+nu)*(7-5.*nu)));
real matd;		matd = 1e-3;
real matw;		matw = 0.8;

// Optimization
int  MaxLoop;	MaxLoop = 50000; 			
// Maximal iteration number
real epsOpt; 	epsOpt = 1.e-3;      
// Criterion for convergence
int  FlagOptMax; FlagOptMax =10; 		
// period for convergence check

// Volume Constraint
real GvMax;       GvMax = 0.45;   
// upper limit of the volume constraint
int  GvLoop;      GvLoop = 15;    
// extrended
real LagGvA;      LagGvA    = 2.0;   
// coefficient of the linear factor
real LagGvinit;   LagGvinit = 1.0;		
// initial value of lagrange multiplar
real LagGvD;
real LagGvC;     LagGvC = 1.0;    
// coefficient for modification of the Lagrange multiplier
real LagGvMax;    LagGvMax  = 5.0;    
// uppper limmit of the Lagrange multiplier
real LagGvMin;    LagGvMin  = 0.1;   
// lower limmit of the Lagrange multiplier

// Level Set Function
bool UseAcceration 	= true;
int  StatIt;StatIt 	= 5;							
// iteration number of starting the acceration
real L;	L = 1.0;						
// characteristic length
real tau; tau = 5e-4;						
// regularization parameter
real CdF; CdF = 1.2;						
// normalization parameter for sensitivity
real dt; dt=0.7;						
// initial time increment

//-----------------------------------------------------------
//  geometry and Mesh
//-----------------------------------------------------------
int[int] 		labs(4);									
// labels on boundaries of squre: (bottom, right, top, left)
int[int]		MeshNum(2);								
// number of mesh: (x cordinate, y cordinate)
real[int] 	Pos0(2),Pos1(2);					
// Positions of diagonal points of squre: (x cordinate, y cordinate)
int  resMesh;	resMesh=4;							
// mesh resolution

mesh 		Shc1;				
// Mesh of the squre 1
labs		= [1,2,3,4];
MeshNum = [50,11]*resMesh;
Pos0		= [0.0,0.0];
Pos1		= [2.0,0.5-0.04];
Shc1 		= square(MeshNum(0),MeshNum(1),
[Pos0(0)+(Pos1(0)-Pos0(0))*x,Pos0(1)+(Pos1(1)-Pos0(1))*y],
label=labs,flags=1,region=1);

mesh 		Shc2;				
// Mesh of the squre  2
labs		= [5,6,7,8];
MeshNum = [50,2]*resMesh;
Pos0		= [0.0,0.5-0.04];
Pos1		= [2.0,0.5+0.04];
Shc2 		= square(MeshNum(0),MeshNum(1),
[Pos0(0)+(Pos1(0)-Pos0(0))*x,Pos0(1)+(Pos1(1)-Pos0(1))*y],
label=labs,flags=1,region=1);

mesh 		Shc3;				
// Mesh of the squre  3
labs		= [9,10,11,12];
MeshNum = [50,11]*resMesh;
Pos0		= [0.0,0.5+0.04];
Pos1		= [2.0,1.0];
Shc3 		= square(MeshNum(0),MeshNum(1),
[Pos0(0)+(Pos1(0)-Pos0(0))*x,Pos0(1)+(Pos1(1)-Pos0(1))*y],
label=labs,flags=1,region=1);

// combine the meshes
mesh 		Sh=Shc1+Shc2+Shc3;

// define Boundary Label
int[int] wall(3);     	wall   = [4,8,12];		
// 	fixed displacement
int[int] traction(1);   traction = [6];				
//	traction
int[int] phiOne(1);     phiOne  = [6];				
//  phi=1.0

//-----------------------------------------------------------
//  Define functional spaces
//-----------------------------------------------------------
// functional space for governing and djoint equation
fespace VhS2(Sh,P2),		VhS1(Sh,P1),	VhS0(Sh,P0);
fespace VhV2(Sh,[P2,P2]),	VhV1(Sh,[P1,P1]);

//-----------------------------------------------------------
//  Definition of variables
//-----------------------------------------------------------
VhS1 		ophi,oophi,LsfDiff;				
// ophi: previous level set function

//-----------------------------------------------------------
//  Define Characteristic function
//-----------------------------------------------------------
macro chiX(phi)  ((0.5+phi/matw*
(15./16-phi^2/matw^2*(5./8-3./16*phi^2/matw^2)))*
(phi>=-matw)*(phi<=matw)+1.*(phi>matw))//
macro chiP(phi)  ( max((1.-matd)*chiX(phi)+matd,matd)) //
macro chiV(phi)  ( (0.5+phi*(15./16-phi^2*(5./8-3./16*phi^2)))*
(phi>=-1.0)*(phi<=1.0)+1.*(phi>1.0))//

//-----------------------------------------------------------
//  Define Governing Equation
//-----------------------------------------------------------
macro u [u1,u2] 			
//	displacement vector
macro tu [tu1,tu2] 			
//	test function for displacement
macro e(u) [dx(u[0]),dy(u[1]),(dx(u[1])+dy(u[0]))] // strain tensor
macro D [[2.*mu+lambda,lambda,0],[lambda,2.*mu+lambda,0],[0,0,mu]] 
//elastic tensor
macro A [[2.*A2+A1,A1,0],[A1,2.*A2+A1,0],[0,0,A2]] 
// for topoligical derivative
macro g [0,-1.0e3] 				
// traction vector
macro Td1 (u,ophi) ((A*e(u))'*e(u)*chiP(ophi)) 
//Topological derivative 1

// define displacement
VhV1 	u,tu;

// governing equation
problem gov(u,tu)
		=int1d(Sh,traction)(g'*tu)
		-int2d(Sh)((D*e(u))'*e(tu)*chiP(ophi))
		+ on(wall,u1=0,u2=0);

//-----------------------------------------------------------
//  Define PDEs and Initialize the LSF
//-----------------------------------------------------------
VhS1 		phi,Vphi;				
// level set function and its test function
real 		TotTime;				
// time
real 		AbsTd1; 				
// absolute value of Td1
real 		LagGv;					
// Lagrange multiplier for volume constraint
real 		LagGvp;					
//Previous Lagrange multiplier

// Accelerated Reaction Diffusion Equation
problem arde(phi,Vphi)
=int2d(Sh)( (phi*Vphi/dt)+(dx(phi)*dx(Vphi)+dy(phi)*dy(Vphi))*tau*L^2)
-int2d(Sh)((CdF*(Td1(u,ophi)/AbsTd1-LagGv)+((2.*ophi-oophi)/dt))*Vphi)
+on(phiOne,phi=1);

// Reaction Diffusion Equation
problem rde(phi,Vphi)
=int2d(Sh)((phi*Vphi/dt)+(dx(phi)*dx(Vphi)+dy(phi)*dy(Vphi))*tau*L^2)
-int2d(Sh)((CdF*(Td1(u,ophi)/AbsTd1-LagGv)+ophi/dt)*Vphi)
+on(phiOne,phi=1);

//  Initialize Level Sets
phi = 4.*(max(cos(10*x*pi),cos(10*y*pi))-0.5); 
//	case (i) of initialize level set function
//phi = 1.0;//2.0*(((y-0.3)>0)-0.5);               
//	case(ii) of initialize level set function
//phi = 2.0*(((y-0.4)>0)-0.5);			
//	case (iii) of initialize level set function
ophi = phi*(phi<=1.0)*(phi>=-1.0)+1.*(phi>1.0)-1.*(phi<-1.0); 
//  ophi is mapped to satisfy the upper and lower constraints.

//-----------------------------------------------------------
//  Define objective functions and constraint functions
//-----------------------------------------------------------
real obj;				
//objective function
real objPre;			
//Previous objective function
real volFDD;		volFDD = int2d(Sh)(1.); 		
// volume of the fixed design domain
real vol; 											 			
// normalized volume
real volInit;													
// initial normalized volume
real Gv;														
// constraint functional for volume constraint
real GvEx;														
// extended constraint functional for volume constraint
real GvMaxEx;													
// extended maximal constraint functional for volume constraint

//compute initial volume
volInit 	= int2d(Sh)(chiV(ophi));

//_/_/_/_/_/_/_/_/_/_/_/_/_/_/_/_/_/_/_/_/_/_/_/_/_/_/_/_/_/
//  Start Optimization Loop
//_/_/_/_/_/_/_/_/_/_/_/_/_/_/_/_/_/_/_/_/_/_/_/_/_/_/_/_/_/
int  FlagOpt;	FlagOpt 	= 0; 					
// flag for check convergence
for (int iter=0;iter<MaxLoop;iter++) {

	//----------------------------------------------------
	//  Solve Governing Equations (Step 2)
	//----------------------------------------------------
	gov;

	//----------------------------------------------------
	//  Compute Objective and Constraint functions (Step 3)
	//----------------------------------------------------
	// compute objective function
	objPre 		= obj;
	obj 		= int1d(Sh,traction)(g'*u);

	// compute volume
	vol 	 	= int2d(Sh)(chiV(ophi));

	// compute extended upper limit
	if (iter<GvLoop) {
		GvMaxEx 	 
		=volInit/volFDD-((volInit/volFDD-GvMax)/GvLoop)*iter;
	}else{
		GvMaxEx = GvMax;
	}
	GvEx = vol/volFDD-GvMaxEx;		
	//extended volume constraint functional
	Gv = vol/volFDD-GvMax;			
	// volume constraint functional

	//--------------------------------------------------------
	//  Plot
	//--------------------------------------------------------
	{//start plot
			plot(ophi,cmm="ophi",fill=1,
			wait=false,boundary=false,WindowIndex=0);
		}//end plot

	//--------------------------------------------------------
	//  Check Convergence (Step 4)
	//--------------------------------------------------------
	if (LsfDiff[].max<=epsOpt && iter>GvLoop && Gv<=epsOpt) {
			FlagOpt++;
			if(FlagOpt>FlagOptMax){ break;}
	}else{ FlagOpt = 0;}

	//--------------------------------------------------------
	//  Compute Lagrange Multipliers and several normalizers (Step 5)
	//--------------------------------------------------------
	//  Compute Normarizers
	AbsTd1 = int2d(Sh)(abs(Td1(u,ophi)))/volFDD;

	//  Compute Lagrange Multiplier for volume constraint
	if(iter==0){
		LagGv 	= LagGvinit;
		LagGvp	= LagGv;
	}
		// compute update value for Lagrange multiplier
		if(GvEx<0){
		LagGvD = -LagGvC*abs((vol/volFDD-GvMaxEx)/GvMaxEx);
		}else{
		LagGvD = LagGvC*abs((vol/volFDD-GvMaxEx)/GvMaxEx);
		}
		// update Lagrange multiplier
		LagGv = LagGvp+LagGvD;

		// set upper and lower limmit of the Lagrange multiplier
		LagGv = max(min(LagGvMax,LagGv),LagGvMin);
		// save update Lagrange multiplier
		LagGvp=LagGv;
		// set extended Lagrange multiplier 
		using the linear modification
		LagGv=LagGvp*(1.0+LagGvA*((vol/volFDD-GvMaxEx)/GvMaxEx));

	//-----------------------------------------------------
	//  Solve Level set function (Step 6)
	//-----------------------------------------------------
	if(iter<=StatIt){rde;}else
	{
		if(UseAcceration){
			arde;
		}else{
			rde;
		}
	}

	//--------------------------------------------------------
	//  map LSF to "ophi" (Step 7)
	//--------------------------------------------------------
	//  ophi is mapped to satisfy the upper and lower constraints.
	oophi =  ophi;
	ophi =  phi*(phi<=1.0)*(phi>=-1.0)+1.*(phi>1.0)-1.*(phi<-1.0);
	LsfDiff=abs(ophi-oophi);
}// End optimization loop

\end{lstlisting}





\end{document}